\documentclass[english,11pt]{article}

\usepackage[T1]{fontenc}
\usepackage{babel,indentfirst}
\usepackage{amsfonts}  
\usepackage{amsmath}  
\usepackage{amssymb} 
\usepackage{amsxtra}
\usepackage{amscd}  
\usepackage{delarray}  
\usepackage{verbatim}
\usepackage{amsthm}
\usepackage[cmtip,arrow,curve,matrix]{xy}
\usepackage{pb-diagram,pb-xy}
\usepackage{graphicx}
\usepackage{ifthen}
\usepackage{syntonly}
\usepackage{supertabular}

\theoremstyle{plain}
\newtheorem{thm*}{Theorem}
\newtheorem{thm}{Theorem}[section]
\newtheorem{prop}[thm]{Proposition}
\newtheorem{cor}[thm]{Corollary}

\newtheorem{lem}[thm]{Lemma}

\theoremstyle{definition}
\newtheorem{defi}[thm]{Definition}
\newtheorem{cor-defi}[thm]{Corollary -- Definition}

\newtheorem{notation}[thm]{Notation}

\begin{document}

\newcommand*{\new}{\noindent \makebox[0cm][c]{\fbox{\small \textbf{Correction récente.}}}}
\newcommand*{\nnew}{\noindent \makebox[0cm][c]{\fbox{\small \textbf{Correction récente.}}}}
\newcommand*{\neww}{\noindent \makebox[0cm][c]{\fbox{\small \textbf{A modifier.}}}}

\newcommand*{\comb}[2]{\begin{pmatrix}#1\\#2\end{pmatrix}}

\newlength{\taille}                                 
\newcommand*{\taillede}[1]{\settowidth{\taille}{$#1$}}  
\newcommand*{\Taillede}[1]{\settowidth{\taille}{#1}}   
\newcommand*{\donnetaille}{|\rule{\taille}{.2pt}|}    
\newcommand*{\tbox}[3][c]{\settowidth{\taille}{$#2$}\makebox[\taille][#1]{$#3$}} 
\newcommand*{\Tbox}[3][c]{\settowidth{\taille}{#2}\makebox[ \taille][#1]{#3}}   
\newcommand*{\mtbox}[2][c]{\makebox[\taille][#1]{$#2$}}   
\newcommand*{\ftbox}[2][c]{\framebox[\taille][#1]{$#2$}}  
\newcommand*{\Ftbox}[2][c]{\framebox[\taille][#1]{#2}}   
\newcommand*{\bbox}[2][c]{\raisebox{0cm}[0cm][0cm]{\makebox[0cm][#1]{$#2$}}} 
\newcommand*{\lbbox}[1]{\raisebox{0cm}[0cm][0cm]{\makebox[0cm][l]{$#1$}}}     
\newcommand*{\Bbox}[2][c]{\raisebox{0cm}[0cm][0cm]{\makebox[0cm][#1]{#2}}}          
\newcommand*{\ssbbox}[2][c]{\raisebox{0cm}[0cm][0cm]{\makebox[0cm][#1]{\tiny$#2$}}}  
\newcommand*{\sbbox}[2][c]{\raisebox{0cm}[0cm][0cm]{\makebox[0cm][#1]{$\scriptstyle #2$}}} 
\newcommand*{\blanc}[1][2em]{\rule{#1}{0pt}}         
\newcommand*{\espace}[1][1.5em]{\rule{0pt}{#1}}         
\newcommand*{\hblanc}[1][2em]{\rule{#1}{0pt}}          
\newcommand*{\vlblanc}[1][1.3em]{\rule{0pt}{#1}}        
\newcommand*{\vrblanc}[1][.7em]{\raisebox{-#1}{\rule{0pt}{#1}}} 
\newcommand*{\vblanc}{\vlblanc\vrblanc}             

\newcommand*{\videbox}[1]{\settowidth{\taille}{$#1$}\makebox[\taille][l]{}} 
\newcommand*{\Videbox}[1]{\settowidth{\taille}{#1}\makebox[\taille][l]{}}  
\newcommand*{\souligner}[1]{\settowidth{\taille}{#1}\makebox[0cm][l]{#1}\rule[-.4ex]{\taille}{.2pt}}
\newcommand*{\ssouligner}[1]{\settowidth{\taille}{#1}\makebox[0cm][l]{#1}\rule[-.8ex]{\taille}{.2pt}}
\newcommand*{\centrer}[2][]{\noindent\mbox{}\hfill\makebox[0cm][c]{#2}\hfill\ifthenelse{\equal{#1}{}}{\makebox{}}{\makebox[0cm][r]{$(#1)$}}}
\newcommand*{\Centrer}[2][]{\noindent\mbox{}\hfill\makebox[0cm][c]{#2}\hfill\makebox[0cm][r]{$#1$}}

\newcommand*{\indice}[1]{\raisebox{-1.5ex}{\scriptsize$#1$}}  
\newcommand*{\under}[1]{_{\bbox{\textstyle #1}}}    
\newcommand*{\Under}[1]{_{\Bbox{\small #1}}}     
\newcommand*{\upper}[1]{^{\bbox{\textstyle #1}}}    
\newcommand*{\Upper}[1]{^{\Bbox{\small #1}}}     
\newcommand*{\er}{\mbox{\textsuperscript{st} }}       
\newcommand*{\ers}{\mbox{\textsuperscript{sts} }}       
\newcommand*{\st}{\mbox{\textsuperscript{st} }}       
\newcommand*{\sts}{\mbox{\textsuperscript{sts} }}      
\newcommand*{\nd}{\mbox{\textsuperscript{nd} }}      
\newcommand*{\nds}{\mbox{\textsuperscript{nds} }}       
\newcommand*{\rd}{\mbox{\textsuperscript{rd} }}       
\newcommand*{\rds}{\mbox{\textsuperscript{rds} }}       
\newcommand*{\eme}{\mbox{\textsuperscript{th} }}      
\newcommand*{\emes}{\mbox{\textsuperscript{thes} }}      
\newcommand*{\ent}{.} 
\newcommand*{\vect}[1]{\mbox{$\overrightarrow{#1}$}}
\newcommand*{\fin}{\rule{0cm}{1ex}\linebreak[3]\rule{0cm}{1ex}\hfill\mbox{$\Box$}} 
\newcommand*{\ouu}{\mbox{ \quad where }}
\newcommand*{\subsetnot}{\makebox[0cm]{\raisebox{-.5ex}[0cm][0cm]{$\;\;\;\scriptstyle\not$}}\subseteq}
\newcommand*{\vide}{\varnothing}
\newcommand*{\point}{\mbox{$\bullet$}}
\newcommand*{\tiret}{\noindent\rule{1ex}{0ex}--\rule{1ex}{0ex}}
\newcommand*{\ttiret}{\noindent\rule{4ex}{0ex}--\rule{1ex}{0ex}}
\newcommand*{\Tiret}{\noindent\rule{1ex}{0ex}\point\rule{1ex}{0ex}}
\newcommand*{\Point}{\noindent\rule{1ex}{0ex}\point\rule{1ex}{0ex}}
\newcommand*{\defin}{\bbox{\hat}=}
\newcommand*{\ou}{\mbox{ or }}
\newcommand*{\et}{\mbox{ and }}


\newcommand*{\A}{\mathcal{A}}
\newcommand*{\B}{\mathcal{B}}
\newcommand*{\C}{\mathcal{C}}
\newcommand*{\D}{\mathcal{D}}
\newcommand*{\E}{\mathcal{E}}
\newcommand*{\F}{\mathcal{F}}
\newcommand*{\G}{\mathcal{G}}
\renewcommand*{\H}{\mathcal{H}}
\newcommand*{\I}{\mathcal{I}}
\newcommand*{\J}{\mathcal{J}}
\newcommand*{\K}{\mathcal{K}}
\renewcommand*{\L}{\mathcal{L}}
\newcommand*{\M}{\mathcal{M}}
\newcommand*{\N}{\mathcal{N}}
\renewcommand*{\O}{\mathcal{O}}
\renewcommand*{\P}{\mathcal{P}}
\newcommand*{\Q}{\mathcal{Q}}
\newcommand*{\R}{\mathcal{R}}
\renewcommand*{\S}{\mathcal{S}}
\newcommand*{\T}{\mathcal{T}}
\newcommand*{\U}{\mathcal{U}}
\newcommand*{\V}{\mathcal{V}}
\newcommand*{\W}{\mathcal{W}}
\newcommand*{\X}{\mathcal{X}}
\newcommand*{\Y}{\mathcal{Y}}
\newcommand*{\Z}{\mathcal{Z}}

\renewcommand*{\AA}{\mathbb{A}}
\newcommand*{\BB}{\mathbb{B}}
\newcommand*{\CC}{\mathbb{C}}
\newcommand*{\DD}{\mathbb{D}}
\newcommand*{\EE}{\mathbb{E}}
\newcommand*{\FF}{\mathbb{F}}
\newcommand*{\GG}{\mathbb{G}}
\newcommand*{\HH}{\mathbb{H}}
\newcommand*{\II}{\mathbb{I}}
\newcommand*{\JJ}{\mathbb{J}}
\newcommand*{\KK}{\mathbb{K}}
\newcommand*{\LL}{\mathbb{L}}
\newcommand*{\MM}{\mathbb{M}}
\newcommand*{\NN}{\mathbb{N}}
\newcommand*{\OO}{\mathbb{O}}
\newcommand*{\PP}{\mathbb{P}}
\newcommand*{\QQ}{\mathbb{Q}}
\newcommand*{\RR}{\mathbb{R}}
\renewcommand*{\SS}{\mathbb{S}}
\newcommand*{\TT}{\mathbb{T}}
\newcommand*{\UU}{\mathbb{U}}
\newcommand*{\VV}{\mathbb{V}}
\newcommand*{\WW}{\mathbb{W}}
\newcommand*{\XX}{\mathbb{X}}
\newcommand*{\YY}{\mathbb{Y}}
\newcommand*{\ZZ}[1][]{\mathbb{Z}\ifthenelse{\equal{#1}{}}{}{/#1\mathbb{Z}}}  
\newcommand*{\Aa}{\mathfrak{A}}
\newcommand*{\Ss}{\mathfrak{S}}

\newcommand*{\Act}{\mathrm{Act}}
\newcommand*{\Ad}{\mathrm{Ad}}
\newcommand*{\An}{\mathrm{An}}
\newcommand*{\Arc}{\mathrm{Arc}}
\newcommand*{\Aut}{\mathrm{Aut}}
\newcommand*{\CPA}{\mathrm{CPA}}
\newcommand*{\CSB}{\mathrm{CSB}}
\newcommand*{\cc}{\mathrm{Comp}}
\newcommand*{\Centr}{\mathrm{C\widetilde{ent}r}}
\newcommand*{\col}{\mathrm{col}}
\newcommand*{\com}{\mathrm{com}\:}
\newcommand*{\Comp}{\mathrm{Comp}}
\newcommand*{\Conf}{\mathrm{Conf}}
\newcommand*{\Couple}{\mathrm{Couple}}
\newcommand*{\Courb}{\C\mathrm{urv}}
\newcommand*{\Curv}{\C\mathrm{urv}}
\newcommand*{\cut}{\mathrm{cut}}
\newcommand*{\Diff}{\mathrm{Diff}}
\newcommand*{\Diffeo}{\mathrm{Diffeo}}
\newcommand*{\Div}{\mathrm{div}}
\newcommand*{\End}{\mathrm{End}}
\newcommand*{\Eq}{\mathrm{$\E$q}}
\newcommand*{\ext}{\mathrm{ext}}
\newcommand*{\Ext}{\mathrm{Ext}}
\newcommand*{\for}{\mathrm{for}}
\newcommand*{\Fix}{\mathrm{Fix}}
\newcommand*{\Gen}{\mathrm{Gen}}
\newcommand*{\genus}{\mathrm{genus}}
\newcommand*{\grad}{\mathrm{grad}\,}
\newcommand*{\Homeo}{\mathrm{Homeo}}
\newcommand*{\Hom}{\mathrm{Hom}}
\newcommand*{\Id}{\mathrm{Id}\,}
\newcommand*{\im}{\mathrm{Im}\:}
\let\Im\im
\newcommand*{\Imp}{\mathrm{Odd}}
\newcommand*{\Ind}{\mathrm{Ind}}
\newcommand*{\Inj}{\mathrm{Inj}}
\newcommand*{\inj}{\mathrm{inj}}
\newcommand*{\Inn}{\mathrm{Inn}}
\newcommand*{\Int}{\mathrm{Int}}
\newcommand*{\Inv}{\mathrm{Inv}}
\newcommand*{\Isom}{\mathrm{Isom}}
\newcommand*{\Jac}{\mathrm{Jac}}
\newcommand*{\Ker}{\mathrm{Ker}\:}
\let\ker\Ker
\newcommand*{\Max}{\mathrm{Max}}
\newcommand*{\MBS}{\mathrm{MBS}}
\newcommand*{\Min}{\mathrm{Min}}
\newcommand*{\Norm}{\mathrm{Norm}}
\newcommand*{\marq}{\mathrm{marq}}
\newcommand*{\Mod}{\M\mathrm{od}}
\newcommand*{\obs}{\mathrm{\emph{obs}}}
\newcommand*{\Odd}{\mathrm{Odd}}
\newcommand*{\Opt}{\O\mathrm{pt}}
\newcommand*{\Or}{\mathrm{Or}}

\newcommand*{\Orb}{\mathrm{Orb}}
\newcommand*{\Out}{\mathrm{Out}}
\newcommand*{\Pair}{\mathrm{Pair}}
\newcommand*{\Per}{\mathrm{Per}}
\newcommand*{\pos}{\mathrm{pos}}
\newcommand*{\Plong}{\mathrm{Plong}}
\newcommand*{\PMod}{\P\M\mathrm{od}}
\newcommand*{\pr}{\mathit{pr}}
\newcommand*{\proj}{proj}
\newcommand*{\PSL}{\mathrm{PSL}}
\newcommand*{\rec}{\mathrm{rec}\:}
\newcommand*{\plug}{\mathrm{plug}\:}
\newcommand*{\rg}{\mathrm{rg}\:}
\newcommand*{\Restr}{\mathrm{restr}}
\newcommand*{\SDiff}{\mathrm{SDiff}}
\newcommand*{\sign}{\mathrm{signature}}
\newcommand*{\SMod}{\S\M\mathrm{od}}
\newcommand*{\SL}{\mathrm{SL}}
\newcommand*{\Sp}{\mathrm{Sp}}
\renewcommand*{\sp}{\mathrm{sp}}
\newcommand*{\sq}{\mathrm{sq}}
\newcommand*{\Sub}{\S\mathrm{ub}}
\newcommand*{\Subet}{\S\mathrm{ub}^{\mathrm{and}}}
\newcommand*{\supp}{\mathrm{supp}}
\newcommand*{\Stab}{\mathrm{Stab}}
\newcommand*{\Sym}{\mathrm{Sym}}
\newcommand*{\tr}{\mbox{\large tr}}
\newcommand*{\Tr}{\mathrm{$\T$\!\emph{r}}}
\newcommand*{\Triv}{\mathrm{$\T$\!\emph{riv}}}
\newcommand*{\Tub}{\mathrm{\emph{Tub}}}
\newcommand*{\Vect}{\mathrm{Vect}}

\newcommand*{\deff}{\overset{\mbox{\scriptsize déf}}{=}}
\newcommand*{\rond}{\circ}
\newcommand*{\DEF}[5]{\mbox{$#1\,:\;
\begin{array}{ccc}#2 &\longrightarrow& #3\\#4 &\longmapsto& #5\end{array}$}}
\newcommand*{\DEFF}[7]{\mbox{$#1\,:\;
\begin{array}{ccc}#2 &\longrightarrow& #3\\#4 &\longmapsto& #5\\#6 &\longmapsto& #7\end{array}$}}
\newcommand*{\sDEF}[5]{\mbox{$\left\{#1\scriptsize\,:\;
\begin{array}{ccc}#2 &\rightarrow& #3\\#4 &\mapsto& #5\end{array}\right\}$}}

\newcommand*{\restr}[2]{\mbox{$#1$}\;\!\rule[-1ex]{.2pt}{2ex}\,\raisebox{-1ex}{$\scriptstyle #2$}}
\newcommand*{\bij}{\mbox{\makebox[0ex][l]{$\hookrightarrow$}$\ \to$}}
\newcommand*{\xtwoheadrightarrow}[2][]{\xrightarrow[\,#1\,]{\,#2\,}\!\!\!\!\!\!\:\to}
\newcommand*{\xhookrightarrow}[2][]{\raisebox{.7ex}{$\scriptstyle \subset$}\!\!\!\!\xrightarrow[\,#1\,]{\,\,#2\,}}
\newcommand*{\xmapsto}[1]{\raisebox{0.12ex}{\tiny$\vdash$}\normalsize\!\!\!\!\xrightarrow{\ #1\ }}
\newcommand*{\homeo}[1][]{\xrightarrow[#1]{\ \cong\ }}
\newcommand*{\longto}{\longrightarrow}


 \renewcommand*{\Tilde}{\makebox[0cm][l]{\raisebox{1.7ex}{\scriptsize $\sim$}}}
 \newcommand*{\ttilde}{\makebox[0cm][l]{\raisebox{1.7ex}{\scriptsize $\:\sim$}}}
 \newcommand*{\ntilde}{\makebox[0cm][l]{\raisebox{1.7ex}{\scriptsize $\sim$}}}
 \newcommand*{\ptilde}{\makebox[0cm][l]{\raisebox{1.4ex}{\scriptsize $\sim$}}}
 \newcommand*{\mtilde}{\makebox[0cm][l]{\raisebox{1ex}{\scriptsize $\sim$}}}

 \newcommand*{\tildeA}{\makebox[0cm][l]{\raisebox{1.7ex}{\scriptsize $\;\sim$}}A}
 \newcommand*{\tildeB}{\makebox[0cm][l]{\raisebox{1.7ex}{\scriptsize $\:\sim$}}B} 
 \newcommand*{\tildeC}{\makebox[0cm][l]{\raisebox{1.7ex}{\scriptsize $\:\sim$}}C}
 \newcommand*{\tildeF}{\makebox[0cm][l]{\raisebox{1.7ex}{\scriptsize $\:\sim$}}F}
 \newcommand*{\tildeG}{\makebox[0cm][l]{\raisebox{1.7ex}{\scriptsize $\:\sim$}}G}
 \newcommand*{\tildeH}{\makebox[0cm][l]{\raisebox{1.7ex}{\scriptsize $\:\sim$}}H}
 \newcommand*{\tildeP}{\makebox[0cm][l]{\raisebox{1.7ex}{\scriptsize $\:\sim$}}P}
 \newcommand*{\tildeQ}{\makebox[0cm][l]{\raisebox{1.7ex}{\scriptsize $\:\sim$}}Q}
 \newcommand*{\tildeR}{\makebox[0cm][l]{\raisebox{1.7ex}{\scriptsize $\:\sim$}}R}
 \newcommand*{\tildeT}{\makebox[0cm][l]{\raisebox{1.7ex}{\scriptsize $\,\sim$}}T}
 \newcommand*{\tildeX}{\makebox[0cm][l]{\raisebox{1.7ex}{\scriptsize $\;\sim$}}\X}
 \newcommand*{\tildeY}{\makebox[0cm][l]{\raisebox{1.7ex}{\scriptsize $\,\sim$}}\Y}
 \newcommand*{\tildeZ}{\makebox[0cm][l]{\raisebox{1.7ex}{\scriptsize $\;\sim$}}\Z}
 \renewcommand*{\tildeY}{\makebox[0cm][l]{\raisebox{1.7ex}{\scriptsize $\,\sim$}}Y}
 \renewcommand*{\tildeZ}{\makebox[0cm][l]{\raisebox{1.7ex}{\scriptsize $\;\sim$}}Z}
 \newcommand*{\tildex}{\makebox[0cm][l]{\raisebox{1.2ex}[0cm]{\scriptsize $\;\!\sim$}}\X}
 \newcommand*{\tildez}{\makebox[0cm][l]{\raisebox{1.2ex}[0cm]{\scriptsize $\;\!\sim$}}Z}
 \newcommand*{\tildeS}{\makebox[0cm][l]{\raisebox{1.7ex}{\scriptsize $\;\!\sim$}}\Sigma}
 \newcommand*{\tildeSigma}{\makebox[0cm][l]{\raisebox{1.7ex}{\scriptsize $\;\!\sim$}}\Sigma}
 \newcommand*{\tildes}{{\makebox[0cm][l]{\raisebox{1.2ex}[0cm]{\scriptsize $\sim$}}\Sigma}}
 \newcommand*{\tildealpha}{\makebox[0cm][l]{\raisebox{1.3ex}{\scriptsize $\;\!\sim$}}\alpha}
 \newcommand*{\tildebeta}{\makebox[0cm][l]{\raisebox{1.5ex}{\scriptsize $\;\!\sim$}}\beta}
 \newcommand*{\tildeb}{\makebox[0cm][l]{\raisebox{1.5ex}{\scriptsize $\sim$}}b}
 \newcommand*{\tildekappa}{\makebox[0cm][l]{\raisebox{1.3ex}{\scriptsize $\;\!\sim$}}\kappa}
 \newcommand*{\tildetau}{\makebox[0cm][l]{\raisebox{1.3ex}{\scriptsize $\;\!\sim$}}\tau}
 \newcommand*{\tildegamma}{\makebox[0cm][l]{\raisebox{1.3ex}{\scriptsize $\;\!\sim$}}\gamma}

 \newcommand*{\rhor}{\rho_{\mathrm{ref}}}
 \newcommand*{\trhor}{\tilde\rho_{\mathrm{ref}}}


\newcommand*{\lllambda}{\makebox[0.08ex][l]{$\lambda$}\makebox[0.08ex][l]{$\lambda$}\lambda}
\newcommand*{\tttau}{\makebox[0.08ex][l]{$\tau$}\makebox[0.08ex][l]{$\tau$}\tau}
\newcommand*{\mmmu}{\makebox[0.04ex][l]{$\mu$}\makebox[0.04ex][l]{$\mu$}
                  \makebox[0.04ex][l]{$\mu$}\makebox[0.04ex][l]{$\mu$}\mu}
\newcommand*{\ooomega}{\makebox[0.04ex][l]{$\omega$}\makebox[0.04ex][l]{$\omega$}
                  \makebox[0.04ex][l]{$\omega$}\makebox[0.04ex][l]{$\omega$}\omega}

\newcommand*{\bbb}[1]{\makebox[0.03ex][l]{$#1$}\makebox[0.03ex][l]{$#1$}
                  \makebox[0.03ex][l]{$#1$}\makebox[0.03ex][l]{$#1$}#1}
\newcommand*{\BBB}[1]{\makebox[0.03ex][l]{#1}\makebox[0.03ex][l]{#1}
                  \makebox[0.03ex][l]{#1}\makebox[0.03ex][l]{#1}#1}

\newcommand*{\rhomon}{\rho_n^{\mathrm{mon}}}
\newcommand*{\bord}{\partial}
\newcommand*{\bordint}{\bord^{\mathrm{inn}}}
\newcommand*{\bordext}{\bord^{\mathrm{ext}}}
\newcommand*{\bordnat}{\bord^{\mathrm{nat}}}

\newcommand*{\Bord}{\mathrm{Bndy}}
\newcommand*{\Bordext}{\Bord^{\mathrm{ext}}}
\newcommand*{\Bordint}{\Bord^{\mathrm{int}}}
\newcommand*{\Bordnat}{\Bord^{\mathrm{nat}}}
\newcommand*{\Sigmaext}{\Sigma^{\mathrm{ext}}}
\newcommand*{\Sigmaint}{\Sigma^{\mathrm{int}}}
\newcommand*{\Sigmanat}{\Sigma^{\mathrm{nat}}}
\newcommand*{\Sigmaij}{\Sigma_{i,\,j}}
\newcommand*{\Sigmaijk}{\Sigma_{i,\,j,\,k}}
\newcommand*{\Sigmapijk}{\Sigma'_{i,\,j,\,k}}

\newcommand*{\dbord}{\mathrm{d}}
\newcommand*{\dbordint}{\dbord^{\mathrm{int}}}
\newcommand*{\dbordnat}{\dbord^{\mathrm{nat}}}
\newcommand*{\ddint}{\dd^{\mathrm{int}}}          
\newcommand*{\ddnat}{\dd^{\mathrm{nat}}}

\newcommand*{\Courbeset}{\C\mathrm{urves}^{\mathrm{ext}}}
\newcommand*{\Curvesxt}{\C\mathrm{urves}^{\mathrm{ext}}}

\newcommand*{\Gaug}{{\G_n^{\mathrm{aug}}}}


\newcommand*{\quot}[2]{\textnormal{{\Large $^{#1}\!/\! _{#2}$}}}
\newcommand*{\lquot}[2]{\textnormal{{\Large $_{#2}\!\backslash\! ^{#1}$}}}
\newcommand*{\rquot}[2]{\textnormal{{\Large $^{#1}\!/\! _{#2}$}}}

\newcommand*{\dessous}[2]{\underset{\scriptstyle #2}{#1}}    
\newcommand*{\dessus}[2]{\overset{\scriptstyle #2}{#1}}     
\newcommand*{\dessousdessus}[3]{\dessus{\dessous{#1}{#2}}{#3}}  
\newcommand*{\ePsi}{\stackrel{\,{}_=}{\Psi}}

\newcommand*{\we}[2]{\raisebox{.7ex}{$#1$}/\!\raisebox{-.7ex}{$#2$}}
\newcommand*{\ssur}[2]{\left(\begin{array}{c}#1\\#2\end{array}\right)}

\newcommand*{\zz}[1][2]{\textnormal{{\large $^{\mathbb{Z}}\!/\! _{#1 \mathbb{Z}}$}}}
\newcommand*{\zzz}[1][2]{\mathbb{Z}/#1 \mathbb{Z}}
\newcommand*{\tc}{\otimes\mathbb{C}}
\newcommand*{\orth}[1]{{#1}^\bot}
\newcommand*{\orthh}{{}^\bot}
\newcommand*{\comp}[1]{{}^c#1}
\newcommand*{\sgn}{\mathrm{\;\!sgn\:}}
\newcommand*{\adh}[1]{\overline{#1}}
\newcommand*{\ouv}[1]{\overset{\circ}{#1}}
\newcommand*{\Ouv}[1]{\overset{\circ}{\overline{#1}}}
\newcommand*{\pt}{\{*\}}
\newcommand*{\transv}{\pitchfork}


\newcommand*{\od}{\partial}
\newcommand*{\dd}{\mathrm{d}}
\newcommand*{\partiel}[3]{\frac{\partial#1}{\partial#2}(#3)}
\newcommand*{\bigpartiel}[3]{{\displaystyle\frac{\partial#1}{\partial#2}}(#3)}
\newcommand*{\smallpartiel}[3]{{\textstyle\frac{\partial#1}{\partial#2}}{\scriptstyle(#3)}}
\newcommand*{\ppartiel}[2]{{\frac{\partial#1}{\partial#2}}}
\newcommand*{\bigppartiel}[2]{{\displaystyle\frac{\partial#1}{\partial#2}}}
\newcommand*{\smallppartiel}[2]{{\textstyle\frac{\partial#1}{\partial#2}}}
\newcommand*{\dpartiel}[4]{\textstyle\frac{\od^2 #1}{\od#2\od#3}\scriptstyle(#4)}
\newcommand*{\derive}[2]{{\displaystyle\frac{\dd#1}{\dd#2}}}
\newcommand*{\matjac}[5][x]{                              
   \begin{pmatrix}                                  
    \smallpartiel{#2_1}{#1_1}{#5} &\dots& \smallpartiel{#2_1}{#1_{#4}}{#5}\\
    \vdots &\ddots & \vdots\\                            
     \smallpartiel{#2_{#3}}{#1_1}{#5} &\dots& \smallpartiel{#2_{#3}}{#1_{#4}}{#5}\\ 
     \end{pmatrix}}                                
\newcommand*{\mathess}[4][x]{                            
   \begin{pmatrix}                                
    \smallpartiel{^2#2}{{#1_1}^2}{#4} &\dots&                   
    \dpartiel{#2}{#1_1}{#1_#3}{#4}\\
     \vdots &\ddots & \vdots\\                          
     \dpartiel{#2}{#1_#3}{#1_1}{#4} &\dots&                   
      \smallppartiel{^2#2}{{#1_#3}^2}{#4}\\                   
      \end{pmatrix}}


\newcommand*{\petito}[2]{\dessous{o}{#1\to 0}\!(#2)}
\newcommand*{\ppetito}[2]{\dessous{o}{#1}\!(#2)}
\newcommand*{\dint}{\int\!\!\!\!\!\;\int}
\newcommand*{\oiint}{\bbox[l]{\!\!\:\int}\bbox[l]{\;\!\int}{\mbox{\small O}}} 
\newcommand*{\bigoiint}{\bbox[l]{\!\!\:\displaystyle\int}\bbox[l]           
         {\displaystyle\,\int}\raisebox{-,2ex}[3,6ex][2,4ex]{\large \:\!O}}  
\newcommand*{\tint}{\int\!\!\!\int\!\!\!\int}
\newcommand*{\vois}[2]{\V^{#1}_{#2}}
\newcommand*{\cvois}[2]{\overline{\V}^{#1}_{#2}}
\newcommand*{\ovois}[2]{\tau^{#1}_{#2}}


\newcommand*{\liste}[3][1]{\mbox{$(#2_{#1},\dots,#2_{#3})$}}            
\newcommand*{\vecteur}[2]{\begin{pmatrix}#1_1\\ \vdots \\ #1_{#2}\end{pmatrix}}
\newcommand*{\matcol}[2]{\begin{pmatrix}#1\\ \vdots \\ #2\end{pmatrix}}

\newcommand*{\bddots}{\ssbbox{\ddots}}
\newcommand*{\bvdots}{\ssbbox{\vdots}}
\newcommand*{\bhdots}{\ssbbox{\dots}}

\newcommand*{\trans}{{}^{\mathrm{\small t}}}


\newcommand*{\TITRE}[1]{\noindent\textbf{#1}}
\newcommand*{\LTITRE}[1]{\noindent\large\textbf{#1}\normalsize}
\newcommand*{\DEM}[1][\!\!]{\noindent\textbf{Proof #1. }}
\newcommand*{\REM}[1][\!\!]{\noindent\textbf{Remark #1. }}
\newcommand*{\REMs}{\noindent\textbf{Remarks. }}
\newcommand*{\CONV}{\noindent\textbf{Convention. }}
\newcommand*{\EX}[1][\!\!]{\noindent\textbf{Example #1. }}
\newcommand*{\EXs}{\noindent\textbf{Examples. }}

\newcounter{etape}[subsection]
\newenvironment{etape}[1][] {\stepcounter{etape}\TITRE{Step
\theetape\ifthenelse{\equal{#1}{}}{.}{: #1.}} \itshape}{}

\newcounter{num}[subsection]
\newenvironment{num}{\stepcounter{num}\TITRE{\thenum.} \itshape}

\newcommand*{\lab}[1]{\newcounter{#1}[section]\setcounter{#1}{\value{etape}}}

\newcounter{remark}[section]
\newenvironment{remark}[1][]
{\stepcounter{remark}\TITRE{Remarque
\theremarke\ifthenelse{\equal{#1}{}}{.}{: #1.}} \itshape}{}

\newcommand*{\encadrer}[1] {
\begin{center}
\fbox{\begin{minipage}{12cm}
 #1
\end{minipage}} \end{center}}

\newcommand*{\Includegraphics}[1] {
\begin{center}
\makebox[0cm][c]{\includegraphics{#1}} \end{center}}

\newcommand*{\Newpage}{\mbox{}\newpage}
\newcommand*{\pageimpaire}{}

\title{Geometric representations of the braid groups}
\author{Fabrice Castel\footnote{Research partially supported by
National University of Singapore Research Grant R-146-000-137-112}}
\bigskip
\date{October 2010}
\maketitle
\bigskip
\bigskip
\bigskip

\TITRE{Abstract:}
We call \emph{geometric representation} any representation of a group in the mapping class group
of a surface. Let $\Sigma_{g,b}$ be the orientable connected compact
surface of genus $g$ with $b$ boundary components, and $\PMod(\Sigma_{g,\,b})$
the associated mapping class group preserving each boundary component.
The main theorem of this paper concerns geometric representations
of the braid group $\B_n$ with $n\geqslant 6$ strands in $\PMod(\Sigma_{g,\,b})$
subject to the only condition that $g\leqslant n/2$.
We prove that under this condition, such representations are either
\emph{cyclic}, that is, their images are cyclic
groups, or are what we call \emph{transvections of monodromy homomorphisms},
defined in the text.

This leads to different results. They will be proved in later papers,
but we explain for each of them how they are deduced from our main theorem.
These corollaries concern
five families of groups: the braid groups $\B_n$ for all $n\geqslant 6$,
the Artin groups of type $D_n$ for all $n\geqslant 6$,
the Artin groups of type $E_n$ for $n\in\{6,\,7,\,8\}$,
the mapping class groups $\PMod(\Sigma_{g,\,b})$
(preserving each boundary component) and the mapping class groups
$\Mod(\Sigma_{g,\,b},\,\bord\Sigma_{g,\,b})$ (preserving
the boundary pointwise), for all $g\geqslant 2$ and all $b\geqslant 0$.
For each of these five families except the Artin groups of type $E_n$,
we are able to describe precisely
the (always remarkable) structure of the endomorphisms.
\bigskip
\bigskip

\TITRE{Keywords:} surface, mapping class group, braid group, rigidity, geometric representation,
Nielsen Thurston's classification, transvection, monodromy morphism.

\bigskip
\bigskip

\TITRE{AMS Subject Classification:}  Primary 20F38, 57M07. Secondary 57M99, 20F36, 20E36, 57M05.
\medskip

            \Newpage

\tableofcontents             

      \Newpage

\section{Introduction}
\bigskip

Some theorems have been established concerning injective endomorphisms of the braid groups or
the mapping class groups. For both families of groups, the outer automorphisms groups were first computed,
in 1981 by Dyer and Grossman, [DyGr], for the braid groups,
and in 1986 - 1988 by Ivanov, [Iv1], and McCarthy, [Mc1], for the mapping class groups.
Later, again for both families of groups, the set of the injective endomorphisms was described,
in 1999 by Ivanov and McCarthy, [IvMc], for the mapping class groups,
and in 2000 by Bell and Margalit, [BeMa], for the braid groups.
For both families of groups, the previous authors
could say a little more on injective endomorphisms
between two different braid groups or the mapping class groups of two distinct surfaces,
with strong restrictions, though.
Our main theorem concerns the homomorphisms from the braid group in the mapping class group,
making the bridge between both families of groups. As corollaries, we are able to prove the
previous theorems and even to strengthen them by getting rid of the injective hypothesis.
\bigskip

\subsection{Presentation of the main objects}
\medskip

\TITRE{The surfaces} in this paper will always be compact, orientable and
oriented 2-manifold whose each connected component
is of negative Euler characteristic. They may have nonempty
boundary. Classically, we denote by $\bord\Sigma$ the topological
boundary of a surface $\Sigma$ and we denote by $\Sigma_{g,\,b}$ the connected surface
of genus $g$ with $b$ boundary components. Implicitly, $g$ and
$b$ will always satisfy $2-2g-b\leqslant-1$.
\bigskip

\TITRE{The mapping class groups} $\Mod(\Sigma)$ of a
surface $\Sigma$ is the group of isotopy classes of
orientation-preserving diffeomorphisms of $\Sigma$.
We denote by $\PMod(\Sigma)$ the finite index subgroup of $\Mod(\Sigma)$
consisting of elements of $\Mod(\Sigma)$ that permute neither the connected components,
nor the boundary components.
\bigskip

\TITRE{A geometric representation} of a group $G$ is a homomorphism from $G$ in the mapping
class group of some surface $\Sigma_{g,\,b}$.
For instance, J. Birman, A. Lubotsky and J. McCarthy have shown in [BiLuMc]
that the maximal rank of the abelian subgroups of the mapping class group $\Mod(\Sigma_{g,\,b})$
is equal to $3g-3+b$.
In other words, they have shown that the group $\ZZ^n$ has faithful geometric representations
in $\Mod(\Sigma_{g,\,b})$ if and only if $n\leqslant  3g-3+b$ holds.
In this paper, we will investigate geometric representations of the braid group.
\bigskip

\TITRE{The braid group} $\B_n$ on $n$ strands is the group defined by the following presentation,
which we will call the \emph{classic presentation}:
\smallskip

 \centrer{$\B_n=\big\langle
\tau_1,\,\tau_2,\dots,\,\tau_{n-1}\ |\ \left\{\begin{array}{l}
 \tau_i \tau_j= \tau_j \tau_i \mbox{ \ if \ } |i-j|\not=1\;\\
 \tau_i \tau_j \tau_i = \tau_j \tau_i \tau_j \mbox{ \ if \ }
 |i-j|=1\ .
\end{array}\right.\,\big\rangle$}
\smallskip

\noindent The generators of this presentation are called the
\emph{standard generators of $\B_n$}.
\bigskip

\subsection{The main theorems:
Theorems \ref{thm:Theorem_principal}, \ref{thm:Theorem_principal_case_à_boundary} and \ref{thm:injectivity_des_homomorphisms}}
\medskip

The aim of this paper is to describe all geometric representations of the braid group $\B_n$
in the mapping class groups $\PMod(\Sigma_{g,\,b})$ and
$\Mod(\Sigma_{g,\,b},\,\bord\Sigma_{g,\,b})$.
The only hypothesis is that the number $n$ of strands of
$\B_n$ and the genus $g$ of the surface $\Sigma_{g,\,b}$ must satisfy:

\centrer{$g\leqslant \frac{n}{2}$,}
\smallskip

\noindent whereas $b$ is any nonnegative integer, as
long as the surface $\Sigma_{g,\,b}$ is of negative Euler characteristic.
Notice that the elements of $\PMod(\Sigma_{g,\,b})$ and
$\Mod(\Sigma,\,\bord\Sigma_{g,\,b})$ do not permute the
boundary components of the surface. Under these assumptions,
we will show that all \emph{non-cyclic} representations
(that is, whose images are not cyclic groups) of the braid group
are ``variations'' of the classical monodromy representations defined below.
As an illustration of this fact, when $n$ is odd, all non-cyclic representations of $\B_n$
are \emph{not far from} being conjugate from each other.
\bigskip

\begin{figure}[!h]
      \Includegraphics{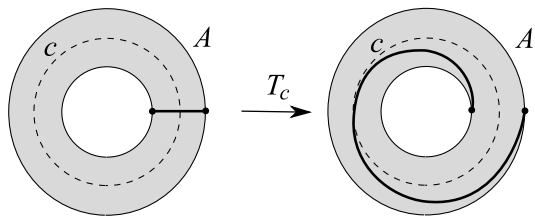}
      \caption{a Dehn twist along a curve $c$.}
      \label{intro-fig:dehnTwist}
     \end{figure}

\TITRE{A monodromy representation} of $\B_n$ will be a geometric representation of $\B_n$ which sends
the different standard generators of $\B_n$ on distinct Dehn twists (or their inverses).
We recall that the Dehn twist $T_c$ along a simple closed curve $c$ in a surface $\Sigma$
is the mapping class of $\Mod(\Sigma)$ defined as follows.
Let $A$ denote an annular neighbourhood of $c$.
Then $T_c$ is the isotopy class of the homeomorphism
which is the identity on $\Sigma \smallsetminus \mathrm{int}(A)$
and transforms the interior of $A$ as shown in Figure \ref{intro-fig:dehnTwist}.
Remember that two Dehn twists $T_a$ and $T_b$ along two distinct curves $a$ and $b$
verify $T_aT_bT_a=T_bT_aT_b$ if and only if the curves $a$ and $b$ meet in one point,
whereas they commute if and only if the curves $a$ and $b$ are disjoint.
Consequently, a
monodromy representation of $\B_n$ can be characterized by the data of
an integer $\varepsilon\in\{\pm 1\}$ and
a \emph{$n$-chain of curves}, that is to say an ordered $(n-1)$-tuple of curves $(a_1,\,a_2,\dots,\,a_{n-1})$
such that for all $i,\,j\in\{1,\dots,\,n-1\}$, the curves $a_i$ and
$a_j$ are disjoint when $|i-j|\not= 1$, and intersect in exactly one
point when $|i-j|=1$. So the monodromy representation $\rho$ associated to the couple
$\big(\,(a_1,\,a_2,\dots,\,a_{n-1}),\,\varepsilon\,\big)$ is defined by setting for all $i\leqslant n-1$:
\smallskip

\centrer{$\rho(\tau_i)={T_{a_i}}^{\varepsilon}$.}
\bigskip

Let $n$ be an integer greater than or equal to 3. Let $G$
be any group, $\rho$ a homomorphism from $\B_n$ to $G$ and $w$ an
element lying in the centralizer of $\rho(\B_n)$ in $G$. The
\emph{transvection $\rho_1$ of $\rho$ with direction $w$} is the homomorphism
defined by setting for all $i\leqslant n-1$:
\smallskip

\centrer{$\rho_1(\tau_i)=\rho(\tau_i)\,w$.}
\smallskip

\TITRE{A transvection of a monodromy representation} $\rho$ of $\B_n$
is thus characterized by a (unique, as will be shown later) triple
$\big(\,(a_1,\,a_2,\dots,\,a_{n-1}),\,\varepsilon,\,W\big)$ where $W$ is a mapping class which preserves each
curve $a_i$, $i\leqslant n-1$, and is defined by setting for all $i\leqslant n-1$:
\smallskip

\centrer{$\rho(\tau_i)={T_{a_i}}^{\varepsilon}\,W$.}
\medskip

\noindent We denote by $\Sigma(\rho)$ the tubular
neighbourhood of the union of the curves $a_i$ where $i$ ranges from 1 to $n-1$.
As will be shown in Proposition \ref{prop:centralizer}, $W$ must preserve
$\Sigma(\rho)$ and induces in $\Mod(\Sigma(\rho))$
either the identity or the \emph{hyper-elliptic element} related to the curves $(a_1,\,a_2,\dots,\,a_{n-1})$,
that is, the only involution of $\Sigma(\rho)$ that preserves $a_i$
and that swaps both sides of its tubular neighbourhood, for all $i\in\{1,\dots,\,n-1\}$.
\medskip

\TITRE{A cyclic representation} of $\B_n$ is a representation whose image is cyclic.
Equivalently, all the standard generators of $\B_n$ have the same image.
\bigskip

Here is the main theorems of this paper.
\smallskip

\begin{thm*}[Representations of $\B_n$ in $\PMod(\Sigma)$ with $n\geqslant 6$]
  \label{thm:Theorem_principal}
\mbox{}\\Let $n$ be an integer greater than or equal to 6 and
$\Sigma$ a surface $\Sigma_{g,\,b}$ with $g\leqslant \frac{n}{2}$ and $b\geqslant 0$.
Then any representation $\rho$ of $\B_n$ in $\PMod(\Sigma)$ is either cyclic,
or is a transvection of monodromy representation. Moreover, such
transvections of monodromy representations exist if and only if
$g\geqslant \frac{n}{2}-1$.
\end{thm*}
\bigskip

Let $\Mod(\Sigma,\,\bord \Sigma)$ be the group of
isotopy classes of all orientation-preserving diffeomorphisms
that coincide with the identity on $\bord\Sigma$. Then:
\smallskip

\begin{thm*}[Representations of $\B_n$ in $\Mod(\Sigma,\,\bord\Sigma)$ with $n\geqslant 6$]
    \label{thm:Theorem_principal_case_à_boundary}
\mbox{}\\ Theorem \ref{thm:Theorem_principal} still holds when $\PMod(\Sigma)$ is replaced by $\Mod(\Sigma,\,\bord\Sigma)$.
\end{thm*}
\bigskip

In addition, we have an easy characterization of the faithfulness of the geometric representations of $\B_n$
when the genus verifies $g\leqslant\frac{n}{2}$.
\medskip

\begin{thm*}[Faithful geometric representation of the braid groups]
    \label{thm:injectivity_des_homomorphisms}
\mbox{}\\Let $n$ be an integer greater than or equal to 6 and
$\Sigma$ a surface $\Sigma_{g,\,b}$ with $g\leqslant \frac{n}{2}$.\\
Let $\rho$ be a homomorphism from $\B_n$ to $\Mod(\Sigma,\,\bord\Sigma)$
or to $\PMod(\Sigma)$.

(i) Case of $\Mod(\Sigma,\,\bord\Sigma)$. The homomorphism $\rho$ is
injective if and only if it is a transvection of monodromy homomorphism
such that the boundary components of $\Sigma(\rho)$ do not bound any
disk in $\Sigma$.

(ii) Case of $\PMod(\Sigma)$. The homomorphism $\rho$ is injective if
and only if it is a transvection of monodromy homomorphism such that
boundary components of $\Sigma(\rho)$ do not bound any disk in
$\Sigma$ and at least one boundary component of $\Sigma(\rho)$ is
not isotopic to any boundary component of $\Sigma$. \end{thm*}
\bigskip

Theorems \ref{thm:Theorem_principal},
\ref{thm:Theorem_principal_case_à_boundary} and \ref{thm:injectivity_des_homomorphisms}
answer the natural question:
``What are the (injective and not injective) homomorphisms between two classical groups?''.
A similar work has been undertaken
by E. Artin (cf. [At3]) and V. Lin (cf. [Ln2]) concerning the homomorphisms from
$\B_n$ to the symmetric group $\Ss_m$ with $m\leqslant n$ (Artin, 1947) and
then with $m\leqslant 2n$ (Lin, 1970 - 2004). However, the importance of
Theorems \ref{thm:Theorem_principal},
\ref{thm:Theorem_principal_case_à_boundary} and \ref{thm:injectivity_des_homomorphisms}
comes above all from their corollaries (see a quick description below,
and a detailed description in Section \ref{sec:corollaries_result}).
\bigskip
\bigskip

\subsection{Summary of the Corollaries}
\medskip

We will deduce the following from Theorems \ref{thm:Theorem_principal},
\ref{thm:Theorem_principal_case_à_boundary} and \ref{thm:injectivity_des_homomorphisms}:
\begin{itemize}
  \item[\;\;(i)]
description of the homomorphisms from $\B_n$ to $\B_m$ with $m\leqslant n+1$.
This generalizes a theorem of Bell and Margalit (see [BeMa]) to the non-injective case;
  \item[\;(ii)]
triviality of $\Mod(\Sigma_{g,\,b})\to\Mod(\Sigma_{g',\,b'})$ as soon as $g'<g$,
which generalizes  a theorem of W. Harvey and M. Korkmaz  (see [HvKo]) to the nonempty boundary case;
  \item[(iii)]
triviality of $\Mod(\Sigma_{g,\,0})\to\Mod(\Sigma_{g+1,\,0})$,
which is a partial answer to a conjecture of J. Berrick and M. Matthey
(see [BkMt]) -- this case can also be derived from the number theoritic conditions of [BkMt];
  \item[(iv)]
description of homomorphisms $\Mod(\Sigma_{g,\,b})\to\Mod(\Sigma_{g+1,\,b})$;
which is a generalization in several ways of a theorem of N. Ivanov and J. McCarthy (see [IvMc]);
  \item[(v)]
non-injectivity of all geometric representations of $A(E_n)$ ($n\in\{6,7,8\}$)
in the expected mapping class groups. This completes a Theorem of B. Wajnryb (see [W]);
  \item[(vi)]
description of the endomorphisms (and automorphisms) of $A(D_n)$
which is a generalization in several ways of a theorem of J. Crisp and L. Paris (see [CrPa], Theorem 4.9).
\end{itemize}
\bigskip

The major improvement of these theorems in comparison with older similar ones comes
from the fact that we deal with homomorphisms instead of injective homomorphisms.
In addition, we gather these various theorems as consequences
of one single main result (essentially Theorem \ref{thm:Theorem_principal_case_à_boundary}).
In Section \ref{sec:corollaries_result}, we state precisely these corollaries and we explain
how to prove them from Theorem \ref{thm:Theorem_principal_case_à_boundary}.
We postpone detailed proofs to further publications. However, the interested reader is referred to
the author's Ph.D thesis for what concerns items (i) - (iv), cf. [Ca].
\bigskip

\subsection{Outline}
\medskip

\noindent The organization of this paper is as follows.
We present in Section \ref{sec:corollaries_result} the corollaries
of Theorems \ref{thm:Theorem_principal},
\ref{thm:Theorem_principal_case_à_boundary} and \ref{thm:injectivity_des_homomorphisms} and explain
briefly how the former can be deduced from the latter.
In Section \ref{sec:conventions_defintions}, we introduce the main definitions and fix the notations
concerning surfaces and their mapping class groups. We also give some general results
on mapping class groups that we shall use later on. Most of them are well-known.
Those ones that are not are proved in the appendix section.
In Section \ref{sec:Proof_Theorem_2_&_3},
we show an important result about lifting braid group geometric representations from $\PMod(\Sigma)$
to $\Mod(\Sigma,\bord\Sigma)$. We then deduce Theorems \ref{thm:Theorem_principal_case_à_boundary}
and \ref{thm:injectivity_des_homomorphisms} from this result and Theorem \ref{thm:Theorem_principal}.
Sections \ref{sec:curves_peripheral} to \ref{sec:end_grande_proof}
are devoted to the proof of Theorem \ref{thm:Theorem_principal}. The major tool will be
Nielsen - Thurston's theory and the canonical reduction system, denoted by $\sigma(F)$
for a mapping class $F$, which were introduced
by Birman, Lubotzky, McCarthy and Ivanov in [BiLuMc] and [Iv2].
The key point of our proof will consist in finding a nice partition of $\sigma(A_i)$
where the $A_i$ are the images of the standard generators of the braid group for a given
geometric representation. This is done in Section \ref{sec:curves_spéciales}. 
\bigskip

\section{Corollaries of Theorems \ref{thm:Theorem_principal},
\ref{thm:Theorem_principal_case_à_boundary} and \ref{thm:injectivity_des_homomorphisms}: Theorems \ref{thm:homomorphisms from Bn_in_Bm} - \ref{thm:automorphisms_A(Dn)}}
        \label{sec:corollaries_result}
\bigskip

\subsection{Homomorphisms between braid groups}
\medskip

The description of the homomorphisms between braid groups
involves the Garside element $\Delta_n$ of the braid group $\B_n$, defined by
$\Delta_n=\tau_1(\tau_2\tau_1)\dots(\tau_{n-1}\tau_{n-2}\dots\tau_1)$,
and the unique involution $\Inv$ of the braid group that sends each standard
generator on its inverse.
\bigskip

Theorems \ref{thm:Theorem_principal_case_à_boundary} and \ref{thm:injectivity_des_homomorphisms} allow us to show the following.
\bigskip

\begin{thm*}[Homomorphisms between braid groups]
            \label{thm:homomorphisms from Bn_in_Bm}
\mbox{}\\ Let $n$ and $m$ be two
integers such that $n\geqslant 6$ and $3\leqslant m\leqslant n+1$.
\begin{itemize}
  \item[\;\;(i)] Case where $m<n$: ([Ln1], 1982) any homomorphism $\varphi$ from $\B_n$ to $\B_m$ is cyclic.

  \item[\;(ii)] Case where $m=n$: any noncyclic endomorphism $\varphi$ of $\B_n$ is a transvection
of inner automorphism possibly precomposed by the involution $\Inv$:
there exist $\gamma\in\B_n$, $\varepsilon=\pm1$ and $k\in\ZZ$ such that for
all $i\leqslant n-1$, we have:

\centrer{$\varphi(\tau_i)=\gamma\,\tau_i^{\,\varepsilon}\,\gamma^{-1}\,{\Delta_n}^{2k}$.}

\noindent Such a homomorphism $\varphi$ is always injective and is an automorphism if and only if $k=0$.

  \item[(iii)] Case where $m=n+1$: let us consider the group $\B_n$ as the subgroup
of $\B_{n+1}$ spanned by the $n-1$ first standard generators of
$\B_{n+1}$. Then, any homomorphism $\varphi$ from $\B_n$ to $\B_{n+1}$
is the restriction to $\B_n$ of a homomorphism from $\B_{n+1}$ in
itself, up to transvection. According to item (ii), if $\varphi$ is
not cyclic, then there exist $\gamma\in\B_{n+1}$,
$\varepsilon=\pm1$ and $k,\,\ell\in\ZZ$ such that for all $i\leqslant n-1$, we have:
\smallskip

\centrer{$\varphi(\tau_i)=\gamma\,\tau_i^{\,\varepsilon}\,{\Delta_n}^{2k}\,\gamma^{-1}\,{\Delta_{n+1}}^{2\ell}$.}
\smallskip

\noindent Such a homomorphism $\varphi$ is always injective (but never surjective). \end{itemize}
\end{thm*}
\bigskip

Historically, the first result in
this direction was found in 1981: Dyer and Grossman computed algebraically the
outer automorphisms group of $\B_n$ and shew that $\Out(\B_n)=\ZZ[2]$
for any $n\geqslant 2$.
As we easily can show that nontrivial transvections
of automorphisms of the braid groups are never onto,
Theorem \ref{thm:homomorphisms from Bn_in_Bm}
provides a new proof of their theorem. In 1982, Lin shew that all homomorphisms from
$\B_n$ to $\B_m$ with $m<n$ are trivial (see [Ln1]). Then, in 2000,
R. Bell and D. Margalit  described all the
injective homomorphisms from $\B_n$ to $\B_n$ and $\B_{n+1}$
(see [BeMa]). Their proof used an algebraic characterization of the braid twists
in which the maximal rank of the abelian subgroups of $\B_n$ played a central
role (that is why they needed the homomorphisms to be injective).
Keeping in mind these results, the new part of
Theorem \ref{thm:homomorphisms from Bn_in_Bm} is the following:
\smallskip

\centrer{\emph{The only non-injective homomorphisms from $\B_n$ to $\B_n$
or $\B_{n+1}$ are cyclic}.}
\smallskip

\noindent In addition Theorem \ref{thm:homomorphisms from Bn_in_Bm} provides a unified proof for
all these results.
\bigskip

\bigskip

\TITRE{Steps of the proof of Theorem \ref{thm:homomorphisms from Bn_in_Bm}.}
\smallskip

\noindent \emph{Step 1.} Let $m\leqslant n+1$ and let $\varphi$ be a non-cyclic homomorphism from $\B_n$ to $\B_m$.
Let $\Sigma$ be the connected surface of genus $\lfloor\frac{m-1}{2}\rfloor$
with one or two boundary components depending on $m$.
Let us denote by $\rho$ a faithful representation of $\B_m$ in $\Mod(\Sigma,\,\bord\Sigma)$.
Thus we have an $(m-1)$-chain of curves $(c_i)_{i\leqslant m-1}$
such that for all $i\leqslant m-1$, we have:
\smallskip

\centrer[1]{$\rho(\tau_i) = T_{c_i}$.}
\smallskip

\noindent Let us consider the following commutative diagram.
\medskip

\centrer{$\begin{diagram}
\node{\B_{n}}\arrow[2]{e,t}{\rho\rond\varphi}\arrow{se,b}{\varphi}
\node[2]{\Mod(\Sigma,\,\bord\Sigma)}\\
\node[2]{\B_{m}}\arrow{ne,b}{\rho} \end{diagram}$}
\smallskip

\noindent According to Theorem \ref{thm:Theorem_principal_case_à_boundary},
$\rho\rond\varphi$ is a transvection of monodromy homomorphism, so there exists a triple
$\big((a_i)_{i\leqslant n-1},\,\varepsilon,\,V\big)$ such that for all $i\in\{1,\,2,\dots,\,n-1\}$, we have:
\smallskip

\centrer[2]{$\rho\rond\varphi(\tau_i) = T_{a_i}^\varepsilon\,V$.}
\medskip

\noindent \emph{Step 2.} The existence of the $(n-1)$-chain of curves $(a_i)_{i\leqslant n-1}$ implies that $n\leqslant m$,
whence Theorem \ref{thm:homomorphisms from Bn_in_Bm}.(i).
\medskip

\noindent \emph{Step 3.} We deduce from (1) and (2) that there exists a mapping class $F$ lying
in the normalizer of $\rho(\B_n)$ in $\Mod(\Sigma,\,\bord\Sigma)$ such that $F(c_i)=a_i$.
\medskip

\noindent \emph{Step 4.} This is the hard part: we have to compute the normalizer of $\rho(\B_n)$ in $\Mod(\Sigma,\,\bord\Sigma)$ and
show that $F$ can be chosen in $\rho(\B_n)$. Then for all $i\leqslant n-1$, we have:
\smallskip

\centrer[3]{$T_{a_i}=\rho(\gamma\tau_i\gamma^{-1})$.}
\medskip

\noindent \emph{Step 5.} It remains to see that $V$ lies in $\rho(\B_n)$. Let $v\in\B_n$ such that $V=\rho(v)$.
We then have:
\smallskip

\centrer[4]{$\rho\rond\varphi(\tau_i)=\rho(\gamma \tau_i^{\,\varepsilon}
\gamma^{-1}\,v)$.}
\smallskip

\noindent But $\rho$ is injective, so (4) implies:
\smallskip

\centrer[5]{$\varphi(\tau_i)=\gamma \tau_i^{\,\varepsilon}
\gamma^{-1}\,v$.}
\smallskip

\noindent Since $V$ lies in the centralizer of $\Im(\rho\rond\varphi)$ and since $\rho$ is injective,
$v$ must belong to the centralizer of $\langle \gamma\tau_i\gamma^{-1},\ i\leqslant n-1\rangle$,
so $\gamma^{-1}v\gamma$ belongs to the centralizer of $\B_n$ in $\B_m$, whence the expression of $v$.
\bigskip


\subsection{Homomorphisms between mapping class groups}
        \label{subsec:cyclic_monodromy_transvection}
\medskip

We introduce $\Mod^\pm(\Sigma)$: the group of isotopy classes of
\emph{all} diffeomorphisms of $\Sigma$ (they may inverse $\Sigma$'s
orientation).
\medskip

\TITRE{Preliminaries.}

\begin{itemize}
\item[\point]\TITRE{Homomorphisms from $\Mod(\Sigma,\,\bord\Sigma)$ to $\Mod(\Sigma',\,\bord\Sigma')$
induced by an embedding.}
Let $\Sigma$ and $\Sigma'$ be two connected oriented
surfaces. Let $F$ be the isotopy class of a possibly
non-orientation-preserving embedding from $\Sigma$ in $\Sigma'$. Then $F$ induces a homomorphism
$\Ad_F$ from $\Mod(\Sigma,\,\bord\Sigma)$ to $\Mod(\Sigma',\,\bord\Sigma')$ as follows.
Let $A\in\Mod(\Sigma,\,\bord\Sigma)$ and let $\bar A$ be a diffeomorphism of $\Sigma$ that represents $A$.
Let $\Sigma''=F(\Sigma)$ and let $\bar F$ be a diffeomorphism from $\Sigma$ into $\Sigma''$ that represents $F$.
Then the isotopy class of $\bar F\bar A\bar F^{-1}$ belongs to
$\Mod(\Sigma'',\,\bord\Sigma'')$ and induces canonically
a mapping class in $\Mod(\Sigma',\,\bord\Sigma')$ which we denote by
$\widetilde\Ad_F(A)$. Such a homomorphism $\widetilde\Ad_F$ will be called
\emph{the homomorphism from $\Mod(\Sigma,\,\bord\Sigma)$ to $\Mod(\Sigma',\,\bord\Sigma')$ induced by the embedding $F$}.
\smallskip

\item[\point]\TITRE{Outer conjugations.} In the above description,
when $\Sigma''=\Sigma'$, we can identify $\Sigma'$ and $\Sigma$ so
that the embedding $F$ becomes an element of
$\Mod^\pm(\Sigma)$. The homomorphism  $\widetilde\Ad_F$ that we get
is then an automorphism of $\Mod(\Sigma,\,\bord\Sigma)$. In this
case, $\widetilde\Ad_F$ will be called an \emph{outer conjugation by
$F$}.\label{page:outer_conjugation}
\smallskip

\item[\point]\TITRE{Hyper-elliptic Involution.}
Given a surface $\Sigma=\Sigma_{g,\,b}$ with $b\in\{0,\,1,\,2\}$ and
a $m$-chain of curves $(a_i)_{i\leqslant m}$ in $\Sigma$ where $m=2g+1$ if $b=1$
and $m=2g+2$ if $b\in\{0,\,2\}$. The \emph{hyper-elliptic involution of $\Sigma$ relative
to $(a_i)_{i\leqslant m}$} is the unique involution of $\Mod(\Sigma)$ which preserves each curve $a_i$, $i\leqslant m$.
It can be represented by an angle $\pi$ rotation over an axis $\delta$ cutting the surface in $m+1$ points.
If $\Sigma=\Sigma_{2,\,0}$, all the hyper-elliptic involutions coincide and will be denoted by $H$; in this case,
$H$ spans the center of $\Mod(\Sigma_{2,\,0})$. See Figure \ref{fig-defi:alphaH}.

\begin{figure}[!h]
 \Includegraphics{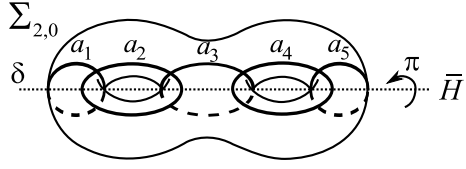}
 \caption{The mapping class $H$ of $\Mod(\Sigma_{2,\,0})$, represented by the rotation $\bar H$.}
 \label{fig-defi:alphaH}
\end{figure}
\smallskip

\item[\point]\TITRE{Cyclic homomorphisms from $\Mod(\Sigma_{2,\,b})$ in any given group.}
A homomorphism from $\Mod(\Sigma_{2,\,b})$ in any given group, with
$b\geqslant 0$, is said to be \emph{cyclic} if its image is cyclic.
\smallskip

\item[\point]\TITRE{Transvection of a homomorphism from the mapping class group in any group}
\mbox{}\\ Let $\Sigma$ be a genus-2 surface, $\M$ one of the mapping
class groups $\PMod(\Sigma)$ or $\Mod(\Sigma,\,\bord\Sigma)$, and
$G$ any group. For any homomorphism $\Psi$ from $\M$ to $G$ and for any
element $g$ belonging to the centralizer of $\Psi(\M)$ in $G$ such
that $g^{10}=1_G$, we will call \emph{transvection of $\Psi$ with
direction $g$} the homomorphism $\Psi'$ that associates $\Psi(T_a)\,g$
to any Dehn twist $T_a$ along a non-separating curve $a$.
\end{itemize}
\bigskip

\TITRE{Previous and new theorems.}
\medskip

So far, the
main result about homomorphisms between mapping class groups was given by N.V. Ivanov and J. McCarthy
in 1999: \smallskip

\begin{thm}[Ivanov, McCarthy, \mbox{[IvMc]}, 1999]
            \label{thm:IvMc}
\mbox{}\\ Let $\Sigma$ be a surface $\Sigma_{g,\,b}$ and $\Sigma'$
be a surface $\Sigma_{g',\,b'}$ with $g\geqslant 2$ and
$(g',\,b')\not=(2,\,0)$, and such that the inequality
\smallskip

\centrer{$|(3g-3+b) - (3g'-3+b')|\leqslant 1$}
\smallskip

\noindent holds. If there exists an injective homomorphism $\rho$ from
$\Mod(\Sigma)$ to $\Mod(\Sigma')$, then $\Sigma'$ is homeomorphic to
$\Sigma$ and $\rho$ is an automorphism induced by a possibly not
orientation-preserving diffeomorphism from $\Sigma$ in $\Sigma'$.
\end{thm}
\smallskip

They have also completed this theorem by dealing with some cases
when $(g',\,b')=(2,\,0)$ or when $g=1$. When $\Sigma'=\Sigma$, this
theorem tells us that the mapping class group is co-Hopfian, that
is, any injective endomorphism of $\Mod(\Sigma)$ is an automorphism.
The computation of $\Out(\Mod(\Sigma))$ also follows from this
theorem.

We turn now to the statements of our theorems. We first deal with
the mapping class group relatively to the boundary. Theorem
\ref{thm:homomorphisms_between_MCGs_complément} is an existence theorem
about nontrivial homomorphisms between mapping class groups. It is
completed by Theorem \ref{thm:homomorphisms_between_MCGs} which
provides a description of these homomorphisms.
\bigskip

\begin{thm*}[Existence of nontrivial homomorphisms from $\Mod(\Sigma,\,\bord\Sigma)$ to $\Mod(\Sigma',\,\bord\Sigma')$]
            \label{thm:homomorphisms_between_MCGs_complément}
\mbox{}\\
Let $\Sigma$ be a surface $\Sigma_{g,\,b}$ and $\Sigma'$ a surface $\Sigma_{g',\,b'}$ with $g\geqslant 2$ and $g'\leqslant
g+1$.
\begin{itemize}
\item When $g=2$, there exist some cyclic nontrivial homomorphisms from $\Mod(\Sigma_{2,\,b},\,\bord\Sigma_{2,\,b})$ in
any mapping class group admitting a subgroup isomorphic to $\ZZ[2]$,
$\ZZ[5]$ or $\ZZ[10]$. When $g\geqslant 3$, there does not exist any
cyclic nontrivial homomorphism.

\item
When $g\geqslant 2$, there exist some noncyclic homomorphisms from
$\Mod(\Sigma,\,\bord\Sigma)$ to $\Mod(\Sigma',\,\bord\Sigma')$ if
and only if one of the two following conditions is satisfied:
$b\not=0$ and $g'\geqslant g$, or $b=0$ and $\Sigma'$ is
homeomorphic to $\Sigma$.
\end{itemize}
\end{thm*}
\smallskip

\begin{thm*}[Homomorphisms from
$\Mod(\Sigma,\,\bord\Sigma)$ to $\Mod(\Sigma',\,\bord\Sigma')$]
            \label{thm:homomorphisms_between_MCGs}
\mbox{}\\
Let $\Sigma$ be a surface $\Sigma_{g,\,b}$ and $\Sigma'$ a surface
$\Sigma_{g',\,b'}$ with $g\geqslant 2$ and $g'\in\{g,\, g+1\}$, and
such that $\Sigma'=\Sigma$ if $b=0$. Any noncyclic homomorphism from
$\Mod(\Sigma,\,\bord\Sigma)$ to $\Mod(\Sigma',\,\bord\Sigma')$ is a
homomorphism induced by the isotopy class of an embedding from $\Sigma$
in $\Sigma'$, or possibly a transvection with direction $H$ (the
hyper-elliptic involution of $\Mod(\Sigma_{2,\,0})$) of such a
homomorphism if $g=2$ and $(g',\,b')=(2,\,0)$. Moreover, if $b=0$, the
homomorphism induced by the isotopy class of an embedding from $\Sigma$
in $\Sigma'$ (up to transvection when
$\Sigma=\Sigma'=\Sigma_{2,\,0}$) is an outer conjugation.
\end{thm*}
\bigskip

Let us now focus on the homomorphisms from $\PMod(\Sigma)$ to $\PMod(\Sigma')$. In most of the cases, they also can simply be
expressed from \emph{homomorphisms induced by an embedding}.
Given an isotopy class of an embedding $F$ from $\Sigma$ in $\Sigma'$ such that
$F$ sends the boundary components of $\Sigma$ on some boundary
components of $\Sigma'$ or on some trivial curves of $\Sigma'$
(\emph{trivial} means here isotopic to a point), we can define a homomorphism $\Ad_F$
from $\PMod(\Sigma)$ to $\PMod(\Sigma')$ exactly like we defined the homomorphism $\widehat\Ad_F$
from $\Mod(\Sigma,\bord\Sigma)$ to $\Mod(\Sigma',\bord\Sigma')$. Such a homomorphism $\Ad_F$ will be called
\emph{the homomorphism from $\PMod(\Sigma)$ to $\PMod(\Sigma')$ induced by the embedding $F$}. When $\Sigma'=\Sigma$,
$\Ad_F$ is an automorphism of $\PMod(\Sigma)$
that we will call the \emph{outer conjugation by $F$}.
The analogues of Theorems \ref{thm:homomorphisms_between_MCGs_complément} and
\ref{thm:homomorphisms_between_MCGs} when each boundary component is fixed only
setwise are Theorems
\ref{thm:homomorphisms_between_PMCGs_complément} and
\ref{thm:homomorphisms_between_PMCGs}.
\bigskip

\begin{thm*}[Existence of noncyclic homomorphisms from $\PMod(\Sigma)$ to $\PMod(\Sigma')$]
            \label{thm:homomorphisms_between_PMCGs_complément}
\mbox{}\\ Let $\Sigma$ be a surface $\Sigma_{g,\,b}$ and $\Sigma'$ a
surface $\Sigma_{g',\,b'}$ with $g\geqslant 2$ and $g'\leqslant
g+1$.
\begin{itemize}
\item When $g=2$ and only in this case,
there exist some cyclic nontrivial homomorphisms from
$\PMod(\Sigma_{2,\,b})$ in any mapping class group that admits a
subgroup isomorphic to $\ZZ[2]$, $\ZZ[5]$ or $\ZZ[10]$.
\item
When $g\geqslant 2$, there exist some noncyclic homomorphisms from
$\PMod(\Sigma)$ to $\PMod(\Sigma')$ if and only if $g'= g$ and
$b'\leqslant b$.
\end{itemize}
\end{thm*}
\medskip

According to Theorem \ref{thm:homomorphisms_between_PMCGs_complément},
studying all noncyclic homomorphisms from
$\PMod(\Sigma)$ to $\PMod(\Sigma')$ when $g'\leqslant g+1$
comes down to studying them when $g'= g$ and $b'\leqslant b$. This is
the aim of Theorem \ref{thm:homomorphisms_between_PMCGs} below.
\medskip

\begin{thm*}[Homomorphisms from
$\PMod(\Sigma)$ to $\PMod(\Sigma')$]
            \label{thm:homomorphisms_between_PMCGs}
\mbox{}\\ Let $\Sigma$ be a surface $\Sigma_{g,\,b}$ and $\Sigma'$ a
surface $\Sigma_{g',\,b'}$ with $g\geqslant 2$, $g'=g$ and
$b'\leqslant b$. Let $\Psi$ be a noncyclic homomorphism from
$\PMod(\Sigma)$ to $\PMod(\Sigma')$. Then there exists an embedding
$F$ from $\Sigma$ in $\Sigma'$ such that $F$ sends the boundary
components of $\Sigma$ on some boundary components of $\Sigma'$ or
on some trivial curves of $\Sigma'$, and such that $\Psi$ is the
homomorphism $\Ad_F$ induced by the embedding $F$, or possibly the
transvection by the hyper-elliptic involution $H$
(see the above Figure \ref{fig-defi:alphaH}) of the homomorphism
$\Ad_F$ if $g=2$ and $(g',\,b')=(2,\,0)$.
\end{thm*}
\bigskip

Among the homomorphisms between mapping class groups provided by
Theorems \ref{thm:homomorphisms_between_MCGs} and
\ref{thm:homomorphisms_between_PMCGs}, let us determine the injective
ones.
\bigskip

\begin{thm*}[Injections from
$\Mod(\Sigma,\,\bord\Sigma)$ into $\Mod(\Sigma',\,\bord\Sigma')$]
            \label{thm:injectivity_des_MCGs}
\mbox{}\\ Let $\Sigma$ be a surface $\Sigma_{g,\,b}$ and $\Sigma'$ a
surface $\Sigma_{g',\,b'}$ with $g\geqslant 2$ and $g'\leqslant
g+1$. Then, a homomorphism from $\Mod(\Sigma,\,\bord\Sigma)$ to $\Mod(\Sigma',\,\bord\Sigma')$ is injective if and only if:
\begin{itemize}
\item[\point] when $b\not=0$: if it is induced, up to transvection when $g=2$,
by an embedding $F$ of $\Sigma$ in $\Sigma'$ such that $F$ sends the
boundary components of $\Sigma$ on pairwise distinct curves in
$\Sigma'$;
\item[\point] when $b=0$ and $\Sigma'=\Sigma$: if it is not cyclic
(it is then an outer conjugation, or possibly a transvection of an
outer conjugation when $g=2$).
\end{itemize}
\end{thm*}
\bigskip

\begin{thm*}[Injections of $\PMod(\Sigma)$ into $\PMod(\Sigma')$]
            \label{thm:injectivity_des_PMCGs}
\mbox{}\\ Let $\Sigma$ be a surface $\Sigma_{g,\,b}$ and $\Sigma'$ a
surface $\Sigma_{g',\,b'}$ with $g\geqslant 2$ and $g'\leqslant
g+1$. Then, a homomorphism from $\PMod(\Sigma)$ to $\PMod(\Sigma')$ is
injective if and only if the two following conditions hold:
\begin{itemize}
\item[\;(i)] the surfaces $\Sigma$ and $\Sigma'$ are homeomorphic,
\item[(ii)] the homomorphism is an outer conjugation (cf definition page \pageref{page:outer_conjugation}), or possibly the transvection with
direction $H$ of an outer conjugation when $\Sigma'$ and $\Sigma$
are homeomorphic to $\Sigma_{2,\,0}$.
\end{itemize}
\end{thm*}
\bigskip

\TITRE{Comparison with Ivanov and McCarthy's Theorem (cf. Theorem \ref{thm:IvMc}), 1999.}
\medskip

The proof of Ivanov and McCarthy's Theorem is based on an algebraic characterization
of the Dehn twists, which is possible only if the maxima of the
ranks of the abelian subgroups of $\Mod(\Sigma)$ and $\Mod(\Sigma')$
differ from at most one. The proof then also requires that the
considered homomorphisms are rank-preserving, hence the considered
homomorphisms have to be injective.
\smallskip

In this paper, instead of using the rank of abelian subgroups
embedded in the mapping class group, we have used homomorphisms from the braid group
to the mapping class group. Since these two groups are quite near
(see for instance the similarities between braid twists and Dehn twists as generating sets)
we get precise results without any injectivity hypothesis.
Let us compare the results of Ivanov and McCarthy with ours. \bigskip

\noindent \emph{Results of Ivanov and McCarthy (1999) that are not
covered in this paper:}
\begin{itemize}
\item For any nonnegative integer $m$ and for any $\varepsilon$ in
$\{0,\,1\}$, there does not exist any injective homomorphism from
$\Mod(\Sigma_{g,\,b+3m})$ to $\Mod(\Sigma_{g+m,\,b+\varepsilon})$,
where $g\geqslant 2$ and $b\geqslant 0$. It is noticeable that the
hypotheses allow the genus of the surface at the target to be
arbitrary large with respect to the genus of the surface at the
source!
\item The elements of the considered mapping class groups can
permute the boundary components.
\end{itemize}
\smallskip

\noindent \emph{Our results (2010) that are not covered by Ivanov
and McCarthy'paper:}
\begin{itemize}
\item[\point] Full description of the homomorphisms from $\Mod(\Sigma_{g,\,b},\,\bord\Sigma_{g,\,b})$
to $\Mod(\Sigma_{g',\,b'},\,\bord\Sigma_{g',\,b'})$ where $g'<g$ and
$g\geqslant 2$, whatever $b$ and $b'$ are. Precisely, all these
homomorphisms are trivial or cyclic.
\item[\point] Full description of the homomorphisms from $\Mod(\Sigma_{g,\,b},\,\bord\Sigma_{g,\,b})$
to $\Mod(\Sigma_{g',\,b'},\,\bord\Sigma_{g',\,b'})$ where $g'=g$ or
$g'=g+1$, and $g\geqslant 2$, whatever $b$ and $b'$ are. In these
cases, there exist noncyclic homomorphisms, and only some of them are
injective.
\item[\point] We also prove these results in a slightly different frame: when the elements
of mapping class group preserve each boundary component setwise
instead of pointwise.
\end{itemize}
\bigskip

More precisely, in this paper, we focus on homomorphisms between two
mapping class groups associated  to the surfaces $\Sigma$ and
$\Sigma'$ with genera $g$ and $g'$ such that $g\geqslant 2$ and
$g'\leqslant g+1$, and whatever their numbers of boundary components
are. We shall thus describe the following sets:
\begin{itemize}
\item all the homomorphisms from $\Mod(\Sigma,\,\bord\Sigma)$ to $\Mod(\Sigma',\,\bord\Sigma')$ (cf. Theorems
\ref{thm:homomorphisms_between_MCGs_complément} and
\ref{thm:homomorphisms_between_MCGs}),
\item all the homomorphisms from $\PMod(\Sigma)$ to $\PMod(\Sigma')$ (cf. Theorems
\ref{thm:homomorphisms_between_PMCGs_complément} and
\ref{thm:homomorphisms_between_PMCGs}),
\item all the injective homomorphisms from $\Mod(\Sigma,\,\bord\Sigma)$
to $\Mod(\Sigma',\,\bord\Sigma')$ (cf. Theorem
\ref{thm:injectivity_des_MCGs}),
\item all the injective homomorphisms from $\PMod(\Sigma)$ to $\PMod(\Sigma')$ (cf. Theorem \ref{thm:injectivity_des_PMCGs}).
\end{itemize}
\noindent Like Ivanov and McCarthy, we will mainly show that the
nontrivial homomorphisms between mapping class groups are \emph{induced
by some embeddings} (the result is however slightly different when
the genus of the surface $\Sigma$ equals 2).
\bigskip

\TITRE{Steps of the proof of Theorems \ref{thm:homomorphisms_between_MCGs_complément} and
\ref{thm:homomorphisms_between_MCGs}.}

\noindent \emph{Step 1.} Let $\Sigma=\Sigma_{g,\,b}$ and $\Sigma'\Sigma_{g',\,b'}$ be two surfaces
with $g,g'\geqslant 2$ and $g'\leqslant g+1$, and let $\Phi$ be a homomorphism from
$\Mod(\Sigma,\bord\Sigma)$ to $\Mod(\Sigma,\bord\Sigma)$.
Let $n=2g+2$, let $\rho_0$ be a monodromy homomorphism from $\B_n$ to
$\Mod(\Sigma,\bord\Sigma)$.
Thus there is a $(2g+1)$-chain of curves $(a_i)_{i\leqslant 2g+1}$ in $\Sigma_{g,\,b}$
such that for all $i\in\{1,\,2,\dots,\,2g+1\}$, we have:
\smallskip

\centrer[1]{$\rho_0(\tau_i) = T_{a_i}$.}
\smallskip

\noindent Now, let $\rho_1$ be the composition $\Phi\rond\rho_0$
(see the following commutative diagram).
\smallskip

\centrer{$\begin{diagram}
\node{\Mod(\Sigma,\,\bord\Sigma)}\arrow[2]{e,t}{\Phi}
\node[2]{\Mod(\Sigma,\,\bord\Sigma)}\\
\node[2]{\B_{n}}\arrow{nw,b}{\rho_0}\arrow{ne,b}{\rho_1}
\end{diagram}$}
\medskip

\noindent Then, according to Theorem \ref{thm:Theorem_principal_case_à_boundary},
$\rho_1$ is either cyclic or is a transvection of monodromy homomorphism.
\smallskip

(i) Suppose first that $\rho_1$ is cyclic. Then $\Phi(T_{a_1})=\Phi(T_{a_2})=\dots$
We can show that the images by $\Phi$ of a generating set of Dehn twists of $\Mod(\Sigma,\,\bord\Sigma)$
agree, so $\Phi$ is cyclic. So $\Phi$ sends the commutator
$[\Mod(\Sigma,\,\bord\Sigma),\,\Mod(\Sigma,\,\bord\Sigma)]$ on 1.
But after a theorem of Korkmaz (cf. [Ko2]), the abelianization of $\Mod(\Sigma,\,\bord\Sigma)$
is trivial when $g\geqslant 3$, so
$[\Mod(\Sigma,\,\bord\Sigma),\,\Mod(\Sigma,\,\bord\Sigma)]=\Mod(\Sigma,\,\bord\Sigma)$
and $\Phi$ is trivial.
\smallskip

(ii) Suppose now that $\rho_1$ is a transvection of monodromy homomorphism. It means that there exists a triple
$\big((a'_i)_{i\leqslant 2g+1},\,\varepsilon,\,V\big)$ such that for all $i\in\{1,\,2,\dots,\,2g+1\}$,
we have:
\smallskip

\centrer[2]{$\rho_1(\tau_i) = T_{a'_i}^{\,\varepsilon}\,V$.}
\smallskip

\noindent Notice that $g'$ has to be at least as great as $g$
to allow the existence of the $2g+1$ chain of curves
$(a'_i)_{i\leqslant 2g+1}$. This proves a large part of Theorem \ref{thm:homomorphisms_between_MCGs_complément}.
\bigskip

\noindent \emph{Step 2.} We assume $\Phi$ is not cyclic, we consider another $(2g+1)$-chain of curves
$(c_i)_{i\leqslant 2g+1}$ in $\Sigma$. Then, there exists a triple
$\big((c'_i)_{i\leqslant 2g+1},\,\kappa,\,W\big)$ such that for all $i\in\{1,\,2,\dots,\,2g+1\}$,
we have:
\smallskip

\centrer[3]{$\rho_1(\tau_i) = T_{c'_i}^{\,\kappa}\,W$.}
\smallskip

\noindent We can show that if a sub-3-chain of curves in $(c_i)_{i\leqslant 2g+1}$ coincides with
a sub-3-chain of curves in $(a_i)_{i\leqslant 2g+1}$, then $\kappa=\varepsilon$ and $V=W$ in (2) and (3).
By showing a result of connectedness related to the set of the $(2g+1)$-chain of curves in $\Sigma$,
we deduce that up to an orientation issue,
$\Phi$ sends each Dehn twist along a non-separating curve $a$ on
a Dehn twist along a non-separating curve composed by the same mapping class $V$ independent from $a$.
The orientation issue can be solved by replacing $\Phi$ by the composition of $\Phi$
by the ``outer-conjugation'' in $\Mod(\Sigma)$ by a mapping class $H$ lying in $\Mod^\pm(\Sigma)$.
\bigskip

\noindent \emph{Step 3.} Using the fact that there exist special relations with Dehn twists coming
from the fact that the abelianization of the mapping class group is trivial, we can deduce that $V$ is trivial.
Hence $\Phi$ is ``twist-preserving'', which means that is sends Dehn twists along non-separating curves
on Dehn twists along non-separating curves. Moreover, if $a$ and $c$ meet in 0 or 1 point,
then $\Phi(T_a)$ and $\Phi(T_c)$ are two Dehn twists along two curves that meet respectively in 0 or 1 point.
\bigskip

\noindent \emph{Step 4.} If $\Sigma$ has no boundary, using a minimal generating set of $\Mod(\Sigma)$
consisting in Dehn twists along non-separating curves, the proof would be almost over, but in the general case,
this is still far from being the case, since $\Phi$ is not supposed
to be injective. Indeed, different Dehn twist can be (and sometimes are) sent by $\Phi$ on the same Dehn twist.
Another big difficulty consists in showing that three Dehn twists cobounding a pair of pants
are sent to three Dehn twists cobounding a pair of pants. These problems are solved by translating algebraically
these topological issues, involving the existence of special extra curves.
\medskip

\noindent \emph{Step 5.} Let $\A$ be a set of curves in $\Sigma$ containing a $(2g+1)$-chain of curves and containing also
additional curves to get a pant decomposition of $\Sigma$. It follows from Step 4 that there exists an embedding $F$
from $\Sigma$ to $\Sigma'$ such that
for any curve $a$ belonging to $\A$, we have $F(a)= a'$ where $a'$ is such that $\Phi(T_a)=T_{a'}$.
We can now compare $\Phi$ with the outer conjugation by $F$ and see that they agree.
\bigskip


\subsection{Endomorphisms and automorphisms of the mapping class group}
\label{par:MCG_dans_lui-même}
\medskip

We complete the previous subsection by focusing on the injective
homomorphisms from \linebreak $\Mod(\Sigma,\,\bord\Sigma)$ to $\Mod(\Sigma',\,\bord\Sigma')$ when $\Sigma'=\Sigma$. In particular,
we obtain Ivanov and McCarthy's theorem stating that $\PMod(\Sigma)$
is co-Hopfian (cf. [IvMc]). Theorems \ref{thm:homomorphisms_between_MCGs_complément}
and \ref{thm:homomorphisms_between_MCGs} allow us to state an
equivalent theorem for $\Mod(\Sigma,\,\bord\Sigma)$ as well (cf. Theorem
\ref{thm:cohopfianité_des_MCGs} below). We will see that when $b=0$, the
group $\Mod(\Sigma,\,\bord\Sigma)$ satisfies a much stronger
property (cf. item (i) of Theorems \ref{thm:cohopfianité_des_MCGs}). As in the previous
subsection, and since the center of $\Mod(\Sigma_{2,\,0})$ is
nontrivial, the case of the surface $\Sigma_{2,\,0}$ is special.
As it is not central in this paper, we will not recall the remarkable structure of $\Aut(\Mod(2,\,0))$
which has been highlighted by McCarthy;
the interested reader is referred to [Mc1].
\medskip

\begin{thm*}[Co-Hopfian property
of $\Mod(\Sigma,\,\bord\Sigma)$ and structure of
$\Aut(\Mod(\Sigma,\,\bord\Sigma))$, where
$\Sigma\not=\Sigma_{2,\,0}$]
            \label{thm:cohopfianité_des_MCGs}
\mbox{}\\Let $\Sigma$ be a surface $\Sigma_{g,\,b}$ where
$g\geqslant 2$.
\begin{itemize}
  \item[\;\;(i)]
The mapping class group $\Mod(\Sigma,\,\bord\Sigma)$ is co-Hopfian,
that is, the injections are automorphisms.
  \item[\;(ii)]
Moreover, when $b=0$ and $g>2$ (resp. $b=0$ and $g=2$), all
the nontrivial (resp. noncyclic) homomorphisms from $\Mod(\Sigma)$ are automorphisms.
\end{itemize}
From now on, we exclude the case $(g,b)=(2,0)$ for which the following would be wrong.
\begin{itemize}
  \item[\;(iii)]
The map $\Ad\;:\;
\Mod^\pm(\Sigma)\to\Aut(\Mod(\Sigma,\,\bord\Sigma))$ is an
isomorphism.
  \item[(iv)]
The outer automorphism group $\Out(\Mod(\Sigma,\,\bord\Sigma))$ of
$\Mod(\Sigma,\,\bord\Sigma)$ is isomorphic to the direct product
$\ZZ[2]\times\Ss_b$, where $\Ss_b$ is the symmetric group on $b$
elements.
\end{itemize}
\end{thm*}
\bigskip


\subsection{Rigidity/triviality of the homomorphisms from $\Mod(\Sigma_{g,\,0})$ to $\Mod(\Sigma_{g',\,0})$}
\label{par:MCG_triviality}
\medskip

Here is a discussion about the exceptional rigidity of the homomorphisms from $\Mod(\Sigma_{g,\,0})$ to $\Mod(\Sigma_{g',\,0})$.
Let us begin by examining the case where $g=1$. As always, the case $g=1$ is special.
\medskip

\begin{figure}[!h]
      \Includegraphics{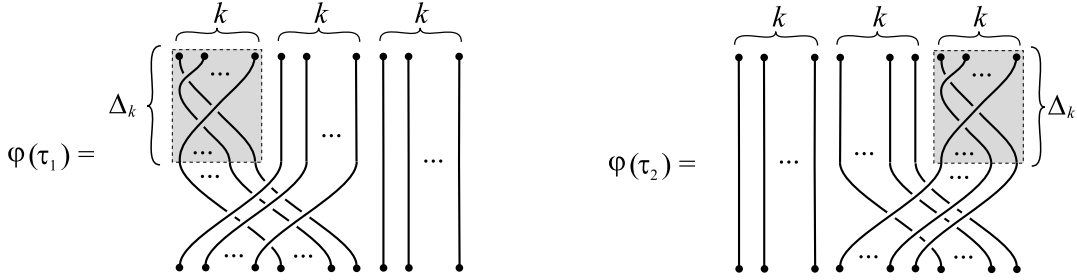}
      \caption{Definition of the homomorphism $\varphi$.}
      \label{fig:morphismeB3B3k}
     \end{figure}

\TITRE{Non-cyclicity of $\Mod(\Sigma_{1,\,0})\to\Mod(\Sigma_{g,\,0})$.}
We can construct non-cyclic homomorphisms from $\Mod(\Sigma_{1,\,0})$ to $\Mod(\Sigma_{g,\,0})$
when $g\not\equiv 0 \mod 3$. Here is an example.
For any $k\in\NN^*$, define an injection $\varphi$ from $\B_3$ to $\B_{3k}$ defined as in Figure \ref{fig:morphismeB3B3k}
where for all integers $i$, ${\Delta_i}$ denotes the Garside element of $\B_i$.
Check that with this definition of $\varphi$, we have:
\smallskip

\centrer{$\varphi({\Delta_3})={\Delta_{3k}}$}
\smallskip

\noindent (Hint: in the braid corresponding to $\varphi({\Delta_3})$,
each strand crosses all the other strands exactly once).
Now, let $\psi$ be the classical monodromy homomorphism
from $\B_{3k}$ to $\Mod(\Sigma_{g,\,b},\,\bord\Sigma_{g,\,b})$
where $(g,\,b)=\big(\frac{3k-1}{2},\,1\big)$ if $k$ is odd and
$(g,\,b)=\big(\frac{3k-2}{2},\,2\big)$ if $k$ is even.
We denote by $d$ or by $d_1$ and $d_2$ the boundary curve(s) of $\Sigma$,
depending on the parity of $k$. Then $\psi({\Delta_{3k}}^4)=T_{d}$ if $k$ is odd
and $\psi({\Delta_{3k}}^2)=T_{d_1}T_{d_2}$ if $k$ is even.
Gluing a disk along $d$, or gluing disks along $d_1$ and $d_2$ induces a canonical homomorphism
$\sq\ :\ \Mod(\Sigma_{g,\,b},\bord\Sigma_{g,\,b})\to\Mod(\Sigma_{g,\,0})$ whose kernel contains
the subgroup $\langle\psi({\Delta_{3k}}^4)\rangle$ of $\psi(\B_{3k})$.
Since $\varphi({\Delta_3})={\Delta_{3k}}$, it follows that ${\Delta_3}^4$ belongs to the kernel
of $\sq\rond\psi\rond\varphi$. But $\Mod(\Sigma_{1,\,0})$ is isomorphic to the quotient of
$\B_3$ by the relation ${\Delta_3}^4=1$, so the homomorphism $\sq\rond\psi\rond\varphi$ induces a well defined homomorphism
from $\Mod(\Sigma_{1,\,0})$ to $\Mod(\Sigma_{g,\,0})$ for any integer
$g\in\{\frac{3k-1}{2},\,k\in2\NN+1\}\cup\{\frac{3k-2}{2},\,k\in2\NN\}$, that is,
for any $g$ coprime with 3. Thus we have shown
the following proposition.
\smallskip

\begin{prop}[Non-cyclicity of $\Mod(\Sigma_{1,\,0})\to\Mod(\Sigma_{g,\,0})$ for all $g$ coprime with 3]
            \label{prop:non-cyclicity}
\mbox{}\\
For any $g$ coprime with 3, there exist homomorphisms from $\Mod(\Sigma_{1,\,0})$ to $\Mod(\Sigma_{g,\,0})$
that are non-cyclic.\fin
\end{prop}
\medskip


However, as far as we know, the existence of non-cyclic (and even nontrivial) homomorphisms from
$\Mod(\Sigma_{g,\,0})$ to $\Mod(\Sigma_{g',\,0})$ seems to be exceptional as the following results
tend to let us think.
\bigskip

\TITRE{Triviality of $\Mod(\Sigma_{g,\,0})\to\Mod(\Sigma_{g',\,0})$ when $g\geqslant 3$ and $g'<g$.}
This is a theorem due to Harvey and Korkmaz [HvKo] which is also included in our Theorem \ref{thm:homomorphisms_between_MCGs}.
\bigskip

Remember that according to Theorem \ref{thm:cohopfianité_des_MCGs}, there are very few homomorphisms from
$\Mod(\Sigma_{g,\,0})\to\Mod(\Sigma_{g',\,0})$ even when $g=g'$:
\medskip

\TITRE{Rigidity of $\End(\Mod(\Sigma_{g,\,0}))$ when $g\geqslant 2$.}
All non-cyclic (and even only nontrivial if $g\geqslant3$) endomorphisms of
$\Mod(\Sigma_{g,\,0})$ are automorphisms.
\bigskip

\TITRE{Triviality of $\Mod(\Sigma_{g,\,0})\to\Mod(\Sigma_{g',\,0})$ when $g\geqslant 3$ and $g'>g$.}
This is an open question outlined in [BkMt] where J. Berrick and M.Mathhey conjecture that for any integer $m>1$,
then for all sufficiently large $g$, the only homomorphism from $\Mod(\Sigma_{g,\,0})$ to $\Mod(\Sigma_{g+m,\,0})$ is
the trivial homomorphism (see Conjecture 4.5). Their conjecture is inspired by a theorem they shew (see Theorem 4.4),
asserting that for any integer $m>1$, there are infinitely many values of $g$ such that the only homomorphism from
$\Mod(\Sigma_{g,\,0})$ to $\Mod(\Sigma_{g+m,\,0})$ is the trivial homomorphism. Our Theorem \ref{thm:homomorphisms_between_MCGs}
brings the answer (``yes'') to their conjecture when $m=1$:
\bigskip

\begin{thm*}{Triviality of $\Mod(\Sigma_{g,\,0})\to\Mod(\Sigma_{g+1,\,0})$ for any $g\geqslant 2$}
\mbox{}\\
For any integer $g\geqslant 2$, the only homomorphism from $\Mod(\Sigma_{g,\,0})$ to $\Mod(\Sigma_{g+1,\,0})$
is the trivial homomorphism.
\fin
\end{thm*}
\bigskip

\REM Notice that according to Proposition \ref{prop:non-cyclicity}, not only do there exist nontrivial homomorphisms from
$\Mod(\Sigma_{1,\,0})$ to $\Mod(\Sigma_{2,\,0})$, but even non-cyclic ones.
\bigskip


\subsection{Geometric representations of the Artin groups of type $D_n$, $E_6$, $E_7$ and $E_8$}
\medskip

Let $S$ be a finite set. A \emph{Coxeter matrix over $S$} is a matrix
$M = \big(m_{\alpha\beta}\big)_{\alpha,\beta\in S}$ indexed by the elements
of $S$ such that $m_{\alpha\alpha} = 1$ for all $\alpha\in S$,
and $m_{\alpha \beta} = m_{\beta \alpha}\in \{2, 3, 4, . . . ,+\infty\}$
for all $\alpha, \beta \in S$, $\alpha \not= \beta$.
A Coxeter matrix 
is usually represented by its \emph{Coxeter graph $\Gamma$}.
This (labeled) graph is defined by the following data.
The vertices of $\Gamma$ are the elements of $S$.
Two vertices $\alpha, \beta$ are joined by an edge if $m_{\alpha \beta} = 3$,
and this edge is labeled by $m_{\alpha \beta}$ if $m_{\alpha \beta} \geqslant 4$.
For $\alpha, \beta \in S$ and $m \in \ZZ_{\geqslant2}$ we denote by
$w(\alpha, \beta\,:\,m)$ the word $\alpha\beta\alpha\dots$ of length $m$.
Define the \emph{Artin group of type $\Gamma$} to be the (abstract) group
$A(\Gamma)$ presented by:
\smallskip

\centrer{$A(\Gamma) = \big\langle\, S \,\,\big|\,\,
w(\alpha, \beta : m_{\alpha \beta}) = w(\beta, \alpha,:\,m_{\alpha \beta})
\mbox{ for  } \alpha \not= \beta \mbox{ and  }m_{\alpha \beta} < +\infty\,\big\rangle$.}
\smallskip

\noindent If the graph $\Gamma$ is connected, then $A(\Gamma)$ is said to be \emph{irreducible}
(it cannot be written as a nontrivial direct product).
A \emph{small-type Artin group} is an Artin group
for which all the entries of the associated Coxeter matrix belongs in $\{2,\,3\}$.
Define the \emph{Coxeter group of type $\Gamma$} to be the (abstract) group
$W(\Gamma)$ presented by:
\smallskip

\centrer{$W(\Gamma) = \big\langle\, S \,\,\big|\,\,
\left\{\begin{array}{l}\alpha^2=1\mbox{ for  all }\alpha\in S\\
w(\alpha, \beta : m_{\alpha \beta}) = w(\beta, \alpha,:\,m_{\alpha \beta})
\mbox{ for  } \alpha \not= \beta \mbox{ and  }m_{\alpha \beta} < +\infty\end{array}\right.
\,\big\rangle$.}
\smallskip

\begin{figure}[!h]
  \Includegraphics{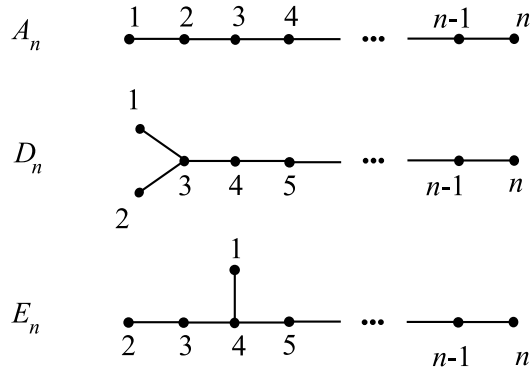}
  \caption{Coxeter Graphs of type $A_{n}$, $D_n$ and $E_n$.}
  \label{fig:grapheCoxeterADE}
\end{figure}

\noindent An Artin group is said to be \emph{of spherical type} if the associated Coxeter group in finite.
In this paper, we are interested by irreducible Artin groups of small spherical type.
After a work of Coxeter (see below),
these consist into two infinite families of Artin groups defined by Coxeter graphs  $A_n$ and $D_n$,
as well as three exceptional Artin groups defined by Coxeter graphs  $E_6$, $E_7$ and $E_8$, see Figure
\ref{fig:grapheCoxeterADE}.
\smallskip

\smallskip

We recognize the group $A(A_n)$ as being the braid group $\B_{n+1}$ with $n+1$ strands.
Artin groups, also called \emph{generalized braid groups}, were first introduced by Tits [Ti] as extensions of
Coxeter groups.
The finite irreducible Coxeter groups, and therefore the irreducible Artin
groups of spherical type, were classified by Coxeter [Cx]. Coxeter groups have been widely studied.
Basic references for them are [Bk] and [Hu].
In contrast, Artin groups are poorly understood in general.
\smallskip

The automorphism groups of the (spherical type) Artin groups are beginning to be explored.
Artin's 1947 paper [At3] was motivated by the problem of determining the automorphism groups of the
braid groups (it is explicit in the introduction). However, the problem itself was only solved 34 years
later by Dyer and Grossman [DyGr].
Until now, except for the braid groups (1981),
Artin groups of type $B_n$ (see [ChCr], 2005),
and Artin groups of rank 2 (see [CrPa2], Theorem 5.1, 2005),
the only known significant result on the automorphism groups
of spherical type Artin groups has been an extension of Artin's results of [At3] to all irreducible Artin groups
of spherical type (see [CoPa], 2003), as well as the computation of the group of the
automorphisms of $A(D_n)$ which leave invariant some normal subgroup
(it was unknown whether these subgroups were characteristic)
(see [CrPa2], Theorem 4.9, 2005).
\bigskip

\begin{figure}[!h]
  \Includegraphics{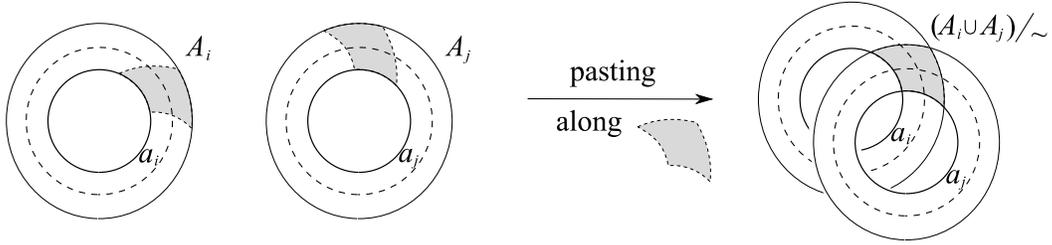}
  \caption{Pasting together two annuli along a portion of each.}
  \label{fig:surfaceCoxeter}
\end{figure}

\begin{figure}[!h]
  \Includegraphics{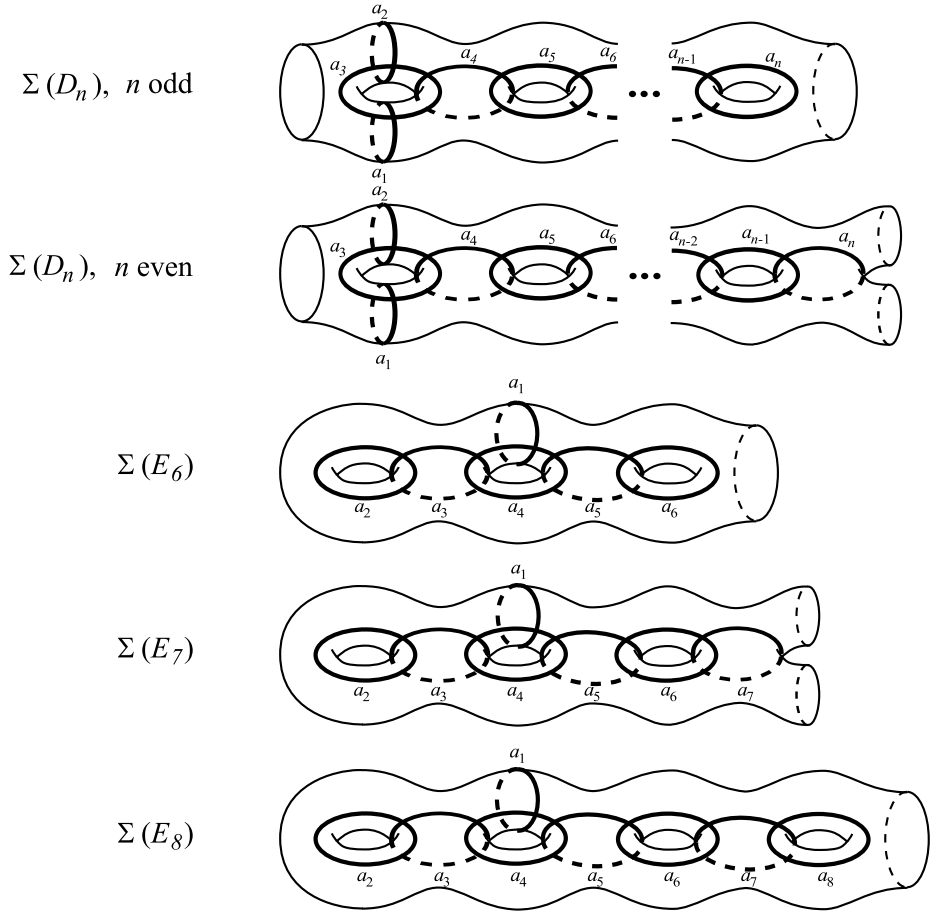}
  \caption{The surface $\Sigma(\Gamma)$ where $\Gamma\in\{\,D_n\, ;, n\geqslant 6\}\cup\{E_6,\,E_7,\,E_8\}$.}
  \label{fig:surfaceCoxeter2}
\end{figure}

\TITRE{Geometric representation of the Artin braid group of finite type.}
From now on, we assume that $A(\Gamma)$ is an small-type Artin group.
Let $n$ be the cardinality of $S$.
We denote by $\tau_1$, $\tau_2$, $\dots$, $\tau_n$ the elements of $S$ in accordance with
the labeled vertices in Figure \ref{fig:grapheCoxeterADE}
and we call them the \emph{standard generators}.
Let us consider a collection of annuli $(A_i)_{i\leqslant n}$.
For each $i\leqslant n$, let $a_i$ be a non contractible curve in $A_i\smallsetminus\bord A_i$.
For every couple $(\tau_i,\,\tau_j)\in S^2$ such that $m_{\tau_i\tau_j}=3$, we paste together a portion of $A_i$
and a portion of $A_j$ as illustrated on Figure \ref{fig:surfaceCoxeter}.
We thus obtain a surface that we denote by $\Sigma(\Gamma)$,
together with a special set of curves $(a_i)_{i\leqslant n}$, see Figure \ref{fig:surfaceCoxeter2}.
Let $n$ be the cardinality of $S$.
We denote the elements of $S$ by $\tau_1$, $\tau_2$, $\dots$, $\tau_n$
where the subscripts are coherent with the graphs in Figure \ref{fig:grapheCoxeterADE}.
Then, we have a reference homomorphism $\rho_\Gamma$ from $A(\Gamma)$
to $\Mod(\Sigma(\Gamma),\,\bord\Sigma(\Gamma))$ given by:
\smallskip

\centrer{$\DEF{\rho_\Gamma}{A(\Gamma)}{\Mod(\Sigma(\Gamma),\,\bord\Sigma(\Gamma))}{\tau_i}{T_{a_i}}$}
\smallskip

\noindent for all $i\leqslant n$.
\bigskip

\noindent

A homomorphism from a small-type Artin group
in any group will be said to be \emph{cyclic} if its image is cyclic.
Equivalently (thanks to the connectedness of the Coxeter graph), a homomorphism is said to be
cyclic if every standard generator is sent on the same element.
A \emph{monodromy homomorphism} from $A(\Gamma)$ in the mapping class group of $\Sigma(\Gamma)$ is a homomorphism
which sends any two distinct standard generators
on two Dehn twists along distinct curves.
A homomorphism from $A(D_n)$ in the mapping class group will be said to be a \emph{degenerate monodromy homomorphism}
if $\tau_1$ and $\tau_2$ are sent on the same Dehn twist,
all other standard generators being sent on Dehn twists on distinct curves.
\emph{Transvections} of such homomorphisms are defined as in the case of the braid group.
\bigskip

The next theorem is redundant with Theorem \ref{thm:Theorem_principal_case_à_boundary} in case where $\Gamma=A_{n-1}$,
and is a direct consequence of Theorem \ref{thm:Theorem_principal_case_à_boundary} when $\Gamma\in\{D_n,\,E_n\}$,
since according to a theorem of Van der Leck (cf. [VL] and [Pa1]), $A(\Gamma)$
contains $A(A_{n-1})$ as a parabolic subgroup (that is, a subgroup spanned by some of the standard generators of $A(\Gamma)$).
\bigskip

\begin{thm*}[Geometric Representations of the irreducible Artin groups of small spherical type]\mbox{}\\
            \label{thm:geometricRepresentationSmalltypeArtinGroup}
Let $\Gamma\in\{A_{n-1},\,D_n\, ;, n\geqslant 6\}\cup\{E_6,\,E_7,\,E_8\}$.
\begin{itemize}
\item[\;(i)]
For any connected surface $\Sigma$ of genus $g$, homomorphisms from $A(\Gamma)$ to $\Mod(\Sigma,\,\bord\Sigma)$
are cyclic if the genus of $\Sigma$ is less than the genus of $\Sigma(\Gamma)$, that is, if $g<\lfloor\frac{n}{2}\rfloor$
for $\Gamma\in\{A_n, \,E_n\}$, or if $g<\lfloor\frac{n-1}{2}\rfloor$ for $\Gamma=D_n$.
\item[(ii)]
For any connected surface $\Sigma$ having the same genus as $\Sigma(\Gamma)$,
non-cyclic homomorphisms from $A(\Gamma)$ to $\Mod(\Sigma,\,\bord\Sigma)$ exist.
They are transvections of (possibly degenerate in the case of $D_n$) monodromy homomorphisms.\fin
\end{itemize}
\end{thm*}
\bigskip

Whereas monodromy representations of $A(\Gamma)$
are known to be faithful for $\Gamma\in\{A_n,\,D_n\}$ (see [BiHi] and [PeVa]),
they are not faithful for $\Gamma\in\{E_6,\,E_7,\, E_8\}$ (see [W]).
Hence Theorem \ref{thm:geometricRepresentationSmalltypeArtinGroup} implies the following.
\bigskip

\begin{thm*}{Geometric representations of $A(E_6)$, $A(E_7)$ and $A(E_8)$ are not faithful}
            \label{thm:monodromy_homomorphism_E6_E7_E8}
\mbox{}\\
Let $\Gamma\in\{E_6,\,E_7,\,E_8\}$.
There is no injective homomorphism from $A(\Gamma)$ to $\Mod(\Sigma(\Gamma),\,\bord\Sigma(\Gamma))$.\fin
\end{thm*}
\bigskip


\subsection{Endomorphisms and automorphisms of the Artin groups of type $D_n$}
\medskip

Remember that the fact the
monodromy representations of $\B_n$ are faithful helped us
to determine $\End(\B_n)$. In the same way, we can use
the fact the
monodromy representations of $A(D_n)$ are faithful to compute $\End(A(D_n))$.
\bigskip

\begin{thm*}[Endomorphisms and automorphisms of $A(D_n)$]
            \label{thm:automorphisms_A(Dn)}
\mbox{}\\
Let $n\geqslant 6$.
\begin{itemize}
\item[\;\;(i)]
 The non-cyclic endomorphisms of $A(D_n)$ are the transvections of (possibly degenerate) monodromy homomorphisms.
\item[\;(ii)]
 The injective endomorphisms of $A(D_n)$ are the transvections of monodromy homomorphisms.
\item[(iii)]
 The automorphisms of $A(D_n)$ are the monodromy homomorphisms.
\item[(iv)]
 $\Out(A(D_n))\cong\left\{\begin{array}{ll}\ZZ[2]&\mbox{ if $n$ is odd},\\
 \ZZ[2]\times\ZZ[2]&\mbox{ if $n$ is even}.\end{array}\right.$
\item[(v)]
The following short exact sequence splits,
\smallskip

\centrer{$1\to\Inn(A(D_n))\to \Aut(A(D_n))\to\Out(A(D_n))\to 1$}
\smallskip

and the image of a possible section of $\Out(A(D_n))$ to $\Aut(A(D_n))$
is the group spanned by $\Inv$ if $n$ is odd, or the group spanned by $\Inv$ and $\theta$
if $n$ is even, where $\Inv$ is the automorphism that sends $\tau_i$ to ${\tau_i}^{-1}$ for all $i\in\{1,\,2,\dots,\,n\}$,
and $\theta$ is the automorphism that swaps $\tau_1$ and $\tau_2$ and that fixes $\tau_i$ for all $i\in\{3,\,4,\dots,\,n\}$.
\end{itemize}
\end{thm*}
\bigskip

\TITRE{Steps of the proof.}
Let $\varphi$ be a non-cyclic endomorphism of $A(D_n)$. Let $\Sigma=\Sigma(D_n)$.
We start from the faithful representation $\rho_{D_n}\ :\ A(D_n)\to\Mod(\Sigma,\,\bord\Sigma)$
(see [PeVa] for the injectivity). Let $\rho=\rho_{D_n}\rond\varphi$ (see the following commutative diagram).
\smallskip

\centrer{$\begin{diagram}
  \node[2]{A(D_n)}\arrow{sw,l}{\varphi}\arrow{se,t}{\rho}\\
  \node[1]{A(D_n)}\arrow[2]{e,t,J}{\rho_{D_n}}
      \node[2]{\Mod(\Sigma,\,\bord\Sigma)}
 \end{diagram}$}
\smallskip

\noindent Let $\B_n$ the braid group included in $A(D_n)$ and spanned by $\tau_i$, with $2\leqslant i\leqslant n$.
According to Theorem \ref{thm:Theorem_principal_case_à_boundary},
the restriction of $\rho$ to $\B_n$ is a transvection of monodromy homomorphism.
\smallskip

\noindent \emph{Step 1.} The case where $\varphi(\tau_1)=\varphi(\tau_2)$ is easy to deal with.
We exclude it, so that $\rho(\tau_1)\not=\rho(\tau_2)$ (since $\rho_{D_n}$ is injective).
\smallskip

\noindent \emph{Step 2.} Using the injectivity of $\rho_{D_n}$, we show that $\varphi$
is a transvection of monodromy homomorphism. We show also that if $\varphi$ is surjective,
then the transvection part has to be trivial.
\smallskip

\noindent \emph{Step 3.} We characterize the elements of $\Mod(\Sigma,\,\bord\Sigma)$
that lie in $\rho_{D_n}(A(D_n))$. They are the ones that lie in $\rho_{A_{n-1}}(A(A_{n-1}))$
after having pasted a disk along the boundary circle of $\Sigma$
that cobounded a pair of pants with $a_1$ and $a_2$.
\smallskip

\noindent \emph{Step 4.} We compute the normalizer $\N_{\Mod(\Sigma,\,\bord\Sigma)}(\rho_{D_n}(A(D_n)))$ of
$\rho_{D_n}(A(D_n))$ in $\Mod(\Sigma,\,\bord\Sigma)$.
\smallskip

\noindent \emph{Step 5.} We finally show that $\rho_{D_n}$ induces the following isomorphs:
\smallskip

\mbox{}
\hfill
$\begin{array}{lll}
\Aut(A(D_n)) &\cong& \N_{\Mod^\pm(\Sigma)}\big(\,\rho_{D_n}(A(D_n))\,\big),\\
\Inn(A(D_n))&\cong& \N_{\PMod(\Sigma)}\big(\,\rho_{D_n}(A(D_n))\,\big).
\end{array}$
\hfill
$\begin{array}{r}
(1)\\
(2)
\end{array}$
\smallskip

\noindent \emph{Step 6.} We then easily compute $A(D_n)$ from step 5 and from the following remark.
\bigskip

\REM From the equality (1) we get a geometric interpretation about the difference between the odd and the even case
for the computation of $\Out(A(D_n))$. The automorphism $\theta$ corresponds in $\Aut(\Mod(\Sigma,\,\bord\Sigma))$
to an obvious inner automorphism when $n$ is odd,
and to an outer automorphism when $n$ is even, since it has to permute two boundary circles.
\bigskip

\REM We deduce from Theorem \ref{thm:automorphisms_A(Dn)} and from Theorem 4.9 in [CrPa2]
the answer to a question asked by J. Crisp and L. Paris in [CrPa2]:
the kernel of the map $\pi_D\ :\ A(D_n)\to\B_n$ sending $\tau_1$ and $\tau_2$ on $\tau_1$ and sending
$\tau_i$ on $\tau_{i-1}$ for all $i\geqslant 3$, which is a free group of rank $n-1$,
happens to be a characteristic subgroup of $A(D(n))$.
\bigskip

\pageimpaire\addcontentsline{toc}{part}{\protect\textsc{\textbf{Proof of Theorems \ref{thm:Theorem_principal_case_à_boundary} and \ref{thm:injectivity_des_homomorphisms}}}}


\part*{Proof of Theorems \ref{thm:Theorem_principal_case_à_boundary} and
\ref{thm:injectivity_des_homomorphisms}}
\bigskip

\section{On monodromy homomorphisms}
        \label{sec:Proof_Theorem_2_&_3}
\bigskip

\subsection{Cyclic homomorphisms, monodromy homomorphisms and transvections}
            \label{par:cyclic_homomorphisms_monodromy_homomorphisms_transvections}
\medskip

Definitions of Cyclic homomorphisms, monodromy homomorphisms and transvections
were given in Subsection \ref{subsec:cyclic_monodromy_transvection}.
Here are some basic facts related to each of them.
\bigskip

\begin{lem}[Criterion 1 for a homomorphism to be cyclic]
      \label{lem:homomorphism_cyclic}
\mbox{}\\
Let $n$ be an integer greater than or equal to $5$, $G$ any
group and $\varphi$ a homomorphism from $\B_n$ to $G$. If there
exist two distinct integers $i$ and $j$ in
$\{1,\,2,\dots,\,n-1\}$ and an nonzero integer $k$ such that
$\varphi({\tau_i}^k)=\varphi({\tau_j}^k)$, then
\smallskip

\centrer{$\varphi({\tau_1}^k)=\varphi({\tau_2}^k)=\dots=\varphi({\tau_{n-1}}^k)$}
\smallskip

\noindent If $k=1$, then $\varphi$ is a cyclic homomorphism.
\end{lem}
\smallskip

\DEM
For the simplicity of the proof, we assume
that $i=1$, but the case $i\not=1$ works in a similar way.
We introduce the element $\delta=\tau_1\tau_2\dots\tau_{n-1}$ lying in $\B_n$.
\smallskip

a) If $j=2$, let us conjugate the equality
$\varphi({\tau_1}^k)=\varphi({\tau_2}^k)$ by $\delta^\ell$ for $\ell\in\{0,\dots,\,n-3\}$.
We get $\varphi({\tau_{1+\ell}}^k)=\varphi({\tau_{2+\ell}}^k)$ for all $\ell$, so the lemma follows.
\smallskip

b) If $j>2$, then let us assume that $j<n-1$. The case $j=n-1$ works in a similar way.
Let us conjugate the equality $\varphi({\tau_1}^k)=\varphi({\tau_j}^k)$ by $\tau_j\tau_{j+1}\tau_j$.
We get $\varphi({\tau_1}^k)=\varphi({\tau_{j+1}}^k)$, and so
$\varphi({\tau_j}^k)=\varphi({\tau_{j+1}}^k)$. Now, let us conjugate this last equality by $\delta^{1-j}$.
We get $\varphi({\tau_1}^k)=\varphi({\tau_2}^k)$ so we have boiled down to the case a).
\fin
\bigskip

\begin{lem}[Criterion 2 for a homomorphism to be cyclic]
        \label{lem:B_n_in_group_abelian}
\mbox{}\\Any homomorphism from $\B_n$ ($n\geqslant 2$) in an abelian group is cyclic.
\end{lem}
\smallskip

\DEM Recall that a homomorphism $\rho$ from a group $G$ in an abelian group $A$
sends $[G,G]$ on $1_A$. Hence $A$ is isomorphic to a quotient of the abelianization of $G$.
Since the abelianization of $\B_n$ is infinite cyclic (cf. [Bi]), the lemma follows.\fin
\bigskip

\begin{lem}[Criterion on the existence of monodromy homomorphisms]
      \label{lem:existence_des_homomorphisms_of_monodromie}
\mbox{}\\Let $n$ be an integer greater than or equal to 3 and
$\Sigma$ a surface $\Sigma_{g,\,b}$. There exist monodromy
homomorphisms from $\B_n$ to $\PMod(\Sigma)$ if and only if
$g\geqslant \frac{n}{2}-1$. \end{lem}
\smallskip

\DEM Notice that the existence of
monodromy homomorphisms only depends on the existence of
$(n-1)$-chains of curves, and such a chain of curves exists if and only if
 $g\geqslant \frac{n}{2}-1$.\fin
\bigskip

The following lemma shows how transvections arise naturally
from central exact sequences of groups. Such sequences are
frequent between mapping class groups.
See for instance Proposition \ref{prop:suites_exactes_entre_MCGs}.(iv)-(vi).
\medskip

\begin{lem}[Criterion for two homomorphisms to be transvections the one of the other]
          \label{lem:transvections}
Let $1\to N\to G\xrightarrow{\psi} \widehat G\to 1$ be a
central exact sequence of groups, let $n$ be an integer greater
than or equal to $3$ and let $\rho$ and $\rho'$ be two
homomorphisms from $\B_n$ to $G$ such that
$\psi\rond\rho=\psi\rond\rho'$. Then
\begin{itemize}
\item[\;(i)] $\rho'$ is a transvection of $\rho$,
\item[(ii)]$\rho$ is cyclic if and only if $\psi\rond\rho$
    is cyclic. \end{itemize} \end{lem}
\smallskip

\DEM Let us prove item (i). For all integers $i$ in
$\{1,\dots,\,n-1\}$, there exists $g_i\in N$ such that
$\rho'(\tau_i)=\rho(\tau_i)g_i$. We have then the following
equalities, true for all integers $i$ in $\{1,\dots,\,n-1\}$:
\smallskip

\centrer{$\begin{array}{rcl}
  \rho'(\tau_i)\,\rho'(\tau_{i+1})\,\rho'(\tau_i) &=&
  \rho(\tau_{i})\,\rho(\tau_{i+1})\,\rho(\tau_{i})\,g_i\,g_{i+1}\,g_i,\\
  \rho'(\tau_{i+1})\,\rho'(\tau_{i})\,\rho'(\tau_{i+1}) &=&
  \rho(\tau_{i+1})\,\rho(\tau_{i})\,\rho(\tau_{i+1})\,g_{i+1}\,g_{i}\,g_{i+1}.
\end{array}$}

\smallskip

\noindent The braid relations in $\B_n$ imply that the four
members in these two equalities must be all equal. Therefore
for any integer $i$ in $\{1,\,2,\dots,\,n-1\}$, we have:
\smallskip

\centrer{$g_i\,g_{i+1}\,g_i=g_{i+1}\,g_{i}\,g_{i+1}$.}
\smallskip

\noindent But for all $i$ and $j$ in $\{1,\,2,\dots,n-1\}$, the
elements $g_i$ and $g_j$ commute, so they all are equal. Hence
$\rho'$ is a transvection of $\rho$.
\smallskip

Let us prove now item (ii). If $\rho$ is cyclic, then
$\psi\rond\rho$ is cyclic. Conversely, if $\psi\rond\rho$ is
cyclic, then it is clear that there exists a cyclic homomorphism
$\rho'$ such that $\psi\rond\rho'=\psi\rond\rho$. According to
item (i), this implies that $\rho$ is a transvection of
$\rho'$. Hence $\rho$ is cyclic. \fin
\bigskip

\TITRE{Homomorphisms of the same nature.}
Let us say that two homomorphisms $\rho_1$ and $\rho_2$ from the braid group to
two possibly different mapping class groups are
\emph{of the same nature} if they both are cyclic or if they both are transvection
of a monodromy homomorphism.
Of course, if two homomorphisms are not on the same nature,
we shall say that they are \emph{of different nature}.
\bigskip

Let $\Sigma$ be a surface and let $\M$ be a subgroup of
$\Mod(\Sigma)$ or of $\Mod(\Sigma,\,\bord\Sigma)$. According to their definition,
any transvection of monodromy homomorphism from $\B_n$ to $\M$ can be
characterized by the data of a $(n-1)$-chain
$(a_1,\dots,\,a_{n-1})$ of curves in $\Sigma$, of an integer
$\varepsilon\in\{\pm1\}$, and of a mapping class $V$ that
commutes with $T_{a_i}$ for all $i\in\{1,\,\dots,\,n-1\}$. We
are going to show that such a triple
$\big(\,(a_1,\dots,\,a_{n-1}),\,\varepsilon,\,V\,\big)$ is
unique.
\bigskip

   \begin{lem}(Uniqueness of the triple representing a transvection)
            \label{lem:uniqueness_triple}
    \mbox{}\\
    Let $n$ be an integer greater than or equal to 5, let $\Sigma$
    be a surface, and let $\M$ be a subgroup of $\Mod(\Sigma)$ or
    of $\Mod(\Sigma,\,\bord\Sigma)$. Let $\rho$ be a transvection
    of monodromy homomorphism from $\B_n$ to $\M$ such that there exist
    two triples
    $\big(\,(a_1,\dots,\,a_{n-1}),\,\varepsilon,\,V\,\big)$ and
    $\big(\,(c_1,\dots,\,c_{n-1}),\,\eta,\,W\,\big)$ satisfying the
    following for all $i\in\{1,\,\dots,\,n-1\}$:
    \smallskip

    \centrer{$\rho(\tau_i)=T_{a_i}^{\,\varepsilon}\,V=T_{c_i}^{\,\eta}\,W.$}
    \smallskip

    \noindent Then, these two triples are equal. \end{lem}
    \medskip

    \DEM Given the properties of $V$ and $W$, the computation of
    $\rho(\tau_1\tau_3^{-1})$ leads to:
    \smallskip

    \centrer{$\left(\,T_{a_1}\, T_{a_3}^{\ -1}\right)^{\varepsilon} =
              \left(\,T_{c_1}\, T_{c_3}^{\ -1}\right)^{\eta}.$}
   \smallskip

   \noindent This is an equality between multitwists (for
    $I(a_1,\,a_3)=I(c_1,\,c_3)=0$), so one of the following cases
    holds:

    \begin{tabular}{llllr}
      -- & either & $\varepsilon=\eta$,  & $a_1=c_1$ and $a_3=c_3$,\\
      -- &  or & $\varepsilon=-\eta$, & $a_1=c_3$ and $a_3=c_1$.
    \end{tabular}
    \hfill
    \begin{tabular}{c}
       \makebox[0cm]{(2)}\\
       \makebox[0cm]{(3)}\\
    \end{tabular}
\smallskip

    \noindent Similarly, the computation of
    $\rho(\tau_1\tau_4^{-1})$ leads to one of the following cases:

    \begin{tabular}{llll}
      -- & either & $\varepsilon=\eta$, & $a_1=c_1$ and $a_4=c_4$,\\
      -- &  or    & $\varepsilon=-\eta$,& $a_1=c_4$ and $a_4=c_1$.
    \end{tabular}
    \hfill
    \begin{tabular}{c}
       \makebox[0cm]{(4)}\\
       \makebox[0cm]{(5)}\\
    \end{tabular}
\smallskip

    \noindent The only compatible choice is that (2) and (4) hold.
    So $\varepsilon = \eta$, and $V=W$.
    Since for all $i\in\{1,\dots,\,n-1\}$, we have
    $T_{a_i}^{\,\varepsilon}\,V=T_{c_i}^{\,\eta}\,W$, we then
    deduce the equality $T_{a_i}=T_{c_i}$, and eventually the
    equality
    $(a_1,\,a_2,\,\dots,\,a_{n-1})=(c_1,\,c_2,\,\dots,\,c_{n-1})$.
    \fin

\bigskip

According to Lemma \ref{lem:uniqueness_triple}, we can set the following definitions.
\bigskip

\begin{defi}[Characteristic elements of a transvection of monodromy homomorphism]
          \label{def:triple}
Let $n$ be an integer greater than or equal to 4, $\Sigma$ a
surface, $\M$ a subgroup of $\Mod(\Sigma)$ or of
$\Mod(\Sigma,\,\bord\Sigma)$. Let $\rho$ be a transvection of
monodromy homomorphism from $\B_n$ to $\M$.
\smallskip

\Point The \emph{characteristic triple of $\rho$} is the
    \emph{unique}  triple
    $\big(\,(a_1,\dots,\,a_{n-1}),\,\varepsilon,\,V\,\big)$
    such that for all $i\leqslant n-1$, we have:
\smallskip

\centrer{$\rho(\tau_i)=T_{a_i}^{\,\varepsilon}\,V$.}
\smallskip

\Point We define:\\
  \ttiret \emph{the characteristic $(n-1)$-chain of the
   transvection $\rho$}, as being the $(n-1)$-chain
   $(a_1,\dots,\,a_{n-1})$,\\
 \ttiret \emph{the characteristic sign of the transvection
   $\rho$}, as being the integer $\varepsilon$,\\
 \ttiret \emph{the direction of the transvection $\rho$}, as
   being the mapping class $V$.
\smallskip

\Point We denote by $\Sigma(\rho)$ the tubular neighbourhood
    of $\dessous{\cup}{{i\leqslant n-1}}a_i$ where
    $(a_i)_{i\leqslant n-1}$ is the characteristic
    $(n-1)$-chain of $\rho$. Notice that $\Sigma(\rho)$ does
    not necessarily belong to $\Sub(\Sigma)$, for the
    boundary components of $\Sigma(\rho)$ can bound some
    disks in $\Sigma$.
\smallskip

\Point The transvection $\rho$ determines a unique pair
    $(\rho^*,\,\varphi)$ of homomorphisms such that for all
    $\xi\in\B_n$, we have:
\smallskip

\centrer{$\rho(\xi)=\rho^*(\xi)\varphi(\xi)$,}
\smallskip

\noindent where, for all $i\leqslant n-1$, we have:

 \centrer{ $\left\{\begin{array}{ccc}
  \rho^*(\tau_i) &=& T_{a_i}^{\,\varepsilon},\\
  \varphi(\tau_i) &=& V.
\end{array}\right.$}

\noindent The monodromy homomorphism $\rho^*$ and the cyclic
homomorphism $\varphi$ will be called respectively the
\emph{monodromy homomorphism} and the \emph{cyclic homomorphism}
\emph{associated to the transvection $\rho$}.
\end{defi}
\bigskip

Thus, the decomposition of a transvection of monodromy homomorphism
gives rise to two homomorphisms: a monodromy homomorphism $\rho^*$
determined by a $(n-1)$-chain of curves and a cyclic homomorphism
determined by the direction of the transvection, which is a
mapping class $V$ belonging to the centralizer of
$\rho^*(\B_n)$ in $\M$. Therefore the computation of this
centralizer is essential.
This will be done in Proposition \ref{prop:centralizer}.
\bigskip






\subsection{Proof of Theorem \ref{thm:Theorem_principal_case_à_boundary}}
\label{sec:énoncés_théorèmes}
\label{par:passage_à_boundary_components}

\medskip

In this subsection, we show the links that exist between the
sets $\Hom(\B_n,\,\PMod(\Sigma))$ and
$\Hom(\B_n,\,\Mod(\Sigma,\,\bord\Sigma))$ cf. Proposition
\ref{prop:passage_PMod_à_ModBord}, and we prove Theorem
\ref{thm:Theorem_principal_case_à_boundary}. Proposition
\ref{prop:passage_PMod_à_ModBord} is actually a corollary
of the following lemma.
\bigskip

\begin{lem}
                        \label{lem:ça_commute}
Let $\Sigma$ be a surface with a  nonempty boundary. Let $F$
and $G$ be mapping classes in $\PMod(\Sigma)$ such that $F$ and
$G$ commute. Let $\tilde F$ and $\tilde G$ be the lifts in
$\Mod(\Sigma,\,\bord\Sigma)$ of respectively $F$ and $G$. Then
$\tilde F$ and $\tilde G$ commute. \end{lem}
\smallskip

\DEM Let us start from the following central exact sequence
linking $\Mod(\Sigma,\,\bord\Sigma)$ and $\PMod(\Sigma)$:
\smallskip

\centrer[*]{$1\to\langle T_d\ ,\
d\in\Bord(\Sigma)\rangle\to\Mod(\Sigma,\,\bord\Sigma)\to\PMod(\Sigma)\to
1$.}
\smallskip

\noindent Then two lifts in $\Mod(\Sigma,\,\bord\Sigma)$ of a
same element of $\PMod(\Sigma)$ differ from a central element.
Therefore, in order to show the lemma, it is enough to show
that, given two elements $F$ and $G$ in $\PMod(\Sigma)$ that
commute, there exist two lifts $\tilde F$ and $\tilde G$ of $F$
and $G$ in $\Mod(\Sigma,\,\bord\Sigma)$ that commute. So we
start from two mapping classes $F$ and $G$ in $\PMod(\Sigma)$
that commute.
\medskip

1. Suppose first that $F$ and $G$ are periodic. Since $\Sigma$
has a nonempty boundary, according to Lemma
\ref{lem:group_cyclic}, $F$ and $G$ span a cyclic group.
Therefore there exists $H\in\PMod(\Sigma)$ such that $F$ and
$G$ are both some powers of $H$. Let $p$ and $q$ be two
integers such that $F=H^p$ and $G=H^q$. Let $\tilde H$ be a
lift of $H$ in $\Mod(\Sigma,\,\bord\Sigma)$. Then $\tilde H^p$
and $\tilde H^q$ are some lifts of $F$ and $G$ that commute.
\smallskip

2. We now turn to the case where $F$ is pseudo-Anosov and $G$
is periodic. Let us denote by $\tilde F$ and $\tilde G$ some
lifts in $\Mod(\Sigma,\,\bord\Sigma)$ of $F$ and $G$
respectively. Since $FGF^{-1}=G$, there exists a central
mapping class $V\in\Mod(\Sigma,\,\bord\Sigma)$ such that
$\tilde F\tilde G\tilde F^{-1}=\tilde G V$. Let $p$ be the
order of $G$ and let $W$ be the central mapping class of
$\Mod(\Sigma,\,\bord\Sigma)$ such that $\tilde G^p=W$. Then, we
have on one hand:
\smallskip

\centrer[1]{$(\tilde F \tilde G \tilde F^{-1})^p=\tilde F
\tilde G^p \tilde F^{-1}=\tilde FW\tilde F^{-1}=W$,}

\noindent and we have on the other hand:

\centrer[2]{$(\tilde F \tilde G \tilde F^{-1})^p= (\tilde G
V)^p=W V^p$.}
\smallskip

\noindent When we compare (1) and (2), it comes out that $V$ is
trivial for $\Mod(\Sigma,\,\bord\Sigma)$ is torsion-free. So
$\tilde F$ and $\tilde G$ commute.
\smallskip

3. If $F$ and $G$ are pseudo-Anosov, according to Proposition
\ref{prop:structure_centralizer}.(iv), there exist two nonzero
integers $p$ and $q$ such that $F^p=G^q$. Let $\ell$ and $k$ be
two integers such that $\ell p+kq=p\wedge q=d$. Let us set
\smallskip

\begin{itemize}
\item[\point] $H=F^kG^\ell$ (so $H$ satisfies $H^p=G^d$ and
    $H^q=F^d$; hence $H$ is pseudo-Anosov),
\item[\point] $P=F(H^{-1})^{(q/d)}$ (so $P^d=1$ and
    $F\in\langle P,\,H\rangle$),
\item[\point] $Q=G(H^{-1})^{(p/d)}$ (so $Q^d=1$ and
    $G\in\langle Q,\,H\rangle$).
\end{itemize}
\smallskip

\noindent According to Lemma \ref{lem:group_cyclic}, since $P$
and $Q$ are two periodic mapping classes that commute, there
exists a mapping class $R\in\langle P,\,Q\rangle$ such that
$\langle P,\,Q\rangle=\langle R\rangle$. Thus $F$ and $G$
belong to the abelian group spanned by $H$ and $R$. Then
according to step 2., two lifts $\tilde H$ and $\tilde R$ of
$H$ and $R$ in $\Mod(\Sigma,\,\bord\Sigma)$ span an abelian
group, too. Moreover, the latter contains two lifts $\tilde F$
and $\tilde G$ of $F$ and $G$ in $\Mod(\Sigma,\,\bord\Sigma)$.
In particular, $F$ and $G$ admit two lifts $\tilde F$ and
$\tilde G$ that commute.
\smallskip

4. Let $F$ and $G$ be any two mapping classes of
$\PMod(\Sigma)$ that commute. Let $\A$ be the set of curves
$\sigma(F)\cup\sigma(G)$. Notice that $\A$ is a simplex
according to Proposition \ref{prop:properties_sigma}.(iii).
\smallskip

4.a) Let us assume that $F$ and $G$ belong to
$\P_\A\Mod(\Sigma)$ (i.e. $F$ and $G$ preserve each curve of
the set $\A=\sigma(F)\cup\sigma(G)$). We are going to describe
for any $H\in \P_\A\Mod(\Sigma)$ a construction of a lift of
$H$ in $\P_\A\Mod(\Sigma,\,\bord\Sigma)$, then we will apply it
to $F$ and $G$. First, let us consider the following
commutative diagram where all the arrows are canonical
($rec_\A$ comes from Proposition
\ref{prop:suites_exactes_entre_MCGs}, the three other homomorphisms
have been introduced in Subsection \ref{par:definitions_MCG}
Definition \ref{défi:MCG}):
\medskip

\centrer{$\begin{diagram}
   \node{\rule{3cm}{0cm}H_3\in}
 \node{\P_\A\Mod(\Sigma,\,\bord\Sigma)} \arrow{s,l}{\for_{\bord\Sigma}}
 \node{\Mod(\Sigma_\A,\,\bord\Sigma_\A)}\arrow{w,t}{rec_\A}\arrow{s,r}{\for_{\bord\Sigma_\A}}
  \node{\ni H_2\rule{3cm}{0cm}}\\
  \node{\rule{3cm}{0cm}H_4,\,H\in}
 \node{\P_\A\Mod(\Sigma)}\arrow{e,b}{\cut_\A}
 \node{\PMod(\Sigma_\A)}
  \node{\ni H_1,\,H_5\rule{3cm}{0cm}}
\end{diagram}$}
\medskip

\noindent For any $H\in\P_\A\Mod(\Sigma)$, let us denote by
$H_i$, $1\leqslant i\leqslant 5$, the following mapping
classes, derived from $H$ when following the diagram above:
\smallskip

$
\begin{array}{lll}
  \point & H_1=\cut_\A(H),                                                          & \mbox{ so } H_1\in\PMod(\Sigma_\A),\\
  \point & H_2 \mbox{ a lift of } H_1 \mbox{ in } \Mod(\Sigma_\A,\,\bord\Sigma_\A), & \mbox{ so } H_2\in\Mod(\Sigma_\A,\,\bord\Sigma_\A),\\
  \point & H_3=rec_\A(H_2),                                                         & \mbox{ so } H_3\in\Mod(\Sigma_\A,\,\bord\Sigma_\A),\\
  \point & H_4=\for_{\bord\Sigma}(H_3),                                             & \mbox{ so } H_4\in\P_\A\Mod(\Sigma,\,\bord\Sigma),\\
  \point & H_5=\cut_\A(H_4),                                                        & \mbox{ so } H_5\in\PMod(\Sigma_\A).
\end{array} $

\smallskip

\noindent The diagram is commutative:
$\for_{\bord\Sigma_\A}=(\cut_\A)(\for_{\bord\Sigma})(rec_\A)$,
so $H_1=H_5$. But we have the following central exact sequence:
\smallskip

\centrer[**]{$1\to\T\to\P_\A\Mod(\Sigma)\xrightarrow{\cut_\A}\PMod(\Sigma_\A)\to
1$ ,}
\smallskip

\noindent where $\T=\langle T_d\ ,\ d\in\A\rangle$. Hence $H$
and $H_4$, the preimages of $H_1$ and $H_5$ by $\cut_\A$,
differ from a multitwist along some curves of $\A$. Hence, up
to elements in $\T$, the mapping class $H_3$ is a lift of $H$.
\medskip

Let us apply this to $F$ and $G$. As $F$ and $G$ commute, $F_1$
and $G_1$ commute. But on each connected component of
$\Sigma_\A$, the restrictions of $F_1$ and $G_1$ are periodic
or pseudo-Anosov, so we can apply what was shown above in steps
1., 2. and 3., and deduce from it that $F_2$ and $G_2$ commute.
Hence $F_3$ and $G_3$ commute as well. Now, as we just saw it
with $H$, there exist $T$ and $T'$ belonging to $\T$ such that
$\tilde F=F_3T$ and $\tilde G=G_3T'$ are some lifts of $F$ and
$G$. Moreover $T$ and $T'$ are central in
$\P_\A\Mod(\Sigma,\,\bord\Sigma)$ and in addition, $F_3$ and
$G_3$ commute, so $\tilde F$ and $\tilde G$ commute.
\medskip

4.b) In the general case, if $F$ and $G$ are any two mapping
classes that commute, let us denote by $\tilde F$ and $\tilde
G$ some lifts of $F$ and $G$ in $\Mod(\Sigma,\,\bord\Sigma)$.
\emph{A priori}, there exists a multitwist $W$ along the
boundary components such that $\tilde F\tilde G\tilde
F^{-1}=\tilde GW$. Once again, let us set
$\A=\sigma(F)\cup\sigma(G)$. Since $F$ and $G$ commute, they
preserve globally $\A$, so there exists a nonzero integer $m$
such that $F^m$ and $G^m$ preserve $\A$ curve-wise. In other
words, $F^m$ and $G^m$ belong to $\P_\A\Mod(\Sigma)$. So,
according to 4.a),
\smallskip

\centrer{$\tilde F^m\tilde G^m\tilde F^{-m}=\tilde G^m$.}
\smallskip

\noindent Now, the equality $\tilde F\tilde G\tilde
F^{-1}=\tilde GW$ implies that $\tilde F\tilde G^m\tilde
F^{-1}=\tilde G^m W^m$, then
\smallskip

\centrer{$\tilde F^m\tilde G^m\tilde F^{-m}=\tilde G^m
W^{m^2}$.}
\smallskip

\noindent  So $W^{m^2}$ is trivial. But
$\Mod(\Sigma,\,\bord\Sigma)$ is torsion-free, hence $W$ is
trivial and $\tilde G$ and $\tilde F$ commute.\fin
\bigskip

\begin{prop}[Lifting from
$\Hom\big(\B_n,\,\PMod(\Sigma)\big)$ in
$\Hom\big(\B_n,\,\Mod(\Sigma,\,\bord\Sigma)\big)$]
          \label{prop:passage_PMod_à_ModBord}
\mbox{}\\ Let $n$ be an integer greater than or equal to 3, let
$\Sigma$ be a surface and $\rho\;:\; \B_n\to\PMod(\Sigma)$ a
homomorphism. Let us recall that we denote by $\for_{\bord\Sigma}$,
or $\for$, the canonical epimorphism from
$\Mod(\Sigma,\,\bord\Sigma)$ to $\PMod(\Sigma)$. Then:
\begin{itemize}
\item[\;\;(i)] There exists a homomorphism $\tilde\rho\;:\;
\B_n\to\Mod(\Sigma,\,\bord\Sigma)$ such that
$\for\rond\tilde\rho=\rho$.

\item[\;(ii)] Such a homomorphism $\tilde\rho$ is unique up to transvection,
that is, if $\tilde\rho_1$ and $\tilde\rho_2$ satisfy
$\for(\,\tilde\rho_1)=\for(\,\tilde\rho_2)=\rho$, then there
exists $V\in\Mod(\Sigma,\,\bord\Sigma)$ such that $V$ is in the
centralizer of $\tilde\rho_1(\B_n)$ and of $\tilde\rho_2(\B_n)$
and satisfies for all $i\in\{1,\dots,\, n-1\}$:
\smallskip

\centrer{$\tilde\rho_2(\tau_i)=\tilde\rho_1(\tau_i)\,V$.}

\item[(iii)] Such a homomorphism $\tilde\rho$ is cyclic if and only if $\rho$
is cyclic.

\item[(iv)] Such a homomorphism $\tilde\rho$ is a transvection of monodromy
homomorphism if and only if $\rho$ is a transvection of monodromy
homomorphism.
\end{itemize}
\end{prop}
\bigskip

\DEM
\smallskip

(i) Let us start from the following central exact sequence:
\smallskip

\centrer[*]{$1\to\langle\,T_d,\,d\in\Bord(\Sigma)\,\rangle\to\Mod(\Sigma,\,\bord\Sigma)\xrightarrow{\for}\PMod(\Sigma)\to
1$.} \smallskip

\noindent For all $i\in\{1,\,\dots,\,n-1\}$, let $A_i$ be a
mapping class of $\Mod(\Sigma,\,\bord\Sigma)$ such that
$\for(A_i)=\rho(\tau_i)$. Then for all
$i\in\{1,\,\dots,\,n-2\}$, we have:
\smallskip

\centrer{$\for(A_iA_{i+1}A_i)=\for(A_{i+1}A_iA_{i+1})$,}
\smallskip

\noindent hence, according to the exact sequence $(*)$, for all
$i\in\{1,\,\dots,\,n-2\}$, there exists a multitwist denoted by
$W_i$ along some boundary components of $\Sigma$ such that
\smallskip

\centrer{$A_iA_{i+1}A_i=A_{i+1}A_iA_{i+1}W_i$.}
\smallskip

\noindent Let us set:
\smallskip

\centrer{$\left\{\begin{array}{l}
  A'_1:=A_1\,,\\
  A'_i:=A_iW_1W_2\cdots W_{i-1} \mbox{ when } 2\leqslant i \leqslant n-1.
\end{array}\right.$}
\smallskip

\noindent Let us recall that the $W_i$ are central. Hence for
all $i\in\{1,\,\dots,\,n-2\}$, we have:
\smallskip

\centrer{$A'_iA'_{i+1}A'_i=A'_{i+1}A'_iA'_{i+1}$.}
\smallskip

\noindent Besides, for all integers $i$ and $j$ smaller than or
equal to $n-1$ such that $|i-j|\geqslant 2$, the mapping
classes $A_i$ and $A_j$ commute, so according to Lemma
\ref{lem:ça_commute}, the mapping classes $A'_i$ and $A'_j$
commute as well. Finally the map $\tilde\rho$ defined by
\smallskip

\centrer{$\tilde\rho(\tau_i)=A'_i$}
\smallskip

\noindent is a homomorphism from $\B_n$ to $\Mod(\Sigma,\,\bord\Sigma)$. Moreover, by construction, we
have $\for(\tilde\rho)=\rho$.
\medskip

(ii) Let $\tilde\rho_1$ and $\tilde\rho_2$ be two homomorphisms from
$\B_n$ to $\Mod(\Sigma,\,\bord\Sigma)$ that satisfy
$\for(\tilde\rho_1)=\for(\tilde\rho_2)$. According to the
central exact sequence ($*$), we can apply Lemma
\ref{lem:transvections}: $\tilde\rho_2$ is a transvection of
$\tilde\rho_1$. \smallskip

(iii) According to Lemma \ref{lem:transvections}, $\tilde\rho$ is
cyclic if and only if $\rho$ is cyclic.
\smallskip

(iv) If $\tilde\rho$ is a transvection of monodromy homomorphism, it
is clear that $\for\rond\tilde\rho$ is still a transvection of
monodromy homomorphism. Conversely, suppose that $\rho$ is a transvection of
monodromy homomorphism characterized by a triple
$\big((a_i)_{1\leqslant i\leqslant n-1},\,\varepsilon,\,V\big)$.
Let $\tilde V$ be a lift of $V$ in
$\Mod(\Sigma,\,\bord\Sigma)$. For all
$i\in\{1,\,\dots,\,n-1\}$, the mapping class $V$ commutes with
$T_{a_i}$, so $V$ fixes $a_i$ and so does $\tilde V$ as well. Hence
$\tilde V$ commutes with $T_{a_i}$ in
$\Mod(\Sigma,\,\bord\Sigma)$ (we also could have quoted Lemma \ref{lem:ça_commute}).
Then we can define $\tilde\rho$ as being the
transvection of monodromy homomorphism characterized by the
triple $\big((a_i)_{1\leqslant i\leqslant n-1},\,\varepsilon,\,\tilde V\big)$.
It is clear that $\tilde\rho$ is a lifting of $\rho$.
According to step (ii), all liftings of $\rho$ are transvections of $\tilde\rho$,
therefore all liftings of $\rho$ are transvections of monodromy homomorphisms.\fin
\bigskip

\TITRE{Theorem \ref{thm:Theorem_principal_case_à_boundary}}
(Homomorphisms from
$\B_n$ to $\Mod(\Sigma,\,\bord\Sigma)$).\\
{\em Let $n$ be any integer greater than or equal to 6. Let
$\Sigma$ be a surface $\Sigma_{g,\,b}$ where $g\leqslant n/2$.
Let $\tilde\rho$ a homomorphism from $\B_n$ to $\Mod(\Sigma,\,\bord\Sigma)$. Then $\tilde\rho$ is cyclic or is
a transvection of monodromy homomorphism. Moreover such
transvections of monodromy homomorphisms exist if and only if
$g\geqslant \frac{n}{2}-1$.}
\smallskip

\DEM The second part of Theorem
\ref{thm:Theorem_principal_case_à_boundary} has been proved in Lemma
\ref{lem:existence_des_homomorphisms_of_monodromie}.
It remains to show that all the homomorphisms from $\B_n$ to $\Mod(\Sigma,\,\bord\Sigma)$ are either cyclic, or
transvections of monodromy homomorphisms. Let $\tilde\rho$ be a
noncyclic homomorphism from $\B_n$ to $\Mod(\Sigma,\,\bord\Sigma)$.
If we compose $\tilde\rho$ with the projection
$\Mod(\Sigma,\,\bord\Sigma)\xrightarrow{\ \for\
}\PMod(\Sigma)$, we get a homomorphism $\rho$ from $\B_n$ to $\PMod(\Sigma)$. According to Proposition
\ref{prop:passage_PMod_à_ModBord}.(iii), $\rho$ is not cyclic, so
according to Theorem \ref{thm:Theorem_principal}, $\rho$ is a
transvection of monodromy homomorphism, then according to Proposition
\ref{prop:passage_PMod_à_ModBord}.(iv), $\tilde\rho$ is a
transvection of monodromy homomorphism.\fin
\bigskip



\subsection{Centralizer of the monodromy homomorphisms}
            \label{par:transvections_monodromy_homomorphisms}
\medskip

\begin{prop}[Centralizer of the image of a monodromy homomorphism]
                \label{prop:centralizer}
\mbox{}\\
Let $n$ be an integer greater than or equal to 5.\\
Let $\Sigma$ be a surface $\Sigma_{g,\,b}$ with $g\geqslant \frac{n}{2}-1$.\\
Let $\M$ be one of the groups $\Mod(\Sigma)$, $\PMod(\Sigma)$, or $\Mod(\Sigma,\,\bord\Sigma)$.\\
Let $\rho\ :\ \B_n\to\M$ be a monodromy homomorphism with $(a_i)_{i\leqslant n-1}$ as characteristic chain of curves.\\
Let $\M^{\Sigma(\rho)}$ be the group of the mapping classes in $\M$ that preserve the subsurface $\Sigma(\rho)$
and that induce the identity in $\Mod\big(\Sigma(\rho)\big)$. Then:
\begin{itemize}
\item[\;(i)] the centralizer of $\rho(\B_n)$ in $\M$ is
    reduced to the group $\M^{\Sigma(\rho)}$, except in the
    below cases (a), (b) or (c) where it is equal to
    the group spanned by $\M^{\Sigma(\rho)}$ and by $Z$,
    where $Z$ is any extension in $\Mod(\Sigma)$ of the
    hyper-elliptic mapping class of $Mod(\Sigma(\rho))$ associated to $(a_i)_{i\leqslant n-1}$.
\smallskip

\noindent Cases (a), (b) and (c) are the following:\\
  \ \ (a) the integer $n$ is odd,\\
  \ \ (b) the curve simplex
      $\{a_1,\,a_3,\,\dots,\,a_{n-1}\}$ is
      non-separating,\\
  \ \ (c) $\Sigma_{\{a_1,\,a_3,\,\dots,\,a_{n-1}\}}$
      consists in two homeomorphic connected components
      and $\M=\Mod(\Sigma)$.

\item[(ii)] for any mapping class $V$ belonging to the
    centralizer of $\rho(\B_n)$ in $\M$, the mapping class
    $V^2$ belongs to $\M^{\Sigma(\rho)}$.
\end{itemize}
\end{prop}

\DEM
\smallskip

\emph{Let us show item (i) when $n$ is odd}.
\smallskip

It is clear that the group
$\M^{\Sigma(\rho)}$ defined in the statement of Proposition
\ref{prop:centralizer} is in the centralizer of
$\rho(\B_n)$. As for the mapping class $Z$, it preserves each curve $a_i$
by definition, so it lies in the centralizer of
$\rho(\B_n)$ as well. Therefore group spanned by $\M^{\Sigma(\rho)}$ and $Z$ is in the centralizer
of $\rho(\B_n)$. Conversely, let us show that any element in the
centralizer of $\rho(\B_n)$ coincides with an element of
$\M^{\Sigma(\rho)}$ possibly composed by $Z$. In this purpose, let
us start from a mapping class $F$ lying in the centralizer of
$\rho(\B_n)$. We will first set some definitions in $\Sigma$, then
we will study $F$.
\bigskip

a) \emph{Definitions of some curves in $\Sigma$.}
\smallskip

\noindent For all integers $i\in\{2,\,3,\dots,\,n-1\}$, let us set:

\centrer{$\Delta_i=\tau_1(\tau_2\tau_1)\dots(\tau_i\tau_{i-1}\dots\tau_1)$.}
\smallskip

\begin{figure}[!h]
 \Includegraphics{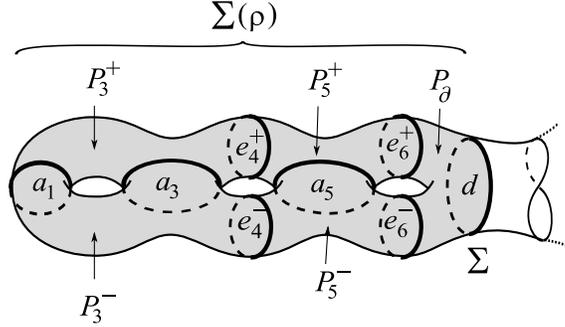}
 \caption{Cutting $\Sigma(\rho)$ in pairs of pants;
 the curve $d$ and the simplex $\A^1$ (case where $n$ is odd).}
 \label{fig:centralisateurCasImpair}
\end{figure}

\noindent For any even $i$ in $\{4,\dots,\,n-1\}$, let $e_i^+$ and
$e_i^-$ be the two curves such that
$\rho(\Delta_{i-1}^2)=T_{e_i^+}T_{e_i^-}$. By induction on the odd
integer $i$ in $\{3,\dots,\,n-2\}$, we define the pairs of pants
$P_i^+$ and $P_i^-$  (cf. Figure \ref{fig:centralisateurCasImpair})
in such a way that:
\smallskip

\begin{itemize}
\item[\point] when $i=3$, let us denote by $P_3^+$ and $P_3^-$
respectively the pairs of pants included in $\Sigma$ whose
boundaries are $\{a_1, a_3, e_4^+\}$ and $\{a_1, a_3, e_4^-\}$
respectively,
\smallskip

\item[\point] when $i$ is an odd integer in $\{5,\dots,\,n-2\}$ and
when $P_{i-2}^+$ and $P_{i-2}^-$ have been defined, even if it means
swapping $e_{i+1}^+$ and $e_{i+1}^-$, we can assume that
$\{e_{i-1}^+, a_i, e_{i+1}^+\}$ and $\{e_{i-1}^-, a_i, e_{i+1}^-\}$
are the boundary components of two pairs of pants that we denote by
$P_i^+$ and $P_i^-$ respectively. \end{itemize}
\smallskip

\noindent Let $d$ be the curve such that
$\rho(\Delta_{n-1}^4)=T_{d}$. We denote by $P_\bord$ the pair of
pants whose boundary is $\{e_{n-2}^+,\,e_{n-2}^-,\,d\}$ (cf. Figure
\ref{fig:centralisateurCasImpair}). We denote by $\A^{0}$ the union
of the curves $a_i$ where $i$ is even in $\{2,\dots,\,n-1\}$, and we
denote by $\A^{1}$ the union of the curves $a_i$ where $i$ is odd in
$\{1,\,2,\dots,\,n-2\}$ and of the curves $e_j^+$ and $e_j^-$ where
$j$ is even in $\{4,\dots,\,n-1\}$.
\medskip

b) \emph{A mapping class $F$ in the centralizer of $\rho(\B_n)$.}
\smallskip

\noindent The mapping class $F$ commutes with $T_{a_i}$ for any odd
$i$ in $\{1,\dots,\,n-2\}$, so $F(a_i)=a_i$. The mapping class $F$
commutes also with $T_{e_j^+}T_{e_j^-}$ for any even $j$ in
$\{4,\dots,\,n-1\}$, so $F(\{e_j^+,\,e_j^-\})=\{e_j^+,\,e_j^-\}$.
Hence $F$ preserves the set of subsurfaces $\{P_i^+,\,P_i^-\}$ for
any odd $i$ in $\{3,\dots,\,n-2\}$. Finally, $F$ commutes with $T_d$
for $T_d$ belongs to $\rho(\B_n)$, so $F$ preserves the curve $d$,
preserves the pair of pants $P_\bord$ and preserves the surface
$\Sigma(\rho)$ included in $\Sigma$ with $d$ as boundary, and
containing the curve $a_1$. Notice that for any odd $i$ in
$\{3,\dots,\,n-4\}$, the pairs of pants $P_{i}^+$ and $P_{i+2}^+$
have $e_{i+1}^+$ as common boundary component, so the pairs of pants
$F(P_{i}^+)$ and $F(P_{i+2}^+)$ have $F(e_{i+1}^+)$ as common
boundary component. Then two situations can happen concerning
$F(P_{i}^+)$ and $F(P_{i}^+)$:
\begin{itemize}
\item either $F(e_{i+1}^+)=e_{i+1}^+$, and then $F(P_{i}^+)=P_{i}^+$ and
$F(P_{i+2}^+)=P_{i+2}^+$,
\item or $F(e_{i+1}^+)=e_{i+1}^-$, and then $F(P_{i}^+)=P_{i}^-$ and
$F(P_{i+2}^+)=P_{i+2}^-$. \end{itemize} Finally, by induction, only
two situations can happen concerning $F$:
\smallskip

\begin{itemize}
\item First alternative: for all odd integers $i\in\{3,\dots,\,n-2\}$,
$F(P_i^+)=P_i^+$. Then $F$ fixes $e_j^+$ and $e_j^-$ for all even
integers $j\in\{4,\dots,\,n-1\}$. We define
$F'\in\Mod(\Sigma(\rho))$ as being the restriction of $F$ to
$\Sigma(\rho)$. \smallskip

\item Second alternative: for any odd $i$
in $\{3,\dots,\,n-2\}$, we have $F(P_i^+)=P_i^-$. Then for any even
$j$ in $\{4,\dots,\,n-1\}$, the mapping class $F$ swaps $e_j^+$ and
$e_j^-$. But for any odd $i$ in $\{1,\dots,\,n-1\}$, the mapping
class $Z$ fixes the curves $a_i$, and for any even $j$ in
$\{4,\dots,\,n-1\}$, the mapping class $Z$ swaps the curves $e_j^+$
and $e_j^-$ . Hence $FZ$ fixes all the curves of $\A^1$. Since $F$
and $Z$ preserve the surface $\Sigma(\rho)$ included in $\Sigma$, we
can define $F'\in\Mod(\Sigma(\rho))$ as being the restriction of
$FZ$ to $\Sigma(\rho)$. \end{itemize}
\smallskip

\noindent Let us examine $F'$. The mapping class $F'$ fixes all the
curves of $\A^1$, hence preserves each subsurface of
$\Sub_{\A^1}(\Sigma(\rho))$, which are pairs of pants, and preserves
each of their boundary components. So $F'$ induces in
$\PMod(\,\big(\Sigma(\rho)\big)_{\A^1}\,)$ a trivial mapping class,
where $\big(\Sigma(\rho)\big)_{\A^1}$ is the surface we get after
having cut $\Sigma(\rho)$ along the curves of $\A^1$. Then,
according to the following exact sequence:
\smallskip

\centrer{$1\to\langle\,T_a,\,a\in\A^1\,\rangle\to\P_{\A^1}\Mod(\Sigma(\rho))
\to\PMod(\,\big(\Sigma(\rho)\big)_{\A^1}\,)\to1$ ,}
\smallskip

\noindent the mapping class $F'$ is a multitwist along the curves of
$\A^1$. However, $F'$ commutes with $T_{a}$, $a\in\A^0$, so
according to Proposition \ref{prop:properties_twists_Dehn}, the
curves of $\A^0$ are reduction curves of $F'$. But each curve of $\A^1$
intersects one of the curves of $\A^0$, so no curve in $\A^1$ can be
an essential reduction curve of $F'$. Since $F'$ was to be a multitwist
along the curves of $\A^1$, it then must be trivial.

Let us come back to the mapping class $F$. The restriction of $F$ to
$\Sigma(\rho)$, or the restriction of $FZ$ to $\Sigma(\rho)$, equals
$\Id$ in $\Mod(\Sigma(\rho))$. Hence the centralizer of $\rho(\B_n)$
in $\M$ is the group spanned by $Z$ and by $\M^{\Sigma(\rho)}$, the
subgroup of $\M$ of the mapping classes inducing the identity
mapping class on $\Sigma(\rho)$.
\bigskip
\bigskip

\emph{Let us show item (i) when $n$ is even}.
\smallskip

\begin{figure}[!h]
 \Includegraphics{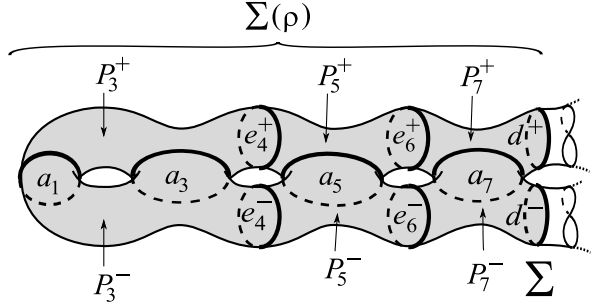}
 \caption{Cutting $\Sigma(\rho)$ in pairs of pants;
 the curves $d^+$, $d^-$ and the simplex $\A^1$ (case where $n$ is even).}
 \label{fig:centralisateurCasPair}
\end{figure}

The proof is very similar to the odd case. Only the inclusion
$\Z_\M(\rho(\B_n))\subset\big\langle \M^{\Sigma(\rho)},\,Z\big\rangle$
Similarly to the odd case, we first define the following
topological objets in $\Sigma(\rho)$, drawn in Figure
\ref{fig:centralisateurCasPair}:
\begin{itemize}
  \item a $(n-1)$-chain of curves $\big(a_i\big)_{1\leqslant i\leqslant n-1}$,
  \item some curves $e_{j}^{+}$ and $e_{j}^{-}$ for any even integer $j\in\{4,\dots,\,n-2\}$,
  \item some curves $d^+$ and $d^-$,
  \item some pairs of pants $P_{i}^{+}$ and $P_{i}^{-}$ for any odd integer $i\in\{3,\dots,\,n-1\}$,
  \item a set $\A$ of the curves $a_i$ where $i$ is odd in $\{1,\dots,\,n-1\}$,
  \item a set $\A^1$ of the curves of $\A$ and of the curves $e_{j}^{+}$ and $e_{j}^{-}$ for any even $j$ in $\{4,\dots,\,n-2\}$.
\end{itemize}

\noindent Let us start from a mapping class $F$ belonging to the
centralizer of $\rho(\B_n)$. Notice that $\Sigma(\rho)=\Sigma(\B_n)$
and let $H$ be the hyper-elliptic mapping class, belonging to
$\Mod(\Sigma(\rho))$). The action of $H$ on the curves of $\A^1$
consists in fixing the curves $a_i$ for any odd $i$ in
$\{1,\dots,\,n-1\}$, in swapping the curves $e_j^+$ and $e_j^-$ for
any even $j$ in $\{4,\dots,\,n-2\}$, and in swapping the curves
$d^+$ and $d^-$. As in the odd case, by considering the action of $F$
on the set of pairs of pants of $\Sub_{\A^1(\Sigma(\rho))}$, we see
that the restriction of $F$ to $\Sigma(\rho)$ coincides either with
$H$ or with the identity of $\Sigma(\rho)$, depending on whether $F$
fixes or swaps the boundary components $d^+$ and $d^-$.

Assume that we are in case (a) or (b). Then the mapping class
$H\in\Mod(\Sigma(\rho))$ can be extended on $\Sigma$. We denote by
$Z$ this extension, which is a mapping class of $\M$. In all the
other cases (different from (a) and (b)), the curves $d^+$ and $d^-$
do not belong to the same orbit under the action of $\PMod(\Sigma)$
on $\Courb(\Sigma,\,\bord\Sigma)$, so in these case, $F$ cannot swap
$d^+$ and $d^-$ and cannot coincide with $H$ on $\Sigma(\rho)$, so
$F$ induces in $\Mod(\Sigma(\rho))$ the identity mapping class. To
conclude,
\begin{itemize}
            \item when one of the conditions (a) or (b) is satisfied,
      the centralizer of $\rho(\B_n)$ in $\M$ is the group spanned by $Z$ and
      $\M^{\Sigma(\rho)}$,
            \item if none of the conditions (a) or (b) is satisfied,
      the centralizer of $\rho(\B_n)$ in $\M$ is the subgroup
      $\M^{\Sigma(\rho)}$.
\end{itemize} This shows item (i).
\bigskip
\bigskip

\emph{Let us show item (ii)}.
\smallskip

Let $V$ belongs to the centralizer of $\rho(\B_n)$ in $\M$.
We assume that the centralizer of $\rho(\B_n)$ in $\M$ is
not reduced to $\M^{\Sigma(\rho)}$ (otherwise, the result is trivial).
There exist $\varepsilon\in\{0,\,1\}$ and $V'\in\M^{\Sigma(\rho)}$ such
that $V=Z^\varepsilon\, V'$. We assume
that $\varepsilon=1$ (otherwise, the result is trivial).
Let us denote by $H$ the mapping class induced by $Z$ on $\Sigma(\rho)$.
Notice that $Z$ and $V'$ preserve $\Sigma(\rho)$ and that
$V'$ induces the identity on $\Sigma(\rho)$, so $ZV'Z$
preserves $\Sigma(\rho)$ and induces the mapping class $H^2=\Id$ in $\Mod(\Sigma(\rho))$. So
$ZV'Z$ belongs to $\M^{\Sigma(\rho)}$. Now, since $V^2=(Z\,V'\,Z)\,V'$,
it follows that $V^2\in\M^{\Sigma(\rho)}$.\fin
\bigskip
\bigskip



\subsection{Proof of Theorem \ref{thm:injectivity_des_homomorphisms}}
        \label{par:Bn_in_MCG}
\medskip

According to Theorems \ref{thm:Theorem_principal} and
\ref{thm:Theorem_principal_case_à_boundary}, when the genus of
$\Sigma$ is bounded by $n/2$, the homomorphisms from $\B_n$ in
the mapping class group associated to the surface $\Sigma$ are
either cyclic, or are some transvections of monodromy homomorphisms.
Consequently, the issue of the injectivity of the homomorphisms from
$\B_n$ in the mapping class group is reduced to the issue of the
injectivity of the monodromy homomorphisms (cf. Proposition
\ref{prop:injectivity_des_homomorphisms_of_monodromie}) and of the
transvections of the monodromy homomorphisms (cf. Proposition
\ref{prop:transvection_of injection}), for the cyclic homomorphisms
obviously cannot be injective. These different results are gathered
in Theorem \ref{thm:injectivity_des_homomorphisms}.
\medskip

Given $\rho$, a transvection of monodromy homomorphism from $\B_n$ in
the mapping class group of a surface $\Sigma$, the surface
$\Sigma(\rho)$ (described in Definition \ref{def:triple}) will help
us in characterizing the injectivity of the transvections of
monodromy homomorphisms.
\bigskip

\begin{prop}[Injectivity of the monodromy homomorphisms]
    \label{prop:injectivity_des_homomorphisms_of_monodromie}
\mbox{}\\ Let $n\geqslant 6$ and let $\Sigma$ be any surface. We distinguish the cases $\Mod(\Sigma,\,\bord\Sigma)$ and $\PMod(\Sigma)$.
\begin{itemize}
\item[\;(i)] A monodromy homomorphism $\tilde\rho\ :\ \B_n\to\Mod(\Sigma,\,\bord\Sigma)$
is injective if and only if we have:

\centrer{$\Bord(\Sigma(\tilde\rho))\subset\Courb(\Sigma,\,\bord\Sigma)$.}
\smallskip

\noindent In
other words, the boundary components of $\Sigma(\tilde\rho)$ do not
bound any disk in $\Sigma$.

\item[\;(ii)] A monodromy homomorphism $\rho\ :\ \B_n\to\PMod(\Sigma)$
is injective if and only if we have:

\centrer{$\left\{\begin{array}{l}
\Bord(\Sigma(\rho))\subset\Courb(\Sigma,\,\bord\Sigma),\\
\Bord(\Sigma(\rho))\not\subset\Bord(\Sigma).\end{array}\right.$ }
\smallskip

\noindent In other words, the
boundary components of $\Sigma(\rho)$ do not bound any disk in
$\Sigma$ and at least one of them is not isotopic to a boundary
component of $\Sigma$.
\end{itemize}
\end{prop}
\bigskip

\DEM

\emph{Let us show item (i).}
\smallskip

Let $\theta$ be the homomorphism induced by $\tilde\rho$ to $\Mod(\Sigma(\tilde\rho),\,\bord\Sigma(\tilde\rho))$. According to
Theorem \ref{thm:Perron_Vannier_1}, $\theta$ is injective. Let
$\iota$ be the inclusion of $\Sigma(\tilde\rho)$ in $\Sigma$ and
$\iota_*$ the homomorphism induced, going from
$\Mod(\Sigma(\tilde\rho),\,\bord\Sigma(\tilde\rho))$ into
$\Mod(\Sigma,\,\bord\Sigma)$, so that
$\tilde\rho=\iota_*\rond\theta$.
\begin{itemize}
\smallskip

\item[\point]
\emph{Necessary condition.} If $\tilde\rho$ is injective, then
$\tilde\rho(\Delta_n^{\,4})$ is not trivial. However when $n$ is odd,
$\tilde\rho(\Delta_n^{\,4})$ coincides with $T_d^{\pm1}$ where $d$ is the
unique boundary component of $\Sigma(\tilde\rho)$, hence $T_d$ must
be nontrivial. In other words, $d\in\Courb(\Sigma,\,\bord\Sigma)$.
When $n$ is even, $\tilde\rho(\Delta_n^{\,4})$ coincides with
$(T_{d_1}T_{d_2})^{\pm2}$ where $d_1$ and $d_2$ are the two boundary
components of $\Sigma(\tilde\rho)$, so at least one of the curves
$d_1$ or $d_2$ has to be nontrivial. Moreover, if one of them is
trivial in $\Sigma$, say $d_1$ for example (cf. Figure
\ref{fig:injectivite}), then
$\iota_*\rond\theta(\big(\tau_{a_1}\dots
\tau_{a_{n-2}}\big)^{2(n-1)})=\iota_*\rond\theta(\big(\tau_{a_1}\dots
\tau_{a_{n-1}}\big)^{n})$, since $\big(T_{a_1}\dots
T_{a_{n-2}}\big)^{2(n-1)}=\big(T_{a_1}\dots
T_{a_{n-1}}\big)^{n}=T_{d_2}$. But this contradicts the injectivity
of $\tilde\rho$, for in $\B_n$, a product of $n(n-1)$ standard
generators can be equal to a product of $2(n-1)(n-2)$ standard
generators only if $n(n-1)=2(n-1)(n-2)$, hence only if
$n\in\{1,\,4\}$. Therefore,
$\{d_1,\,d_2\}\subset\Courb(\Sigma,\,\bord\Sigma)$.
\begin{figure}[!h]
 \Includegraphics{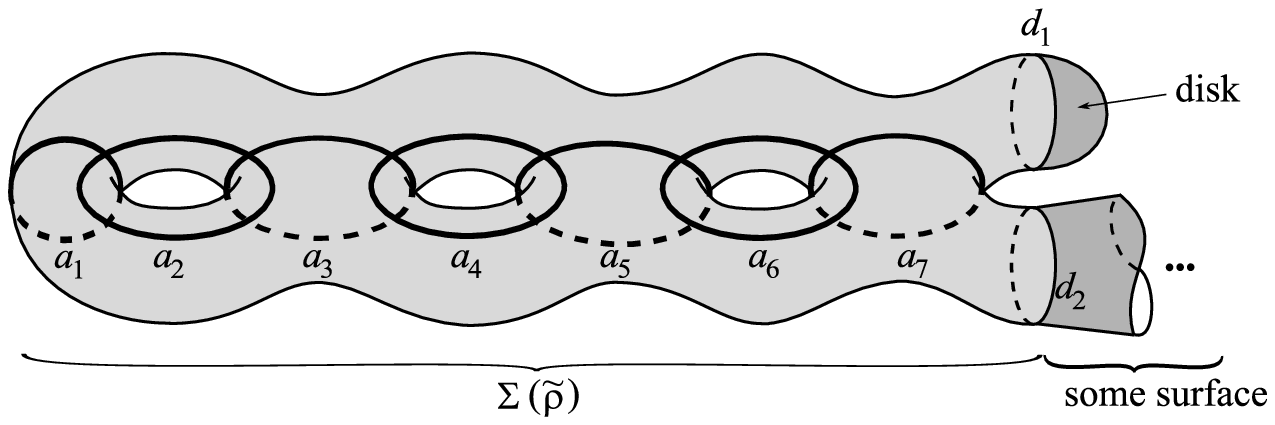}
 \caption{Example with $n=8$.}
 \label{fig:injectivite}
\end{figure}
\smallskip

\item[\point]
\emph{Sufficient condition.} We assume that
$\Bord(\Sigma(\tilde\rho))\subset\Courb(\Sigma,\,\bord\Sigma)$. Then
according to Theorem \ref{thm:inclusionMCG}, in the case where $n$
is odd, or in the one where $n$ is even and where the two boundary
components of $\Sigma(\tilde\rho)$ are not isotopic in $\Sigma$,
$\iota_*$ is injective. In the case where $n$ is even and where the
two boundary components $d_1$ and $d_2$ of $\Sigma(\tilde\rho)$ are
isotopic in $\Sigma$, according to Theorem \ref{thm:inclusionMCG}
again, we have $\Ker(\iota_*)=\langle T_{d_1}T_{d_2}^{\,-1}\rangle$.
Now, according to Theorem \ref{thm:Perron_Vannier_1}, $\theta(\B_n)$
coincides with $\SMod(\Sigma(\tilde\rho),
\,\bord\Sigma(\tilde\rho))$ whereas $T_{d_1}T_{d_2}^{-1}$ does not
belong to $\SMod(\Sigma(\tilde\rho), \,\bord\Sigma(\tilde\rho))$.
Hence $\Ker(\iota_*)\cap\Im(\theta)=\{1\}$, so $\iota_*(\theta)$ is
injective. Finally, in all the cases, $\iota_*\rond\theta$ is
injective, and hence $\tilde\rho$, too.
\end{itemize}
\bigskip

\emph{Let us show item (ii).}
\smallskip

Again, let $\theta$ be the homomorphism induced by $\rho$ to $\Mod(\Sigma(\rho),\,\bord\Sigma(\rho))$. According to Theorem
\ref{thm:Perron_Vannier_1}, $\theta$ is injective. Let $\iota$ be
the inclusion of $\Sigma(\rho)$ in $\Sigma$ and $\iota_*$ the
induced homomorphism from $\Mod(\Sigma(\rho),\,\bord\Sigma(\rho))$ to $\PMod(\Sigma)$, so that $\rho=\iota_*\rond\theta$. The homomorphism
$\iota_*$ is not necessarily injective.
\smallskip

\begin{itemize}

\item[\point]
\emph{Necessary condition.} As in the case of item (i) with
$\tilde\rho$, it is necessary that $\Bord(\Sigma(\rho))\subset$
$\Courb(\Sigma,\,\bord\Sigma)$, but since the Dehn twists along
boundary components are trivial in $\PMod(\Sigma)$, it is necessary
that $\Bord(\Sigma(\rho))\not\subset\Bord(\Sigma)$.
\smallskip

\item[\point] \emph{Sufficient condition.} Let us assume
    that
    $\Bord(\Sigma(\rho))\subset\Courb(\Sigma,\,\bord\Sigma)$
    and that \linebreak
    $\Bord(\Sigma(\rho))\not\subset\Bord(\Sigma)$,
    and let us check that $\rho$ is injective. Let us
    denote by $\Sigma'$ the complement of $\Sigma(\rho)$ in
    $\Sigma$; we assume that if a boundary component of
    $\Sigma(\rho)$ is isotopic to a boundary component of
    $\Sigma$, these two boundary components coincide. With
    this assumption, all the connected components of
    $\Sigma'$ are of negative Euler characteristic. Let
$\bordint\big(\Sigma(\rho)\big)=\bord\big(\Sigma(\rho)\big)\smallsetminus\bord\Sigma$.
Since $\Bord(\Sigma(\rho))\not\subset\Bord(\Sigma)$, the
surface $\Sigma'$ is nonempty and
$\bordint\big(\Sigma(\rho)\big)\not=\vide$. The image of
$\rho$ lies in $\PMod(\Sigma,\,\Sigma')$, which is
isomorphic to
$\Mod\big(\,\Sigma(\rho),\,\bordint\big(\Sigma(\rho)\big)\,\big)$.
If
$\bordint\big(\Sigma(\rho)\big)=\bord\big(\Sigma(\rho)\big)$,
Theorem \ref{thm:Perron_Vannier_1} can be applied and
$\rho$ is then injective. This is always what happens in
the case when $n$ is odd, but when $n$ is even, it can
happen that
$\bordint\big(\Sigma(\rho)\big)=\{d\}\not=\bord\big(\Sigma(\rho)\big)$,
where $d$ is one of the two boundary components of
$\Sigma(\rho)$. In this case, $\PMod(\Sigma,\,\Sigma')$ is
isomorphic to $\Mod(\Sigma(\rho),\,d)$, and $\rho$ induces
a homomorphism $\varsigma$ from $\B_n$ to $\Mod(\Sigma(\rho),\,d)$, which is injective if and only if
$\rho$ is injective itself. Let us denote by $\pr$ the
canonical projection of
$\Mod(\Sigma(\rho),\,\bord\Sigma(\rho))$ to
$\Mod(\Sigma(\rho),\,d)$. Then $\varsigma=\pr\rond\theta$.
Moreover, we have $\Ker(\varsigma) =
\Ker(\pr)\cap\theta(\B_n) =\{1\}$, for $\Ker(\pr)=\langle
T_{d'}\rangle$ where $d'$ is the boundary component of
$\Sigma(\rho)$ different from $d$, but $T_{d'}$ does not
belong to $\theta(\B_n)$. Hence $\varsigma$ is injective,
so $\rho$ is injective.\fin \end{itemize}
\bigskip

\begin{prop}[Injectivity of the transvections]
          \label{prop:transvection_of injection}
\mbox{}\\ \indent (i) Let $n$ be an integer greater than or equal to
3, $G$ any group, $\rho$ a homomorphism from $\B_n$ to $G$ and $\rho_1$
a transvection of $\rho$. If $G$ is torsion-free (for example, if
$G=\Mod(\Sigma,\,\bord\Sigma)$ with $\bord\Sigma\not=\vide$), then
$\rho$ is injective if and only if $\rho_1$ is injective.

(ii) Let $n$ be an integer greater than or equal to 6 and $\Sigma$ a
surface of genus $g\leqslant \frac{n}{2}$. Then for any monodromy
homomorphism $\rho$ from $\B_n$ to $\PMod(\Sigma)$, any transvection of
$\rho$ is injective if and only if $\rho$ is injective. \end{prop}
\medskip

\DEM
Let $\rho$ be an injective homomorphism from $\B_n$ to $G$ and let
$\rho_1$ be a transvection of $\rho$. We show by contradiction that
$\rho_1$ is injective, then we will have shown item (i) of
Proposition \ref{prop:transvection_of injection}, for $\rho$ can
also be seen as a transvection of $\rho_1$. Let us then assume that
there exists $\xi$, a nontrivial element of $\B_n$ such that
$\rho_1(\xi)=1$. Let $\varphi\;:\; \B_n\to G$ the cyclic homomorphism associated
to the transvection $\rho_1$ such that for all
$\zeta\in\B_n$, we have $\rho_1(\zeta)=\rho(\zeta)\varphi(\zeta)$.
Hence, we have $\rho(\xi)=\varphi(\xi)^{-1}$,
but since $\varphi(\xi)$ lies in the centralizer of $\rho(\B_n)$ in $G$,
so does $\rho(\xi)$. But $\rho$ is injective, so $\xi$
belongs to the center of $\B_n$, that is, $\xi$ is a power of $\Delta_n^{\,2}$.
Therefore, there exists a nonzero integer $k$ such that
\smallskip

\centrer[1]{$\big(\rho(\Delta_n)\varphi(\Delta_n)\big)^{2k}=1$.}
\smallskip

\noindent Now, if $G$ is torsion-free, we have
$\rho(\Delta_n)\varphi(\Delta_n)=\rho_1(\Delta_n)= 1$.
Then, as above, $\varphi (\Delta_n)$ lies in the
centralizer of $\rho(\B_n)$ in $G$, hence so does $\rho(\Delta_n)$.
Now, since $\rho$ is injective, $\Delta_n$ must belong to the center
of $\B_n$, hence must be a power of $\Delta_n^{\,2}$, which is absurd.
This proves item (i) of Proposition \ref{prop:transvection_of injection}.
\bigskip
\medskip

\Point \emph{Let us show the direct implication of item (ii): if $\rho$ is an injective monodromy homomorphism, then
the transvection $\rho_1$ of $\rho$ is injective.}
\medskip

We are going to show by contradiction that $\rho_1$ is injective.
As above, if $\rho_1$ is not injective, there exists a
nonzero integer $k$ such that (1) holds.
Now, $\Delta_n$ is a
product of $\frac{n(n-1)}{2}$ generators, so $\varphi(\Delta_n^{\,2})$ is
the $n(n-1)$-th power of a mapping class $V$ belonging to the
centralizer of $\rho(\B_n)$ in $\PMod(\Sigma)$. Let us denote by $W$
the mapping class $V^2$. According
to Proposition \ref{prop:centralizer}, $W$ lies in $\M^{\Sigma(\rho)}$,
where $\M^{\Sigma(\rho)}$ is the group of the mapping classes in $\PMod(\Sigma)$
that preserve $\Sigma(\rho)$ and that induce the identity in $\PMod\big(\Sigma(\rho)\big)$.
So we have $\varphi(\Delta_n^{\,2})=W^{\frac{n(n-1)}{2}}$, and hence, with (1), we get:
\smallskip

\centrer[2]{$W^{k\frac{n(n-1)}{2}}=\rho(\Delta_n^{\,-2k})$.}
\smallskip

\noindent Now, if $n$ is odd, then $\rho(\Delta_n^{\,4})=T_d^{\,\pm1}$ where
$d$ is the unique boundary component of $\Sigma(\rho)$ , whereas if
$n$ is even, then $\rho(\Delta_n^{\,2})=\big(T_{d_1}T_{d_2}\big)^{\,\pm1}$
where $d_1$ and $d_2$ are the two boundary components of
$\Sigma(\rho)$. So the mapping class $W$ lying in
$\M^{\Sigma(\rho)}$ satisfies:
\smallskip

$\left\{\begin{array}{llll}
W^{kn(n-1)}=T_d^{\mp k}, & \mbox{if $n$ is odd,} & \mbox{where }\{d\}=\Bord(\Sigma(\rho));\\
W^{\frac{kn(n-1)}{2}}=\big(T_{d_1}T_{d_2}\big)^{\mp k}, & \mbox{if $n$ is
even,} & \mbox{where }\{d_1,\,d_2\}=\Bord(\Sigma(\rho)).\vblanc
\end{array}\right.$
\hfill
\begin{tabular}{l}
\makebox[0cm]{(3)}\\
\makebox[0cm]{(4)}
\end{tabular}
\smallskip

\noindent Let us recall that since $\rho$ is injective by
assumption, then according to Proposition
\ref{prop:injectivity_des_homomorphisms_of_monodromie}, at least one of
the boundary components of $\Sigma(\rho)$ is not trivial in
$\Sigma$, hence the mapping classes $T_d$ in case a) and
$T_{d_1}T_{d_2}$ in case b) are not trivial.
\medskip

Let us that (3) leads to a contradiction. The curve $d$ is
a separating curve of $\Sigma$. Let us call $\Sigma'$ the connected
component of $\Sigma_d$ different from $\Sigma(\rho)$. According to
Proposition \ref{prop:links_between_MCGs}.(i),
$\PMod\big(\Sigma,\,\Sigma(\rho)\big)$ is isomorphic to
$\PMod(\Sigma',\,d)$. The mapping class $W$ can thus be seen as a
periodic mapping class in $\PMod(\Sigma',\,d)$. Let us call $m$ its
period. According to Lemma $\ref{lem:without_torsion}$, there exists
an integer $p$ coprime with $m$ such that $W^m=T_d^{\,p}$, so the
equality $W^{kn(n-1)}=T_d^{\mp k}$ implies $W^{n(n-1)}=T_d^{\mp 1}$.
Therefore the period of $W$ is $n(n-1)$, so it is greater than or
equal to 42 since $n\geqslant 7$. But $\Sigma(\rho)$ is of genus
$\frac{n}{2}-1$ and $\Sigma$ is of genus at most $\frac{n}{2}$ hence
$\Sigma'$ is of genus at most 1. But there does not exist any
nontrivial periodic mapping class in a genus-0 surface whose
boundary components are fixed, according to Corollary
\ref{cor:periodic_we_a_sphere}, and the order of periodic mapping
classes on a genus-1 surface with a nonempty boundary and whose
boundary components are fixed is bounded by 6, according to
Corollary \ref{cor:order_periodic_à_boundary_genre_1}. This is a
contradiction.
\smallskip

Let us show that (4) leads also to a contradiction. If $\Sigma$
is the gluing of $\Sigma(\rho)$ on itself by identifying both of its
boundary components $d_1$ and $d_2$ and if we call $d$ the image of
$d_1$ in $\Sigma$, then $\M^{\Sigma(\rho)}$ is the cyclic group
spanned by $T_{d}$, so there exists an integer $m$ such that
$W=T_{d}^{\,m}$. But on the other hand,
$W^{\frac{kn(n-1)}{2}}=\big(T_{d_1}T_{d_2}\big)^{\mp k}=T_{d}^{\mp
2k}$, whence $(T_{d}^m)^{\frac{kn(n-1)}{2}}=T_{d}^{\mp2k}$. This is
absurd, for $T_{d}$ is not a torsion element. Now, if we assume that
$d_1$ and $d_2$ are two distinct curves in
$\Courb(\Sigma,\,\bord\Sigma)$, we know that at least one  is not a
boundary component of $\Sigma$. Notice that $V$ used to preserve
$\{d_1,\,d_2\}$, so $W$ (which is equal to $V^2$) preserves $d_1$
and $d_2$. Finally, as previously in case a), $W$ can be seen as a
periodic mapping class of period at least 15 (for $\frac{n(n-1)}{2}$
equals at least 15 when $n\geqslant 6$), on a surface of genus zero
or one. As explained above, this is absurd.
\bigskip
\medskip

\Point \emph{Let us show the reverse implication of item (ii): if $\rho$ is a non-injective monodromy homomorphism, then
the transvection $\rho_1$ of $\rho$ is not injective either.}
\medskip

Let $\rho$ be a non-injective monodromy homomorphism from $\B_n$ to $\PMod(\Sigma)$ and let $\rho_1$ be a transvection of $\rho$. Let us
show that $\rho_1$ is not injective. According to Proposition
\ref{prop:injectivity_des_homomorphisms_of_monodromie}.(ii), it can
exist several reasons for $\rho$ not being injective. We distinguish
two cases, whether $\Bord(\Sigma(\rho))\cap\Courb(\Sigma)$ is empty
or not.
\smallskip

\begin{itemize}
\item[] If $\Bord(\Sigma(\rho))\cap\Courb(\Sigma)$ is empty,
then $\M^{\Sigma(\rho)}$ is trivial. Hence according to Proposition
\ref{prop:centralizer}, when $n$ is odd, the centralizer of
$\rho(\B_n)$ is spanned by $\rho(\Delta_n^{\,2})$ which is of order 2, and
when $n$ is even, this centralizer is trivial. So, in both cases, we
have $\rho_1(\Delta_n^{\,4})=1$, and $\rho_1$ is not injective.
\smallskip

\item[] If $\Bord(\Sigma(\rho))\cap\Courb(\Sigma)$ is not empty,
whereas $\rho$ is not injective, then necessarily, $n$ is even and
the boundary component of $\Bord(\Sigma(\rho))$ that is not in
$\Courb(\Sigma)$ bounds  a disk. In this case, we are going to
exhibit two elements of $\B_n$ that do not commute although their
images do. This will show that $\rho_1$ is not injective. Let us set
$\Delta_{n-2}=$
$\tau_1(\tau_2\tau_1)\dots(\tau_{n-2}\tau_{n-3}\dots\tau_1)$. We
have seen in the proof of Proposition
\ref{prop:injectivity_des_homomorphisms_of_monodromie} that in our
situation, we have $\rho(\Delta_n^{\,2})=\rho(\Delta_{n-2}^{\ \ 4})$ (cf.
Figure \ref{fig:injectivite}). Hence in particular
$\rho(\Delta_{n-2}^{\ \ 4})$ commutes with $\rho(\tau_{n-1})$, so
$\rho_1(\Delta_{n-2}^{\ \ 4})$ commutes with $\rho_1(\tau_{n-1})$.
But $\Delta_{n-2}^{\ \ 4}$ do not commute with $\tau_{n-1}$.
Therefore, $\rho_1$ is not injective. \fin \end{itemize}
\bigskip

Now,  Theorem \ref{thm:injectivity_des_homomorphisms}  is a direct corollary from
Theorems \ref{thm:Theorem_principal} and
\ref{thm:Theorem_principal_case_à_boundary}, and from Propositions
\ref{prop:injectivity_des_homomorphisms_of_monodromie} and
\ref{prop:transvection_of injection}.
\bigskip

\pageimpaire\addcontentsline{toc}{part}{\protect\textsc{\textbf{Proof of Theorem \ref{thm:Theorem_principal}}}}

\part*{Proof of Theorem \ref{thm:Theorem_principal}}
\bigskip
\label{titreIII}

We turn now to the proof Theorem \ref{thm:Theorem_principal}.
The proof ends in Section \ref{sec:théorèmes}.
We first prove it when $n$ is an even number.
We will deduce from it in Subsection \ref{par:théorèmes} the case when $n$ is odd.
Throughout Sections \ref{sec:curves_peripheral} - \ref{sec:théorèmes},
we will introduce the following pieces of notation.
We gather them here for later reference.
\medskip

\pageimpaire\addcontentsline{toc}{subsection}{\protect{Summary of the upcoming notation}}

\LTITRE{Summary of the upcoming notation}

\begin{itemize}
  \item[\point]
$n$ is an \emph{even} integer greater than or equal to 6,

$g$ and $b$ are nonnegative integers such that
$2-2g-b\leqslant -1$ and $g\leqslant \frac{n}{2}$,

$\rho$ is a given homomorphism from $\B_n$ to $\PMod(\Sigma_{g,\,b})$,

from Section \ref{sec:curves_spéciales} on, $\rho$ will be assumed to be noncyclic;
\smallskip

  \item[\point]
$\tau_i$ where $i\in\{1,\,2,\,\dots,\,n-1\}$ is the $i$\eme standard
generator of $\B_n$.

$\Delta_n$ is the \emph{Garside-element} of $\B_n$, and
$\delta$ is the ``$\frac{1}{n}$-flip'' of $\B_n$. They are defined by:
\smallskip

\centrer{$\begin{array}{lll}
\Delta_n &=&
\tau_1(\tau_2\tau_1)(\tau_3\tau_2\tau_1)\dots(\tau_{n-1}\tau_{n-2}\dots
\tau_1),\\
\delta&=&\tau_1\tau_2\dots\tau_{n-1},
\end{array}$}

\noindent and satisfy the following well-known properties:
\begin{itemize}
  \item[\point] for all integers $i\in\{1,\,2,\dots,\,n-1\}$, we have $\Delta_n\, \tau_i = \tau_{n-i}\, \Delta_n$;
  \item[\point] for all integers $i\in\{1,\,2,\dots,\,n-2\}$, we have $\delta\, \tau_i = \tau_{i+1}\, \delta$,
                        and $\delta^2\, \tau_{n-1} = \tau_{1}\, \delta^2$;
  \item[\point] the element $\Delta_n^{\,2} = \delta^n$ span the center of $\B_n$.
\end{itemize}

\noindent In section \ref{sec:n_th_generator}, we define an $n$\eme ``standard'' generator of $\B_n$,
namely $\tau_0=\delta\tau_{n-1}\delta^{-1}$.
When the subscript $i$ satisfies $i<0$ or $i>n-1$,
$\tau_i$ is the $j$\eme standard generator of $\B_n$ where $j$ is
the remainder  of the euclidian division of $i$ by $n$;

\smallskip

  \item[\point]
$A_i$ is defined by  $A_i=\rho(\tau_i)$ for all integers $i$,

$\G$ is the set $\{A_1,\,A_2,\dots,\,A_{n-1}\}$ and $\G_0$ is the set $\{A_0,\,A_1,\dots,\,A_{n-1}\}$,

$\Imp(n)$ is the set $\{1,\,3,\dots,\,n-1\}$ of the $\frac{n}{2}$ first odd integers,

$\X$ the set $\{A_i,\,i\in\Imp(n)\}$ and $\Y$ the set $\{A_i,\,i\in\{1,\,2,\dots,\,n-1\}\smallsetminus\Imp(n)\}$;
\smallskip

  \item[\point]
$\sigma(\G_0)$ is the union $\dessous{\cup}{i\leqslant
n-1} \sigma(A_i)$ and $\sigma(\X)$ is the union
$\dessous{\cup}{i\in \Imp(n-1)} \sigma(A_i)$.

  \item[\point]
We have a natural action from $\B_n$ on $\PMod(\Sigma)$ by conjugation via $\rho$ given by:
\smallskip

\centrer{$\xi.A=\rho(\xi)A\rho(\xi)^{-1}$.}
\smallskip

\noindent We will introduce the subgroups
$\J=\langle\delta\rangle$ and $\H=\langle\,\gamma_i\,,\ i\in\Imp(n)\,\rangle$ of $\B_n$
where $\gamma_i=\tau_i\tau_{i+1}\tau_i\tau_{i+2}\tau_{i+1}\tau_i$ for all $i\in\Imp(n)$.
Then the action of $\B_n$ on $\PMod(\Sigma)$ induces an action of $\J$ on $\G_0$ and an action of $\H$ on $\X$.
\end{itemize}
\bigskip

\pageimpaire\addcontentsline{toc}{subsection}
{\protect{Outline of the proof of Theorem \ref{thm:Theorem_principal}}}

\LTITRE{Outline of the proof of Theorem \ref{thm:Theorem_principal}}
\bigskip

\blanc  \TITRE{Section \ref{sec:curves_peripheral}:}

\noindent A curve $a$ of $\sigma(\G)$ is a \emph{peripheral curve}
if $a$ is a separating curve with the following property: one of the
connected components of $\Sigma_a$ is of genus 0. Let $\sigma_p(\G)$
be the set of \emph{peripheral curves}. We show that up to
transvection, we can assume that $\sigma_p(\G)$ is empty, which is a
first way to simplify the study of $\sigma(\G)$. From Section
\ref{sec:homomorphisms_irreducible} on, we will assume that
$\sigma_p(\G)$ is empty. We also show in section
\ref{sec:curves_peripheral} that in many cases, we can assume
without loss of generality that $\Sigma$ is a surface without
boundary.
\bigskip

\blanc  \TITRE{Section \ref{sec:homomorphisms_irreducible}:}

\noindent We show that, although $\sigma_p(\G)$ is assumed to be
empty, if $\rho$ is not cyclic, then $\sigma(\G)$ is nonempty. From
Section \ref{sec:curves_spéciales} on, we will assume that $\rho$ is
not cyclic, so that $\sigma(\G)$ is nonempty.
\bigskip

\blanc  \TITRE{Section \ref{par:group_of_tresses_Artin}:}

\noindent We introduce an $n$\eme generator $\tau_0=\delta\tau_{n-1}\delta^{-1}$ to the standard presentation of
the braid group in order to get a cyclic action of the subgroup of $\B_n$ spanned by $\delta$
on the set of the generators $\{\tau_0,\,\tau_1,\dots,\,\tau_{n-1}\}$. This action will be fundamental
in the next sections.
\bigskip

\blanc  \TITRE{Section \ref{sec:curves_spéciales}:}

\noindent We prove that there exists a partition of $\sigma(\G_0)$ in
two sets of curves: $\sigma_s(\G_0)$ and $\sigma_n(\G_0)$, both of them
satisfying each interesting properties. For instance, for any curve
$a$ of $\sigma_s(\G_0)$, there exists a unique
$i\in\{0,\,1,\dots,\,n-1\}$ such that $a\in\sigma(A_i)$. As for the
set of curves $\sigma_n(\G_0)$, it is stable by the action of $\B_n$
induced by $\rho$ on $\Courb(\Sigma)$.
\bigskip

\blanc  \TITRE{Section
\ref{sec:the_curves_spéciales_are_non-separating}:}

\noindent We show that $\sigma(\G_0)$ contains only non-separating
curves.
\bigskip

\blanc  \TITRE{Section \ref{sec:sigmaX}:}

\noindent We describe the set of curves $\sigma(\X)$ in the surface
$\Sigma$.
\bigskip

\blanc  \TITRE{Section \ref{sec:end_grande_proof}:}

\noindent We gather the results of the previous sections and we show
that $\rho$ is a transvection of monodromy homomorphism: Theorem
\ref{thm:Theorem_principal} when $n$ is even is shown. Finally, we deduce from it
Theorem \ref{thm:Theorem_principal} in the general case.
\bigskip

\section{On geometric representations of $\B_n$ in nonempty-boundary surfaces}
        \label{sec:curves_peripheral}
\bigskip

This section is divided into two parts.
\smallskip

\Tiret In Subsection \ref{par:BndansMCG0b}, we are interested by homomorphisms
from the braid group in the mapping class group of genus-0 surfaces.
We are going to prove that all homomorphisms from $\B_n$ to $\PMod(\Sigma_{0,\,b})$
are cyclic, see Theorem \ref{thm:BndansMCG0b}.
\smallskip

\Tiret In Subsection \ref{par:curves_peripheral}, we are interested by
homomorphisms from the braid group to the mapping class group of a
surface of genus $g\geqslant 1$ with $b\geqslant 2$ boundary
components. Given such a homomorphism $\rho$, we focus on some special separating curves
related to $\rho$ which will be called \emph{peripheral curves}.
Three propositions will be useful for the remainder of the paper, namely:\\
\tiret Proposition \ref{prop:stability_curves_peripheral} (on the stability of peripheral curves),\\
\tiret Proposition \ref{prop:pas_of_curve_peripheral} (on ``getting rid of'' peripheral curves),\\
\tiret Proposition \ref{prop:as_if_surface_without_boundary} (on ``getting rid of'' the boundary).\\
The two last propositions
will follows from the first one. Moreover, Proposition
\ref{prop:as_if_surface_without_boundary} utilizes Theorem
\ref{thm:BndansMCG0b}.
\bigskip
















\subsection{Geometric representations of $\B_n$ in genus-0 surfaces are cyclic}
                  \label{par:BndansMCG0b}
\medskip

We are going to use
the fact that some mapping class groups are \emph{bi-orderable}.
\medskip

\begin{defi}[Bi-orderable group]
                  \label{defi:group_bi-orderable}
A group $G$ is \emph{bi-orderable} if there exists a linear ordering
$\leqslant$ on $G$ invariant by left and right multiplications
(namely, if $f \leqslant g$, then $h_1fh_2 \leqslant h_1gh_2$ for
all $f,\,g,\,h_1,\,h_2\in G$). In what follows, we will denote by
$\leqslant$ the ordering of all the bi-orderable groups that we are
going to meet, and by $<$ the strict order associated to
$\leqslant$.
\end{defi}
\medskip

\begin{prop}
          \label{prop:bi-orderable_cyclic}
Any homomorphism from $\B_n$ in a bi-orderable group is cyclic.
\end{prop}
\smallskip

\DEM Let $G$ be an orderable group and let $\varphi$ be a homomorphism
from $B_n$ to $G$. Let us assume that
$\varphi(\tau_1)<\varphi(\tau_2)$. Let $\gamma$ be the element
$\tau_1\tau_2\tau_1$. Then
$\varphi(\gamma\tau_1\gamma^{-1})<\varphi(\gamma\tau_2\gamma^{-1})$.
Since $\gamma\tau_1\gamma^{-1}=\tau_2$ and
$\gamma\tau_2\gamma^{-1}=\tau_1$ we have
$\varphi(\tau_2)<\varphi(\tau_1)$, which is absurd. In the same way,
assuming that $\varphi(\tau_2)<\varphi(\tau_1)$ leads to a
contradiction. Hence $\varphi(\tau_1)=\varphi(\tau_2)$, so $\varphi$
is cyclic. \fin
\bigskip

Thanks to Proposition \ref{prop:ModS_bi-orderable} below, we will be
able to apply Proposition \ref{prop:bi-orderable_cyclic} to the
mapping class groups and so prove Theorem \ref{thm:BndansMCG0b}.
\smallskip

\begin{prop}[Bonatti,\,Paris]
        \label{prop:ModS_bi-orderable}
\mbox{}\\
For any genus-0 surface $\Sigma$, the mapping class group
$\Mod(\Sigma,\,\bord\Sigma)$ is bi-orderable.\fin
\end{prop}
\smallskip

As a corollary, we get:
\smallskip

\begin{thm}[Homomorphisms from $\B_n$ to $\Mod(\Sigma_{0,\,b}, \bord\Sigma_{0,\,b})$ and to $\PMod(\Sigma_{0,\,b})$]
                \label{thm:BndansMCG0b}
\mbox{}\\Let $\Sigma$ be a genus-0 surface. For all integers
$n$ greater than or equal to 3, any homomorphism from $\B_n$ to
$\Mod(\Sigma,\,\bord\Sigma)$, respectively to $\PMod(\Sigma)$,
is cyclic. \end{thm}
\bigskip

\DEM Any homomorphism from $\B_n$ to $\Mod(\Sigma,\,\bord\Sigma)$ is
cyclic, for according to Proposition \ref{prop:ModS_bi-orderable},
$\Mod(\Sigma,\,\bord\Sigma)$ is bi-orderable. As for the homomorphisms
from $\B_n$ to $\PMod(\Sigma)$, according to Proposition
\ref{prop:passage_PMod_à_ModBord}, they can be lifted in homomorphisms
from $\B_n$ to $\Mod(\Sigma,\,\bord\Sigma)$, which are cyclic. Hence
again, according to Proposition \ref{prop:passage_PMod_à_ModBord},
the homomorphisms from $\B_n$ to $\PMod(\Sigma)$ are cyclic.\fin
\bigskip

\subsection{Peripheral curves}
                  \label{par:curves_peripheral}
\medskip

\TITRE{Starting point for the next sections.}

\tiret Let $n$ be an even number greater than or equal to 6, let $\Sigma=\Sigma_{g,\,b}$ where $g\geqslant 1$ and
$b\geqslant 2$.

\tiret Let $\rho$ be a homomorphism from $\B_n$ to $\PMod(\Sigma)$.

\tiret For all integers $i\leqslant n-1$, $\rho(\tau_i)$ will be denoted by $A_i$.

\tiret We denote by $\G$ the set $\{A_1,\,A_2,\dots,\,A_{n-1}\}$.

\noindent We aim to get as much information as possible concerning $\sigma(A)$ for all $A\in\G$.
\bigskip

\begin{figure}[!h]
 \Includegraphics{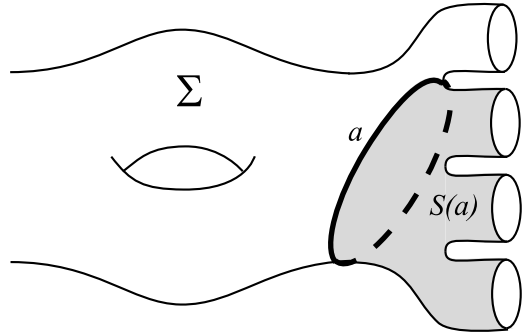}
 \caption{Example of peripheral curve.}
 \label{fig:courbeSep}
\end{figure}

\TITRE{Peripheral curves, $\sigma_p(\G)$.}
\mbox{}\\A curve $a$ of $\sigma(\G)$ is said to be
\emph{peripheral} if it separates $\Sigma$ in two connected
components and if the genus of one of them is zero (cf. Figure
\ref{fig:courbeSep}). The set of peripheral curves will be denoted
by $\sigma_p(\G)$. We will denote by $\sigma_p(A)$ the set of curves
$\sigma_p(\G)\cap\sigma(A)$.
\medskip

\TITRE{The mapping class group $\P_\A\Mod(\Sigma)$}
\mbox{}\\
Remember that if $\A$ is a curve simplex in $\Curv(\Sigma)$, we define
$\P_\A\Mod(\Sigma)$ as the subgroup of $\PMod(\Sigma)$
consisting on all mapping classes that fixes each curve of $\A$
and that preserve each subsurface of $\Sub_{\A}(\Sigma)$). If
$\A=\{a\}$, then we write $\P_a\Mod(\Sigma)$ instead of
$\P_{\{a\}}\Mod(\Sigma)$.
\medskip

\begin{prop}[Stability of peripheral curves]
          \label{prop:stability_curves_peripheral}
\mbox{}
\begin{itemize}
\item[\;(i)] We have the equalities: $\sigma_p(A_1)=\sigma_p(A_2)=\dots=\sigma_p(A_{n-1})=\sigma_p(\G)$.
\item[(ii)] We have the inclusion:
    $\rho(\B_n)\subset\P_{\sigma_p(\G)}\Mod(\Sigma)$.
    \end{itemize}
\end{prop}
\medskip

\begin{lem}
        \label{lem:stability_curves_peripheral}
Let $x$ be a peripheral curve. If a mapping class $F$ in
$\PMod(\Sigma)$ verifies $I\big(F(x),\,x\big)=0$, then $F$ lies
in $\P_x\Mod(\Sigma)$.
\end{lem}

\DEM[of Lemma \ref{lem:stability_curves_peripheral}] Let $S$ be the
holed sphere bounded by $x$. Since $I\big(F(x),\,x\big)=0$, if
$F(x)$ is distinct from $x$, then either $F(x)$ is in $S$ and $F(S)$
is included in $S$, or $F(x)$ is outside of $S$ and $S$ is included
in $F(S)$ (recall that
$\Bord(F(S))\cap\Bord(\Sigma)=\Bord(S)\cap\Bord(\Sigma)$ for $F$
belongs to $\PMod(\Sigma)$). Since $S$ and $F(S)$ are homeomorphic,
these two hypotheses are absurd. Finally, $F(x)=x$. Moreover, $S$
and $F(S)$ are located on the same side of $x$, so $F$ does not swap
the two side-neighbourhoods of $x$.\fin
\medskip

\DEM[of Proposition \ref{prop:stability_curves_peripheral}]

(i) Let $A$ and $C$ be two mapping classes of $\G$ that commute, $x$
a curve of $\sigma_p(A)$ and $Z$ a mapping class such that
$ZAZ^{-1}=C$. Then $Z(x)\in\sigma_p(C)$, and since $AC=CA$, we have
$I(\sigma(A),\,\sigma(C))=0$, so $I(x,\,Z(x))=0$. Now, according to
Lemma \ref{lem:stability_curves_peripheral} which we apply to the
curve $x$ and the mapping class $Z$, we get $Z(x)=x$, and hence
$x\in\sigma(C)$. This shows that for all
$i,j\in\{1,\,2,\dots,\,n-1\}$ such that $|i-j|_n>1$, we have
$\sigma_p(A_i)=\sigma_p(A_j)$. We then easily can deduce that
\smallskip

\centrer{$\sigma_p(A_1)=\sigma_p(A_2)=\dots=\sigma_p(A_{n-1})=\sigma_p(\G)$.}
\smallskip

(ii) For all $x\in\sigma_p(\G)$, all $A\in\G$, we have just seen
that $x\in\sigma(A)$, so $I(x,\,A(x))=0$, hence according to Lemma
\ref{lem:stability_curves_peripheral}, $A(x)=x$ and $A$ does not
swap the two connected components of $\Sigma_x$. Hence $A$ belongs
to $\P_x\Mod(\Sigma)$. Since $\G$ span $\rho(\B_n)$, this proves the
second part of Proposition
\ref{prop:stability_curves_peripheral}.\fin
\bigskip

\begin{prop}[Killing the peripheral curves]
    \label{prop:pas_of_curve_peripheral}
\mbox{}\\Let $\Sigma'$ be the connected component of nonzero genus
of $\Sigma_{\sigma_p(\G)}$. Then:
\begin{itemize}
        \item[\;(i)]
For all $\xi\in\B_n$, $\rho(\xi)$ induced a mapping class in
$\PMod(\Sigma')$ that we denote by $\rho'(\xi)$. The obtained map
$\rho'\;:\; \B_n\to\PMod(\Sigma')$ is a homomorphism.
\smallskip

        \item[(ii)]
The homomorphisms $\rho$ and $\rho'$ are of the same nature: one is
cyclic (respectively is a transvection of monodromy homomorphism) if and
only if the other is.
\end{itemize}
\end{prop}
\medskip

\DEM

(i) Let $\Sigma'$ be the connected component of nonzero genus of
$\Sigma_{\sigma_p(\G)}$ and let $\U$ be the subset of curves of
$\sigma_p(\G)$ that bound the subsurface of $\Sigma$ isomorphic to
$\Sigma'$. According to Proposition
\ref{prop:stability_curves_peripheral} on the stability of the
peripheral curves, $\rho(\B_n)$ is included in $\P_\U\Mod(\Sigma)$.
Let us denote by $\pi'$ the homomorphism from $\PMod(\Sigma_\U)$ to $\PMod(\Sigma')$. Then $\rho'=\pi'\rond \cut_\U\rond\rho$, so
$\rho'$ is indeed a homomorphism.
\smallskip

(ii) According to Proposition
\ref{prop:passage_PMod_à_ModBord}, there exists a lift
$\tilde\rho'$ of $\rho'$ in
$\Hom(\B_n,\,\Mod(\Sigma',\,\bord\Sigma'))$, which is of the
same nature as $\rho'$. For all $\xi\in\B_n$, if we extend the
mapping class $\tilde\rho'(\xi)$ by the identity on $\Sigma$,
and if we then postcompose it by the canonical ``forget''
homomorphism $\for_{\bord\Sigma}\;:\;
\Mod(\Sigma,\,\bord\Sigma)\to\PMod(\Sigma)$, we get a homomorphism
$\rho_1$ from $\B_n$ to $\PMod(\Sigma)$ such that $\pi'\rond
\cut_\U\rond\rho_1=\rho'$. By construction, $\rho_1$ is of the
same nature as $\rho'$.
\smallskip

Let $\Sigma''$ be the union of the subsurfaces of
$\Sub_{\U}(\Sigma)$ distinct from $\Sigma'$. Let us denote by
$\pi''$ the homomorphism from $\PMod(\Sigma_\U)$ to $\PMod(\Sigma'')$. According to Theorem \ref{thm:BndansMCG0b},
$\pi''\rond \cut_\U\rond\rho$ is a cyclic homomorphism. Let $W$ be
the mapping class $\pi''\rond \cut_\U\rond\rho(\tau_1)$ of
$\PMod(\Sigma'')$ and $\widetilde W$ the mapping class
$\rho(\tau_1)\big(\rho_1(\tau_1)\big)^{-1}$ of
$\PMod(\Sigma,\,\Sigma')$. Notice that $\widetilde W$ induces
$W$ on $\PMod(\Sigma'')$. Let $\rho_2$ be the transvection of
$\rho_1$ with direction $\widetilde W$ (i.e. for all integers
$i$ in $\{1,\,2,\dots,\,n-1\}$, we have
$\rho_2(\tau_i)=\rho_1(\tau_i)\widetilde W$). Notice that
$\rho_2$ and $\rho_1$ are of the same nature, so $\rho_2$ and
$\rho'$ are of the same nature.

On the other hand, we have the following central exact sequence:
\smallskip

\centrer{$1\to\langle
T_u,\,u\in\U\rangle\to\P_\U\Mod(\Sigma)\xrightarrow{\cut_\U}\PMod(\Sigma_\U)\to1$.}
\smallskip

\noindent Since $\cut_\U\rond\rho_2=\cut_\U\rond\rho$, it comes
that, according to Lemma \ref{lem:transvections}, $\rho_2$ is a
transvection of $\rho$, hence $\rho_2$ and $\rho$ are of the same
nature.
\smallskip

Finally, $\rho'$ and $\rho$ are of the same nature. \fin
\bigskip

\begin{defi}[The ``squeeze map'' of a surface with a nonempty boundary]
          \label{defi:projection_canonical_without_boundary}
\mbox{}\\ Starting from the surface $\Sigma_{g,\,b}$ with $b>0$, let
$\Sigma_{g,\,0}$ be the surface without boundary obtained from
$\Sigma_{g,\,b}$ by squeezing each boundary component to a point. We
get a surface which exceptionally can be a sphere or a torus. There
is a canonical surjective continuous map from $\Sigma_{g,\,b}$ to
$\Sigma_{g,\,0}$ that we will denote by $\sq\;:\;
\Sigma_{g,\,b}\to\Sigma_{g,\,0}$. The map $\sq$ induces a canonical
homomorphism between mapping class groups that we denote by:
\smallskip

\centrer{$\sq^*\;:\; \PMod(\Sigma_{g,\,b})\to\Mod(\Sigma_{g,\,0})$.}
\end{defi}
\medskip

\begin{prop}
          \label{prop:as_if_surface_without_boundary}
Let $\A$ be a curve simplex and $\K$ a subgroup of $\PMod(\Sigma)$
such that $\A$ is $\K$-stable and such that the cardinality of any
curve orbit in $\A$ under the action of $\K$ is at least 3. Then the
canonical surjective continuous map $\sq\;:\;
\Sigma\to\Sigma_{g,\,0}$ induces an
isomorphism from the graph $\Gamma(\Sigma,\,\A)$ (cf. Definition
\ref{defi:SubASigma}) to the graph
$\Gamma(\Sigma_{g,\,0},\,\sq(\A))$. In particular, the cardinality
of $\Sub_\A(\Sigma)$ is smaller than or equal to $2g-2$, and the
cardinality of  $\A$ is smaller than or equal to $3g-3$. Moreover,
for any mapping class $F\in\K$, for any curve $a\in\A$ and for any
subsurface $S\in\Sub_\A(\Sigma)$, we have:
\smallskip

\centrer{$\sq(F(a))=\sq^*(F)(\sq(a))$,}

\centrer{$\sq(F(S))=\sq^*(F)(\sq(S))$.}
\end{prop}
\medskip

\DEM

 1. Let us show that
no subsurface of $\Sub_\A(\Sigma)$ can be sent by $\sq$ on a sphere
minus one or two disks. Hence no curve of $\A$ is sent on a
contractible curve by $\sq$, and for any two curves of $\A$, they
cannot be sent on the same isotopy class in $\Sigma'$. Once we have
shown this, we have shown that the sets $\A$ and $\sq(\A)$ of curves
have the same cardinality.
\smallskip

a) By assumption, the set $\A$ does not contain any fixed point
under the action of $\K$. So $\A$ does not contain any
peripheral curve, according to Lemma
\ref{lem:stability_curves_peripheral}. Hence no subsurface of
$\Sub_\A(\Sigma)$ can be sent by $\sq$ on a sphere minus a
disk. Therefore, no curve of $\A$ is sent on a contractible
curve.
\smallskip

b) Let us show that no subsurface of $\Sub_\A(\Sigma)$ can be
sent by $\sq$ on a sphere minus two disks, which is equivalent
to say that for any two curves of $\A$, they cannot be sent on
the same isotopy class in $\Sigma'$. If there did exist two
distinct curves $a$ and $a'$ of $\A$ such that
$\sq(a)=\sq(a')$, then it would exist in $\Sub_\A(\Sigma)$ a
genus-0 subsurface $S$ whose boundary would consists in some
boundary components included in $\bord\Sigma$ and exactly two
boundary components $a$ and $a'$ that do not belong to
$\Bord(\Sigma)$. But any mapping class of $\PMod(\Sigma)$ which
globally preserves $\A$ should preserve the surface $S$, since
$\Bord(S)\cap\Bord(\Sigma)$ is not empty, and hence should
preserve the pair $\{a,\,a'\}$. This is in contradiction with
our hypotheses, since the cardinality of the orbit of $a$ under
the action of $\K$ must be greater than or equal to 3.
\smallskip

According to a) and b), the cardinalities of the sets $\A$ and
$\sq(\A)$ are equal.
\smallskip

2. Since $\A$ and $\sq(\A)$ have the same cardinality, the map $\sq$
induces a graph isomorphism $\Psi$ from the graph
$\Gamma(\Sigma,\,\A)$ to the graph
$\Gamma(\Sigma_{g,\,0},\,\sq(\A))$. Moreover, as the map $\sq$ and
the homomorphism $\sq^*$ are canonical, the action of $\K$ on
$\Gamma(\Sigma,\,\A)$ induces an action of $\sq^*(\K)$ on
$\Gamma(\Sigma_{g,\,0},\,\sq(\A))$ and the expected commutation
properties hold. \fin
\bigskip
\bigskip

\section{Irreducible geometric representations of $\B_n$}
\label{sec:homomorphisms_irreducible}
\bigskip

\TITRE{Hypotheses.}\\
Let $n\geqslant 6$ an even number, let $\Sigma=\Sigma_{g,\,b}$ with $g\leqslant\frac{n}{2}$, and let
$\rho\ :\ \B_n\to\PMod(\Sigma)$ such that:

\begin{tabular}{lll}
\point & $\rho$ is non-cyclic      & by assumption,\\
\point & $\sigma_p(\G)=\vide$      & by assumption, inspired by Proposition \ref{prop:pas_of_curve_peripheral}.\\
\end{tabular}
\smallskip

\noindent
Remember that for all integers $i$, $\rho(\tau_i)$ will be denoted by $A_i$;
and the set $\{A_1,\,A_2,\dots,\,A_{n-1}\}$ will be denoted by $\G$.
\medskip

\TITRE{Irreducible homomorphisms, periodic homomorphisms, pseudo-Anosov homomorphisms}
\mbox{}\\
In $\B_n$, the standard generators are conjugate,
so their images by $\rho$ are all reducible, all periodic or all pseudo-Anosov.
We will say that $\rho$ is an \emph{irreducible homomorphism
from $\B_n$} if $\rho(\tau_1)$ is an \emph{irreducible mapping class},
that is, if $\sigma\big(\rho(\tau_1)\big)=\vide$. If
$\rho(\tau_1)$ is periodic (respectively pseudo-Anosov), we will say
that $\rho$ is \emph{periodic} (respectively \emph{pseudo-Anosov}).
\medskip

To prove Theorem \ref{thm:Theorem_principal}, we need to prove
that, up to an element in the centralizer of $\G$, the elements
of $\G$ are Dehn Twists. In this purpose, focusing on
$\sigma(\G)$ will be efficient, but we first need to prove that
$\sigma(\G)$ is not empty! This is precisely the aim of this
section, whose main theorem is the following.
\smallskip

\begin{thm}
              \label{thm:homomorphisms_irreducible}
Any irreducible homomorphism from $\B_n$ to $\PMod(\Sigma)$ is cyclic.
\end{thm}
\smallskip

We will distinguish the case of the periodic homomorphisms (cf.
propositions \ref{prop:periodic_implies_cyclic_b>0} and
\ref{prop:periodic_implies_cyclic_b=0}) from the one of the
pseudo-Anosov homomorphisms (cf. Proposition
\ref{prop:pseudo-Anosov_implies_cyclic}). The proof of this theorem
is short when $\Sigma$ has a nonempty boundary but the involved
methods are inefficient when the boundary of $\Sigma$ is empty. When
$\bord\Sigma=\vide$, we argue by contradiction: we assume that
$\rho$ is not cyclic, we exhibit a finite subgroup of $\rho(\B_n)$
and we show that its cardinality exceeds the theoretical maximal
cardinality of a finite subgroup of $\Mod(\Sigma)$.
\smallskip

In the first subsection, we present some results on the relations in
$\rho(\B_n)$, which will be useful to fix a lower bound to the
cardinality of some subgroups of $\rho(\B_n)$. The second and third
subsections (\ref{par:homomorphisms_periodic} and
\ref{par:homomorphisms_Anosov}) are devoted to the proof of Theorem
\ref{thm:homomorphisms_irreducible} in the case of the periodic
homomorphisms and of the pseudo-Anosov homomorphisms, respectively.


\subsection{Cardinalities of some abelian subgroups of $\rho(\B_n)$}
\medskip

\Tiret Let $\varphi$ be a homomorphism from $\B_n$ in some group. For any element $\xi$ in $\B_n$, we
write $\bar\xi$ instead of $\varphi(\xi)$.
\smallskip

\Tiret For all even positive integers positif $N$, let us denote by
$\Imp(N)=\{1,\,3,\dots,\,N-1\}$ the set of odd positive integers
smaller than $N$. Let $r$ be the integer $\frac{n}{2}$. Thus
$\Imp(n)$ contains the $r$ first odd positive integers.
\smallskip

\Tiret Let
$\L_n$ be the free abelian subgroup
$\big\langle\,\tau_i\,,\,i\in\Imp(n)\,\big\rangle_{\B_n}$ of $\B_n$.
The group $\L_n$ is isomorphic to $\ZZ^r$.
\medskip

The aim of this
subsection is to study the algebraic structure of the abelian group
$\varphi(\L_n)$, that is to say to study the structure of the
quotients of $\L_n$, cf. Lemma \ref{lem:relations}, and then to
compute the cardinality of $\varphi(\L_n)$, cf. Lemma
\ref{lem:cardinality}. Let us begin by stating an elementary case that
was already proved, see Lemma \ref{lem:homomorphism_cyclic}.
\smallskip

\begin{lem}
        \label{lem:powers_égales}
If there exist two distinct integers $i$ and $j$ smaller than or
equal to $n-1$ and a nonzero integer $\ell$ such that
$\bar\tau_i^{\,\ell}=\bar\tau_j^{\,\ell}$, then we have:
\smallskip

\Centrer[\Box]{$\bar\tau_1^{\,\ell}= \bar\tau_2^{\,\ell}=\dots=
\bar\tau_{n-1}^{\,\ell}$.}
\end{lem}

\newcommand*{\Run}{\mathrm{R1}}
\newcommand*{\Rdeux}{\mathrm{R2}}
\newcommand*{\Rtrois}{\mathrm{R3}}

\begin{lem}
        \label{lem:relations}
There exist four nonnegative integers $M$, $m$, $d$, $s$ such that
the group $\varphi(\L_n)$ is isomorphic to the quotient of $\L_n$ by
the three following relations:
\smallskip

\Centrer[\Run(M)]{$\tau_1^M=\tau_3^M=\dots=\tau_{n-1}^M=1$,}
\smallskip

\Centrer[\Rdeux(m)]{$\tau_1^m=\tau_3^m=\dots=\tau_{n-1}^m$,}
\smallskip

\Centrer[\Rtrois(d,\,s)]{$(\tau_1\tau_3\dots
\tau_{2r-1})^d=\tau_1^{s}$.}
\smallskip

\noindent When $M$ is nonzero, $m$ and $d$ are also nonzero and the
integers $M$, $m$, $d$, $s$ satisfy the following divisibility
relations:

\centrer{\begin{tabular}{ll}
\Point & $m$ divides $M$;\\
\Point & $d$ divides $m$ and $m$ divides $s$;\\
\Point & $M$ divides $(r-\frac{s}{d})m$.\\
\end{tabular}}
\smallskip

\noindent Finally, $d=1$ if and only if $\varphi$ is cyclic (that is
to say, if $m=1$).
\end{lem}
\medskip

\REM In the following subsection, we will apply this lemma to the
image of $\rho$ where $\rho$ will be a homomorphism from $\B_n$ in a
mapping class group. The integer $M$ will be then the order of
$\rho(\tau_1)$ and the integer $m$ will be the order of
$\rho(\tau_1\tau_3^{-1})$.
\medskip

\DEM

\noindent 0. Before starting, notice that if we choose $M=m=d=s=0$,
the quotient of $\L_n$ by $\big(\,\Run(0), \Rdeux(0),
\Rtrois(0,\,0)\,\big)$ is equal to $\L_n$. \medskip

\noindent 1. \emph{Let us show that any relation that holds in
$\varphi(\L_n)$ is equivalent to a set of relations of type $\Run$,
$\Rdeux$ and $\Rtrois$.}
\smallskip

Since $\L_n$ is abelian, any relation in $\varphi(\L_n)$ can be
written as follows:
\smallskip

\centrer[1]{$\bar\tau_1^{k_1}\bar\tau_3^{k_3}\dots
\bar\tau_{n-1}^{k_{n-1}}=1$,}
\smallskip

\noindent where the $k_i$, $i\in \Imp(n)$, are $r$ not all zero. We
are going to show that relation (1) is equivalent to a set of
relations of type $\Run(M)$, $\Rdeux(m)$, $\Rtrois(d,s)$ where $M$,
$m$, $d$, $s$ are some integers. We distinguish three cases: a), b)
and c) below.
\smallskip

a) If the $k_i$ are all equal, (1) is exactly the relation
$\Rtrois(k_1,\,0)$.
\smallskip

b) Suppose now that the $k_i$ are not all equal. Let us consider the
differences $|k_i-k_j|$, $i, j\in\Imp(n)$. They are not all zero.
Let us assume for example that $k_1-k_3\not=0$. Then after having
conjugated (1) by
$\bar\tau_1\bar\tau_2\bar\tau_1\bar\tau_3\bar\tau_2\bar\tau_1$, we
get:
\smallskip

\centrer[2]{$\bar\tau_3^{k_1}\bar\tau_1^{k_3}\big(\bar\tau_{5}^{k_5}\dots
\bar\tau_{n-1}^{k_{n-1}}\big)=1$.}
\smallskip

\noindent Now, if we compare (1) and (2), we have:
\smallskip

\centrer[3]{$\bar\tau_1^{k_1}\bar\tau_3^{k_3}=\bar\tau_3^{k_1}\bar\tau_1^{k_3}$,}
\smallskip

\noindent whence $\bar\tau_1^{k_1-k_3}=\bar\tau_3^{k_1-k_3}$, and
hence, according to Lemma \ref{lem:powers_égales}:
\smallskip

\centrer[4]{$\bar\tau_1^{|k_1-k_3|}=\bar\tau_2^{|k_1-k_3|}=\dots=\bar\tau_{n-1}^{|k_1-k_3|}$.}

\noindent We repeat this argument for all the pairs
$(i,\,j)\in\Imp(n)^2$, $i\not=j$. Let $p$ be the greatest common
divisor  of $\{|k_i-k_j|\,,\ i,j\in\Imp(n)\}$. We get relation (5),
which is equal to $\Rdeux(p)$:
\smallskip

\centrer[5]{$\bar\tau_1^{p}=\bar\tau_2^{p}=\dots=\bar\tau_{n-1}^{p}$.}
\smallskip

\noindent For all $i\in\Imp(n)$, the euclidian division of $k_i$ by
$p$ provides two integers $q_i$ and $k'_i$ such that
\smallskip

\centrer{$k_i=q_ip+k'_i$ where $0\leqslant k'_i<p$.}
\smallskip

\noindent Since the $k_i$, $i\in\Imp(n)$, differ one from the other
by a multiple of $p$, the $k'_i$, $i\in\Imp(n)$, are all equal. Let
us call $k'$ this integer. Thanks to relation (5), relation (1)
implies:
\smallskip

\centrer[6]{$\bar\tau_1^{k'}\bar\tau_3^{k'}\dots
\bar\tau_{n-1}^{k'}=\bar\tau_1^{\,\big(-p\sum q_i\big)}$.}
\smallskip

\noindent In $\varphi(\L_n)$, relation (6) is equivalent to
R1($-p\sum q_i$) if $k'=0$, and is equivalent to R3($k',\,-p\sum
q_i$) if $k'\not=0$. Hence if the $k_i$ are not all equal, then
relation (1) implies R2($p$) and R1($-p\sum q_i$), or R2($p$) and
R3($k',\,-p\sum q_i$). And conversely, the set of the relations
R2($p$) and R1($-p\sum q_i$) implies relation (1), as does the set
of the relations R2($p$) and R3($k',\,-p\sum q_i$).
\smallskip

Finally, any relation in $\varphi(\L_n)$ is equivalent to a set of
relations of type R1, R2 and R3. This terminates the proof of step
1..\bigskip

\noindent 2. \emph{Let us now show that there exist four integers
$M$, $m$, $d$, $s$ such that $\varphi(\L_n)$ is isomorphic to the
quotient of $\L_n$ by the three relations R1($M$), R2($m$) and
R3($d$,$s$).}
\smallskip

Let us define $M$, $m$, $d$, $s$ as follows:

\centrer{$E_1=\big\{\,k\in\NN^*\ |\ \bar\tau_1^k=1\,\big\}$ and
$M=\left\{\begin{array}{l}\min(E_1)\mbox{ if }E_1\not=\vide\\0
\mbox{ if }E_1=\vide,\end{array}\right.$}

\centrer{$E_2=\big\{\,k\in\NN^*\ |\
\bar\tau_1^k=\bar\tau_3^k\,\big\}$ and
$m=\left\{\begin{array}{l}\min(E_2)\mbox{ if }E_2\not=\vide\\0
\mbox{ if }E_2=\vide,\end{array}\right.$}

\centrer{$E_3=\big\{\,k\in\NN^*\ |\
(\bar\tau_1\bar\tau_3\dots\bar\tau_{n-1})^k\in\langle\bar\tau_1\rangle\,\big\}$
and $d=\left\{\begin{array}{l}\min(E_3)\mbox{ if }E_3\not=\vide\\0
\mbox{ if }E_3=\vide,\end{array}\right.$}

\centrer{$\vlblanc s \mbox{ is chosen arbitrarily in }
\big\{\,k\in\NN^*\ |\
(\bar\tau_1\bar\tau_3\dots\bar\tau_{n-1})^d=\bar\tau_1^k\,\big\}.$}
\medskip

\noindent Then by definition of $M$, $m$, $d$ and $s$, the three
relations R1($M$), R2($m$) and R3($d$,$s$) hold in $\varphi(\L_n)$.
According to step 1., any relation like (1) comes from some
relations of type R1($M$), $M\in\ZZ$, R2($m$), $m\in\ZZ$ and
R3($d,s$), $d,\,s\in\ZZ$, that take place in $\varphi(\L_n)$. Let us
then show that any relation R of type R1, R2 or R3 that holds in
$\varphi(\L_n)$ comes from the three relations R1($M$), R2($m$) and
R3($d$,$s$).
\smallskip

\begin{itemize}
  \item If R is of type R1:
Let $M'$ be a nonzero integer such that the relation R1($M'$) is
satisfied in $\varphi(\L_n)$. Then $E_1$ is nonempty, hence $M$ is
nonzero. Notice that the union of both relations R1($M$) and
R1($M'$) is equivalent to the relation R1($M\wedge M'$), where
$a\wedge b$ is the greatest common divisor of $a$ and $b$. However,
by definition of $M$, we have $M\leqslant M\wedge M'$, so $M=M\wedge
M'$, hence $M$ divides $M'$. Consequently, R1($M$) implies R1($M'$).
\smallskip
  \item If R is of type R2:
Similarly, any relation of type R2($m'$) where $m'$ is a nonzero
integer is induced by $R2(m)$. \smallskip
  \item If R is of type R3:
If there exist two nonzero integers $d'$ and $s'$ with $d'\not=0$
such that R3($d',\,s'$) takes place in $\varphi(\L_n)$, then $E_3$
is nonempty. So $d$ is nonzero and the conjunction of R3($d',\,s'$)
and R3($d,\,s$) induces R3($kd'+\ell d,\,ks'+\ell s$), for all
integers $k$ and $\ell$. Let us choose $k$ and $\ell$ such that
$kd'+\ell d=d' \wedge d$. By definition of $d$, we have $d\leqslant
(d'\wedge d)$, so $d=(d'\wedge d)$ and $d$ divides $d'$. Let $p$ be
the integer $d'/d$. We have:
\smallskip

\centrer{$ \left\{\begin{array}{l}
 \Rtrois(d,\,s)\\
 \Rtrois(d',\,s')
\end{array}\right. \iff \left\{\begin{array}{l}
 \Rtrois(d,\,s)\\
 \Rtrois(pd,\,ps)\\
 \Rtrois(pd,\,s')
\end{array}\right. \iff \left\{\begin{array}{l}
 \Rtrois(d,\,s)\\
 \Rtrois(pd,\,ps)\\
 \Run(|ps-s'|)
\end{array}\right. \iff \left\{\begin{array}{l}
 \Rtrois(d,\,s)\\
 \Run(|ps-s'|)
\end{array}\right. $}\smallskip

\noindent Again, the definition of $M$ implies that $|ps-s'|$ is a
multiple of $M$. Hence R3($d',\,s'$) comes from the three relations
R1($M$), R2($m$) and R3($d$,$s$).
\end{itemize}

 \bigskip

3. \emph{Let $M$, $m$, $d$ and $s$ be the integers defined in step
2.. According to step 2., $\varphi(\L_n)$ is isomorphic to the group
$\L_n$ quotiented by the three relations R1($M$), R2($m$), and
R3($d$,$s$). Let us show that if $M$ is nonzero, then $m$ and $d$
are also nonzero. Then let us determine the divisibility relations
that link these four integers.}
\smallskip
\begin{itemize}
  \item[\point]
If $M$ is nonzero, then R1($M$) implies R2($M$) and R3($M$,0), so by
definition of $m$ and $d$, we have that $m$ and $d$ are nonzero.
  \item[\point]
Since R1($M$) implies R2($M$), the relations R2($M$) and R2($m$)
coexist in $\varphi(\L_n)$, so R2($M\wedge m$) is satisfied, too.
Then, by definition of $m$, $m$ is smaller than or equal to $M\wedge
m$, so we have $m=M\wedge m$. Hence $m$ divides $M$.
  \item[\point]
Similarly, the relation R2($m$) implies R3($m$,$rm$). Now,
R3($m$,$rm$) and R3($d,\,s$) imply a third relation R3($u,\,v$)
where $u=m\wedge d$, and $v$ is an integer determined by $r$, $m$,
$d$ and $s$. But by definition of $d$, $d$ is smaller than or equal
to $u$. Hence $d$ divides $m$.
  \item[\point]
As for the integer $s$ in R3($d,\,s$), we have seen in step 1. that
R3($d,\,s$) implies R2($s$), so $m$ divides $s$.
  \item[\point]
We still have to show that $M$ divides $(r-\frac{s}{d})m$. Let us
start from the relation R3($d$,$s$) in which we replace
$\bar\tau_1^s$ by $\bar\tau_1^{kd}$ where $k=\frac{s}{d}$. We get:
\smallskip

\centrer[7]{$\big[\big(\bar\tau_1\bar\tau_3\bar\tau_5\dots\bar\tau_{n-1})\bar\tau_1^{-k}\big]^d=1$,}

\noindent then:

\centrer[8]{$\displaystyle\big[(\bar\tau_3\bar\tau_{1}^{-1})(\bar\tau_5\bar\tau_{1}^{-1})
\dots(\bar\tau_{n-1}\bar\tau_{1}^{-1})\,\bar\tau_1^{(r-k)}\big]^d=1$.}
\smallskip

\noindent Since $m$ is a multiple of $d$, we get:
\smallskip

\centrer[9]{$\big[(\bar\tau_3\bar\tau_{1}^{-1})(\bar\tau_5\bar\tau_{1}^{-1})\dots(\bar\tau_{n-1}\bar\tau_{1}^{-1})\,
\bar\tau_1^{(r-k)}\big]^m=1$.}
\smallskip

\noindent Now, according to R2($m$), for all
$i\in\Imp(n)\smallsetminus\{1\}$, we have
$(\bar\tau_i\bar\tau_{1}^{-1})^m=1$. Hence (9) implies:
\smallskip

\centrer[]{$\bar\tau_1^{(r-k)m}=1$.}
\smallskip

\noindent In other words, R1$\big((r-k)m\big)$ takes place in
$\varphi(\L_n)$. Then, as before, we deduce from definition of $M$
that $M$ divides $(r-k)m$. \end{itemize}
\bigskip

4. \emph{Let us show that $\varphi$ is cyclic if and only if $d=1$.}
\smallskip

If $\varphi$ is cyclic, then R2(1) holds, and so does R3(1, $r$).
Conversely, if $d=1$, let us show that $m=1$. Let $\xi$ be the
element $\tau_1^{s}$. If $d=1$, we have:
\smallskip

\centrer[10]{$\bar\tau_1\bar\tau_3\bar\tau_5\dots
\bar\tau_{n-1}=\bar\xi$,}

\noindent whence:

\centrer[11]{$\bar\tau_1^{-1}\bar\tau_3^{-1}=\bar\tau_5\dots
\bar\tau_{n-1}\bar\xi^{-1}$.}

\smallskip

\noindent Since $\bar\xi=\bar\tau_1^s$ and $m$ divides $s$, then
$\bar\xi$ is a multiple of $\bar\tau_1^m$. But according to the
relation R2($m$) and Lemma \ref{lem:powers_égales}, we have
$\bar\tau_1^m=\bar\tau_2^m=\dots=\bar\tau_{n-1}^m$, so
$\bar\tau_1^m$ is central in $\varphi(\B_n)$, so $\bar\xi$ is
central in $\varphi(\B_n)$. According to equality (11), it follows
that $\bar\tau_2$ commutes with the right-hand side, hence
$\bar\tau_2$ commutes with the left-hand side. So we get:
\smallskip

\centrer[]{$\bar\tau_2\bar\tau_1^{-1}\bar\tau_3^{-1}=\bar\tau_1^{-1}\bar\tau_3^{-1}\bar\tau_2$,}

\noindent whence

\centrer[]{$\bar\tau_1\bar\tau_2\bar\tau_1^{-1}=\bar\tau_3^{-1}\bar\tau_2\bar\tau_3$,}

\noindent but

\centrer[]{$\bar\tau_3^{-1}\bar\tau_2\bar\tau_3=\bar\tau_2\bar\tau_3\bar\tau_2^{-1}$,}

\noindent so

\centrer[]{$\bar\tau_1\bar\tau_2\bar\tau_1^{-1}=\bar\tau_2\bar\tau_3\bar\tau_2^{-1}$,}

\noindent and by conjugating by $\bar\tau_1\bar\tau_2\bar\tau_3\bar\tau_4$:

\centrer[]{$\bar\tau_2\bar\tau_3\bar\tau_2^{-1}=\bar\tau_3\bar\tau_4\bar\tau_3^{-1}$,}

\noindent whence

\centrer[]{$\bar\tau_1\bar\tau_2\bar\tau_1^{-1}=\bar\tau_3\bar\tau_4\bar\tau_3^{-1}$,}

\noindent then

\centrer[]{$\bar\tau_2=
  \big(\bar\tau_1\bar\tau_2\bar\tau_1^{-1})\bar\tau_1(\bar\tau_1\bar\tau_2\bar\tau_1^{-1}\big)^{-1}=
  \big(\bar\tau_3\bar\tau_4\bar\tau_3^{-1})\bar\tau_1(\bar\tau_3\bar\tau_4\bar\tau_3^{-1}\big)^{-1}=\bar\tau_1$,}

\noindent so $\bar\tau_2=\bar\tau_1$. Hence $\varphi$ is cyclic.\fin
\bigskip

\begin{defi}[$\L_n(M,m,d,s)$]
\mbox{}\\For all quadruples of integers ($M$, $m$, $d$, $s$), as
soon as this definition makes sense (i.e. when $m$ is a multiple of
$d$ and $M$ is a multiple of $m$, according to Lemma
\ref{lem:relations}), let us denote by $\L_n(M,m,d,s)$ the group
$\langle\tau_i,\ i\in\Imp(n)\rangle$ quotiented by the relations
$(\Run(M))$, $(\Rdeux(m))$, $(\Rtrois(d,s))$. For example,
$\L_n(0,0,0,0)\cong\ZZ^r$ and $\L_n(M,M,M,rM)\cong(\ZZ[M])^r$.
\end{defi}
\medskip

\begin{lem}[Cardinality of $\L_n(M,m,d,s)$]
        \label{lem:cardinality}
\mbox{}\\For all $M>0$, $m\geqslant 2$, $d$ and $s$, the cardinality
of $\L_n(M,m,d,s)$ is equal to $qdm ^{r-1}$ where $q=\frac{M}{m}$
and $r=\frac{n}{2}$. \end{lem}
\medskip

\DEM The group $\L_n(0,0,0,0)$ is spanned by:
\smallskip

\centrer[1]{$\tau_1,_,\tau_3,\,\tau_5,\,\dots,\,\tau_{n-3},\,\tau_{n-1}$.}
\smallskip

\noindent Let us set $ u_i=\tau_i\tau_1^{-1}$ for all
$i\in\Imp(n)\smallsetminus\{1\}$. We set $k=\frac{s}{d}$ ($k$ is an
integer for, according to Lemma \ref{lem:relations}, $d$ divides $m$
which divides $s$). Then we set $ w= \big(u_3 u_5\dots u_{n-3}
u_{n-1}\big)\tau_1^{(r-k)}$. Thanks to a change of variables, we go
from the set (1) spanning $\L_n(0,0,0,0)$ to the below set (2) still
spanning $\L_n(0,0,0,0)$:
\smallskip

\centrer[2]{$\tau_1,\, u_3,\, u_5,\,\dots,\, u_{n-3},\, w$.}
\smallskip

\noindent With this change of variables, the relation R1($M$) is now
equivalent to:
\smallskip

\centrer[3]{$\tau_1^M=1,\, u_3^M=1,\, u_5^M=1,\,\dots,\,
u_{n-3}^M=1,\, w^M=1$.}
\smallskip

\noindent Let us denote by $\xi\mapsto\bar\xi$ the canonical
homomorphism from $\L_n(0,0,0,0)$ to $\L_n(M,M,M,rM)$, which is the
quotient of $\L_n(0,0,0,0)$ by R1($M$). According to Lemma
\ref{lem:relations}, $M$ divides $(r-k)m$, so in $\L_n(M,M,M,rM)$,
we have:
\smallskip

\centrer{$\bar w^m=\big(\bar u_3 \bar u_5\dots \bar u_{n-3} \bar
u_{n-1}\big)^m\bar \tau_1^{(r-k)m}=\big(\bar u_3 \bar u_5\dots \bar
u_{n-3} \bar u_{n-1}\big)^m$.}
\smallskip

\noindent Hence the relation R2($m$) in $\L_n(M,M,M,rM)$ is
equivalent to:
\smallskip

\centrer[4]{$\bar u_3^m=1,\,\bar u_5^m=1,\,\dots,\,\bar
u_{n-3}^m=1,\,\bar w^m=1$.}
\smallskip

\noindent Finally, in $\L_n(M,M,M,rM)$, the relation R3($d$,$s$) is:
$(\bar\tau_1\bar\tau_3\dots\bar\tau_{n-1})^d=\bar\tau_1^s$. Let us
replace $\tau_1^s$ by $\tau_1^{kd}$, the relation $R3(d,s)$ is
equivalent to
$\big[\big(\tau_1\bar\tau_3\bar\tau_5\dots\bar\tau_{n-1})\bar\tau_1^{-k}\big]^d=1$,
and then to
\smallskip

\centrer{$\displaystyle\big[(\bar\tau_3\bar\tau_{1}^{-1})(\bar\tau_5\bar\tau_{1}^{-1})
\dots(\bar\tau_{n-1}\bar\tau_{1}^{-1})\,\bar\tau_1^{(r-k)}\big]^d=1$,}
\smallskip

\noindent so the relation $R3(d,s)$ in $\L_n(M,M,M,rM)$ is
equivalent to:
\smallskip

\centrer[5]{$\bar w^d=1$.}
\smallskip

\noindent Finally, since $m$ divides $M$ and since $d$ divides $m$,
the set of relations R1($M$), R3($m$) and R3($d$,$s$) is equivalent
in $\L_n(0,0,0,0)$ to:
\smallskip

\centrer[6]{$\tau_1^M=1,\,u_3^m=1,\,u_5^m=1,\,\dots,\,u_{n-3}^m=1,\,w^d=1$.}
\smallskip

\noindent Therefore a presentation by generators and relations of
the group $\L_n(M,m,d,s)$ can be obtained from the lines (2) and
(6). Therefore $\L_n(M,m,d,s)$ is isomorphic to
$\ZZ[M]\times\big(\ZZ[m]\big)^{r-2}\times\ZZ[d]$. So its cardinality
is $Mm^{r-2}d=qdm^{r-1}$ where $q=\frac{M}{m}$.\fin
\bigskip


\subsection{Periodic geometric representations of $\B_n$}
\label{par:homomorphisms_periodic}
\medskip

In the first proposition, we deal with the case where the boundary
of $\Sigma$ is nontrivial ($\Sigma=\Sigma_{g,\,b}$ with $b>0$). The
remainder of this subsection is devoted to the case without boundary
($b=0$), which is harder. For all homomorphisms $\rho$ from $\B_n$ to $\PMod(\Sigma)$ and for all $i\leqslant n-1$, we denote by $A_i$ the
mapping class $\rho(\tau_i)$.
\bigskip

\begin{prop}[The periodic homomorphisms from $\B_n$ to $\PMod(\Sigma_{g,\,b})$, $b>0$, are cyclic]
                \label{prop:periodic_implies_cyclic_b>0}
\mbox{}\\
Any periodic homomorphism $\rho$ from $\B_n$ to $\PMod(\Sigma_{g,\,b})$
where $g\leqslant \frac{n}{2}$ and $b>0$ is cyclic. \end{prop}
\bigskip

\DEM Let $\rho$ be a periodic homomorphism from $\B_n$ to $\PMod(\Sigma)$. Notice that the mapping classes $A_i$ for all
$i\leqslant n-1$ are conjugate. So they are periodic and have the
same order. Let us call $m$ this order.

Since the boundary of $\Sigma$ is nonempty, according to Lemma
\ref{lem:group_cyclic}, $A_1$ and $A_4$ span a cyclic group that we
denote by $\Gamma$. Any generator of $\Gamma$ is a product of powers
of $A_1$ and $A_4$, so $\Gamma$ is a cyclic group of order $m$. Now,
the subgroups of $\Gamma$ spanned on one hand by $A_1$ and on the
other hand by $A_4$ have the same order $m$, so each of $A_1$ and
$A_4$ spans independently $\Gamma$. Thus, in $\G$, any two standard
generators that commute span the same cyclic group. Hence $A_2$ span
the same cyclic group as $A_4$, that is, the same cyclic group as
$A_1$. In particular $A_2$ and $A_1$ commute. But $A_2$ and $A_1$
satisfy a braid relation, so they have to be equal. Then, according
to Lemma \ref{lem:homomorphism_cyclic}, $\rho$ is cyclic. \fin
\bigskip

\begin{prop}[The periodic homomorphisms from $\B_n$ to $\Mod(\Sigma_{g,\,0})$ are cyclic]
                \label{prop:periodic_implies_cyclic_b=0}
\mbox{}\\Let $n$ be an integer greater than or equal to 6 and
$\Sigma$ the surface $\Sigma_{g,\,0}$ such that $g\leqslant
\frac{n}{2}$. Any periodic homomorphism $\rho$ from $\B_n$ to $\Mod(\Sigma)$ is cyclic.
\end{prop} \medskip

\DEM We argue by contradiction. Let $\rho$ be a periodic
homomorphism. We assume that $\rho$ is not cyclic. We separate the
cases according to the orders of $A_1$ and $A_3A_1^{-1}$ (the
order of $A_3A_1^{-1}$ is different from 1 since $\rho$ is not
cyclic).
\smallskip

1. \textbf{When $A_1$ is of order 2}.
\smallskip

\noindent If for all $i\leqslant n-1$, the mapping class $A_i$
is of order 2, then $\rho(\B_n)$ is isomorphic to a quotient of
the symmetric group $\Ss_n$. It follows easily from the
simplicity of $\Aa_n$ (remember that $n\geqslant 6$) that the
only quotients of $\Ss_n$ are $\{1\}$, $\ZZ[2]$ and $\Ss_n$. So
the only nontrivial quotient of $\Ss_n$ is $\ZZ[2]$. However
we assume that $\rho$ is not cyclic, so $\rho(\B_n)$ has to be
isomorphic to the group $\Ss_n$. In particular, $\rho(\B_n)$ is
finite, and its cardinality is:
\smallskip

\centrer{$n!=(n-1)!\times n\ \geqslant\ 5!\times n\ =\ 120n\
\geqslant\ 240g$.}
\smallskip

\noindent Now, according to the ``84($g-1$)'' theorem (see Corollary
\ref{cor:borne_cardinality_Mod}), the cardinality of a finite
subgroup of $\Mod(\Sigma)$ is bounded by $42|\chi(\Sigma)|=84g-84$,
whence a contradiction.
\bigskip

2. \textbf{When $A_1$ is not of order 2, but $A_3A_1^{-1}$ is of
order 2}.
\smallskip

\noindent If $(A_3A_1^{-1})^2=1$, then $A_3^2=A_1^2$, so according
to Lemma \ref{lem:powers_égales}, we have
$A_1^2=A_2^2=\dots=A_{n-1}^2$. Let $\Z$ be the centralizer
$\Z_{\Mod(\Sigma)}(A_1^2)$ of $A_1^2$ in $\Mod(\Sigma)$ and let $p$
be the canonical homomorphism from $\Z$ to $\Z/\langle A_1^2\rangle$.
Notice that $\rho(\B_n)$ is included in $\Z$, so we can consider the
homomorphism $p\rond\rho$ from $\B_n$ to $\Z/\langle A_1^2\rangle$. It
is not cyclic, for $A_3A_1^{-1}$ is not a power of $A_1^2$: indeed,
if it existed an integer $k$ such that $A_3=A_1^{(1+2k)}$, by
conjugation, we would have $A_5=A_1^{(1+2k)}$, and so $A_3=A_5$. But
this is absurd for $\rho$ is not cyclic. Thus, the homomorphism
$p\rond\rho$ is not cyclic, but $p\rond\rho(\tau_1^2)=p(A_1^2)=1$.
Then, as we have seen it in step 1., $p\rond\rho(\B_n)$ is
isomorphic to the group $\Ss_n$ and hence contains at least 240g
elements. But $\rho(\B_n)$ is a (central) extension of
$p\rond\rho(\B_n)$ by the finite group $\langle A_1^2\rangle$, in
other words, the following sequence is exact:
\smallskip

\centrer{$1\to\langle
A_1^2\rangle\to\rho(\B_n)\xrightarrow{\:p\:}p\rond\rho(\B_n)\to1$.}
\smallskip

\noindent Hence $\rho(\B_n)$ is a finite group that contains at
least 480g elements. As in step 1., this is absurd.
\bigskip

3. \textbf{Where $A_1$ is of order $M\geqslant 3$ and $A_3A_1^{-1}$
is of order $m$ with $3\leqslant m\leqslant M$}.
\smallskip

\noindent According to Kerckhoff's Theorem (cf.
\ref{thm:Kerckhoff}), the abelian group
$\langle\,A_i,\,i\in\Imp(n)\,\rangle$ being finite, there exist a
hyperbolic metric $g$ on $\Sigma$ and an injective homomorphism from
$\langle\,A_i,\,i\in\Imp(n)\,\rangle$ to $\Isom(\Sigma\,;\,g)$. Let
us denote by $\F$ its image and $\bar A_i$ the image of $A_i$ for
all $i\in\Imp(n)$. Let us recall that we assume that $\rho$ is not
cyclic, so $\bar A_1\not=\bar A_3$. We will show that the action of
$\F$ on the points of $\Sigma$ is free, for if an element of $\F$
had a fixed point in $\Sigma$, it would automatically have many,
actually too much compared with Corollary \ref{cor:points_fixes}. We
will conclude by showing that if the elements of $\F$ do not have
any fixed point, the inequality linking $\chi(\Sigma)$ and
$\chi(\Sigma/\langle\F\rangle)$ given by Lemma
\ref{lem:periodic_Euler} cannot be satisfied, whence the
contradiction.
\medskip

a) Let us show that the action of $\F$ on $\Sigma$ is free.
\smallskip

Let $x$ be a point of $\Sigma$ and let $\Stab(x)$ be the subgroup of
$\F$ that fixes the point $x$. Let us assume that $\Stab(x)$ is not
reduced to $\{1\}$. Let us recall that two isometries that fixes a
same point and that have the same differential in this point are
equal (cf. Lemma \ref{lem:rotation}). But the differential of an
isometry in a fixed point is a rotation. Therefore $\Stab(x)$ is a
cyclic group. Let $G$ be an isometry spanning $\Stab(x)$, let $M'$
be its order, with $2\leqslant M'\leqslant M$, for on one hand $G$
is not the identity, on the other hand $G$ belongs to the abelian
group $\F$ spanned by elements of order $M$. We are going to count
the number $\ell$ of fixed points of $G$. On one hand, according to
Corollary \ref{cor:points_fixes}, we have:
\smallskip

\centrer[1]{$\displaystyle \ell\leqslant 2+\frac{2g}{M'-1}$.}
\smallskip

\noindent On the other hand, if $G$ commutes with another isometry
$G'$, then the images by $G'$ of all fixed points of $G$ are again
fixed points of $G$. Since the group $\F$ is abelian, the set of
fixed points of $G$ contains the orbit of $x$ by the group $\F$, so:
\smallskip

\centrer[2]{$\displaystyle |\Orb(x)|\leqslant \ell$,}
\smallskip

\noindent where $\Orb(x)$ is the orbit of $x$. By definition of $\F$
and according to Lemma \ref{lem:relations}, there exist four
integers $M'$, $m'$, $d$ and $s$ such that $\F$ is isomorphic to
$\L_n(M',m',d,s)$. Now, since $M$ is the order of $A_1$ and $m$ is
the order of $A_3A_1^{-1}$, we have $M'=M$ and $m'=m$. Hence $\F$ is
isomorphic to $\L_n(M,m,d,s)$ and according to Lemma
\ref{lem:cardinality}, the cardinality of $\F$ is $qd{(m)}^{r-1}$
where $q=\frac{M}{m}$ and $r=\frac{n}{2}$. We can then compute the
cardinality of the orbit of $x$:
\smallskip

\centrer[3]{$\displaystyle |\Orb(x)|=\frac{|\F|}{|\Stab(x)|}=
\frac{qd{(m)}^{r-1}}{M'}$.}
\smallskip

\noindent From (1), (2) and (3), we get:
\smallskip

\centrer[4]{$\displaystyle
\frac{qd{(m)}^{r-1}}{M'}=|\Orb(x)|\leqslant \ell\leqslant
2+\frac{2g}{M'-1}$.}
\smallskip

\noindent By multiplying all by $\frac{M'}{q}$, we get:
\smallskip

\centrer[5]{$\displaystyle d{(m)}^{r-1}\leqslant
2\frac{M'}{q}+2g\,\frac{1}{q}\,\frac{M'}{M'-1}$.}
\smallskip

\noindent Since $\frac{M'}{q}\leqslant\frac{M}{q}= m$, we can bound
$\frac{M'}{q}$ by $m$. We bound $g$ by $r$, $\frac{1}{q}$ by $1$,
and $\frac{M'}{M'-1}$ by 2. Then (5) becomes:
\smallskip

\centrer[6]{$d(m)^{r-1}\leqslant2m+4r$, with
$\left\{\mbox{\begin{tabular}{l}
  $r=\frac{n}{2}\geqslant 3$,\\
  $d\geqslant2$ and $d$ divides $m$, according to Lemma \ref{lem:relations},\\
  $m\geqslant 3$ by hypothesis.\\
\end{tabular}}\right.$}
\smallskip

\noindent When $m=3$, we have $d=3$, so (6) becomes:
\smallskip

\centrer[7]{$3^r\leqslant 6+4r$,}
\smallskip

\noindent but this equation is never satisfied for $r\geqslant 3$
(for $r=3$, we get $27\leqslant 6+12$ which is absurd, and for
$r>3$, this is even more flagrant). When $m\geqslant 4$, let us
consider equation (6), we bound  $4r$ by $mr$ in the right-hand
side, we divide the left-hand side and the right-hand side by $m$,
then in the left-hand side, we replace $m$ by its lower bound: 4,
and $d$ by its lower bound: 2. We get:
\smallskip

\centrer[8]{$2\times4^{(r-2)}\leqslant2+r$,}
\smallskip

\noindent that is not satisfied for $r=3$ and certainly not for
$r>3$. Thus, it was absurd to assume that $\Stab(x)\not=\{1\}$.
Hence the action of $\F$ on $\Sigma$ is free.
\medskip

b) Let us apply the Riemann-Hurwitz' formula (cf. Lemma
\ref{lem:periodic_Euler}) to the finite group $\F$:
\smallskip

\centrer[9]{$\chi(\Sigma)+\sum (|\F|-o(Q_i))=|\F|.\chi(\Sigma/\F)$.}
\smallskip

\noindent The surface $\Sigma$ satisfies $\chi(\Sigma)=2-2g$.
Besides, as the action of $\F$ on $\Sigma$ is free, there is no
point of ramification $Q_i$ in the surface $\Sigma/\F$ hence $\sum
(|\F|-o(Q_i))=0$. So the two terms of equality (9) are negative.
Since the elements of $\F$ preserve the orientation, $\Sigma/\F$ is
an orientable closed surface with $\chi(\Sigma/\F)\leqslant-2$. But
the order of $\F$ is $qd(m)^{r-1}$ with $q\geqslant 1$, $d\geqslant
2$, $m\geqslant 2$ and $r\geqslant g$, so $|\F|\geqslant 2^g$, so
the equality $(9)$ implies $2-2g\leqslant2^g(-2)$, i.e.:
\smallskip

\centrer[10]{$g\geqslant1+2^g$, with $g\geqslant 0$,}
\smallskip

\noindent which is absurd. \fin
\bigskip


\subsection{Pseudo-Anosov geometric representations of $\B_n$}
\label{par:homomorphisms_Anosov}
\medskip

\begin{prop}[The pseudo-Anosov homomorphisms from $\B_n$ to $\PMod(\Sigma)$ are cyclic]
            \label{prop:pseudo-Anosov_implies_cyclic}
\mbox{}\\Let $n$ be an integer greater than or equal to 6 and let
$\Sigma$ be a surface $\Sigma_{g,\,b}$ where $g\leqslant
\frac{n}{2}$. Any pseudo-Anosov homomorphism from $\B_n$ to $\PMod(\Sigma)$ is cyclic. \end{prop}
\bigskip

\DEM Let $\rho\;:\; \B_n\to\PMod(\Sigma)$ be a pseudo-Anosov
homomorphism. For all $i\leqslant n-1$, we set again $A_i=\rho(\tau_i)$.
\smallskip

1. The mapping class $A_1$ is pseudo-Anosov, so according to
Proposition \ref{prop:structure_centralizer}.(iv), its
centralizer is virtually infinite cyclic. Since the mapping
class $A_3$ commutes with $A_1$, there exist two nonzero
integers $p$ and $p'$ such that $A_1^{p'}=A_3^{p}$. By
conjugating this equality by $A_3A_4A_3$, we get
$A_1^{p'}=A_4^{p}$. Hence $A_3^p=A_4^p$, so according to Lemma
\ref{lem:powers_égales}:
\smallskip

\centrer[1]{$A_1^p=A_2^p=A_3^p\dots=A_{n-1}^p$.}
\smallskip

Let us exploit this. We separate the cases whether $b>0$ (cf. 2.) or
$b=0$ (cf. 3. - 5.).
\smallskip

2. When $b>0$, we produce a direct proof. According to Proposition
\ref{prop:passage_PMod_à_ModBord}, there exists $\tilde\rho$, a lift
of $\rho\in\Hom(\B_n,\,\PMod(\Sigma))$ in
$\Hom(\B_n,\,\Mod(\Sigma,\,\bord\Sigma))$. For all $i\leqslant n-1$,
let us denote by $\tildeA_i$ the mapping class $\tilde\rho(\tau_i)$,
so that $\tildeA_i$ is a lift of $A_i$. Let us set then
$W=\tildeA_3\tildeA_1^{-1}$. Since $A_1^p=A_3^p$, the mapping class
$W^p$ is a multitwist along the boundary components. Let $Z$ be the
mapping class $(\tildeA_1\tildeA_2\tildeA_3)^2$. Then
$Z\tildeA_1Z^{-1}=\tildeA_3$ and $Z\tildeA_3Z^{-1}=\tildeA_1$, so
$ZW^pZ^{-1}=W^{-p}$. Since $W$ is central in
$\Mod(\Sigma,\,\bord\Sigma)$, $Z$ and $W^p$ commute, so we have
$W^{2p}=\Id$. But $\Mod(\Sigma,\bord\Sigma)$ is torsion-free
according to Lemma \ref{lem:without_torsion}, so $W=\Id$. Hence
$\tildeA_1=\tildeA_3$ and $\tilde\rho$ is cyclic, hence $\rho$ is
cyclic.
\medskip

3. When $b=0$, we argue by contradiction and we assume that $\rho$
is not cyclic. Then, according to Lemma \ref{lem:powers_égales}, the
$A_i$, $1\leqslant i\leqslant n-1$, are pairwise distinct. Let us
consider the group $\Centr(A_1^p)$ (cf. Definition \ref{defi:centr}
of Subsection \ref{par:mapping_classes_ps-An}). According to (1),
$\rho(\B_n)\subset\Centr(A_1^p)$. Let $\ell$ be the homomorphism
associated to $A_1^p$ defined by Proposition
\ref{prop:structure_centralizer}. According to this proposition, the
cardinality of $\Ker(\ell)$ satisfies:
\smallskip

\centrer[2]{$|\Ker(\ell)|\leqslant 6|\chi(\Sigma)|$.}
\smallskip

\noindent Since all the $A_i$, $i\leqslant n-1$, are conjugate in
$\Centr(A_1^p)$, $\ell(A_i)$ is independent of the index $i$ when
$i$ ranges from $1$ to $n-1$. Hence the group spanned by
$A_jA_k^{-1}$ where $j,\, k\leqslant n-1$ is included in
$\Ker(\ell)$. Yet, we are going to show that its cardinality is
greater than $6|\chi(\Sigma)|$, whence the contradiction.

\medskip

4. Let us assume that $(A_1A_3^{-1})$ is of order $p=2$. Recall that
$\rho$ is not cyclic, so $A_1$ and $A_3$ are different. Then the
subgroup $\F$ of $\Ker(\ell)$ defined by:
\smallskip

\centrer{$\F:=\big\langle\,A_iA_{n-1}^{-1}\,;\,1\leqslant i\leqslant
n-3\,\big\rangle$}
\smallskip

\noindent is isomorphic to a quotient $\Ss_{n-2}$ by the
homomorphism: $(12)\mapsto A_1A_{n-1}^{-1}$, $(23)\mapsto
A_2A_{n-1}^{-1}$, \dots, $(n-3,n-2)\mapsto
A_{n-3}A_{n-1}^{-1}$. However, $A_1\not=A_3$, so this quotient
is neither $\{1\}$ nor $\ZZ[2]$. Then, when $n\geqslant 8$,
because of the simplicity of $\Aa_6$ this quotient has to be
$\Ss_{n-2}$. When $n=6$, the only quotient of $\Ss_4$ different
from $\{1\}$, $\ZZ[2]$ and $\Ss_4$ is the quotient of $\Ss_4$
by the normal closure of the element (12)(34). The image of
(12)(34) in $\F$ is $(A_1A_5^{-1})(A_3A_5^{-1})$, which is
equal to $A_1A_3A_5^{-2}$, hence equal to $A_1A_3^{-1}$, for
$A_3^{-2}=A_5^{-2}$ according to (1). Since $A_1\not=A_3$,
$A_1A_3^{-1}$ is not trivial. Hence $\F$ is not isomorphic to
the above quotient of $\Ss_4$. Then, even when $n=6$, $\F$ is
isomorphic to $\Ss_{n-2}$. Hence $\Ker(\ell)$ contains $\F$
that owns $(n-2)!$ elements. When $n\geqslant 8$, we get:
\smallskip

\centrer[3]{$|\Ker(\ell)|\geqslant|\F|=(n-2)!\geqslant 5! (n-2) >
6(n-2)\geqslant 6(2g-2) = 6|\chi(\Sigma)|$.} \medskip

\noindent  But (2) and (3) lead to a contradiction, this is the
expected contradiction.
\smallskip

The case $n=6$ implies $|\F|=|\Ss_4|=4!=24\geqslant6|\chi(\Sigma)|$,
since $n=6$ implies that $g\leqslant \frac{n}{2}=3$, and then
$|\chi(\Sigma)|\leqslant 4$. However $\Ker(\ell)$ contains the
element $A_{n-2}A_{n-1}^{-1}$, too, which is different from any
element of $\F$, for $\F$ is in the centralizer of $A_{n-1}$,
whereas $A_{n-2}A_{n-1}^{-1}$ is not. Indeed, if
$A_{n-2}A_{n-1}^{-1}$ was in the centralizer of $A_{n-1}$, then
$A_{n-2}$ and $A_{n-1}$ would commute. However they satisfy a braid
relation, so they would be equal and $\rho$ would be cyclic: this is
absurd. Thus $\Ker(\ell)$ contains $\F$ and the element
$A_{n-2}A_{n-1}^{-1}$, which does not belong to $\F$. Since the
cardinality of $\F$ satisfies $|\F|\geqslant6|\chi(\Sigma)|$, then
the cardinality of $\Ker(\ell)$ satisfies
$|\Ker(\ell)|>6|\chi(\Sigma)|$, which contradicts (2). This is the
expected contradiction.
\medskip

5. Let us assume that $(A_1A_3^{-1})$ is of order $p\geqslant 3$ and
let us consider the abelian groups $\H$ and $\H'$ defined by:
\smallskip

\centrer{$\H:=\big\langle\, A_i\ ,\
i\in\Imp(n)\,\big\rangle_{\Mod(\Sigma)}$ }
\smallskip

\centrer{and \quad $\H':=\big\langle\, A_iA_{n-1}^{-1}\ ,\
i\in\Imp(n-2)\,\big\rangle_{\Mod(\Sigma)}$.}
\smallskip

\noindent Let us apply Lemma \ref{lem:relations} to these two
groups:
\smallskip

\begin{itemize}
   \item \emph{Concerning the group $\H$.}
It is clear that there exist two integers $d$ and $s$ such that $\H$
is isomorphic to $\L_n(0,p,d,s)$, where $d\not=1$, $d$ divides $p$,
and $p$ divides $s$, according to Lemma \ref{lem:relations}.
Moreover, since $\ell(A_1)=\dots=\ell(A_{n-1})>0$, all these
relations have to be homogeneous, so by considering the relation
R3($d,s$), it follows $s=rd$. The relation R3($d,s$) then becomes:
\smallskip

\centrer[4]{$\big(\dessous{\prod}{i\in\Imp(n)}A_i\big)^{d}=A_1^{rd}$,}
\smallskip

\noindent and it implies:
\smallskip

\centrer[5]{$\big(\dessous{\prod}{i\in\Imp(n-3)}A_iA_{n-1}^{-1}\big)^{d}=(A_1A_{n-1}^{-1})^{rd}$,}
\smallskip

    \item \emph{Concerning the group $\H'$.}
It is clear that $\H'$ is isomorphic to $\L_{n-2}(p,p,d',s')$ where
$d'$ and $s'$ are to be determine. The relation $(\Rtrois(d',s'))$
is equivalent to:
\smallskip

\centrer[6]{$\big(\dessous{\prod}{i\in\Imp(n-3)}A_iA_{n-1}^{-1}\big)^{d'}=(A_1A_{n-1}^{-1})^{s'}$,}
\smallskip

\noindent and implies:
\smallskip

\centrer[7]{$\big(\dessous{\prod}{i\in\Imp(n)}A_i\big)^{d'}=A_1^{s'}A_{n-1}^{\,(rd'-s')}$.}
\smallskip

\noindent Since $p$ divides $s'$, according to Lemma
\ref{lem:relations}, we have $A_1^{s'}=A_{n-1}^{s'}$ and (7)
becomes:
\smallskip

\centrer[8]{$\big(\dessous{\prod}{i\in\Imp(n)}A_i\big)^{d'}=A_{n-1}^{\,rd'}$.}
\end{itemize}

\noindent In (8), by conjugation, we can replace $A_{n-1}^{\,rd'}$
by $A_{1}^{\,rd'}$. Let us compare the equalities (4) and (8). By
definition of $d$, it follows from that comparison that $d$ divides
$d'$. Moreover, by comparing (5) and (6), it follows by definition
of $d'$ that $d'$ divides $d$. Thus $d'=d$, so $\H'$ is isomorphic
to $\L_{n-2}(p,p,d,s')$, and according to Lemma \ref{lem:relations},
the following holds:
\smallskip

\centrer[9]{$d\geqslant 2$, $d$ divides $p$.}
\smallskip

\noindent Then according to Lemma \ref{lem:cardinality}, $|\H'|=
dp^{r-2}$. The only pairs $(p,\,d)$ that respect (9) and such that
$p<6$ are (3,3), (4,2), (4,4) and (5,5). However, if
$(p,\,d)=(4,\,2)$, then $r$ is even. Indeed, as we saw it in the
lines preceding (4), $p$ divides $s$ and $s$ is equal to $rd$. We
check in the following table all the possible values of $dp^{r-2}$
for the pairs $(p,d)$ where $p<6$, as a function of $r$, and we give
a lower bound to the values $dp^{r-2}$ for the pairs $(p,d)$ with
$p\geqslant 6$, as a function of $r$.

\smallskip

\centrer{$\begin{array}{|c|c|c|c|c|c|} \hline
  \vblanc
  \raisebox{-1ex}{$r$}\backslash\raisebox{1ex}{$(p,\,d)$} & (3,3) & (4,2) & (4,4) & (5,5) & \mbox{$(p,d)$) with $p\geqslant6$} \\
\hline
  r=3 & 9 & - & 16 & 25 & dp^{r-2} \geqslant 2\times 6=12\\
\hline
  r=4 & 27 & 32 & 64 & 125 & dp^{r-2}\geqslant 2\times6^2=72\\
\hline
  r\geqslant 5 & 27\times3^{r-4} & 32\times4^{r-4} & 64\times4^{r-4} & 125\times5^{r-4} & dp^{r-2}\geqslant 2p^{r-2}=72\times 6^{r-4}\\
\hline \end{array}$}
\medskip

\centrer{Table 1 -- Computation of $dp^{r-2}$ as a function of $d$,
$p$ and $r$.}
\bigskip

\noindent According to Table 1, for all $r\geqslant 3$, the
expression $dp^{r-2}$ achieves its lower bound when $p=d=3$, hence:
\smallskip

\centrer[11]{$|\H'|\geqslant 3^{r-1}$.}
\smallskip

\noindent However, $\Ker(\ell)$ contains also the element
$A_2A_{n-1}^{-1}$, which is of order $p$, too. And $A_2A_{n-1}^{-1}$
does not commute with $A_1$ (otherwise, $A_2$ would commute with
$A_1$, we would have $A_1=A_2$ and $\rho$ would be cyclic). Since
$\H'$ is in the centralizer of $A_1$, the mapping class
$A_2A_{n-1}^{-1}$ cannot belong to $\H'$. Similarly
$(A_2A_{n-1}^{-1})^{-1}$ cannot belong to $\H'$. But
$A_2A_{n-1}^{-1}$ and its inverse are distinct, for $p\geqslant 3$.
Then the group $\big\langle\,\H'\cup A_2A_{n-1}^{-1}\,\big\rangle$
contains the following set:
\smallskip

\centrer{$\big\{\,H(A_2A_{n-1}^{-1})^k,\ H\in\H',\
k\in\{-1,\,0,\,1\}\,\big\}$.}
\smallskip

\noindent Its cardinality is $3|\H'|$. Hence $|\Ker(\ell)|\geqslant
3^r$. But for all integers of $r\geqslant 3$, the number $3^r$ is
greater than $6(2r-2)$, which is greater than or equal to
$6|\chi(\Sigma)|$. Thus:

\smallskip

\centrer[12]{$|\Ker(\ell)|>6|\chi(\Sigma)|$.} \medskip

\noindent This contradicts assertion (2). This is the expected
contradiction and the end of the proof.\fin
\bigskip

\section{Adding an $n$\eme generator to the standard presentation of $\B_n$}
    \label{par:group_of_tresses_Artin}
    \label{sec:n_th_generator}
\bigskip

\TITRE{Completing the standard presentation of the braid group}
\mbox{}\\
We introduce the following element of the braid group $\B_n$:
\smallskip

\centrer{$\delta=\tau_1\tau_2\dots\tau_{n-1}$.}
\smallskip

\noindent We then define an $n$\eme generator $\tau_0$ by setting:
\smallskip

\centrer{$\tau_0=\delta\tau_{n-1}\delta^{-1}$.}

\noindent  We adopt the following convention:
for all integers $k\in\ZZ$, $\tau_k$ is the standard generator
$\tau_\ell$ where $\ell$ is the remainder of the euclidian division
of $k$ by $n$. Then, for all $i\in\{0,\dots,\,n-1\}$, we have:

\centrer{ $\delta\,\tau_{i}\,\delta^{-1}=\tau_{i+1}$.}
\smallskip

\noindent Moreover, for all pairs of integers $(i,\,j)$, we
denote by $|i-j|_n$ the integer $\min(\{|i-j+kn|, k\in\ZZ\})$.
Then, the braid group with $n$ strands $\B_n$ admits the following
presentation:
\begin{itemize}
 \item generators: $\tau_i$, $i\in\{0,\,1,\dots,\,n-1\}$,
 \item relations: for all $i,\,j\in\{0,\,1,\dots,\,n-1\}$:
\quad $\left\{\begin{array}{l}
 \tau_{i+1}\tau_{i+2}\dots\tau_{i+n-1}=\tau_{j+1}\tau_{j+2}\dots\tau_{j+n-1}\;\\
 \tau_i \tau_j= \tau_j \tau_i \mbox{ when } |i-j|_n\not=1\;\\
 \tau_i \tau_j \tau_i = \tau_j \tau_i \tau_j \mbox{ when }
 |i-j|_n=1\ .
\end{array}\right.$ \end{itemize}

\noindent From now on, the \emph{standard generators} of $\B_n$ will refer to $\tau_i$, $i\in\{0,\dots,\,n-1\}$.
Such modifications of the classical presentation of the braid group $\B_n$ have
been generalized by V. Sergiescu (cf. [S]). We do not need $\tau_0$ to define a braid homomorphism, since
$\tau_0$ can be obtained from the other generators,
but we need $\tau_0$ to consider the action of $\langle\delta\rangle$
by conjugation on the set of the $n$ standard generators. Depending on our need,
$\tau_0$ will appear explicitly in the statements or not.
\medskip

Let $D=\rho(\delta)=A_1A_2\dots A_{n-1}$. We set of course $A_0=\rho(\tau_0)=DA_{n-1}D^{-1}$
and more generally, we set $A_k=\rho(\tau_k)$.
Let $\G_0=\{A_0,\,A_1,\dots,\,A_{n-1}\}=\G\cup\{A_0\}$.
We introduce $A_0$ only now but all the statements that were true for $\G$ are also true for $\G_0$.

\bigskip

\section{The special curves $\sigma_s(\G_0)$}
    \label{sec:curves_spéciales}
\bigskip

\TITRE{Hypotheses.}\\
Let $n\geqslant 6$ an even number, let $\Sigma=\Sigma_{g,\,b}$ with $g\leqslant\frac{n}{2}$, and let
$\rho\ :\ \B_n\to\PMod(\Sigma)$ such that:

\begin{tabular}{lll}
\point & $\rho$ is non-cyclic      & by assumption,\\
\point & $\sigma_p(\G_0)=\vide$      & by assumption, inspired by Proposition \ref{prop:pas_of_curve_peripheral},\\
\point & $\sigma(\G_0)\not=\vide$    & as a consequence of the non-cyclicity of $\rho$, after Theorem \ref{thm:homomorphisms_irreducible}.
\end{tabular}
\medskip

\noindent For all $i$, we have set $A_i=\rho(\tau_i)$ .
We will study $A_1,\,A_2,\dots,\,A_{n-1}$ via their
canonical reduction systems $\sigma(A_1),\
\sigma(A_2),\dots,$ $\sigma(A_{n-1})$.
Recall that we aim to prove that $\rho$ is a transvection of monodromy homomorphism, that is,
there exists a triple $\big((a_i)_{i\leqslant n-1},\,\varepsilon,\,V\big)$
such that for all
integers $i\in\{1,\,2,\dots,\,n-1\}$, we have:
\smallskip

\centrer{$\begin{array}{rcl}
A_i&=&T_{a_i}^{\,\varepsilon}\,V,\\
\sigma(A_i)&=&\{a_i\}\cup\sigma(V).
\end{array}$}
\smallskip

\noindent The curves $a_i$ and the curves in $\sigma(V)$
play two very different roles in the description of $\rho$.
For instance, $I(a_i,\,a_{i+1})\not=0$ whereas $I(\sigma(V),\,\sigma(\G_0))=0$.
In order to stress this matter of fact, we set the following definition.
\bigskip

\TITRE{Normal curves, special curves}
\mbox{}\\We say that a curve $a$
belonging to $\sigma(\G_0)$ is \emph{special} if it satisfies
$I(a,\,\sigma(\G_0))\not=0$, and that it is \emph{normal} if it satisfies
$I(a,\,\sigma(\G_0))=0$. We denote by $\sigma_s(\G_0)$ the set of
special curves, and by $\sigma_n(\G_0)$ the set of normal curves.
Moreover, for all $A\in\G_0$, we set\smallskip

\centrer{$\sigma_n(A)=\sigma_n(\G_0)\cap\sigma(A)$,}

\centrer{$\sigma_s(A)=\sigma_s(\G_0)\cap\sigma(A)$.}
\smallskip

\noindent Thus, when $\rho$ is a transvection of monodromy homomorphism as above,
the curves $a_i$ are special whereas
the curves of $\sigma(V)$ are normal.
\bigskip

The aim of this section is to prove Proposition \ref{prop:propiétés_normales_spéciales}
which constitutes both the hardest and the most important step in the proof of
Theorem \ref{thm:Theorem_principal}. We mention in Proposition \ref{prop:propiétés_normales_spéciales}
an action of $\B_n$ on $\Curv(\Sigma)$ which will be introduced in Subsection
\ref{par:action_we_simplex_of_curves}.
\medskip

\begin{prop}[Properties of the normal curves and the special curves]
      \label{prop:propiétés_normales_spéciales}
Let $a$ be a curve of $\sigma(\G_0)$. Then:
\smallskip
\begin{itemize}
\item[\;\;(i)] The following statements are equivalent:
\begin{itemize}
  \item[\ \point] $a$ is special,
  \item[\ \point] there exists a unique $i\in\{0,\,1,\dots,\,n-1\}$ such that $a\in\sigma(A_i)$,
  \item[\ \point] $I(a,\,\delta.a)=I(a,\,\delta^{-1}.a)\not=0$,
  \item[\ \point] $|\langle\delta\rangle.a|=n$.
\end{itemize}
\smallskip

\item[\;(ii)] The following statements are equivalent:
\begin{itemize}
  \item[\ \point] $a$ is normal,
  \item[\ \point] for all $i\in\{0,\dots,\,n-1\}$, we have $a\in\sigma(A_i)$,
  \item[\ \point] $I(a,\,\sigma(\G_0))=0$,
  \item[\ \point] $|\langle\delta\rangle.a|<n$.
\end{itemize}
\smallskip

\item[(iii)] The curves of $\sigma(\G_0)$ split as follows:
\begin{itemize}
  \item[\ \point] $\sigma(\G_0)$ admits the partition:
  $\sigma(\G_0)=\sigma_s(\G_0)\sqcup\sigma_n(\G_0)$,
  \item[\ \point] $\sigma_s(\G_0)$ is nonempty and contains $n$ or $2n$
  curves, depending on whether $|\sigma_s(A_1)|=1$ or $|\sigma_s(A_1)|=2$.
\end{itemize}
\smallskip

\item[(iv)] The curves of $\sigma(\G_0)$ satisfy the following properties of stability:
\begin{itemize}
  \item[\ \point] $\sigma_n(\G_0)$ is stable under the action of $\B_n$ on $\Courb(\Sigma)$
          and the restriction of this action on $\sigma_n(\G_0)$ is cyclic:
          for all $a\in\sigma_n(\G_0)$ and all integers $i,\,j\in\{0,\dots,\,n-1\}$, we have $A_i(a)=A_j(a)$.
  \item[\ \point] $\sigma_s(\G_0)$ is stable under the action of $\J$ on $\Courb(\Sigma)$.
\end{itemize}
\end{itemize}
\end{prop}
\medskip







\subsection{Action of $\B_n$ on the simplexes of curves}
\label{par:action_we_simplex_of_curves}
\medskip

\TITRE{Cyclic actions.}
\begin{itemize}
\item An action of a group $G$ on a set $\E$ will be said
\emph{cyclic} if its structural homomorphism $\varphi\ :\ G\to\Ss(\E)$
is such that the quotient $G/\Ker(\varphi)$ is a cyclic group.

\item A \emph{cyclic action on a set $\E$}
will be an action of a cyclic group on $\E$.

\item A \emph{cyclic action on a graph $\Gamma$} will be the data of a homomorphism from
a cyclic group in $\Aut(\Gamma)$. Recall that an element in $\Aut(\Gamma)$,
the \emph{automorphisms group of $\Gamma$}, is a pair of bijections, one
on the set of the vertices, the other on the set of the edges, such that the
images of the extremities of an edge are the extremities of the
image of this edge.
\end{itemize}

\TITRE{The subgroups $\F_n$ and $\F_n^*$ of $\B_n$}
\begin{itemize}
\item There exists a unique homomorphism $\lambda\ :\ \B_n\to\ZZ$ such that $\lambda(\tau_1)=1$.
It is called sometimes \emph{the degree} or  \emph{the exponent}.
Let $[\B_n,\,\B_n]$ be the commutators group of $\B_n$. Let $\big\langle\!\langle\ \
\rangle\!\big\rangle_{\B_n}$ denote the normal closure in $\B_n$. Then we have
$\Ker(\lambda)=[\B_n,\,\B_n]=\big\langle\!\langle\,\tau_3\tau_1^{-1}\,
\rangle\!\big\rangle_{\B_n}$. We denote by $\F_n$ this group. We also define
\smallskip

\centrer{$\F_n^* = \big\langle\,\tau_i\tau_{1}^{-1}\,,\ 3\leqslant i\leqslant n-1\,\big\rangle_{\B_n}
       \subsetnot \Ker(\lambda).$}
\smallskip

\item  Let us make some remarks about these two subgroups.
\smallskip

1. Let $\varphi$ be a homomorphism from $\B_n$ in a group $G$.
The homomorphism $\varphi$ is cyclic if and only if $\varphi(\tau_1)=\varphi(\tau_2)=\dots=\varphi(\tau_{n-1})$,
if and only if the induced homomorphism by $\varphi$ on $\F_n$ is trivial. More conceptually, we could have said that
$\varphi$ is cyclic if and only if $\varphi(\B_n)$ is cyclic, if and only if
$\B_n/\Ker(\varphi)$ is cyclic, if and only if $[\B_n,\,\B_n]\subset \Ker(\varphi)$
(since $\B_n/[\B_n,\,\B_n]$ is cyclic), if and only if the induced homomorphism by $\varphi$ on $\F_n=[\B_n,\,\B_n]$ is trivial.
\smallskip

2. A homomorphism from $\F_n$ in a group $G$ is trivial if and only if the induced homomorphism by $\varphi$ on $\F_n^*$
is trivial, since $\big\langle\!\langle\,\F_n^*\,\rangle\!\big\rangle_{\B_n}=\F_n$.
\smallskip

3. The group $\F_n^*$ is clearly isomorphic to $\B_{n-2}$.
A homomorphism $\varphi$ from $\F_n^*$ will of course be said \emph{cyclic}
if $\varphi(\F_n^*)$ is cyclic, which is equivalent to saying that
$\varphi(\tau_3\tau_1^{-1})=\varphi(\tau_4\tau_1^{-1})=\dots=\varphi(\tau_{n-1}\tau_1^{-1})$.
\smallskip

4. All homomorphisms from $\B_n$ that induce by restriction to $\F_n^*$ a cyclic homomorphism are cyclic themselves.
Indeed if $\varphi$ is such a homomorphism, then the equality
$\varphi(\tau_3\tau_1^{-1})=\varphi(\tau_4\tau_1^{-1})$ holds. Consequently,
we have $\varphi(\tau_3)=\varphi(\tau_4)$ and hence $\varphi$ is cyclic.
\end{itemize}
\medskip

We sum up these remarks in the following Lemma.
\medskip

\begin{lem}[Cyclic homomorphisms/actions and the subgroups $\F_n^*$ and $\F_n$]
                \label{lem:cyclic_and_F_n}
\mbox{}\\
Let $G$ be any group, let $\varphi\ :\ \B_n\to G$ be a homomorphism and $\bar\varphi\ :\ \F_n^*\to G$ the induced homomorphism
by restriction. Then :
\begin{itemize}
\item[\;\;(i)] If $\varphi\ :\ \B_n\to G$ is cyclic, then $\bar\varphi\ :\ \F_n^*\to G$ is trivial.
\item[\;(ii)] If $\bar\varphi\ :\ \F_n^*\to G$ is cyclic, then $\varphi\ :\ \B_n\to G$ is cyclic.
\end{itemize}
Similarly, if we are given an action of $\B_n$ on any set $\E$ and its induced action of $\F_n^*$ on $\E$, then:
\begin{itemize}
\item[(iii)] If the action of $\B_n$ is cyclic, then  the action of $\F_n^*$ is trivial.
\item[(iv)] If the action of $\F_n^*$ is cyclic, then the action of $\B_n$ is cyclic.\fin
\end{itemize}
\end{lem}
\bigskip

In order to determine whether an action of $\B_n$ is cyclic or not, we have at our disposal
the following proposition due to Artin. We will use it several times to prove Proposition \ref{prop:action_we_simplex_of_curves}.
\medskip

\begin{prop}[Artin, cf. \mbox{[At3]}]
              \label{prop:fait_2}
\mbox{}\\
Let $n\geqslant5$ and $\E$ a set of $n-1$ elements or less. Then any action of
$\B_n$ on $\E$ is \emph{cyclic}.
\fin
\end{prop}
\medskip

\begin{prop}[All action of $\B_n$ on a $\B_n$-stable curve simplex is cyclic]
            \label{prop:action_we_simplex_of_curves}
\mbox{}\\Let $\A$ be a curve simplex in $\Courb(\Sigma)$
stable by the action of $\B_n$ via $\rho$ on $\Courb(\Sigma)$. Then
the actions of $\B_n$ induced by $\rho$ on $\A$, on
$\Sub_\A(\Sigma)$ and on $\Bord(\Sigma_\A)$ are cyclic.
\end{prop}
\medskip

\DEM Let us recall the
statement. Let $n$ be an even integer greater than or equal to
6, $\Sigma$ a surface $\Sigma_{g,\,b}$ where $g\leqslant
\frac{n}{2}$ and $\rho$ a homomorphism from $\B_n$ to $\PMod(\Sigma)$. Let $\A$ be a curve simplex in
$\Courb(\Sigma)$ stable by the action of $\B_n$ via $\rho$ on
$\Courb(\Sigma)$. We want to show that the actions induced by
$\B_n$ on $\A$, on $\Sub_\A(\Sigma)$ and on
$\Bord(\Sigma_{\A})$ are cyclic.
\smallskip

1. \emph{Let us show that the action of $\B_n$ on $\Sub_\A(\Sigma)$
is cyclic.}
\smallskip

\noindent Since the action of $\B_n$ on $\Courb(\Sigma)$ preserves
$\A$, the action of $\B_n$ on $\Sub(\Sigma)$ preserves
$\Sub_\A(\Sigma)$. Let us then consider the action of $\B_n$ on
$\Sub_\A(\Sigma)$.

The subsurfaces that have some common boundary components with
$\Sigma$ are fixed points of the action of $\B_n$. Let $\C$ be
the set of the subsurfaces of $\Sub_\A(\Sigma)$ that have no
common boundary component with $\Sigma$. For all $S\in\C$, we
have $\chi(S)=\chi(\for_{\bord\Sigma}(S))$. Now, the sum
$\sum_{S\in\C}\,\chi(\for_{\bord\Sigma}(S))$ is greater than or
equal to $\chi(\Sigma_{g,\,0})=2-2g$. Moreover, for all
$S\in\C$, $\chi(\for_{\bord\Sigma}(S))\leqslant -1$, hence the
cardinality of $\C$ satisfies $|\C|\leqslant 2g-2$, and finally
$|\C|\leqslant n-2$. Hence according to Proposition
\ref{prop:fait_2}, $\B_n$ acts cyclicly on the surfaces of
$\C$. Finally, $\B_n$ acts cyclicly on $\Sub_\A(\Sigma)$.
\medskip

2. \emph{Let us show that the action of  $\B_n$ on $\A$ is cyclic.}
\smallskip

\noindent Let $a$ be a curve of $\A$. We are going to study the
action of $\B_n$ on $\B_n.a$. We distinguish two cases, whether
$|\F_n.a|<n-2$ or not. In both cases, we will show that the action
of $\B_n$ on $\B_n.a$ is cyclic. This will be enough for $\A$ is a
union of orbits of curves under the action of $\B_n$.
\medskip

2.a) \emph{Case where the curve $a$ satisfies $|\F_n.a|<n-2$.}
\smallskip

\noindent Let $c$ be a curve in $\B_n.a$ and let $\gamma$ be an
element of $\B_n$ such that $c=\gamma.a$. Since $\F_n$ is normal in
$\B_n$, it follows that
$\F_n=\{\gamma\varphi\gamma^{-1},\,\varphi\in\F_n\}$. Then we have:
\smallskip

\centrer[*]{$|\F_n.c|=\big|\{\gamma\varphi\gamma^{-1}.c\,,\,\varphi\in\F_n\}\big|=
\big|\{\gamma\varphi.a\,,\,\varphi\in\F_n\}\big|=\big|\{\varphi.a\,,\,\varphi\in\F_n\}\big|=|\F_n.a|<n-2$.}
\smallskip

\noindent Let us distinguish two sub-cases, whether $n\geqslant 8$
or $n=6$:
\begin{itemize}
  \item When $n\geqslant 8$, we can apply Proposition \ref{prop:fait_2} to the action of
$\F_n^*$ (isomorphic to $\B_{n-2}$) on $\F_n.c$. Then the action of $\F_n^*$ on
$\F_n.c$ is cyclic. \smallskip

  \item When $n=6$, we cannot apply Proposition \ref{prop:fait_2}
to $\F_n^*$ for $\F_n^*$ is isomorphic to the braid group on $4$ strands
only. We have seen that the orbit of $c$ under $\F_n$ contained at
most three elements according to ($*$). Hence the action of $\F_n^*$
on $\F_n.c$ is described by a homomorphism from $\B_4$ to $\Ss_3$. Since
such a homomorphism sends the standard generators $\tau_1$, $\tau_2$ and
$\tau_3$ of $\B_4$ on three conjugate elements in $\Ss_3$, they must
be three transpositions, three 3-cycles, or three times the
identity.
\begin{itemize}
\item
If $\tau_1$, $\tau_2$ and $\tau_3$ are sent on three transpositions,
then $\tau_1$ and $\tau_3$ are sent on the same
element since they commute;
\item
If $\tau_1$, $\tau_2$ and $\tau_3$ are sent on three 3-cycles, then
the homomorphism is cyclic for the set of 3-cycles span in $\Ss_3$ a
subgroup isomorphic to $\ZZ[3]$;
\item
If $\tau_1$, $\tau_2$ and $\tau_3$ are sent on three times the
identity, then the homomorphism is trivial.
\end{itemize}
So whatever this homomorphism from $\B_4$ to $\Ss_3$ is, the elements
$\tau_1$ and $\tau_3$ have the same image. This means that in the
group $\F_n^*$, the elements $\tau_3\tau_1^{-1}$ and
$\tau_5\tau_1^{-1}$ have the same action on $\F_n.c$.
\end{itemize}
\smallskip

\noindent Finally, for any $c\in\B_n.a$ and for any even integer $n$
greater than or equal to $6$, the action of $\F_n^*$ on $\F_n.c$ is cyclic.
Hence the action of $\F_n^*$ on $\B_n.a$ is cyclic. Then, according to Lemma
\ref{lem:cyclic_and_F_n}.(iv), the action of $\B_n$ on $\B_n.a$ is cyclic.
\bigskip

2.b) \emph{Case where the curve $a$ satisfies $|\F_n.a|=m\geqslant
n-2$. }\smallskip

\noindent Let $S$ and $S'$ be the two subsurfaces (possibly equal)
of $\Sub_\A(\Sigma)$ containing the curve $a$ in their boundary.
Since the action of $\F_n$ is trivial on $\Sub_\A(\Sigma)$, the set
of curves $\F_n.a$ is included in $\Bord(S)\cap\Bord(S')$, so:
\begin{itemize}
  \item if $S\not=S'$, $S$ and $S'$ are two subsurfaces glued together along at least $m$ curves in $\Sigma$,
  so the surface resulting from the gluing is a subsurface of $\Sigma$ of genus at least $m-1$,
  so $g\geqslant m-1$,
  \item and if $S=S'$, then $S$ is a marked surface and its mark
   contains $m$ curves, so the surface resulting from the gluing
   is a subsurface of $\Sigma$ of genus at last $m$, so $g\geqslant
   m$.
\end{itemize}
Hence in both cases, $g\geqslant m-1$. However on one hand
$g\leqslant \frac{n}{2}$, on the other hand, $m\geqslant n-2$. So we
get $\frac{n}{2}\geqslant n-3$. The only possible integer
$n\geqslant 6$ that satisfies this condition is $n=6$. Then $g=3$,
$m=n-2=4$, $S\not=S'$, and $\Bord(S)\cap\Bord(S')$ is reduced to
$\F_n.a$. Moreover, the whole genus of $\Sigma$ comes from the
gluing of $S$ and $S'$ along the curves of $\F_n.a$, so it cannot
exist in $\Sub_\A(\Sigma)$ another pair of subsurfaces $(T,\,T')$
such that $\Bord(T)\cap\Bord(T')$ contains $m$ curves of $\A$. Hence
$\B_n.a=\F_n.a$. We can then apply Proposition \ref{prop:fait_2} to
the action of $\B_n$ (a braid group of order 6) on $\B_n.a$ (a
simplex of 4 curves): the action of $\B_n$ on
$\B_n.a$ is trivial.
\bigskip

3. \emph{Let us show that the action of $\B_n$ on $\Bord(\Sigma_\A)$ is cyclic.}
\smallskip

\noindent Let $a$ be a curve of $\A$ and let $a^+$ and $a^-$ be the
two boundary components of $\Sigma_\A$ coming from the cut along of
the curve $a$. According to step 2., the action of $\F_n$ on $\A$ is
trivial, so the action of $\F^*_n$ on $\Courb(\Sigma)$ via $\rho$
fixes the curve $a$, so the action of $\F^*_n$ on $\Bord(\Sigma_\A)$
via $\rho$ preserves $\{a^+,\,a^-\}$. But $\F_n^*$ is isomorphic to
$\B_{n-2}$ and $n-2>2$, so according to Proposition
\ref{prop:fait_2}, the action of $\F_n^*$ on $\{a^+,\,a^-\}$ is
cyclic. Since this is true for all curve $a$
of the set $\A$, the action of $\F_n^*$ on $\Bord(\Sigma_\A)$ is cyclic.
So according to Lemma \ref{lem:cyclic_and_F_n}.(iv),
the action of $\B_n$ on $\Bord(\Sigma_\A)$ is cyclic.
\fin
\medskip





\subsection{Curve simplexes which are orbits under $\J$ in $\sigma(\G_0)$ have less than $n$ curves}
\label{par:action_of_J}
\medskip

Recall that $n\geqslant 6$ is an even number, $\Sigma=\Sigma_{g,\,b}$ where
$g\leqslant\frac{n}{2}$, and $\rho\ :\ \B_n\to\PMod(\Sigma)$ is a noncyclic
homomorphism. Let $\delta$ be the element $\tau_1\tau_2\dots \tau_{n-1}$ of $\B_n$
and let $\J$ be the subgroup of $\B_n$ spanned by $\delta$. For all
integers $i\in\{0,\dots,\,n-1\}$, we have:
\smallskip

\centrer{$\delta\tau_i\delta^{-1}=\tau_{i+1}$.}
\smallskip

\noindent We deduce from it an action of $\J$ on $\G_0$ defined as
follows: for all integers $i\in\{0,\dots,\,n-1\}$, we set:
\smallskip

\centrer{$\delta.A_i=\rho(\delta)A_i\,\rho(\delta)^{-1}=A_{i+1}$.}
\smallskip

\noindent We define also an action of $\J$ on $\P(\G_0)$, the power
set of $\G_0$, by setting for all subsets $\K$ of $\G_0$:
\smallskip

\centrer{$\delta.\K=\{\delta.A,\,A\in\K\}=\{\rho(\delta)A\rho(\delta)^{-1},\,A\in\K\}.$}
\medskip

\TITRE{The $\J$-coloration $\sigma$ and the associated spectrum
$\sp$.}\\
A \emph{$\J$-coloration} on a $\J$-set $\E$ (i.e. a set $\E$
together with an action of $\J$ on $\E$) is a map col of $\G_0$ in
$\P(\E)$ (the power set of $\E$) compatible with the actions of $\J$
on $\G_0$ and on $\E$, i.e. such that for any $A\in\G_0$, we have:
\smallskip

\centrer{$col(\delta.A)=\delta.col(A)$.}
\smallskip

\noindent Given a $\J$-coloration col, we call \emph{spectrum}
associated to col the map of $\E$ in $\P(\G_0)$ that associates to any
element $e\in\E$ the following set $\{A\in\G_0,\ |\ e\in col(A)\}$.
\smallskip

\begin{prop}
The restriction of the map $\sigma$ from $\G_0$ in $\Courb(\Sigma)$,
which associates to any mapping class $A\in\G_0$ its canonical
reduction system $\sigma(A)\subset\Courb(\Sigma)$, is a
$\J$-coloration.
\end{prop}
\smallskip

\DEM The map $\sigma$ is a $\J$-coloration, since, for all $A\in\G_0$,
we have:

\centrer{$\begin{array}{lll}
  \delta.\sigma(A) & = & \{\delta.a,\,a\in\sigma(A)\}\\
          & = & \{\rho(\delta)(a),\,a\in\sigma(A)\}\\
          & = & \{a',\,a'\in\sigma(\rho(\delta)A\rho(\delta)^{-1})\}\\
          & = & \sigma(\delta.A).
\end{array}$} \fin
\smallskip

\begin{notation}
In this section, we will denote by $\sp$ the spectrum associated to
the $\J$-coloration $\sigma$. Thus by definition, for all
$a\in\Courb(\Sigma)$,
\smallskip

\centrer{$\sp(a)=\{A\in\G_0\ |\ a\in\sigma(A)\}$.}
\end{notation}

By considering the action of $\J$ on $\G_0$, on $\sigma(\G_0)$ and on
$\Sub_{\sigma(\G_0)}(\Sigma)$, we will show the following result.
\medskip

\TITRE{Proposition \ref{prop:cardinaln}.} {\em Let $a$ be a curve of
$\sigma(\G_0)$. Then $\J.a$ contains at most $n$ curves. The limit
case $|\J.a|=n$ can be achieved only when $\J.a$ is not a simplex.}
\medskip

\TITRE{Steps of the proof.} The proof of this proposition calls for:
\begin{itemize}
\item
 a lemma on the graphs together with a cyclic action (cf. Lemma
\ref{Z-graph}),
\item
 a lemma proposing a first version of Proposition
\ref{prop:cardinaln} (cf. Lemma \ref{lem:simplex_maximum}),
\item
 a lemma treating a special case (cf. Lemma
\ref{lem:simplex_maximum_case_particular}),
\item
 a corollary proposing a second version of Proposition
\ref{prop:cardinaln} (cf. Corollary \ref{cor:simplex_maximum2}).
\end{itemize}
The proof of Proposition \ref{prop:cardinaln} is on page
\pageref{demo:cardinaln}.
\bigskip

In the remainder of this section, we will use the following notation.
\medskip

\TITRE{Notation.}
\mbox{}
For all integers $k$ and all nonzero integer $d$, let us denote by
$[k]_d$ the remainder of the euclidian division of $k$ by $d$.
For instance, recall that we write $A_k$
instead of $A_{[k]_n}$, and $\tau_k$ instead of $\tau_{[k]_n}$.
\bigskip

\begin{figure}[!h]
 \Includegraphics{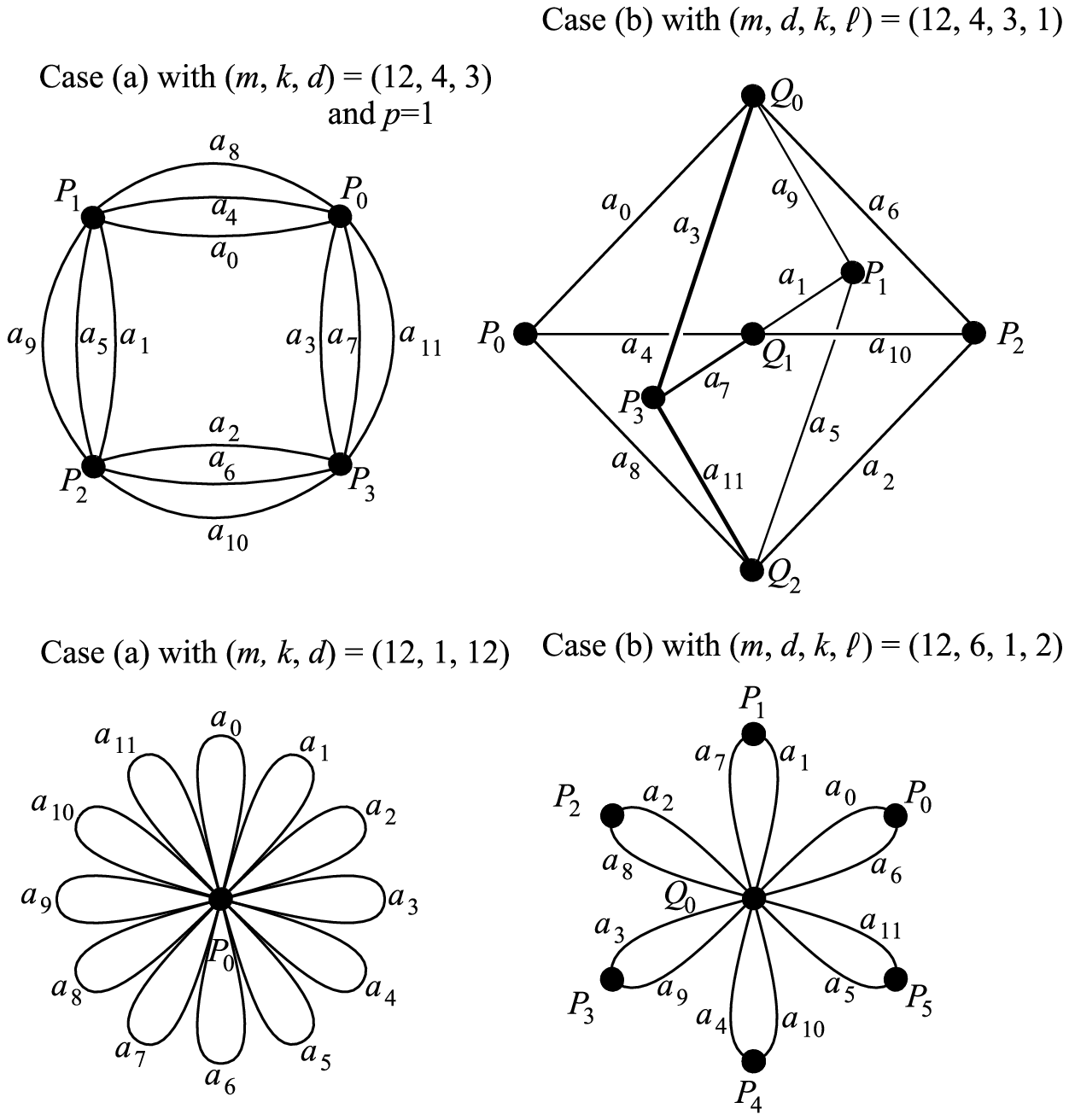}
 \caption{Four examples of graphs with 12 edges, together with a transitive $\ZZ$-action on the edges.}
 \label{fig:grapheLimite}
\end{figure}

\begin{lem}
 \label{Z-graph}
Let $\Gamma$ be a connected nonoriented graph whose number of edges
is $m$. We assume that there exists an action of $\ZZ$ on $\Gamma$,
which is compatible with its structure of graph, and which is
transitive on the set of edges. Then the pair graph-action
$(\Gamma,\,.)$ is one of the following pairs: \smallskip

\begin{itemize}
\item[(a)] The graph $\Gamma$ consists in $k$ vertices and $m$ edges
where $k$ is equal to 2 or to a divisor of $m$ (cf. Figure
\ref{fig:grapheLimite}):
\smallskip

  $\left\{
  \begin{array}{cl}
  \mbox{vertices:} & \S=\{P_i\,;\ 0\leqslant i\leqslant k-1\}\,,\\
  \mbox{edges:} & \parbox[t]{12cm}
           {$\A=\{a_0,\dots,\,a_{m-1}\}$ and there exists an integer $p$
           coprime
           with $k$ ($p=1$ if $k=1$) such that for all integers $i\in\{0,\dots,\,m-1\}$,
           the edge $a_i$ joins the vertices
           $P_{[i]_k}$ and $Q_{[i+p]_k}$. Thus, for all integers
           $i\in\{0,\dots,\,k-1\}$,
           the vertices $P_i$ and $P_{[i+p]_k}$ are
           joined by the $d$ edges
           $a_{[i]_m},\,a_{[i+k]_m},\dots,\,a_{[i+(d-1)k]_m}$,
           where $d=m$ if $k=2$ and $d=\frac{m}{k}$ otherwise.}
  \end{array}\right.$
\smallskip

The action associated to this graph $\Gamma$ is given by
$1.P_i=P_{[i+1]_k}$ and $1.a_i=a_{[i+1]_m}$ for all integers $i$
$i\in\{0,\dots,\,k-1\}$.
\medskip

\item[(b)] The graph $\Gamma$ consists in $k+\ell$ vertices and $m$ edges
where $k$ and $\ell$ are two integers greater than or equal to 1 and
coprime and $m$ is a multiple of $k\ell$ (cf. Figure
\ref{fig:grapheLimite}):
\smallskip

  $\left\{
  \begin{array}{cl}
  \mbox{vertices:} & \S=\S_1\sqcup\S_2 \mbox{ where } \S_1=\{P_0,\dots,\,P_{k-1}\} \mbox{ and } \S_2=\{Q_0,\dots,\,Q_{\ell-1}\},\\
  \mbox{edges:} & \parbox[t]{12cm}
           {$\A=\{a_0,\dots,\,a_{m-1}\}$ such that for all integers
           $i\in\{0,\dots,\,m-1\}$, the edge $a_i$ joins the vertices
           $P_{[i]_k}$ and $Q_{[i]_\ell}$. Thus, for all integers $i\in\{0,\dots,\,k-1\}$ and $j\in\{0,\dots,\,m-1\}$,
           the vertices $P_i$ and $Q_j$ are joint by $d$ edges where
           $d=\frac{m}{k\ell}$.}
  \end{array}
  \right.$
\smallskip

The action associated to this graph $\Gamma$ is given by
$1.P_i=P_{[i+1]_k}$, by $1.Q_i=Q_{[i+1]_\ell}$ and by
$1.a_i=a_{[i+1]_m}$ for all integers $i\in\{0,\dots,\,k-1\}$ and
$j\in\{0,\dots,\,m-1\}$.
\end{itemize}
\end{lem}
\bigskip

\DEM We check easily that the proposed graphs together with the
$\ZZ$-actions described in the statement exist (cf. Figure
\ref{fig:grapheLimite}) and that the actions are transitive on the
edges. Conversely, let us show that under these assumptions,
$\Gamma$ is necessarily one of the announced graphs. We begin by the
graphs with just one or two vertices:
\begin{itemize}
\item if $\Gamma$ has only one vertex and consists in a bouquet of
$m$ circles, then this is a special of case (a) with $k=1$;

\item if $\Gamma$ has exactly two vertices and if they are swapped by the $\ZZ$-action,
then we are in case (a) with $k=2$;

\item if $\Gamma$ has exactly two vertices and if they are fixed by the $\ZZ$-action,
then we are in case (b) with $(k,\,\ell)=(1,\,1)$.
\end{itemize}
\smallskip

Let us focus on the graphs together with a transitive $\ZZ$-action
on the vertices, having at least 3 vertices. Let $m$ be the number
of edges. Since the action is cyclic, $m\ZZ$ acts trivially on the
edges. Notice that any vertex $P$ can be identified by the set of
edges ending in $P$. Indeed, if two distinct vertices were the
extremities of the same edges, then by connectedness, the set of
vertices of the graph would be reduced to these two vertices, which
contradicts our hypotheses. Hence any trivial action on the edges
induces a trivial action on the vertices. Thus $m\ZZ$ acts trivially
on the set of vertices. Hence, We can quotient the action of $\ZZ$
by $m\ZZ$ and thus get an action of $\ZZ[m]$ on $\Gamma$ that acts
freely and transitively on the edges. The action of $\ZZ[m]$ on the
non ordered pairs of vertices $\{p,\,q\}$ where $p$ and $q$ are the
extremities of a same edge is hence transitive as well. We deduce
that there exist one or two orbits of vertices under the action of
$\ZZ[m]$, whether the extremities of a same edge belong to a same
orbit or not.
\smallskip

\begin{figure}[!h]
 \Includegraphics{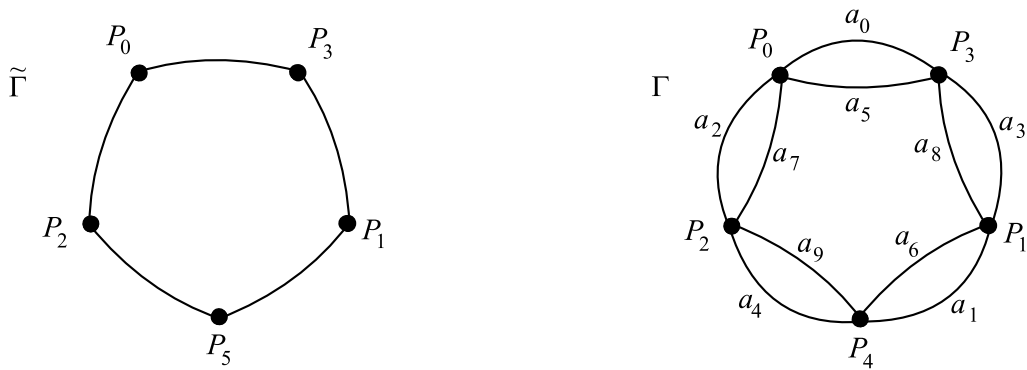}
 \caption{Example of a graph of type (a) where $k=5$, $\ell=3$, $d=2$,
$m=10$.}
 \label{fig:grapheTypeA}
\end{figure}

\TITRE{Case (a): one single orbit of vertices.} Let $k$ be the
number of vertices with $k\geqslant 3$. Since the vertices form a
single orbit under the action of $\ZZ[m]$, $k$ must divide $m$, so
the $k\ZZ[m]$-action on $\Gamma$ must fix the vertices; and for each
pair of vertices $(S_1,\,S_2)$ linked by some edge, $k\ZZ[m]$ acts
freely and transitively on the $d=m/k$ edges whose extremities are
$S_1$ and $S_2$. Let $\widetilde\Gamma$ be the graph obtained from
the graph $\Gamma$ when we identify the edges having the same
extremities. The quotient of $\ZZ[m]$ by $k\ZZ[m]$, isomorphic to
$\ZZ[k]$, acts on $\widetilde\Gamma$ and acts transitively on the
$k$ edges and the $k$ vertices of the graph $\widetilde\Gamma$. Let
us call $P_0$, $P_1$,\dots, $P_{k-1}$ the $k$ vertices of
$\widetilde\Gamma$ so that for all $\ell\in\ZZ[k]$, we have
$\ell.P_0=P_\ell$. Let $p$ be an integer in $\{1,\dots,\,k-1\}$ such
that the vertices $P_0$ and $P_p$ are joined by an edge. We obtain
the left-hand side graph in Figure \ref{fig:grapheTypeA}. Notice
that $k$ and $p$ are coprime, because the graph $\widetilde\Gamma$
(and consequently the graph $\Gamma$) would not be connected. Let us
come back to the graph $\Gamma$, we denote its vertices in the same
way: we denote by $a_0$ one of the $d$ edges ending in $P_0$ and
$P_p$, and for all $\ell\in\{1,\,2,\dots,\,m-1\}$, we denote by
$a_\ell$ the edge $\ell.a_0$ of extremities $P_{[\ell]_k}$ and
$P_{[\ell+p]_k}$. We get the right-hand side graph in Figure
\ref{fig:grapheTypeA}.
\smallskip
\smallskip

\TITRE{Case (b): two orbits of vertices.}

\noindent Let $\A$ be the set of edges. We assume now that, for all
edges $a\in\A$:
\smallskip

\centrer[1]{the extremities of $a$ belong to two different orbits.}
\smallskip

\noindent Let $a_0$ be an edge, $P_0$ and $Q_0$ the extremities of
$a_0$. Let $\S_1$ be the orbit of $P_0$ and let $k$ be the
cardinality of $\S_1$. Similarly, let $\S_2$ be the orbit of $Q_0$
and let $\ell$ be the cardinality of $\S_2$. We name the vertices of
$\S_1$ and of $\S_2$ so that for all $i\in\{0,\dots,k-1\}$,
$j\in\{0,\dots,\ell-1\}$, and $p\in\ZZ$,
\smallskip

\centrer[2]{$p.P_i=P_{[i+p]_k}$, $p.Q_i=Q_{[i+p]_\ell}$, and
$p.a_i=a_{[i+p]_m}$.}
\smallskip

\noindent We have then $\S_1=\{P_0,\dots,\,P_{k-1}\}$ and
$\S_2=\{Q_0,\dots,\,Q_{\ell-1}\}$. The integers $k$ and $\ell$ must
divide $m$, for the action of the integer $m$ on the vertices is
trivial. Since the cardinality of the orbit of $P_0$ is $k$, the
stabilizer of $P_0$ is $k\ZZ[m]$. Similarly the stabilizer of $Q_0$
is $\ell\ZZ[m]$. We deduce the following equalities between sets:
\smallskip

\centrer[3]{ $\left\{\begin{array}{l}
\{a\in\A\ |\ P_0 \mbox{ is an extremity of } a\}=(k\ZZ[m]).a_0,\\
\{a\in\A\ |\ Q_0 \mbox{ is an extremity of } a\}=(\ell\ZZ[m]).a_0,
\end{array}\right.$}
\smallskip

\noindent We will say that two edges are \emph{adjacent} if they
share at least one extremity in common. Then:\smallskip

\centrer[4]{$\{a\in\A\ |\ a \mbox{ is adjacent to }
a_0\}=\big(k\ZZ[m]\cup\ell\ZZ[m]\big).a_0$.}
\smallskip

\noindent Since the action of $\ZZ[m]$ is transitive on the edges,
equality (4) holds not only for $a_0$, but for all the edges in
$\Gamma$. Then, given a path of edges starting with the edge $a_0$,
namely a finite sequence of edges $(a'_0=a_0,\, a'_1,\dots,\,a'_r)$,
$r\geqslant 1$, such that $a'_i\cap a'_{i+1}\not=\vide$ for all
$i\leqslant r-1$, the last edge $a'_r$ must satisfy:
\smallskip

\centrer{$a'_r\in \big(k\ZZ[m] + \ell\ZZ[m]\big).a'_0$.}
\smallskip

\noindent But $\Gamma$ is connected, hence any edge of $\Gamma$ can
be seen as the last edge of some path of edges starting with $a$,
hence $k\ZZ[m] + \ell\ZZ[m]=\ZZ[m]$. Since $k<m$ and $\ell<m$ , it
follows that:
\smallskip

\centrer[5]{$k$ and $\ell$ are coprime.}
\smallskip

\begin{figure}[!h]
 \Includegraphics{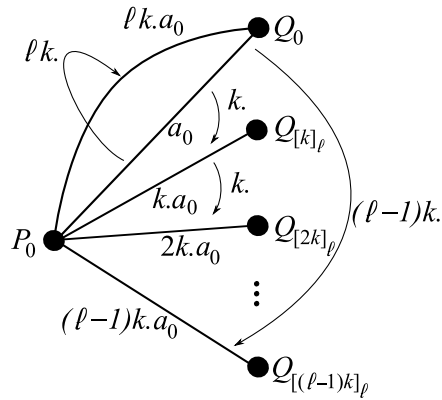}
 \caption{Action of $\ZZ$ on $a$. The points $Q_1,\,Q_2,\dots,\Q_\ell$ form the orbit of $Q_1$.}
 \label{fig:deuxOrbitesDeSommets}
\end{figure}

\noindent Let us determine $d$, the number of edges having the same
extremities as $a_0$ (cf. Figure \ref{fig:deuxOrbitesDeSommets}).
According to (3), the set of edges is
$(k\ZZ[m]).a\,\cap\,(\ell\ZZ[m]).a$, so, using (5) we get:
\smallskip

\centrer[6]{$\{\mbox{edges of extremities $P_0$ and
$Q_0$}\}=\big((k\ell)\ZZ[m]\big).a_0$,}
\smallskip

\noindent Hence we count exactly $d=\frac{m}{k\ell}$ edges
(including $a_0$)  having the same extremities as $a_0$.
\smallskip

\noindent From $m$, $k$, $\ell$, we can now describe completely
$\Gamma$ and the action of $\ZZ[m]$ on $\Gamma$. For all
$i\in\{0,\dots,\,k-1\}$ and all $j\in\{0,\dots,\,\ell-1\}$,
according to the Chinese theorem, there exists a unique integer
$u\in\{0,\dots,\,k\ell-1\}$ such that $u$ is congruent to $i$ modulo
$k$ and to $j$ modulo $\ell$. Then according to (6):
\smallskip

\centrer[7]{$\begin{array}{lll} \{\mbox{edges of extremities $P_i$
and $Q_j$}\}&=&\{(u+pk\ell).a_0,\,0\leqslant p\leqslant
d-1\}\\&=&\{a_u,\, a_{u+k\ell},\dots,\,a_{u+(d-1)k\ell}\}.
\end{array}$}
\smallskip

\noindent Thus we get the $m=dk\ell$ edges of $\Gamma$. Such a graph
is characterized by the triple $(m,\,k,\,\ell)$ or equivalently by
the triple $(d,\,k,\,\ell)$. \fin\bigskip

\begin{figure}[!h]
 \Includegraphics{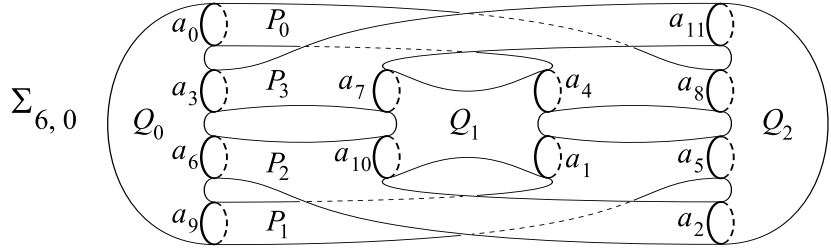}
 \caption{A simplex of $12$ curves $\{a_i, i\in\ZZ[12]\}$ in $\Sigma_{6,\,0}$,
      and a  $\ZZ$-action in $\Mod(\Sigma_{6,\,0})$
      such that for all $k\in\ZZ$ and all $i\in\{0,\dots,\,11\}$
:}

  \centrer{$k.a_i=a_{[i+k]_{12}}$,}\\
  \centrer{$k.P_i=P_{[i+k]_{4}}$,}\\
  \centrer{$k.Q_i=Q_{[i+k]_{3}}$.}
 \label{fig:case_exceptional}
\end{figure}

\begin{lem}
            \label{lem:simplex_maximum}
Let $\Sigma$ be a surface $\Sigma_{g,\,b}$. Let $\A$ be a simplex of
at least three curves in $\Sigma$. We assume that there exists a
homomorphism $\ZZ\to\PMod(\Sigma)$ whose image preserves $\A$ and that
induces a transitive action  on the curves of $\A$. Then the
cardinality of $\A$ is smaller than or equal to $2g$. The equality
$|\A|=2g$ can happen only if $g=6$ and $b=0$. In this case, the
position of the curves of $\A$ and the action of $\ZZ$ on these
curves are unique, up to homeomorphism; this case is represented in
Figure \ref{fig:case_exceptional}, where the action of $\ZZ$ is
given as follows: let us denote by $\{a_i, 0\leqslant i\leqslant
11\}$ the set of curves $\A$, and by $\{P_0,\,P_1,\,P_2,\,P_3\}$ and
$\{Q_0,\,Q_1,\,Q_2\}$ the seven subsurfaces of $\Sub_\A(\Sigma)$,
then for all $k\in\ZZ$ and all $i\in\{0,\dots,\,11\}$, we have
$k.a_i=a_{[i+k]_{12}}$, $k.P_i=P_{[i+k]_{4}}$,
$k.Q_i=Q_{[i+k]_{3}}$. \end{lem}
\medskip

\DEM
\smallskip

1. \emph{Let us show Lemma \ref{lem:simplex_maximum} in the case
where $b=0$.}\smallskip

Let $\Gamma$ be the graph $\Gamma(\Sigma\,;\,\A)$. The action of
$\ZZ$ on $\Sigma$ induces an action of $\ZZ$ on $\Gamma$ that is
transitive on the edges. Then $\Gamma$ is one of the graphs
described by Lemma \ref{Z-graph}. We are going to bound the
cardinality $|\A|$ as a function of $g$, the genus of $\Sigma$. To
do so, according to Lemma \ref{Z-graph}, we denote by:
\begin{itemize}
 \item[\point] $m=|\A|$ the number of curves of $\A$, also equal to the number
 of edges of $\Gamma$,
 \item[\point] $\S$ the set of vertices of $\Gamma$,
 \item[\point] $c$ the number of independent cycles of $\Gamma$, we
 have: $c = 1 + |\A| - |\S|$,
 \item[\point] $h$ the number of vertices in $\Gamma$ of degree 1 or
 2.
\end{itemize} Let us recall that $b=0$, hence the vertices of
degree 1 or 2 correspond in $\Sigma$ to connected components having
only one or two boundary components. Therefore, these connected
components must be of nonzero genus, whence $g\geqslant c+h$.
Moreover, by hypothesis, $m\geqslant 3$.
\smallskip

\noindent In the case (a) of Lemma \ref{Z-graph}, we set $k=|\S|$.
Then:
\begin{itemize}
  \item if $k=1$, then $c=m$, but $g\geqslant c$, so $m\leqslant g$. Therefore, $m<2g$;
  \item if $k=2$, then $c=1+m-2=m-1$, whence $g\geqslant c=m-1>\frac{m}{2}$, so $m<2g$;
  \item if $k\geqslant 3$, we then set $d=\frac{m}{k}$. We have
$c=1+m-k=1+\frac{(d-1)m}{d}$. Hence:
\begin{itemize}
\item if $d=1$, all the vertices are of degree 2, hence $h=k=m$ and we have: $g\geqslant c+k = 1+m$,
so $m<2g$;
\item if $d\geqslant 2$, then: $g\geqslant c\geqslant 1+\frac{m}{2}$,
so $m<2g$. \end{itemize} \end{itemize}
\medskip

\noindent In the case (b) of Lemma \ref{Z-graph}, let $k$ and $\ell$
be the cardinalities of the two orbits $\S_1$ and $\S_2$ of vertices
and $d=\frac{m}{k\ell}$. Even if it means swapping $\S_1$ and
$\S_2$, we can assume that $k\leqslant \ell$. Let us recall that $k$
and $\ell$ are coprime. For each triple $(d,k,\ell)$ that respects
these two conditions, let us compare $m$ and $2g$.
\begin{itemize}
\item if $(d,k,\ell)=(1,1,m)$, then the $\ell$ vertices of $\S_2$ are of degree 1, so $h=\ell=m$, so $g\geqslant h=\ell=m$ , so $m<2g$;
\item if $(d,k,\ell)=(1,2,\frac{m}{2})$, then
$c=1+m-(2+\frac{m}{2})=\frac{m}{2}-1$. Now, the $\ell=\frac{m}{2}$
vertices of $\S_2$ are of degree 2, so $h=\frac{m}{2}$ and
$g\geqslant c+h=m-1$, so $m<2g$;
\item if $(d,k,\ell)=(1,3,4)$, then $c=6$ and $m=12$, so $g\geqslant\frac{m}{2}$, so $m\leqslant2g$;
\item if $(d,k,\ell)=(1,3,5)$, then $c=8$ and $m=15$, so $g\geqslant\frac{m+1}{2}$, so $m<2g$;
\item if $(d,k,\ell)=(1,k,\ell)$ with $k=3$ and $\ell\geqslant 7$, or $k\geqslant 4$ and $\ell\geqslant 5$
(recall that $k$ and $\ell$ are coprime), then
$c=1+m(1-\frac{1}{\ell}-\frac{1}{k})$. Then
$\frac{1}{\ell}+\frac{1}{k}\leqslant \frac{1}{2}$, so $g\geqslant
c\geqslant 1+\frac{m}{2}$, so $m<2g$;
\item if $(d,k,\ell)=(2,1,\frac{m}{2})$, then a vertex is of degree $m$ and $\frac{m}{2}$
vertices are of degree 2, so $h=\frac{m}{2}$ and
$c=1+m-(1+\frac{m}{2})$, so $g\geqslant h+c=m$, so $m<2g$;
\item if $(d,k,\ell)=(d,1,\frac{m}{d})$ with $d\geqslant 3$, then we have $d+1$ vertices and $dm$
edges. So $c=1-(1+\frac{m}{d})+(dm)=\frac{d-1}{d}m$. But
$\frac{d-1}{d}>\frac{1}{2}$, so $g\geqslant c\geqslant \frac{m}{2}$
and finally, $m<2g$;
\item if $(d,k,\ell)=(m,1,1)$, then $g\geqslant c=1+m-2=m-1$. But $m\geqslant 3$, so $m<2g$;
\item if $(d,k,\ell)$ satisfies $d\geqslant 2$, $k\geqslant 2$ and $\ell\geqslant 2$, we have $c=1+m-(k+\ell)$. But
$k\ell=\frac{m}{d}\leqslant\frac{m}{2}$, so $k+\ell\leqslant
\frac{m}{2}$ (indeed, a sum of integers is always smaller than or
equal to a product of these two integers as soon as they are greater
than or equal to 2), so $g\geqslant c\geqslant 1+\frac{m}{2}$, so
$m<2g$.
\end{itemize}
\smallskip

\noindent Finally, in all the cases, $m\leqslant 2g$. The equality
case comes only in the case (b), when the triple $(d,\,k,\,\ell)$
equals $(1,\,3,\,4)$, cf. Figure \ref{fig:grapheLimite}, top-right-hand graph. The corresponding surface together with the curves of
$\A$, indexed in a self-understanding way with respect to the action
of $\ZZ$, is the surface $\Sigma_{6,\,0}$ depicted in Figure
\ref{fig:case_exceptional}.
\medskip

2. \emph{Let us show Lemma \ref{lem:simplex_maximum} in the case
where $b>0$.}\smallskip

Let us assume that we have a surface $\Sigma$ together with a
simplex $\A$ of at least three curves, and a homomorphism $\varphi$ of
$\ZZ$ to $\PMod(\Sigma)$ whose image preserves $\A$ and induces on
its curves a transitive action. Then, after having applied the map
$\sq$ to $\Sigma$ and $\A$, and after having replaced $\varphi$ by
$\sq^*\rond\varphi$, we have boiled down to the case without
boundary. Since the simplex $\A$ contains at least three curves and
since the action of $\ZZ$ induced by $\varphi$ is transitive on
$\A$, we can apply Proposition
\ref{prop:as_if_surface_without_boundary}: $\A$ and $\sq(\A)$
consist in the same number of curves. Hence $|\A|\leqslant 2g$, and
if $|\A|=2g$, then $|\sq(\A)|=2g$, so according to what we just saw
in the case without boundary, the only pair
$\big(\sq(\Sigma),\,\sq(\A)\big)$ satisfying $|\sq(\A)|=2g$ is the
one of Figure \ref{fig:case_exceptional}. Now, according to
Proposition \ref{prop:as_if_surface_without_boundary}, for all
subsurface $S\in\Sub_\A(\Sigma)$, we have:
\smallskip

\centrer{$\sq(\,\varphi(1)(S)\,)=\sq^*(\varphi(1))\,(\sq(S))$.}
\smallskip

\noindent In the case without boundary, no subsurface is
preserved by $\sq^*(\varphi(1))$, so no subsurface is preserved
by $\varphi(1)$ in the case with boundary. This is absurd for a
subsurface of $\Sub_\A(\Sigma)$ having some common boundary
components with $\Sigma$ is preserved by any mapping class of
$\PMod_\A(\Sigma)$. Hence in a surface of genus 6 with
boundary, the curves of $\A$ cannot be arranged as in Figure
\ref{fig:case_exceptional}. Therefore, when the boundary is not
empty, we have $|\A|< 2g$.\fin
\bigskip
\bigskip

In order to better understand the situation depicted in Figure
\ref{fig:case_exceptional}, when the role of $\A$ is played by the
set of curves $\sigma(\G_0)$, coming from a homomorphism $\rho$ from
$\B_{12}$ to $\Mod(\Sigma_{6,\,0})$, we prove the following lemma:
\medskip

\begin{lem}
          \label{lem:simplex_maximum_case_particular}
Let $\Sigma$ be the surface $\Sigma_{6,\,0}$ and let $\rho$ be a
homomorphism from $\B_{12}$ to $\Mod(\Sigma)$. We assume that there
exists a simplex $\A$ of 12 curves in $\sigma(\G_0)$ such that a
subgroup of $\J$ acts transitively on $\A$. Then $\rho$ is cyclic.
\end{lem}
\medskip

\DEM Let $\A$ be a simplex of 12 curves in $\sigma(\G_0)$ such that a
subgroup $\K$ of $\J=\langle\delta\rangle_{\B_{12}}$ acts
transitively on it. We want to show that $\rho$ is then cyclic. To
do so, we assume that $\rho$ is not cyclic and we look for a
contradiction (actually, we will use the fact that $\rho$ is not
cyclic to show that $\rho$ is cyclic! This is the expected
contradiction).
\smallskip

Let $\gamma$ be a generator of $\K$ and let $k$ be an integer such
that $\gamma=\delta^k$. Even if it means replacing $\A$ by
$\delta^\ell.\A$ where $\ell$ is an integer, we can assume that
$\A\cap\sigma(A_0)$ is not empty. Let $a_0$ be one of the curves of
$\A\cap\sigma(A_0)$. For all $i\in\{1,\dots,\,11\}$, let us denote
by $a_{i}$ the curve $\gamma^i.a_{0}$. It belongs to
$\sigma(A_{ki})$. Then $\A$ is the set $\{a_j,\,0\leqslant
j\leqslant 11\}$, and the surface $\Sigma$ together with the curves
of $\A$ and with the action of $\K$ on $\Sigma$ and on $\A$ is (up
to homeomorphism) the surface together with the 12 curves and with
the $\ZZ$-action depicted in Figure \ref{fig:case_exceptional}.
\medskip

1. \emph{Let us show that $k$ is coprime with 3.}
\smallskip

Let us argue by contradiction. Let us assume that $k$ is a multiple
of $3$ and let us set $k'=\frac{k}{3}$. Let us then set $\L=\langle
\delta^{2k'}\rangle$. We are going to show that $\L.a_0$ contains 18
curves, then to show that this is absurd. First, we check that:
\begin{itemize}
\item[\point]$a_6=\gamma^6.a_0\not=a_0$, so $\delta^{6k}.a_0\not=a_0$, so
$(\delta^{2k'})^{9}.a_0\not=a_0$;
\item[\point]$a_8=\gamma^8.a_0\not=a_0$, so $\delta^{8k}.a_0\not=a_0$, so
$(\delta^{2k'})^{12}.a_0\not=a_0$;
\item[\point]$\gamma^{12}.a_0=a_0$, so $\delta^{12k}.a_0=a_0$, so
$(\delta^{2k'})^{18}.a_0=a_0$.
\end{itemize}

\noindent Consequently, $\L.a_0$ contains 18 curves.
Furthermore, $\L.a_0$ is included in
$\sigma(A_0)\cup\sigma(A_2)\cup\dots\cup\sigma(A_{n-2})$, which
is a simplex, for the group
$\langle\,A_0,\,A_2,\dots,\,A_{n-2}\,\rangle$ is abelian.
Finally, $\L.a_0$ is a simplex of 18 curves in $\Sigma$, but
this is absurd because as it is well known, the greatest
simplex in $\Sigma$ contains $3g-3+b=15$ curves. Hence $k$ is
coprime with $3$.\medskip

2. \emph{Let us show that $I\big(\,\A,\,\sigma(\G_0)\,\big)=0$.}
\smallskip

\begin{figure}[!h]
 \Includegraphics{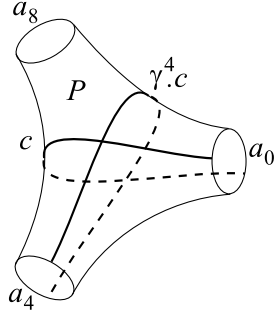}
 \caption{The curves $c$ and $\gamma^4.c$ must intersect in $P$.}
 \label{fig:intersectionDansPantalon}
\end{figure}

Let us argue by contradiction. We assume that there exists a curve
$c$ of $\sigma(\G_0)$ that intersects a curve of $\A$. Since
$\sigma(\G_0)$ is stable by $\K$, we can assume without loss of
generality that $c$ belongs to $\sigma(A_{\pm1})$ and intersects
$a_0$. According to Figure \ref{fig:case_exceptional}, there exists
a pair of pants $P$ in $\Sigma$ whose boundary components are $a_0$,
$a_4=\gamma^4.a_0$ and $a_8=\gamma^8.a_0$. But the curves $a_4$ and
$a_8$ belong to $\sigma(\gamma^{4}.A_0)\cup\sigma(\gamma^{8k}.A_0)$,
hence belong to
$\sigma(\delta^{4k}.A_0)\cup\sigma(\delta^{8k}.A_0)$. Now, $k$ is
coprime with 3 according to step 1., so $4k$ and $8k$ belong to
$4\ZZ\smallsetminus 12\ZZ$. Hence
$\sigma(\delta^{4k}.A_0)\cup\sigma(\delta^{8k}.A_0)\subset\sigma(A_4)\cup\sigma(A_8)$.
Moreover $A_{\pm1}$ commutes with $A_4$ and $A_8$, so the only
boundary component of $P$ that $c$ intersects is $a_0$, and neither
$a_4$ nor $a_8$. Let us consider now the image of this situation by
$\gamma^4$. The pair of pants $P$ is stable by $\gamma^4$, but $c$
is sent on a curve $\gamma^4.c$ that intersects only $a_4$, cf.
Figure \ref{fig:intersectionDansPantalon}. Because of a lack of room
in the pair of pants $P$, these two curves $c$ and $\gamma^4.c$ must
intersect. Yet $c$ and $\gamma^4.c$ belong respectively to
$\sigma(A_{\pm1})$ and $\sigma(A_{4k\pm1})\cup\sigma(A_{8k\pm1})$.
Since $A_{\pm1}$ commutes with $A_{4k\pm1}$ and $A_{8k\pm1}$, the
curves $c$ and $\gamma^4.c$ cannot intersect: this is a
contradiction.

\medskip

3. \emph{Let us show that $\A=\J.a_0$; in other words, everything
happens as if $\K=\J$.}\smallskip

Let us denote by $\delta\K$ the set $\{\delta\xi,\,\xi\in\K\}$ and
by $\delta.\A$ the set $\{\delta.x,\,x\in\A\}$. Let us assume that
$\K\not=\J$. Then it is clear that $\K$ and $\delta\K$ do not share
any element in common. Since $\A$ and $\delta.\A$ are the orbits of
curves under the action of $\K$, we have $\A=\delta.\A$ or
$\A\cap\delta.\A=\vide$. If $\A\cap\delta.\A=\vide$, then
$\A\cup\delta.\A$ contains 24 curves. In addition, according to step
2., we have $I(\sigma(\G_0),\,\A)=0$, so $\delta.\A\subset\sigma(\G_0)$,
so $I(\delta.\A,\,\A)=0$, so $\A\cup\delta.\A$ is a simplex.
Finally, $\A\cup\delta.\A$ is a simplex of 24 curves. But this is
absurd for in $\Sigma$, the largest simplex contains $3g-3+b=15$
curves. Hence $\A=\delta.\A$. Then $\A$ is stable by $\J$, so $\A$
(namely $\K.a_0$) is equal to $\J.a_0$. So $\K.a_0=\J.a_0$. So
everything happens as if $\K=\J$ (actually we have proven that $k$
is coprime with 12).
\medskip

\TITRE{Notation.}
\begin{itemize}
  \item[\point]
Let $D$ be the mapping class $\rho(\delta)$.

  \item[\point]
For all integers $i\in\{1,\dots,\,11\}$, let $a_i$ be the curve
$D^i(a_0)$ so that $a_i\in\sigma(A_i)$. According to step 3., we
have the equality $\A=\{a_i,\,0\leqslant i\leqslant 11\}$.

  \item[\point]
Let $Q_0$ (respectively $Q_1$, resp. $Q_2$) be the subsurface of
$\Sub_\A(\Sigma)$ that is bounded by the curves $a_0$, $a_3$, $a_6$
and $a_9$ (resp. $a_1$, $a_4$, $a_7$ and $a_{10}$, resp. $a_2$,
$a_5$, $a_8$ and $a_{11}$), cf. Figure \ref{fig:case_exceptional}.
\end{itemize}
\medskip

4.a) \emph{Let us show that $\sigma(\G_0)\subset\A\cup\Courb(Q_0\cup
Q_1 \cup Q_2)$ and for all $x\in\sigma(\G_0)\smallsetminus \A$, we
have $D^3(x)=x$.}
\smallskip

If there exists a curve $x$ belonging to $\sigma(\G_0)\smallsetminus
\A$, the curve $x$ does not intersect any curve of $\A$ according to
step 2.. Now, all the subsurfaces of $\Sub_\A(\Sigma)$, except
$Q_0$, $Q_1$ and $Q_2$, are some pairs of pants, so $x$ belongs to
$\Courb(Q_0)$, $\Courb(Q_1)$, or $\Courb(Q_2)$. Even if it means
considering $D(x)$ or $D^2(x)$ instead of $x$, we can assume that
$x$ belongs to $\Courb(Q_0)$. Then $D^3(x)$ belongs to $\Courb(Q_0)$
as well. Let $i$ be an integer of $\{0,\dots,\,11\}$ such that
$x\in\sigma(A_i)$. Then $D^3(x)$ belongs to $\sigma(A_{i+3})$. Since
$A_i$ and $A_{i+3}$ commute, $x$ and $D^3(x)$ cannot intersect. But
$\Courb(Q_0)$ does not contain any simplex of two curves, so
$D^3(x)=x$.
\medskip

4.b) \emph{Let us show that $\sigma(D)\subset\A\cup\Courb(Q_0\cup
Q_1 \cup Q_2)$ and for all $x\in\sigma(D)\smallsetminus \A$, we have
$D^3(x)=x$.}
\smallskip

According to step 3., $\A=\J.a_0$, so $\A$ is stable by the mapping
class $D$, so the curves of $\A$ are reduction curves of $D$. Hence
the curves of $\sigma(D)$ do not intersect the curves of $\A$. Since
all the subsurfaces of $\Sub_\A(\Sigma)$ except $Q_0$, $Q_1$ and
$Q_2$ are pairs of pants, the curves of $\sigma(D)$ belong to
$\A\cup\Courb(Q_0)\cup\Courb(Q_1)\cup\Courb(Q_2)$. Moreover, for all
$i\in\{0,\,1,\,2\}$, the surface $Q_i$ does not contain any simplex
of more than one curve, since in the case of $Q_i$, we have
$3g-3+b=1$. Hence $\sigma(D)\cap\Courb(Q_i)$ is empty or is reduced
to one curve. Now, $D^3$ preserves $\sigma(D)\cap\Courb(Q_i)$, so
$D^3$ must preserves this curve.
\medskip

5. \emph{Let us show that $\sigma(\G_0)\smallsetminus\A \subset
\sigma(D)\smallsetminus \A$, and let us describe
$\sigma(D)\smallsetminus \A$: if $\sigma(D)\smallsetminus \A$ is not
empty, then $\sigma(D)\smallsetminus \A$ contains three curves, one
included in $Q_0$ which we denote by $c_0$, and two other curves
$c_1=D(c_0)$ and $c_2=D^2(c_0)$ included respectively in $Q_1$ and
$Q_2$.}
\smallskip

According to the action of $D$ on $\A$, $D$ permutes $Q_0$, $Q_1$
and $Q_2$, whereas $D^3$ preserves each of them. Let us denote by
$\widehat{D^3}$ the restriction of $D^3$ on $Q_0$: $\widehat{D^3}$
belongs to $\Mod(Q_0)$. Let us focus on $\sigma(D)\smallsetminus
\A$, depending on the nature of $\widehat{D^3}$, and let us show
that $\sigma(\G_0)\smallsetminus \A\subset\sigma(D)\smallsetminus \A$.
\smallskip

\begin{itemize}
\item
  If $\widehat{D^3}$ is pseudo-Anosov, then so is $(\widehat{D^{3}})^4$.
  But $D^{12}$ is in the center of $\rho(\B_{12})$,
  so all the curves of $\sigma(\rho(\B_{12}))$ are some reduction curves of
  $D^{12}$. Therefore $Q_0$ do not contain any curve of
  $\sigma(\rho(\B_{12}))$. Thus,
\smallskip

\centrer{$\sigma(\G_0)\smallsetminus\A\,\subset\,\sigma(D)\smallsetminus\A\,=\,\vide$.}
\smallskip

\item
  If $\widehat{D^3}$ is periodic, then $\widehat{D^3}$ would be the isotopy class of a
  positive diffeomorphism of finite order according to
  Kerckhoff's Theorem. Now, according to Kerékj\`art\`o's Theorem (cf. [Kj]), such a
  diffeomorphism is conjugate to a rotation of the sphere. But if
  such a rotation, of order 4 here, preserves a curve $c$, it
  preserves also each of both hemispheres bounded by this
  curve. Hence one of these two hemispheres contains the orbit of $a_0$,
  that is to say the four boundary components $a_0$, $a_3$, $a_6$,
  $a_9$, so the other hemisphere is homeomorphic to a disk.
  Hence the curve $c$ bounds a disk. Hence $\widehat{D^3}$ does not preserve
  any curve of $\Courb(Q_0)$. Hence according to step 4.a) step 4.b),
\smallskip

\centrer{$\sigma(\G_0)\smallsetminus\A\,\subset\,\sigma(D)\smallsetminus\A\,=\,\vide$.}
\smallskip

\item
  If $\widehat{D^3}$ is reducible, let us denote by $c_0$ an essential
  reduction curve of $\widehat{D^3}$. Then
  $\sigma(\widehat{D^3})=\{c_0\}$ according to step 4.b). Hence $\J.c_0$
  is a set of three curves, one in $Q_0$, one in $Q_1$
  and one in $Q_2$. Now, any curve of $\Courb(Q_0)$ different from $c_0$ intersects $c_0$,
  so by definition of $\sigma(\widehat{D^3})$,
  this curve is not a reduction curve of $\widehat{D^3}$,
  hence it cannot be preserved by $\widehat{D^3}$.
  So according to step 4., $\sigma(\G_0)\cap\Courb(Q_0)\subset\{c_0\}$.
  Hence if $\widehat{D^3}$ is reducible, we have:
\smallskip

\centrer{$\sigma(\G_0)\smallsetminus\A\,\subset\,\sigma(D)\smallsetminus\A\,=\,\J.c_0\,=\,\{c_0,\,c_1,\,c_2\}$,}
\smallskip

\noindent where $c_1=D(c_0)$ and $c_2=D^2(c_0)$. \end{itemize}
\bigskip

In order to discuss later (in steps 7. and 8.) about the stability
of $\A$ under the action of $\B_{12}$, we are going to study in step
6. the stability of $\J.c_0$ under the action of $\B_{12}$. Of
course the set $\J.c_0$ has some meaning if the curve $c_0$ is
defined, that is to say when $\sigma(D)\smallsetminus\A\not=\vide$
according to step 5..
\medskip

6. \emph{Let us show that if $\sigma(D)\smallsetminus\A\not=\vide$,
then the action of $\B_{12}$ via $\rho$ on $\Courb(\Sigma)$
preserves $\J.c_0$.}
\smallskip

By definition of $c_0$ according to step 5., $c_0$ belongs to
$\sigma(D)$. Since $\sigma(D)$ is stable by $D$, $\sigma(D)$
contains $\J.c_0$, and if $a_0$ belongs to $\sigma(D)$, then
$\sigma(D)$ contains $\J.a_0$ that is equal to $\A$, whereas if
$a_0$ does not belong to $\sigma(D)$, then $\sigma(D)\cap\A=\vide$.
Now, according to step 5., $\sigma(D)$ is included in
$\J.c_0\cup\A$. Hence $\sigma(D)=\J.c_0$ or
$\sigma(D)=\J.c_0\cup\A$. But, for all $\xi\in\B_{12}$, the mapping
class $\rho(\xi)$ commutes with $D^{12}$, so
$\rho(\xi)(\sigma(D^{12}))=\sigma(D^{12})$, and so
$\rho(\xi)(\sigma(D))=\sigma(D)$, since $\sigma(D^{12})=\sigma(D)$.
Hence the action of $\B_{12}$ via $\rho$ on $\Courb(\Sigma)$
preserves the curves of $\sigma(D)$. Then if $\sigma(D)=\J.c_0$, we
have shown that the action of $\B_{12}$ via $\rho$ on
$\Courb(\Sigma)$ preserves $\J.c_0$. We still have to study the case
where $\sigma(D)=\J.c_0\cup\A$. We are going to show that:
\smallskip

\centrer[$1$]{Any mapping class that preserves $\J.c_0\cup\A$
preserves $\J.c_0$ and preserves $\A$.}
\smallskip

\noindent Since the action of $\B_{12}$ via $\rho$ on
$\Courb(\Sigma)$ preserves $\sigma(D)=\J.c_0\cup\A$, it preserves
$\J.c_0$. So, proving (1) is enough to show step 6..
\medskip

\begin{figure}[!h]
 \Includegraphics{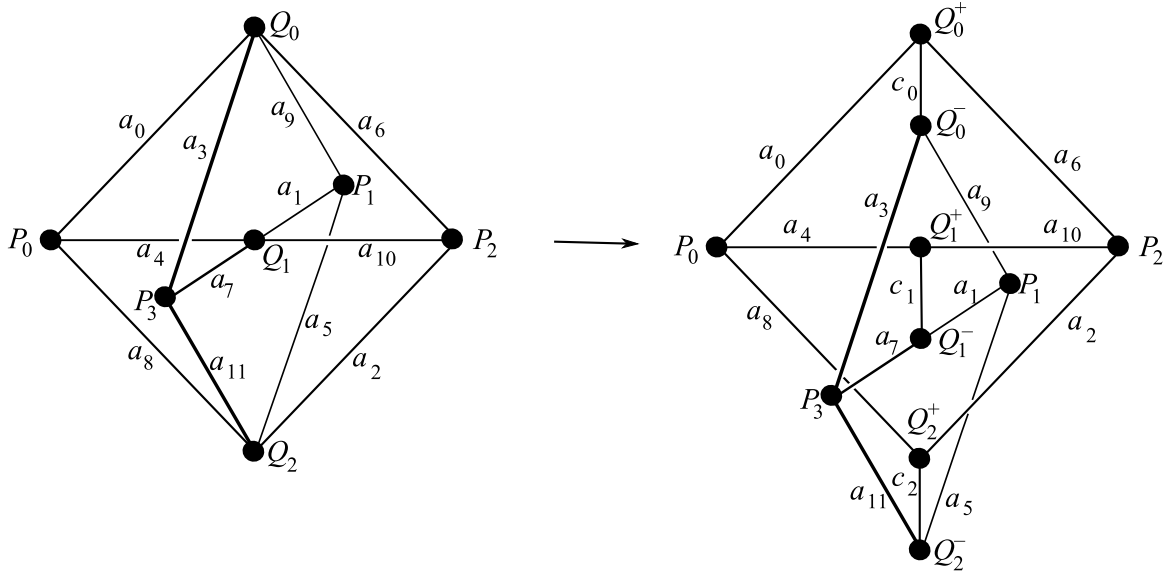}
 \caption{The graph $\Gamma(\Sigma,\,\A)$ on the left-hand side, the graph $\Gamma(\Sigma,\,\J.c_0\cup\A)$ on the right-hand side.}
 \label{fig:apparitionDeC0}
\end{figure}

Let us show assertion ($1$). The curve $c_0$ lies in $Q_0$, is
separating in $Q_0$ and so induces a partition of $\Bord(Q_0)$ in
two subsets: the boundary components located on an edge of $c_0$ and
the boundary components located on the other one. Since the curve
$c_0$ is stable by $D^3$, this partition must be stable by $D^3$.
The boundary components of $Q_0$ are the curves $a_0$, $a_3$, $a_6$,
$a_9$, and their images by $D^3$ are respectively $a_3$, $a_6$,
$a_9$, $a_0$, so this partition can only be
$\{a_0,\,a_6\}\sqcup\{a_3,\,a_9\}$. Indeed, it is clear that the two
other partitions $\{a_0,\,a_3\}\sqcup\{a_6,\,a_9\}$ and
$\{a_0,\,a_9\}\sqcup\{a_3,\,a_6\}$ are not stable by $D^3$. Let us
consider the graph $\Gamma(\Sigma,\,\sigma(D))$, cf. Figure
\ref{fig:apparitionDeC0} on the right-hand side.
 We see that the smallest injective cycle of
edges containing $a_0$ contains four edges, for example the cycle
$a_0$, $a_4$, $a_{10}$, $a_6$, whereas the smallest injective cycles
of edges containing $c_0$ contains six edges, for example the cycle
$a_0$, $c_0$, $a_3$, $a_{11}$, $c_2$, $a_8$. Then, any mapping class
preserving $\J.c_0\cup\A$ induces an action on the graph
$\Gamma(\Sigma,\,\J.c_0\cup\A)$, but cannot swap in this graph two
edges such that $a_0$ and $c_0$, since they have different
combinatoric properties, as we have just seen it. Hence any mapping
class that preserves $\J.c_0\cup\A$ preserves $\J.c_0$ and preserves
$\A$: the statement ($1$) is proved.
\bigskip

7. \emph{Let us show that if the action of $\B_{12}$ on
$\Courb(\Sigma)$ preserves $\A$, then $\rho$ is cyclic.}
\smallskip

Since we assume that the action of $\B_{12}$ on $\Courb(\Sigma)$
preserves the curve simplex $\A$, we can apply Proposition
\ref{prop:action_we_simplex_of_curves} and conclude that $\B_{12}$
acts cyclicly on $\A$, and so $\F_{12}$ acts trivially on $\A$. Let
us consider the subgroup $\F_{12}^*$ of $\F_{12}$ which is
isomorphic to $\B_{10}$. Then for any surface $S\in\Sub_\A(\Sigma)$,
the homomorphism $\rho$ restricted to $\F_{12}^*$ induces a homomorphism in
$\PMod(S)$. Since $S$ is of genus zero and $\F_{12}^*$ is isomorphic
to $\B_{10}$, we can apply Theorem \ref{thm:BndansMCG0b} and
conclude that the restriction of $\rho$ to $\F_{12}^*$ induces in
$\PMod(S)$ a cyclic homomorphism. Hence the restriction of $\rho$ to
$\F_{12}^*$ induces in $\PMod(\Sigma_\A)$ a cyclic homomorphism. Thus
the group $\rho(\F^*_{12})$ included in $\P_\A\Mod(\Sigma)$ is sent
in $\PMod(\Sigma_\A)$ on an abelian subgroup $\widehat G$. Let us
consider the below diagram.
\smallskip

\centrer{$\begin{diagram}
   \node{\rule{3cm}{0cm}\rho(\F^*_{12})\ \subset}
 \node{\P_\A\Mod(\Sigma)} \arrow{e,t}{\cut_\A}
 \node{\PMod(\Sigma_\A)}
  \node{\supset \widehat G\rule{3cm}{0cm}}\\
  \node{\rule{3cm}{0cm}G\ \subset}
 \node{\P_\A\Mod(\Sigma)}\arrow{ne,t}{\cut_\A}
 \node{\Mod(\Sigma_\A,\,\bord\Sigma_\A)}\arrow{w,b}{rec_\A}\arrow{n,r}{\for_{\bord\Sigma_\A}}
  \node{\supset \widetilde G=\for_{\bord\Sigma_\A}^{\:\:-1}(\widehat G)}
\end{diagram}$}
\smallskip

\noindent Starting from $\widehat G$, we successively define the
groups:
\smallskip

\begin{itemize}
\item[\point] $\widetilde G=\for_{\bord\Sigma_\A}^{\:\:-1}(\widehat
G)$ where $\for_{\bord\Sigma_\A}\;:\;
\Mod(\Sigma_\A,\,\bord\Sigma_\A)\to\PMod(\Sigma_\A)$ is the
canonical ``forget'' homomorphism. According to Lemma
\ref{lem:ça_commute}, $\widetilde G$ is abelian,
\smallskip

\item[\point] $G=rec_\A(\widetilde G)$ where $rec_\A\;:\; \Mod(\Sigma_\A,\,\bord\Sigma)\to\P_\A\Mod(\Sigma)$
is the gluing homomorphism along the curves of $\A$. The group $G$ is
abelian since $\widetilde G$ is. \end{itemize}
\smallskip

\noindent By construction, $\rho(\F_{12}^*)$ is included in $G$, an
abelian group, so $\rho(\F_{12}^*)$ is abelian. Since $\F_{12}^*$ is
isomorphic to $\B_{10}$ and since the abelianization of $\B_{10}$ is
cyclic, if follows that $\rho(\F_{12}^*)$ is cyclic. So the mapping
class
$\rho(\tau_3\tau_4^{-1})=\rho(\tau_3\tau_1^{-1})\big(\rho(\tau_4\tau_1^{-1})\big)^{-1}$
coincides with the identity. Hence $\rho(\tau_3)=\rho(\tau_4)$, so
according to Lemma \ref{lem:homomorphism_cyclic}, the homomorphism $\rho$ is
cyclic.
\bigskip

8. \emph{Let us show that if $\rho(\B_{12})$ does not preserve $\A$,
then, again, $\rho$ is cyclic.}
\smallskip

If $\rho(\B_{12})$ does not preserve $\A$, it is clear that there
exists a curve $a'\in\A$ and a mapping class $F\in\G_0$ such that
$F(a')$ does not belong to $\A$. Even if it means conjugating $F$ by
a power of $D$, we can assume without loss of generality that
$a'=a_0$. For all $i\in\{0\}\cup\{2,\,3,\dots,\,10\}$, the mapping
class $A_i$ commutes with $A_0$, hence preserves $\sigma(A_0)$,
hence sends $a_0$ in $\sigma(\G_0)$. But, according to step 5.,
$\sigma(\G_0)$ is included in $\A\cup\sigma(D)$, that is to say that
$\sigma(\G_0)$ is included in $\A$ if
$\sigma(D)\smallsetminus\A=\vide$, or in $\A\cup\J.c_0$ if
$\sigma(D)\smallsetminus\A\not=\vide$. If
$\sigma(D)\smallsetminus\A\not=\vide$, then $\J.c_0$ is stable by
$\rho(\B_{12})$ according to step 6., so $A_i$ cannot send $a_0$ in
$\J.c_0$. Hence in all the cases, $A_i$ sends $a_0$ in $\A$. But we
have seen that $F(a_0)$ did not belong to $\A$, so $F$ must belong
to $\{A_1,\,A_{11}\}$. Even if it means conjugating $F$ by
$A_{11}(A_0A_{11})(A_1A_0A_{11})$, we can assume without loss of
generality that $F=A_1$. Let us sum up: we have shown that
\smallskip

\centrer[2]{if $\rho(\B_{12})$ does not preserve $\A$, then
$A_1(a_0)$ does not belong to $\A$.} \smallskip

\noindent This implies that
\smallskip

\centrer[3]{$\sp(a_0)\subset\{A_0,\,A_2\}$,}
\smallskip

\noindent for otherwise, it would exist $j\in\{3,\,4,\dots,\,11\}$
such that $a_0\in\sigma(A_j)$ and as we have just seen it, we would
deduce that $A_1(a_0)$ would belong to $\A$, which would contradict
(2). Remember that by definition of $a_0$, $A_0$ belongs to
$\sp(a_0)$. Hence (3) implies that
$\sp(a_0)\in\big\{\{A_0\},\,\{A_0,\,A_2\}\big\}$, whence by
conjugation,
\smallskip

\centrer[4]{$or$ $\left.\begin{array}{l}
  \mbox{for all } i\in\{0,\,1,\dots,\,n-1\},\ \mbox{ we have } \sp(a_i)=\{A_i\}\;\\
  \mbox{for all } i\in\{0,\,1,\dots,\,n-1\},\ \mbox{ we have } \sp(a_i)=\{A_i,\,A_{i+2}\}.
\end{array}\right\}$}

\noindent For all $\ell\in\{4,\,5,\dots,\,11\}$ and all
$\varepsilon\in\{1,\,2\}$, $A_{\varepsilon}$ and $A_\ell$ commute,
so $A_{\varepsilon}(a_\ell)$ belongs to $\sigma(A_\ell)$, which is
itself included in $\sigma(\G_0)$. According to step 5., $\sigma(\G_0)$
is included in $\A\cup\sigma(D)$, so $A_{\varepsilon}(a_\ell)$
belongs to $\A$ or possibly to $\J.c_0$ if
$\sigma(D)\smallsetminus\A\not=\vide$. But we have seen that
$\J.c_0$ (if it exists) is stable by the action of $\B_{12}$
according to step 6., hence $A_{\varepsilon}(a_\ell)$ belongs to
$\A$. Let us show that $A_{\varepsilon}(a_\ell)=a_\ell$. We argue
differently, depending on the two cases mentioned by (4).
\smallskip

\begin{itemize}
\item
If $\sp(a_\ell)=\{A_\ell\}$, then the spectrum of each curve of $\A$
is reduced to a singleton, and the canonical reduction system of
each mapping class of $\G_0$ contains only one element of $\A$. Since
$A_{\varepsilon}(a_\ell)\in\A\cap\sigma(A_\ell)$, it follows that
$A_{\varepsilon}(a_\ell)=a_\ell$.
\smallskip

\item
If $\sp(a_\ell)=\{A_\ell,\,A_{\ell+2}\}$, then $a_\ell$ and
$a_{\ell-2}$ are the only curves of $\A$ in $\sigma(A_\ell)$. Hence
$A_{\varepsilon}(a_\ell)=a_\ell$ or
$A_{\varepsilon}(a_\ell)=a_{\ell-2}$. If
$A_{\varepsilon}(a_\ell)=a_{\ell-2}$, we would have
$\sp(A_{\varepsilon}(a_\ell))=\{A_{\ell-2},\,A_\ell\}$. But
$A_{\varepsilon}$ commutes with at least one of the two mapping
classes $A_{\ell-2}$ or $A_{\ell+2}$. If $A_{\varepsilon}$ commutes
with $A_{\ell-2}$, the fact that $A_{\ell-2}\not\in\sp(a_\ell)$
implies that $A_{\ell-2}\not\in\sp(A_{\varepsilon}(a_\ell))$;
whereas if $A_{\varepsilon}$ commutes with $A_{\ell+2}$, the fact
that $A_{\ell+2}\in\sp(a_\ell)$ implies that
$A_{\ell+2}\in\sp(A_{\varepsilon}(a_\ell))$. In the two cases, the
fact that $\sp(A_{\varepsilon}(a_\ell))=\{A_{\ell-2},\,A_\ell\}$ is
contradicted, so $A_{\varepsilon}(a_\ell)\not=a_{\ell-2}$, so
$A_{\varepsilon}(a_\ell)=a_{\ell}$. \end{itemize}\smallskip

\noindent Hence $A_{\varepsilon}$ preserves each curve $a_\ell$,
$4\leqslant \ell\leqslant11$. Let us set\smallskip

\centrer{$\A'=\{\,a_\ell,\ 4\leqslant \ell\leqslant11\,\}$.}
\smallskip

\begin{figure}[!h]
 \Includegraphics{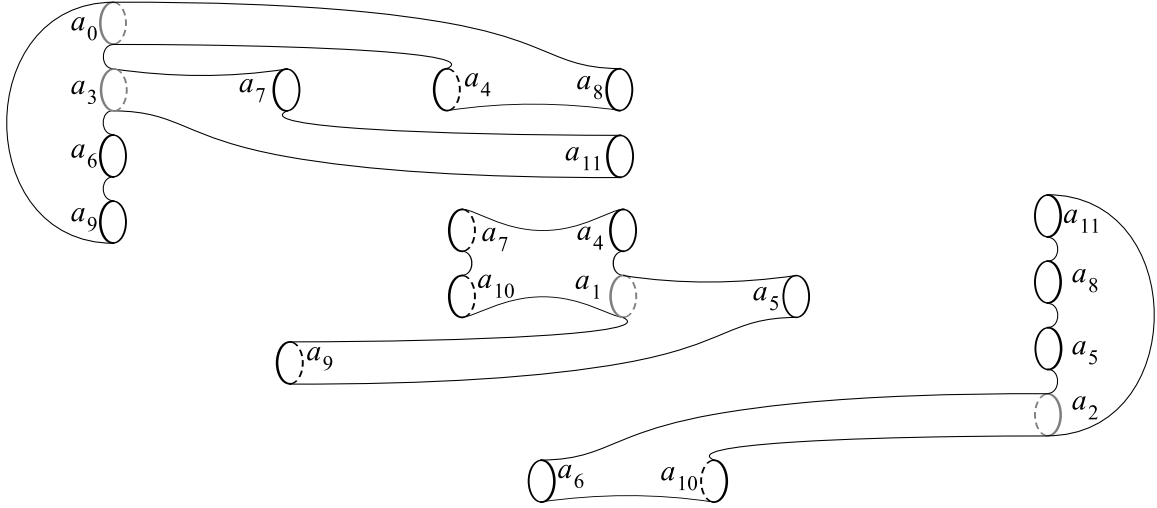}
 \caption{The surface $\Sigma_{\A'}$ is a disjoint union of genus-0 surfaces.}
 \label{fig:composantesGenreNul}
\end{figure}

\noindent Let us consider the surface $\Sigma_{\A'}$ Figure
\ref{fig:composantesGenreNul}. Since $A_{\varepsilon}$ preserves
each curve of $\A'$, $A_{\varepsilon}$ cannot permute the connected
components of $\Sigma_{A'}$. Hence $A_{\varepsilon}$ (recall that
$\varepsilon$ ranges over $\{1,\,2\}$) induces a mapping class
$\tildeA_\varepsilon$ in $\PMod(\Sigma_{\A'})$. But we have the
following canonical isomorphism:
\smallskip

\centrer{$\PMod(\Sigma_{\A'})=\dessous{\prod}{S\in\Comp(\Sigma_{\A'})}\PMod(S)$,}
\smallskip

\noindent where $\Comp(\Sigma_{\sigma(\G_0)})$ denotes the set of all
connected components of $\Sigma_{\sigma(\G_0)}$. In other words,
$\PMod(\Sigma_{\A'})$ is a direct product of mapping class groups of
surfaces of genus zero (cf. Figure \ref{fig:composantesGenreNul}).
 For any connected component $S$ of $\Sigma_{\A'}$, the
images of $\tildeA_1$ and $\tildeA_2$ induced in $\PMod(S)$ satisfy
a braid relation. But $S$ is of genus zero, so we can apply Theorem
\ref{thm:BndansMCG0b}. Hence $\tildeA_1=\tildeA_2$. Finally, $A_1$
and $A_2$ induce the same mapping class in $\PMod(\Sigma_{\A'})$, so
according to the following central exact sequence:
\smallskip

\centrer{$1\to\langle T_{a_i},\,0\leqslant i\leqslant
n-1\rangle\to\PMod_{\A'}(\Sigma)\to\PMod(\Sigma_{\A'})\to1$,}
\smallskip

\noindent the mapping classes $A_1$ and $A_2$ differ from a
multitwist that is central in $\PMod_{\A'}(\Sigma)$, so $A_1$ and
$A_2$ commute in $\Sigma$. But $A_1$ and $A_2$ satisfy also a braid
relation, so $A_1$ and $A_2$ have to be equal. In other words,
$\rho$ is cyclic. \fin
\bigskip
\bigskip

\begin{cor}
      \label{cor:simplex_maximum2}
Let $n$ be an integer greater than or equal to 6, $\Sigma$ a
surface $\Sigma_{g,\,b}$ where $g\leqslant \frac{n}{2}$ and
$\rho$ a noncyclic homomorphism from $\B_n$ to $\PMod(\Sigma)$.
Then, any curve simplex on which a subgroup of $\J$ acts
transitively contains strictly less than $2g$ curves. \end{cor}
\medskip

\DEM According to Lemma \ref{lem:simplex_maximum}, such a simplex
$\A$ contains strictly less than $2g$ curves, except if $g=6$ and
$b=0$ where $\A$ can contain 12 curves. But according to Lemma
\ref{lem:simplex_maximum_case_particular}, the fact that $\rho$ is
not cyclic forbids that $\A$ contains 12 curves. This proves the
corollary.\fin
\bigskip

We can now prove Proposition \ref{prop:cardinaln}. Let us recall
that $n\geqslant 6$, $\Sigma=\Sigma_{g,\,b}$ where $g$ and $b$ are
some integers such that $g\leqslant n/2$, and $\rho$ is a noncyclic
homomorphism from $\B_n$ to $\PMod(\Sigma)$.
\medskip

\begin{prop}
                \label{prop:cardinaln}
Let $a$ be
a curve of $\sigma(\G_0)$. Then $\J.a$ contains at most $n$ curves.
The limit case $|\J.a|=n$ takes place if and only if $\J.a$ is not a
simplex.
\end{prop}
\bigskip

\label{demo:cardinaln}

\DEM Let $a$ be a curve of $\sigma(\G_0)$. If the orbit $\J.a$ is a
curve simplex, Corollary \ref{cor:simplex_maximum2} can be applied:
$\J.a$ contains strictly less than $2g$ curves, hence  strictly less
than $n$ curves.

If $\J.a$ is not a simplex, let us show that $|\J.a|\leqslant n$. We
can assume without loss of generality that $a\in\sigma(A_0)$, in
other words, $\sp(a)\supset\{A_0\}$. For all
$k\in\{0\}\cup\{2,\,3,\dots,\,n-2\}$, the mapping classes $A_0$ and
$A_{k}$ commute, so according to Proposition
\ref{prop:properties_sigma}.(iii), $I(\sigma(A_0),\,\sigma(A_k))=0$, so
$I(a,\,\delta^k.a)=0$. But $\J.a$ is not a simplex. Hence $a$
intersects one of the curves of $\J.a$ that belongs necessarily to
$\sigma(A_{n-1})$ or to $\sigma(A_1)$. These two cases are symmetric
and we can assume without loss of generality that
$I(a,\,\sigma(A_1))\not=0$. Then, since $A_0$ and $A_2$ are the only
mapping classes that do not commute with $A_1$, and since
$I(a,\,\sigma(A_1))\not=0$, we deduce that
$\sp(a)\subset\{A_0,\,A_2\}$. Let $k$ be the least positive integer
such that $\delta^k.a=a$. The integer $k$ satisfies
$\delta^k.\sp(a)=\sp(a)$, but since $\sp(a)\subset\{A_0,\,A_2\}$,
then $k$ must be a multiple of $n$. Let $p$ be the integer
$\frac{k}{n}$, so that $|\J.a|=pn$.

Let us denote by $\K=\langle\delta^2\rangle$. Since $|\J.a|=pn$, we
have $|\K.a|=p\frac{n}{2}\geqslant pg$. But the set
$\K.A_0=\{A_0,\,A_2,\,A_4,\dots\,A_{n-2}\}$ consists on elements that commute,
so $\sigma(\K.A_0)$ is a simplex. Hence $\K.a$, that is included in
$\sigma(\K.A_0)$, is a simplex, too. In other words, $\K.a$ is a
curve simplex on which $\K$, a subgroup of $\J$, acts transitively.
Hence according to Corollary \ref{cor:simplex_maximum2},
$|\K.a|<2g\leqslant n$. Hence $p=1$ and $|\J.a|=n$.\fin
\bigskip





\subsection{Partition of $\sigma(\G_0)$: $\sigma(\G_0)=\sigma_s(\G_0)\sqcup\sigma_n(\G_0)$}
                \label{par:proof_spéciales_normales}
\medskip

\begin{prop}[Partition of $\sigma(\G_0)$]
            \label{prop:spéciales_normales}
\mbox{}\\Any curve $a$ belonging to
$\sigma(\G_0)$ satisfies either all the left-hand side properties (1g)
- (6g),
or all the right-hand side properties (1d) - (6d).\\
\begin{tabular}{lcccr}
(1g)\vlblanc  &\rule{2cm}{0cm} $I(a,\,\delta.a)=0$ \rule{2cm}{0cm} &; & \rule{2cm}{0cm}$I(a,\,\delta.a)\not=0$ \rule{2cm}{0cm} & (1d)\\
(2g)\vlblanc                 & $|\sp(a)|\geqslant 2$ &; & $|\sp(a)|=1$ &                    (2d)\\
(3g)\vlblanc            & $I(a,\,\sigma(\G_0))=0$ &; & $I(a,\,\sigma(\G_0))\not=0$ &      (3d)\\
(4g)\vlblanc & $\forall k,\ \sp(a)\not\subset\{A_k,\,A_{k+2}\}$ &; & $\exists k\ |\ \sp(a)\subset\{A_k,\,A_{k+2}\}$ &(4d)\\
(5g)\vlblanc                      & $|\J.a|<n$ &; & $|\J.a|= n$ &                       (5d)\\
(6g)\vlblanc                & $\J.a$ is a simplex &; & $\J.a$ is not a simplex &               (6d)\\
\end{tabular}
\end{prop}
\bigskip

Notice that each of the six lines of this table contains two
opposite assertions (knowing that $|\J.a|\leqslant n$, as seen in
Proposition \ref{prop:cardinaln}), so that any curve $a$ satisfies
exactly one assertion per line. We are going to show that all the
right-hand side assertions are equivalent (and so all the
left-hand side assertions are equivalent).
\bigskip

\DEM

\Point Let us first show the
cycle of implications (1d) $\Rightarrow$ (3d) $\Rightarrow$ (4d)
$\Rightarrow$ (5d) $\Rightarrow$ (6d) $\Rightarrow$ (1d).
\smallskip

\TITRE{(1d) $\Rightarrow$ (3d).} This first implication is trivial.
\smallskip

\TITRE{(3d) $\Rightarrow$ (4d).} Since $I(a,\,\sigma(\G_0))\not=0$, there exists an integer $k$
such that $I(a,\,\sigma(A_{k+1}))\not=0$. But for all
$i\in\{0,\dots,\,n-1\}\smallsetminus\{k,\,k+2\}$, the mapping
classes $A_{k+1}$ and $A_i$ commute, so according to Proposition
\ref{prop:properties_sigma}.(iii), $I(\sigma(A_1),\,\sigma(A_i))=0$. Hence
$a\not\in\sigma(A_i)$. Thus
$\sp(a)\subset\{A_k,\,A_{k+2}\}$.
\medskip

\TITRE{(4d) $\Rightarrow$ (5d)}
Let us assume for example that there exists an integer $k$ such
that $\sp(\delta^k.a)=\{A_0,\,A_2\}$ (the situation is even more
simple if $\sp(\delta^k.a)$ is a singleton). Then for all integers
$i$, we have $\sp(\delta^{k+i}.a)=\{A_i,\,A_{i+2}\}$. But
$\delta^i.a=a$ only if $\sp(\delta^{k+i}.a)=\sp(\delta^k.a)$, hence
only if $i$ is a multiple of $n$. So $|\J.a|\geqslant n$. So,
according to Proposition \ref{prop:cardinaln},
$|\J.a|=n$.\smallskip

\TITRE{(5d) $\Rightarrow$ (6d)}
This is the second part of Proposition \ref{prop:cardinaln}.
\smallskip

\TITRE{(6d) $\Rightarrow$ (1d)}
Let us start from a curve $a\in\sigma(\G_0)$, and let us denote
by $\A$ its orbit $\J.a$. We have:\smallskip

\centrer{$\displaystyle\begin{array}{lll}
 I(\A,\,\A) & = & \vlblanc\displaystyle  \sum_{0\leqslant i,j\leqslant n-1} I(\delta^i.a\,,\, \delta^j.a)\\
      & = & \vlblanc\displaystyle  \sum_{\scriptsize\begin{array}{c}0\leqslant i\leqslant n-1\\j\in\{i-1,\,i+1\}\end{array}}
                     I(\delta^i.a\,,\, \delta^{j}.a)\\
      & = & \vlblanc\displaystyle n\,I(a\,,\, \delta.a) + n\,I(a\,,\,\delta^{-1}.a)\\
      & = &     \displaystyle 2n\,I(a\,,\, \delta.a).
\end{array}$}
\smallskip

\noindent Hence if $\A$ is not a simplex, in other words if
$I(\A\,,\, \A)\not=0$, then $I(a\,,\, \delta.a)\not=0$.
\bigskip

\Point We terminate the proof of Proposition \ref{prop:spéciales_normales}
by showing the implications (1d) $\Rightarrow$ (2d) $\Rightarrow$ (4d).
\smallskip

\TITRE{(1d) $\Rightarrow$ (2d)}
If $I(a,\,\delta.a)\not=0$, then the action of $\delta^{-1}$ on
the pair $(a,\,\delta.a)$ implies that $I(\delta^{-1}.a,\,a)\not=0$.
Let $i$ be an integer in $\{0,\,1,\dots,\,n-1\}$ such that
$a\in\sigma(A_i)$. Then we have $\delta.a\in\sigma(A_{i+1})$ and
$\delta^{-1}.a\in\sigma(A_{i-1})$. Hence according to Proposition
\ref{prop:properties_sigma}.(iii), the inequality $I(a,\,\delta.a)\not=0$
implies the inclusion $\sp(a)\subset\{A_i,\,A_{i+2}\}$ and
similarly, the inequality $I(\delta^{-1}.a,\,a)\not=0$ implies the
inclusion $\sp(a)\subset\{A_{i-2},\,A_i\}$. Finally, we have
$\sp(a)=\{A_i\}$.\smallskip

\TITRE{(2d) $\Rightarrow$ (4d)} This last implication is trivial.
\fin
\bigskip

\TITRE{Fundamental remark.} Recall that we have defined the special curves as being the curves $a\in\sigma(\G_0)$
satisfying $I(a,\,\sigma(\G_0))\not=0$, and the normal curves as being the curves $a\in\sigma(\G_0)$
satisfying $I(a,\,\sigma(\G_0))=0$. In other words,
the special curves are those satisfying the right-and side assertions of Proposition
\ref{prop:spéciales_normales} and the normal curves are those
satisfying the left-hand side assertions of Proposition
\ref{prop:spéciales_normales}.
\medskip

\subsection{Stability and existence results} \label{par:proof_stability_existence}
\medskip

This subsection is devoted to Propositions
\ref{prop:stability_des_curves_spéciales}, and
\ref{prop:stability_des_curves_normales} concerning the stability of
the normal and the special curves , and to Propositions
\ref{prop:existence_des_curves_spéciales} and
\ref{prop:cardinal_de_sigmas} concerning the existence of special
curves, all being crucial, and more especially Proposition
\ref{prop:existence_des_curves_spéciales}.\bigskip

\begin{prop}[Stability of the special curves]
        \label{prop:stability_des_curves_spéciales}
\mbox{}
\begin{itemize}
\item[\;(i)] The set $\sigma_s(\G_0)$ is $\J$-stable.

\item[(ii)] For any $i\leqslant n-1$, the set
$\sigma_s(A_i)$ is stable by all the elements of
$\G_0\smallsetminus\{A_{i-1},\,A_{i+1}\}$.
\end{itemize}
\end{prop}
\medskip

\DEM
\smallskip

\noindent \emph{Let us show item (i).}
\smallskip

For all $a\in\sigma_s(\G_0)$, we have $I(a,\,\delta.a)\not=0$. From
the point of view of the curve $\delta.a$, we have
$I\big(\delta.a,\,\sigma(\G_0)\big)\not=0$. As $\G_0$ is stable by the
action of $\delta$ via $\rho$ by conjugation on $\PMod(\Sigma)$, the
curve $\delta.a$ belongs to $\sigma(\G_0)$. Finally, $\delta.a$
belongs to $\sigma_s(\G_0)$. Hence $\sigma_s(\G_0)$ is $\J$-stable.
\medskip

\noindent \emph{Let us show item (ii).}
\smallskip

To simplify the proof, we set $i=0$. Let $a$ be a special curve of
$\sigma_s(A_0)$. Let $j$ be an integer in $\{3,\dots,\,n-2\}$ so
that $A_j$ commutes with $A_0$ and $A_1$. Then
$A_j(a)\in\sigma(A_{0})$ and $A_j(\delta.a)\in\sigma(A_{1})$. By
hypothesis, $a$ is special, so $I(a,\,\delta.a)\not=0$. Then, when
we apply $A_j$ to the pair $(a,\,\delta.a)$, we get
$I(A_j(a),\,\sigma(A_1))\not=0$, so
$A_j(a)\in\sigma_s(A_0)$.\smallskip

Symmetrically, when $j=2$, if we  replace $\delta.a$ by
$\delta^{-1}.a$, we show that $I(A_2(a),\,\sigma(A_{n-1}))\not=0$,
so $A_2(a)\in\sigma_s(A_0)$.\smallskip

We have one more case to deal with, when $j=0$. Let us start again
from $a\in\sigma_s(A_0)$. Then $A_0(a)\in\sigma(A_0)$. We want to
show that $A_0(a)$ is a special curve. We are going to show that
$\big|\,\J.(A_0(a))\,\big|=n$, which is enough according to
Proposition \ref{prop:spéciales_normales}. Since $a$ is special, we
have $\sp(a)=\{A_0\}$. Now, for all $\ell\in\{2,\dots,\,n-2\}$, the
mapping class $A_\ell$ commutes with $A_0$ and
$a\not\in\sigma(A_\ell)$, so $A_0(a)\not\in\sigma(A_\ell)$. Hence we
have:
\smallskip

\centrer{$\{A_0\}\subset\sp(A_0(a))\subset\{A_{n-1},\,A_0,\,A_1\}$.}
\smallskip

\noindent So for all integers $k$,
$\delta^k.\sp(A_0(a))=\sp(A_0(a))$ if and only if $k$ is a multiple
of $n$. So, if $\delta^k.A_0(a)=A_0(a)$ then $k$ is a multiple of
$n$, so $|\J.(A_0(a))|\geqslant n$. Hence according to Proposition
\ref{prop:cardinaln}, we have $|\J.(A_0(a))|=n$. So $A_0(a)$ is a
special curve.\fin\bigskip

\begin{prop}[Stability of the normal curves]
        \label{prop:stability_des_curves_normales}
\mbox{}\\The set
$\sigma_n(\G_0)$ is $\B_n$-stable and the actions of $\B_n$ via $\rho$
on $\sigma_n(\G_0)$, on $\Sub_{\sigma_n(\G_0)}(\Sigma)$ and on
$\Bord(\Sigma_{\sigma_n(\G_0)})$ are cyclic.
\end{prop}
\medskip

\DEM Let us recall that according to Proposition
\ref{prop:spéciales_normales}, $\sigma_n(\G_0)$ is a simplex. Then if
we show that $\sigma_n(\G_0)$ is $\B_n$-stable, we can apply
Proposition \ref{prop:action_we_simplex_of_curves} (according to
which any action of $\B_n$ on a curve simplex $\B_n$-stable is
cyclic) and deduce from it Proposition
\ref{prop:stability_des_curves_normales}. Let $a$ be a normal curve.
We proceed as follows:
\begin{itemize}
\item[1.] We show that $A_1(a)$ belongs to $\sigma(\G_0)$.
\item[2.] We show that $A_1(a)$ is not special, hence is normal.
\item[3.] Therefore $\sigma(\G_0)$ is $\B_n$-stable.
\end{itemize}
\smallskip

1. Let $a$ be a normal curve. According to assertion (2g) of
Proposition \ref{prop:spéciales_normales}, $\sp(a)$ contains at
least two elements $A_i$ and $A_j$ with $0\leqslant i<j\leqslant
n-1$, and according to assertion (4g) of the same proposition, we
can assume that $j\not\in\{i+2,\,i+n-2\}$. In particular
$\{i,\,j\}\not=\{0,\,2\}$. Therefore $A_1$ commutes with at least
one of the two mapping classes $A_i$ and $A_j$, so
$A_1(a)\in\sigma(\G_0)$.
\smallskip

2. Let us assume that $A_1(a)$ is special. If $A_1(a)$ did not
belong to $\sigma_s(A_0)$ or to $\sigma_s(A_2)$, then according to
Proposition \ref{prop:stability_des_curves_spéciales}.(ii),
$A_1^{-1}\big(\,A_1(a)\,\big)$ would still be a special curve. But
$A_1^{-1}(A_1(a))=a$, which is a normal curve. Therefore $A_1(a)$
belongs to $\sigma_s(A_0)$ or $\sigma_s(A_2)$. The situation being
symmetric, we can assume tha:t
\smallskip

\centrer[1]{$A_1(a)\in\sigma_s(A_0)$.}
\smallskip

\noindent Having assumed that $A_1(a)$ was special, we have
$I(A_1(a),\,\delta.A_1(a))\not=0$, so
\smallskip

\centrer[2]{$I\big(\,a,\,A_1^{-1}(\delta.A_1(a))\,\big)\not=0$.}
\smallskip

\noindent Now, according to (1), $A_1(a)\in\sigma(A_0)$, so
$\delta.A_1(a)\in\sigma(A_1)$ and
$A_1^{-1}(\delta.A_1(a))\in\sigma(A_1)$. Then (2) implies:
\smallskip

\centrer{$I\big(\,a,\,\sigma(\G_0)\,\big)\not=0$,}
\smallskip

\noindent which is absurd since $a$ is a normal curve.
\medskip

3. Let us conclude: we have just shown that
$A_1(\sigma_n(\G_0))=\sigma_n(\G_0)$. But
$\sigma_n(\G_0)=\sigma(\G_0)\smallsetminus\sigma_s(\G_0)$ and the two sets
$\sigma(\G_0)$ and $\sigma_s(\G_0)$ are $\J$-stable (for $\sigma_s(\G_0)$,
this comes from Proposition
\ref{prop:stability_des_curves_spéciales}), so $\sigma_n(\G_0)$ is
$\J$-stable. Hence for all integers $i\in\{1,\dots,\,n-1\}$, we have
\linebreak $(\rho(\delta^i)A_1\rho(\delta^i)^{-1})(\sigma_n(\G_0))=\sigma_n(\G_0)$,
in other words, $A_i(\sigma_n(\G_0))=\sigma_n(\G_0)$. Since $\G_0$ spans
$\rho(\B_n)$, is follows that $\sigma_n(\G_0)$ is stable by
$\rho(\B_n)$.\fin\bigskip

\begin{prop}[Spectrum of the normal curves]
            \label{prop:spectrum_des_curves_normales}
\mbox{}\\The spectrum of a normal curve is always equal to $\G_0$.
\end{prop}
\medskip

\DEM Let $a$ be a normal curve. There exists an integer $k$ such
that the curve $a'=\rho(\delta^k)(a)$ belongs to $\sigma(A_0)$. Then
the curve $(A_0A_1A_0A_3^{-3})(a')$ belongs to $\sigma(A_1)$ for
\smallskip

\centrer{$(A_0A_1A_0A_3^{-3})A_0(A_0A_1A_0A_3^{-3})^{-1}=A_1$.}
\smallskip

\noindent But since the action of $\B_n$ is cyclic on the normal curves
according to Proposition \ref{prop:stability_des_curves_normales},
the action of $(A_0A_1A_0A_3^{-3})$ is trivial and
$(A_0A_1A_0A_3^{-3})(a')=a'$, so the curve $a'$ belongs to
$\sigma(A_1)$. For all $i\in\{0,\,1,\dots,\,n-1\}$, the same
argument can be repeated, so $a'$ belongs to $\sigma(A_i)$ for all
$i\in\{0,\,1,\dots,\,n-1\}$. By conjugating this by
$\rho(\delta^{-k})$, it follows that the curve
$a=\rho(\delta^{-k})(a')$ belongs to $\sigma(A_i)$ for all
$i\in\{0,\,1,\dots,\,n-1\}$. The proposition is proved. \fin
\bigskip

\begin{prop}[Existence of the special curves]
        \label{prop:existence_des_curves_spéciales}
The set $\sigma_s(\G_0)$ is not empty.
\end{prop}
\medskip

\DEM Let us recall that the homomorphism $\rho$ is assumed to be
noncyclic. We argue by contradiction: we assume that all the curves
of $\sigma(\G_0)$ are normal.

First, since $\rho$ is not cyclic and according to Theorem
\ref{thm:homomorphisms_irreducible}, $\sigma(\G_0)$ is not empty. Then the
set $\sigma(\G_0)$ of curves (being all normal by assumption) is a
simplex, according to Proposition \ref{prop:spéciales_normales}.
Moreover, this simplex is $\B_n$-stable according to Proposition
\ref{prop:stability_des_curves_normales}. The Proposition
\ref{prop:action_we_simplex_of_curves} can be applied to the simplex
$\sigma(\G_0)$, so the action of $\B_n$ on
$\Bord(\Sigma_{\sigma(\G_0)})$ is cyclic and the one of $\F_n$ on
$\Bord(\Sigma_{\sigma(\G_0)})$ is trivial. Hence the homomorphism $\rho$
induces a homomorphism $\bar\rho$ from $\F_n$ to $\PMod(\Sigma_{\sigma(\G_0)})$. Recall that we have the following
canonical isomorphism:
$$\PMod(\Sigma_{\sigma(\G_0)}) =
\prod_{S\in\Comp(\Sigma_{\sigma(\G_0)})}\PMod(S),$$ where
$\Comp(\Sigma_{\sigma(\G_0)})$ denotes the set of all connected
components of $\Sigma_{\sigma(\G_0)}$. Let $S$ be one of these
components and let $\bar\rho_S$ be the homomorphism induced by
$\rho$ on $\PMod(S)$. The homomorphism $\bar\rho_S$ is irreducible,
that is, for all $i\in\{3,\dots,\,n-1\}$, the element
$\tau_i\tau_1^{-1}$ is sent on a irreducible mapping class.
Indeed, according to Proposition
\ref{prop:properties_sigma}.(v),
$\sigma(\rho(\tau_i\tau_1^{-1}))\subset\sigma(\rho(\tau_i))\cup\sigma(\rho(\tau_1^{-1}))$,
but these two sets $\sigma(\rho(\tau_i))$ and
$\sigma(\rho(\tau_1^{-1}))$ do not contain any curve in $S$, so
$\sigma(\bar\rho_S(\tau_i\tau_1^{-1}))=\vide$. We will say that
$\bar\rho_S$ is \emph{periodic} or \emph{pseudo-Anosov} whether
$\bar\rho_S(\tau_3\tau_1^{-1})$ is periodic or pseudo-Anosov.
Let us denote by $A=\bar\rho_S(\tau_3\tau_1^{-1})$,
$B=\bar\rho_S(\tau_4\tau_1^{-1})$,
$C=\bar\rho_S(\tau_5\tau_1^{-1})$, and $Z=ABACBA$. \medskip

\textbf{Case where $\bar\rho_S$ is pseudo-Anosov\footnote{We
cannot apply the results of Section
\ref{sec:homomorphisms_irreducible} for $\F_n^*$ can be isomorphic
to $\B_4$ and this case is not treated in Section
\ref{sec:homomorphisms_irreducible}. On the other hand, the
techniques involved here would have appeared quite complicated
in Section \ref{sec:homomorphisms_irreducible} when the boundary of
$\Sigma$ is nonempty, and simply do not work when the boundary
is empty.}:} According to Proposition
\ref{prop:structure_centralizer}.(iv), the centralizer of a
pseudo-Anosov mapping class is virtually infinite cyclic. Then,
since $C$ is in the centralizer of $A$, there exist two nonzero
integers $p$ and $q$ such that $A^p=C^q$. By conjugating this
equality by $Z$, we get $C^p=A^q$. Hence
$A^{p^2}=C^{qp}=A^{q^2}$. We deduce that $p=q$ or $p=-q$. We
are going to show that both of these two equalities are absurd.
For this purpose, we set
$\gamma_{13}=\tau_1\tau_2\tau_1\tau_3\tau_2\tau_1$ and
$\gamma_{35}=\tau_3\tau_4\tau_3\tau_5\tau_4\tau_3$. Notice
that:
\smallskip

\centrer{$\begin{array}{c}
  \gamma_{13}\,\tau_1\,{\gamma_{13}}^{-1}=\tau_3 \mbox{ and } \gamma_{35}\,\tau_1\,{\gamma_{35}}^{-1}=\tau_1,\\
  \gamma_{13}\,\tau_3\,{\gamma_{13}}^{-1}=\tau_1 \mbox{ and } \gamma_{35}\,\tau_3\,{\gamma_{35}}^{-1}=\tau_5,\\
  \gamma_{13}\,\tau_5\,{\gamma_{13}}^{-1}=\tau_5 \mbox{ and } \gamma_{35}\,\tau_5\,{\gamma_{35}}^{-1}=\tau_3.\\
\end{array}$}
\smallskip

\noindent Hence, if we set $\upsilon=\gamma_{13}{\gamma_{35}}^{-1}$,
then the element $\upsilon$ belongs to $\F_n$ and satisfies:
\smallskip

\centrer{$\begin{array}{lll}
  \upsilon\,\tau_1\,\upsilon^{-1}=\tau_3, & \upsilon\,\tau_3\,\upsilon^{-1}=\tau_5, & \upsilon\,\tau_5\,\upsilon^{-1}=\tau_1,\\
  \upsilon\,(\tau_3\tau_1^{-1})\,\upsilon^{-1}=\tau_5\tau_3^{-1}, &
  \upsilon\,(\tau_5\tau_3^{-1})\,\upsilon^{-1}=\tau_1\tau_5^{-1}, &
  \upsilon\,(\tau_1\tau_5^{-1})\,\upsilon^{-1}=\tau_3\tau_1^{-1}.
\end{array}$}
\smallskip

\noindent Then if we set $U=\rho_S(\upsilon)$, then:
\smallskip

\centrer{$\begin{array}{llll}
  \mbox{we have: } & UAU^{-1}=CA^{-1}, & U(CA^{-1})U^{-1}=C^{-1}, & UC^{-1}U^{-1}=A,\\
  \mbox{so if $C^p=A^p$: } & UA^pU^{-1}=\Id, & U(\Id)U^{-1}=A^{-p}, & UA^{-p}U^{-1}=A^p,\\
  \mbox{and if $C^p=A^{-p}$: } & UA^pU^{-1}=A^{-2p}, & U(A^{-2p})U^{-1}=A^{p}, & UA^{p}U^{-1}=A^p.
\end{array}$}
\smallskip

\noindent So if $C^p=A^p$, then $A^p$ is equal to $\Id$, which is
absurd for $A$ is pseudo-Anosov; whereas if $C^p=A^{-p}$, then we
have $UA^pU^{-1}=A^{p}$ and $UA^pU^{-1}=A^{-2p}$ whence
$A^p=A^{-2p}$ and so $A^{3p}=\Id$, which is also absurd. Hence
$\rho_S$ is not pseudo-Anosov.
\medskip

\textbf{Case where $\bar\rho_S$ is periodic.} Notice that the
mapping classes $A$, $B$, $C$ are conjugate. There are then periodic
of same order. Let us call $m$ this order. We restrict the domain of
the homomorphism $\bar\rho_S$ to $\F_n^*$, which is isomorphic to
$\B_{n-2}$. According to Proposition
\ref{prop:passage_PMod_à_ModBord}, we can lift this new homomorphism
from $\F_n^*$ to $\PMod(S)$ into a homomorphism $\tilde \rho_S$ from
$\F^*_n$ to $\Mod(S,\,\bord S)$. Let us denote by $\tildeA$,
$\tildeB$, $\tildeC$, the images by $\tilde \rho_S$ of
$\tau_3\tau_1^{-1}$, $\tau_4\tau_1^{-1}$, $\tau_5\tau_1^{-1}$. Then
$\tildeA^m$ and $\tildeC^m$ are multitwists along some curves of
$\Bord(S)$. Since these multitwists are in the center of
$\Mod(S,\,\bord S)$ and since
$\tildeZ\tildeA^m\tildeZ^{-1}=\tildeC^m$ (where
$\tildeZ=(\tildeA\tildeB\tildeC)^2$), we have $\tildeA^m=\tildeC^m$.
Therefore, $\tildeA\tildeC^{-1}$ satisfies
$(\tildeA\tildeC^{-1})^m=1$ in $\Mod( S,\,\bord S)$. But
$\Mod(S,\,\bord S)$ is torsion-free, so $\tildeA\tildeC^{-1}$ is
trivial and $\tildeA=\tildeC$. This implies that in $\PMod(S)$, we
have $\bar\rho_S(\tau_3\tau_1^{-1})=\bar\rho_S(\tau_5\tau_1^{-1})$.
\smallskip

This last equality holds for any connected components $S$ of
$\Sigma_{\sigma(\G_0)}$. Hence, by considering the homomorphism
$\bar\rho\;:\; \F_n\to\PMod(\Sigma_{\sigma(\G_0)})$, we have shown
that $\bar\rho(\tau_3\tau_1^{-1})=\rho(\tau_5\tau_1^{-1})$. Then
$\bar\rho(\tau_5\tau_3^{-1})$ coincides with the identity of
$\PMod(\Sigma_{\sigma(\G_0)})$. Let us recall that $\B_n$ stabilizes
$\Sigma_{\sigma(\G_0)}$. Then by conjugation in $\B_n$, we deduce that
the homomorphism $\bar\rho\;:\;\F_n\to\PMod(\Sigma_{\sigma(\G_0)})$ is
trivial. Hence the image of the restriction of the homomorphism $\rho$
to $\F_n^*$ to $\PMod(\Sigma)$ is included in the abelian group
spanned by the Dehn twists along the curves of $\sigma(\G_0)$. Hence,
according to Lemma \ref{lem:B_n_in_group_abelian}, the restriction
of $\rho$ to $\F_n^*$ is a cyclic homomorphism. Therefore
$\rho(\tau_5\tau_3^{-1})$ is the identity, so
$\rho(\tau_3)=\rho(\tau_5)$. Hence, according to Lemma
\ref{lem:homomorphism_cyclic}, $\rho$ is a cyclic homomorphism . This is
contradicts our hypotheses, so the proposition is proved.\fin
\bigskip

\begin{prop}[Cardinality of $\sigma_s(\G_0)$]
            \label{prop:cardinal_de_sigmas}
\mbox{}\\The set $\sigma_s(\G_0)$ contains $n$ or $2n$ curves, depending
on whether $|\sigma_s(A_1)|=1$ or $|\sigma_s(A_1)|=2$.
\end{prop}
\medskip

\DEM First, according to Proposition
\ref{prop:stability_des_curves_spéciales}, we have
$\delta.\sigma_s(A_i)=\sigma_s(A_{i+1})$, so the cardinality of
$\sigma_s(\G_0)$ is equal to $n$ times the one of $\sigma_s(A_1)$. But
$\sigma_s(A_1)$ contains one or two curves for on one hand, it
cannot be empty since $\sigma_s(\G_0)$ is not empty, and on the other
hand, it cannot contain three curves or more as we are going to show
it. This prove the proposition. \smallskip

Let us then show that $|\sigma_s(A_1)|<3$. Notice that the elements
$A_1,\,A_3,\,A_5,\,\dots,\,A_{n-1}$ pairwise commute, so the set of
curves
$\sigma_s(A_1)\cup\sigma_s(A_3)\cup\dots\cup\sigma_s(A_{n-1})$ is a
simplex. If $\sigma_s(A_1)$ contained at least three curves, the set
of curves
$\sigma_s(A_1)\cup\sigma_s(A_3)\cup\dots\cup\sigma_s(A_{n-1})$ would
be a simplex $\A$ of at least $\frac{3n}{2}$ curves that would be
stable by the action of $\langle\delta^2\rangle$ on
$\Courb(\Sigma)$. However, the orbits included in $\A$ under the
action of $\langle\delta^2\rangle$ contain at least $\frac{n}{2}$
curves, which is greater than or equal to 3, so we can apply
Proposition \ref{prop:as_if_surface_without_boundary}: after having
squeezed the boundary components of $\Sigma$, the simplex $\A$ still
contains at least $\frac{3n}{2}$ distinct curves. But
$\frac{3n}{2}\geqslant3g$ whereas the cardinality of all simplex in
a surface without boundary of genus $g$ is bounded by $3g-3$. This
is the expected contradiction. \fin
\bigskip

\section{The special curves are not separating}
    \label{sec:the_curves_spéciales_are_non-separating}
\bigskip

\TITRE{Hypotheses.}\\
Let $n\geqslant 6$ an even number, let $\Sigma=\Sigma_{g,\,b}$ with $g\leqslant\frac{n}{2}$, and let
$\rho\ :\ \B_n\to\PMod(\Sigma)$ such that:

\begin{tabular}{lll}
\point & $\rho$ is non-cyclic      & by assumption,\\
\point & $\sigma_p(\G_0)=\vide$      & by assumption, inspired by Proposition \ref{prop:pas_of_curve_peripheral},\\
\point & $\sigma_s(\G_0)\not=\vide$    & as a consequence of the non-cyclicity of $\rho$, after Proposition
\ref{prop:existence_des_curves_spéciales}.
\end{tabular}
\bigskip

\noindent In this section, we are going to show the following proposition:
\smallskip

\begin{prop}
                \label{prop:curves_spéciales_separating}
The curves of $\sigma_s(\G_0)$ are not separating. \end{prop}
\medskip

\noindent \point\mbox{} In Subsection \ref{par:H_and_X}, we present
a subset $\X$ of $\G_0$, stable under the action of $\H$ via $\rho$,
where $\H$ is a subgroup of $\B_n$. The set $\X$ is smaller than
$\G_0$, but the action of $\H$ on $\X$ is $r$ times transitive, where
$r$ is the cardinality of $\X$. Moreover, $\X$ will consists in
elements of $\G_0$ which pairwise commute, so the union of their
canonical reduction systems will be a simplex. These aspects will be
very useful.
\medskip

\noindent \point\mbox{} Subsection
\ref{par:proof_prop_curves_spéciales_separating} is devoted to the
proof of Proposition \ref{prop:curves_spéciales_separating}. The
proof will be topological and the bound of $\frac{n}{2}$ on the
genus of $\Sigma$ is essential here. If we wanted to replace the
bound $\frac{n}{2}$ by $\frac{n}{2}+1$, our method would fail.
However, we conjecture that the bound $\frac{n}{2}$ is not the best one.
\bigskip


\subsection{The subgroup $\H$ of $\B_n$ and its action on the subset $\X$ of $\G_0$}
      \label{par:H_and_X}
\medskip

\TITRE{The subset $\Imp(n)$ of $\{0,\,1,\dots,\,n-1\}$ and the subset $\X$ of $\G_0$}.\\
For all positive integers $m$, let $\Imp(m)$ be the set of
the first odd integers smaller than or equal to $m$. Let $\X=\{A_i,\
i\in\Imp(n)\}$ be the subset of $\G_0$. We set
$\sigma(\X)=\cup_{i\in\Imp(n)}\sigma(A_i)$.
The elements of $\X$ commute pairwise, so the curves in
$\sigma(\X)$ cannot intersect each other. Thus $\sigma(\X)$ is a
curve simplex.
\medskip

\TITRE{The subgroup $\H$ of $\B_n$}
\mbox{}\\For all integers $i$ belonging to $\Imp(n)$, we set
\smallskip

\centrer{$\gamma_i=\tau_i\tau_{i+1}\tau_i\tau_{i+2}\tau_{i+1}\tau_i$,}
\smallskip

\noindent where for all integers $k$, we denote by $\tau_k$ the
standard generator $\tau_{\ell}$ where $\ell$ is the remainder of
the euclidian division of $k$ by $n$. The group $\H$ is the subgroup
of $\B_n$ defined by
\smallskip

\centrer{$\H:=\langle\,\gamma_i\,,\ i\in\Imp(n)\,\rangle$.}
\medskip

\begin{prop}[Properties of the group $\H$]\mbox{}
           \label{prop:properties_of_H}
\begin{itemize}
\item[\;\;(i)] The action of $\H$ by conjugation via $\rho$ on $\PMod(\Sigma)$ preserves $\X$.
\item[\;(ii)] The homomorphism $\H\to\Ss(\X)$ of the action of $\H$ on $\X$, where
$\Ss(\X)$ is the symmetric group on the elements of $\X$, is
surjective. Consequently, this action is $\frac{n}{2}$ times
transitive.
\item[(iii)] The action of $\H$ on $\Courb(\Sigma)$ preserves $\sigma_s(\X)$.
\end{itemize} \end{prop}
\smallskip

\DEM Let us first remark that for all $i\in\Imp(n)$, the element
$\gamma_i$ acts by conjugation on the subset
$\{\tau_j,\,j\in\Imp(n)\}$ of $\B_n$ in the following way:
\smallskip

\centrer{$\gamma_i\tau_j\gamma_i^{-1}= \left\{\begin{array}{lll}
  \tau_{j} &\mbox{ if }& j\not\in\{i,\,i+2\}\;\\
  \tau_{j-2} &\mbox{ if }& j=i+2\;\\
  \tau_{j+2} &\mbox{ if }& j=i\ .
\end{array}\right.$}
\smallskip

\noindent We deduce the following:

\begin{itemize}
\item[\;\;(i)]
The group $\H$ acts by conjugation on $\B_n$ and preserves the set
$\{\tau_j,\,j\in\Imp(n)\}$, so the group $\H$ acts via $\rho$ on
$\PMod(\Sigma)$ and preserves $\X$.
\item[\;\;(ii)]
For all $i\in\Imp(n)$, the homomorphism $\phi\;:\; \H\to\Ss(\X)$, where
$\Ss(\X)$ is the symmetric group on the elements of $\X$, sends
$\gamma_i$ on the transposition that swaps $A_i$ and $A_{i+2}$. Then
$\phi(\H)$ is a subgroup of $\Ss(\X)$ containing $\frac{n}{2}-1$
transpositions with disjoint supports, so $\phi$ is surjective.
\item[(iii)]
The action of $\H$ on $\Courb(\Sigma)$ preserves $\sigma(\X)$, since
the action of $\H$ on $\PMod(\Sigma)$ preserves $\X$. Now, according
to Proposition \ref{prop:stability_des_curves_normales}, the action
of $\B_n$ on $\Courb(\Sigma)$ preserves $\sigma_n(\G_0)$, so the
action of $\H$ on $\Courb(\Sigma)$ preserves $\sigma_n(\G_0)$. So the
action of $\H$ on $\Courb(\Sigma)$ also preserves the complement of
$\sigma_n(\G_0)$ in $\sigma(\X)$, which is $\sigma_s(\X)$. \fin
\end{itemize}
\bigskip


\subsection{Proof of Proposition \ref{prop:curves_spéciales_separating}}
      \label{par:proof_prop_curves_spéciales_separating}

\medskip

We will need the definition of \emph{natural boundary} and
\emph{special boundary} of a subsurface.
\smallskip

\TITRE{Natural boundary and special boundary.}
\mbox{}\\Let $n$ be an even integer greater than or equal to 6, let
$\Sigma$ be a surface and $\rho$ a homomorphism from $\B_n$ to $\PMod(\Sigma)$. Let $\A$ be a curve simplex included in
$\Courb(\Sigma)$. For any subsurface $S$ of $\Sub_\A(\Sigma)$, a
boundary component $d$ of $S$ will be said to be \emph{natural} if
it belongs to $\Bord(\Sigma)$, and will be said to be \emph{special}
if it belongs to $\sigma_s(\G_0)$. The union of the natural boundary
components will be called the natural boundary, and the union of the
special boundary components will be called the special boundary.
\medskip

To prove Proposition \ref{prop:curves_spéciales_separating}, we proceed in five steps.\bigskip

\begin{figure}[!h]
 \Includegraphics{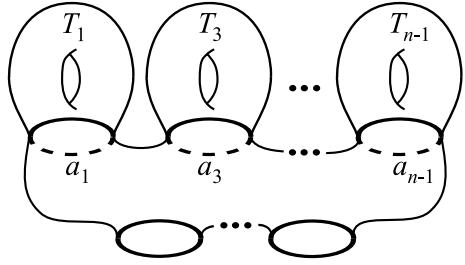}
 \caption{The surface $\Sigma$ and the separating curves $a_i$, $i\in\Imp(n)$ of $\sigma(\X)$.}
 \label{fig:gTores}
\end{figure}

\begin{etape}
If there exists in $\sigma_s(\G_0)$ a separating curve, then for all
$i\in\{0,\,1,\dots,\,n-1\}$, the set of curves $\sigma_s(A_i)$
contains exactly a separating curve that bounds a torus with one
hole (cf. Figure \ref{fig:gTores}). Thus, $\Sigma$ is a surface of
genus $g=\frac{n}{2}$. \end{etape}
\medskip

\DEM[of step 1] If there exists a separating curve in
$\sigma_s(\G_0)$, then there exists at least one separating curve in
$\sigma_s(A_1)$. Let us call it  $a_1$. Let $\A$ be the set of
curves $\H.a_1$. Since we have $a_1\in\sigma(\X)$ and since $\X$ is
$\H$-stable, the set $\A$ is included in $\sigma(\X)$, hence is a
simplex. Let us consider the graph $\Gamma(\Sigma,\,\A)$. Since the
curves of $\A$ are separating, if we remove from
$\Gamma(\Sigma,\,\A)$ one of its edges, we get a disconnected graph.
Hence the graph $\Gamma(\Sigma,\,\A)$ contains no cycle: this graph
is a finite tree. So it contains leaves (vertices of degree 1). Let
$T$ be a subsurface of $\Sub_\A(\Sigma)$ corresponding to a leaf in
$\Gamma(\Sigma,\,\A)$. Since exactly one curve of $\A$ bounds $T$,
$T$ is not stable by the action of $\H$, so $T$ contains no natural
boundary component. Hence $T$ has only one boundary component, so
$T$ is of nonzero genus. Each subsurface in the orbit $\H.T$ can be
identified by the curve $a\in\A$ that bounds it. Since the number of
curves in $\A$ is at least equal to the cardinality of $\X$, there
exists at least $\frac{n}{2}$ disjoint subsurfaces homeomorphic to
$T$ in $\Sub_\A(\Sigma)$. But $\Sigma$ is of genus $g\leqslant
\frac{n}{2}$, so:
\begin{itemize}
\item  there exist exactly $\frac{n}{2}$ such subsurfaces,
\item  these subsurfaces are tori with one hole,
\item  and $\Sigma$ is of genus $g=\frac{n}{2}$.
\end{itemize}
So, there exist exactly $\frac{n}{2}$ curves in $\A$, hence one
separating curve in each set $\sigma_s(A_i)$, $i\in\Imp(n)$.
Moreover, the complement of these $\frac{n}{2}$ tori is a genus-0
surface having $\frac{n}{2}$ special boundary components and $b$
natural boundary components (cf. Figure \ref{fig:gTores}).\fin
\bigskip

From now on, we are in the situation described by step 1 and we
adopt the following notation.
\medskip

\TITRE{Situation described by step 1}\mbox{}\\
For all $i\in\{0,\,1,\dots,\,n-1\}$, let us denote by $a_i$ the
unique separating curve of $\sigma_s(A_i)$ and by $T_i$ the torus
with one hole, included in $\Sigma$ and bounded by $a_i$. Let us
denote by $S$ the genus-0 surface obtained from $\Sigma$ minus the
tori $T_i$ where $i$ ranges over $\Imp(n)$ (cf. Figure
\ref{fig:gTores}).
Let us choose some representatives $\bar a_i$,
$i\in\{0,\,1,\dots,\,n-1\}$ in tight position, of the curves $a_i$,
$i\in\{0,\,1,\dots,\,n-1\}$ (such a system of representatives is
unique up to isotopy, according to Proposition
\ref{prop:existence_uniqueness_representatives_curves}). From these
representatives of curves, we deduce the representatives $\bar S$
and $\bar T_i$ for all $i\in\{0,\,1,\dots,\,n-1\}$, of the
subsurfaces $S$ and $T_i$ for all $i\in\{0,\,1,\dots,\,n-1\}$.
\medskip


\begin{etape} \label{symmétrie}
There exists an orientation preserving diffeomorphism $\bar F$ of
$\Sigma$ that preserves the boundary components of $\Sigma$ such
that
\begin{itemize}
  \item
 for all $i\in\{1,\,2,\,3\}$, we have $\bar F(\bar a_i)=\bar a_{4-i}$,
  \item
 and for all $i\in\{5,\,6,\dots,\,n-1\}$, we have $\bar F(\bar a_i)=\bar a_{i}$.
\end{itemize} \end{etape}
\medskip

\DEM[of step 2] Let $\bar G$ be a diffeomorphism representing the
mapping class $\rho(\gamma_1)$ where $\gamma_1$ is defined by:
\smallskip

\centrer{$\gamma_1:=\tau_1\tau_2\tau_1\tau_3\tau_2\tau_1$.}
\smallskip

\noindent Let us denote by $\bar\A$ the set of representatives of
curves $\{\bar a_2\}\cup\{\bar a_i,\,i\in\Imp(n)\}$ and $\bar\A'$
the set $\{\bar G(\bar a_2)\}\cup\{\bar G(\bar
a_i),\,i\in\Imp(n)\}$. According to Proposition
\ref{prop:properties_of_H},
\begin{itemize}
  \item
 for any $i\in\{1,\,2,\,3\}$, we have
 $\rho(\gamma_1)A_i\rho(\gamma_1)^{-1}=A_{4-i}$,
  \item
 and for any $i\in\{5,\,6,\,\dots,\,n-1\}$, we have $\rho(\gamma_1)
 A_i\rho(\gamma_1)^{-1}=A_{i}$.
\end{itemize}
\smallskip

\noindent Besides, for any $i\in\{0,\,1,\dots,\,n-1\}$, the curve
$a_i$ is the unique separating curve belonging to $\sigma_s(A_i)$,
so:
\smallskip

\begin{itemize}
  \item
 for any $i\in\{1,\,2,\,3\}$, the curve $\bar G(\bar a_i)$ is isotopic to
 $\bar a_{4-i}$,
  \item
 and for any $i\in\{5,\,6,\,\dots,\,n-1\}$, the curve $\bar G(\bar a_i)$ is isotopic to $\bar
 a_{i}$.
\end{itemize}
\smallskip

\noindent Then, the sets of curves $\bar\A$ and $\bar\A'$ are both
weakly isotopic. Let us recall that $\bar\A$ is a set of
representatives of curves in tight position. Hence the
representatives of curves of $\bar\A$ do not bound any bigon. But
$\bar\A'$ is the image of $\bar\A$ by $\bar G$. So the
representatives of curves of $\bar\A'$ do not bound any bigon
either. Hence $\bar\A'$ is a set of representatives of curves in
tight position. Moreover $\bar\A$ is without triple intersection,
hence so is $\bar\A'$. Then, according to Proposition
\ref{prop:uniqueness_PMI}, $\bar\A$ and $\bar\A'$ are in the same
isotopy class. In other words, there exists a diffeomorphism
isotopic to the identity $\bar H$ such that $\bar
H(\bar\A)=\bar\A'$. Then the diffeomorphism $\bar F$ defined by
$\bar F:=\bar H^{-1}\bar G$ satisfies the assertions of the
statement.\fin
\bigskip

\TITRE{Arcs.}\mbox{}\\
For all $i,\,j,\,k$ in $\{1,\,2,\,3\}$, let us denote by
$\Arc_j^k(i)$ the set of closures of the connected components of
$\bar a_i\smallsetminus\big(\bar a_{j}\cup \bar a_{k}\big)$. We will
say that an element of $\Arc_j^k(i)$ is an \emph{arc included in
$\bar a_i$ with extremities in $\bar a_j$ and $\bar a_k$}.
\bigskip


\begin{etape}
                                \label{pas_of_lacet}
The arcs of $\Arc_1^1(2)$ and $\Arc_3^3(2)$ are included
respectively in $\bar T_1$ and $\bar T_3$. In other words, the only
arcs included in $\bar a_2\cap \bar S$ belong to $\Arc_1^3(2)$.
\end{etape}
\smallskip

\begin{figure}[!h]
 \Includegraphics{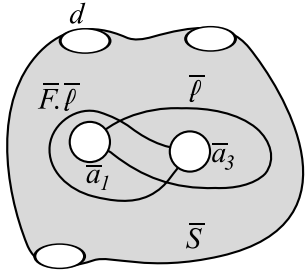}
 \caption{The existence of an arc $\bar\ell$ in $\Arc_1^1(2)$
 induces
 the existence of an arc $\bar F(\bar\ell)$ in $\Arc_3^3(2)$ intersecting $\bar\ell$.}
 \label{fig:intersectionArcs}
\end{figure}

\DEM[of step 3] Let us argue by contradiction. Let us consider an
arc $\bar\ell$ belonging to $\Arc_1^1(2)$ and included in $\bar S$.
Since $\bar S$ is of genus zero, $\bar\ell$ separates $\bar S$ in
two connected components. One of them contains the boundary $\bar
a_3$. Since the curves $\bar a_1$ and $\bar a_2$ do not cobound any
bigon, the other component is not a disk, so it contains a special
boundary component or a natural boundary component. In both cases,
let us call $d$ this boundary component. Finally, we have a path
$\bar\ell$ which separates $\bar S$ in two components, one
containing $\bar a_3$, the other containing $d$. Similarly $\bar
F(\bar\ell)$ belongs to $\Arc_3^3(2)$ (for $\bar\ell$ belongs to
$\Arc_1^1(2)$), is included in $\bar S$ and separates $\bar S$ in
two components, one containing $\bar F(\bar a_3)$ that is equal to
$\bar a_1$, the other containing $\bar F(d)$ that is equal to $d$.
We deduce from it Figure \ref{fig:intersectionArcs} where it is
clear that $\bar\ell$ and $\bar F(\bar\ell)$ intersect, which is
absurd: $\bar F(\bar\ell)$ and $\bar\ell$ cannot intersect, for they
are both included in the same curve.\fin
\bigskip

\TITRE{Arc paths, rectangles, hexagons, octogons.}\mbox{}\\
We call \emph{arc path} a union of arcs such that this union is
homeomorphic to a segment or to a circle. In the first case, the arc
path will be said to be \emph{open}, in the second it will be said
to be \emph{closed}.

If a connected component $D$ of $\bar T_2\cap \bar T_1$, of $\bar
T_2\cap \bar S$, or of $\bar T_2\cap \bar T_3$ is homeomorphic to a
disk whose boundary is a closed arc path, each of the arcs of this
arc path will be called \emph{edge of $D$}. Such a connected
component with four edges will be called a \emph{rectangle}, with
six edges a \emph{hexagon}, and with eight edges an \emph{octogon}.
\bigskip


\begin{etape}
                                                    \label{découpage}
The connected components of $\bar T_2\cap \bar S$ are rectangles.
The connected components of $\bar T_2\cap \bar T_1$ (respectively
$\bar T_2\cap \bar T_3$) consist in exactly one hexagon and some
rectangles. \end{etape}
\medskip

\DEM[of step 4]
\smallskip

1. \emph{connected components of $\bar T_2\cap \bar S$.}
\smallskip

\begin{figure}[!h]
 \Includegraphics{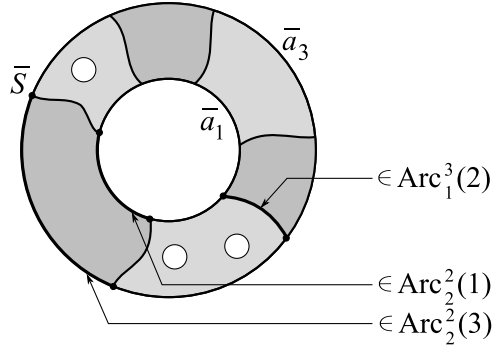}
 \caption{The surface $\bar S$ is seen as an annulus,
 the dark grey parties are the connected components of $\bar T_2\cap \bar S$.}
 \label{fig:rectangles}
\end{figure}

We can see $\bar S$ as an annulus whose boundary components are
$\bar a_1$ and $\bar a_3$, minus $\frac{n}{2}-2+b$ open disks,
corresponding to the curves $a_i$,
$i\in\Imp(n)\smallsetminus\{1,\,3\}$ and to the boundary components
of $\Sigma$ (cf. Figure \ref{fig:rectangles}). The torus $\bar T_2$
contains none of these curves and none of these boundary components
so the connected components of $\bar T_2\cap\bar S$ are simply
connected and of genus zero, hence are homeomorphic to disks. All
the boundary components of a component $\bar C$ of $\bar T_2\cap
\bar S$ are some arc paths leaning on the curves $\bar a_1$, $\bar
a_2$ and $\bar a_3$. But we have seen that the arcs included in
$\bar a_2\cap \bar S$ belong to $\Arc_1^3(2)$. It is easy to see
that in such an annulus, the only injective arc paths that contain
some arcs of $\Arc_1^3(2)$ and that bound disks are rectangles: two
edges belong to $\Arc_1^3(2)$, one edge to $\Arc_2^2(1)$ and one
edge to $\Arc_2^2(3)$.
\medskip

2. \emph{Connected components of $\bar T_2\cap (\bar T_1\cup \bar
T_3)$.}
\smallskip

\begin{figure}[!h]
 \Includegraphics{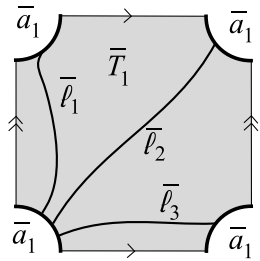}
 \caption{An example of three disjoint, non-isotopic arcs belonging to $\Arc_1^1(2)$, in $\bar T_1$.}
 \label{fig:troisChemins}
\end{figure}

We can see $\bar T_2$ as the gluing along arcs included in $\bar
a_1$ and $\bar a_3$ of the connected components of $\bar T_2\cap
\bar T_1$, $\bar T_2\cap \bar S$ and $\bar T_2\cap \bar T_3$. In
this proof, we call \emph{domains} these connected components. A
domain of $\bar T_2\cap \bar T_1$ is bounded by some arcs belonging
to $\Arc_2^2(1)$ and to $\Arc_1^1(2)$. Notice that in $\bar T_1$ (as
in any such torus with one hole), there exist at most three pairwise
disjoint, non-isotopic arcs, whose extremities belong to $\bord\bar
T_1$. So the set of arcs in $\Arc_1^1(2)$ contains at most three
isotopy classes (cf. Figure \ref{fig:troisChemins}).

\begin{figure}[!h]
 \Includegraphics{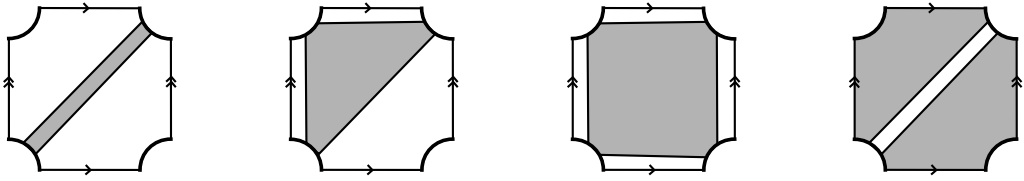}
 \caption{The four types of domains: from left to right:
 the rectangle, the hexagon, the octogon, the cylinder with bigonal boundary components.}
 \label{fig:quatreDomaines}
\end{figure}

But the arcs of $\Arc_1^1(2)$ and of $\Arc_2^2(1)$ constitute the
boundary components of the domains of $\bar T_1\cap\bar T_2$. We
deduce that there exist only four possible types of domains in $\bar
T_1$: rectangles, hexagons, octogons and \emph{cylinders with
bigonal boundary components}, that is to say spheres with two
boundary components, such that each is a path of two arcs (cf.
Figure \ref{fig:quatreDomaines}).

The connected components of $\bar T_3\cap\bar T_2$ satisfy the same
properties, since they are the images by the diffeomorphism $\bar F$
of step 2 of the connected components of $\bar T_1\cap\bar T_2$.
\smallskip

3. \emph{Euler characteristic computation.}
\smallskip

\begin{figure}[!h]
 \Includegraphics{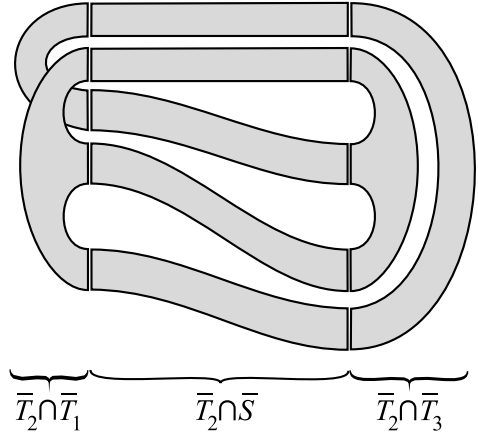}
 \caption{The torus $\bar T_2$, built up from two hexagons and some rectangles.}
 \label{fig:decoupage}
\end{figure}

We are going to determine the contribution of each domain to the
Euler characteristic of $\bar T_2$ (equal to $-1$ since $\bar T_2$
is a torus with one hole). Let us recall that $\bar T_2$ is the
gluing of the domains along the arcs of $\Arc_2^2(1)\cup\Arc_2^2(3)$
included in $\bar T_2$. Hence if a domain $D$ has exactly $k$ edges
belonging to $\Arc_2^2(1)\cup\Arc_2^2(3)$, then its contribution to
$\chi(\bar T_2)$ amounts to $\chi(D)-\frac{k}{2}$. Indeed, when we
add up the Euler characteristics of all the domains, each of the
gluing arcs has been counted twice. To compute the contribution of a
domain to the Euler characteristic of $\bar T_2$, we hence need to
add to its own Euler characteristic $-\frac{1}{2}$ as a corrective
term for each gluing arc included in the boundary of $D$. Thus,
\begin{itemize}
\item the rectangles have a 0-contribution;
\item the hexagons have a $(-\frac{1}{2})$-contribution;
\item the cylinders with bigonal boundary components have a $(-1)$-contribution;
\item the octogons have a $(-1)$-contribution.
\end{itemize} But according to step 2, the domains of $\bar
T_2\cap \bar T_3$ are diffeomorphic to the ones of $\bar T_2\cap
\bar T_1$. Besides, the contribution of the domains of $\bar T_2\cap
\bar S$ is zero, since all the domains of $\bar T_2\cap \bar S$ are
rectangles (as we have seen it in step 1.). Hence the global
contribution of the domains of $\bar T_2\cap \bar T_1$ and the
global contribution of the domains of $\bar T_2\cap \bar T_3$ must
both equal $-\frac{1}{2}$. Therefore $\bar T_2\cap \bar T_1$
contains exactly one hexagon and some rectangles. Same thing for
$\bar T_2\cap \bar T_3$. (An example of the torus $\bar T_2$ built up from two hexagons and
some rectangles, according to the conclusion of step 4, is given
Figure \ref{fig:decoupage}.)\fin
\bigskip

\begin{etape}
We end in a contradiction. \end{etape}
\smallskip

\begin{figure}[!h]
 \Includegraphics{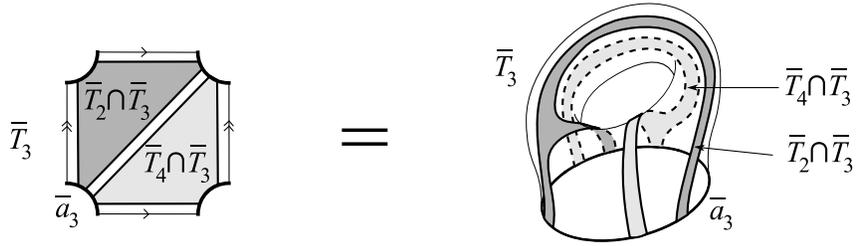}
 \caption{Two disjoint hexagons, one included in $\bar T_2\cap\bar
 T_3$, the other included in $\bar T_4\cap\bar T_3$.}
 \label{fig:deuxHexagones}
\end{figure}

\DEM[of step 5 and end of the proof of Proposition
\ref{prop:curves_spéciales_separating}]\\ The torus $\bar T_3$
contains the connected components of $\bar T_2\cap \bar T_3$ and of
$\bar T_4\cap \bar T_3$, which are pairwise disjoint since $\bar
T_2\cap \bar T_4=\vide$. There are two hexagons among them, one
included in $\bar T_2\cap \bar T_3$, the other included in $\bar
T_4\cap \bar T_3$. Each of them contains three edges included in
$\bar a_3$ (cf. Figure \ref{fig:deuxHexagones}).

\smallskip

Let us recall that it is possible to include in a torus with one
hole only three pairwise disjoint and non-isotopic arcs with
extremities in the boundary. Therefore our two hexagons are arranged
as in Figure \ref{fig:deuxHexagones}. In particular, let us remember
that:
\smallskip

\centrer[1]{\begin{tabular}{l}the six edges included in the boundary
of the torus $\bar T_3$ belong\\ alternatively to one and to the
other of the two tori $\bar T_2$ and $\bar T_4$.\end{tabular}}
\medskip

\begin{figure}[!h]
 \Includegraphics{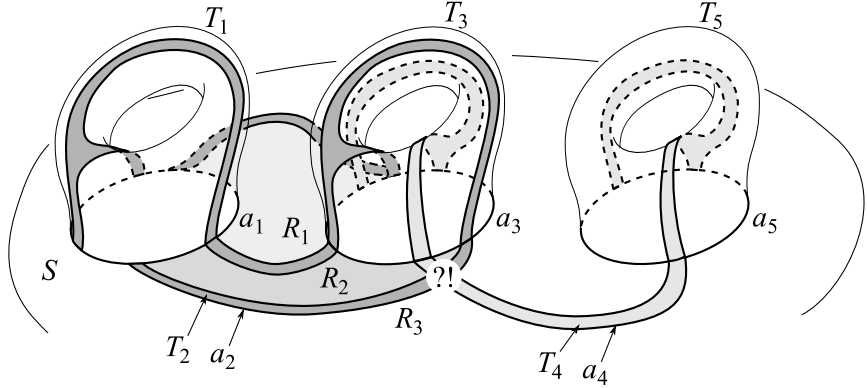}
 \caption{In this configuration where the hexagons $\bar T_2\cap\bar  T_3$ and $\bar T_4\cap\bar  T_3$ are ``nested'' in $\bar T_3$,
    the tori $\bar T_2$ and $\bar T_4$ intersect, which should not happen yet.}
 \label{fig:arcsProlonges}
\end{figure}

Let us describe how the tori $\bar T_2$ and $\bar T_4$ are embedded
in $\Sigma$. The following description is depicted in Figure
\ref{fig:arcsProlonges}.

1. Since $\bar T_2\cap \bar T_3$ contains a hexagon, $\bar T_2\cap
\bar a_3$ contains at least three connected components.

2. These at least three connected components extend in $\bar
T_2\cap\bar S$ in at least three rectangles. Indeed, let us recall
that the rectangles of $\bar T_2\cap\bar S$ have only one edge in
$\bar a_3$, so two distinct connected components of $\bar T_2\cap
\bar a_3$ are the edges of two distinct rectangles of $\bar
T_2\cap\bar S$. Now, there are at least three such rectangles in
$\bar T_2\cap \bar S$. Each of them has an edge in $\bar a_1$ and an
edge in $\bar a_3$. Since $\bar S$ is of genus zero, we deduce that
$\bar S\smallsetminus \bar T_2$ contains at least three connected
components. We will name by \emph{region} each of these connected
components.

3. For instance, consider the regions $R_1$, $R_2$, $R_3$, in Figure
\ref{fig:arcsProlonges}. Since the two tori $\bar T_2$ and $\bar
T_4$ are disjoint, the rectangles of $\bar T_4\cap\bar S$ are inside
some of these regions. However, all the rectangles of $\bar
T_4\cap\bar S$ have an edge included in $\bar a_5$, so they should
all be located in the region containing the curve $\bar a_5$.

4. But this is impossible, for according to statement (1), there
exist rectangles of $\bar T_4\cap\bar S$ in at least three distinct
regions (cf. Figure \ref{fig:arcsProlonges}).\fin

\bigskip

\section{Description of $\sigma(\X)$ in $\Sigma$}
    \label{sec:sigmaX}
\bigskip

\TITRE{Hypotheses.}\\
Let $n\geqslant 6$ an even number, let $\Sigma=\Sigma_{g,\,b}$ with $g\leqslant\frac{n}{2}$, and let
$\rho\ :\ \B_n\to\PMod(\Sigma)$ such that:

\begin{tabular}{lll}
\point & $\rho$ is non-cyclic      & by assumption,\\
\point & $\sigma_p(\G_0)=\vide$      & by assumption, inspired by Proposition \ref{prop:pas_of_curve_peripheral},\\
\point & $|\sigma_s(A_1)|\in\{1,\,2\}$ & according to \ref{prop:cardinal_de_sigmas}\\
\point & curves of $\sigma_s(\G_0)$ are not separating & according to Proposition \ref{prop:curves_spéciales_separating}\\
\end{tabular}
\smallskip

\noindent Recall that
$\X=\{A_1,\,A_3,\,A_5,\dots,\,A_{n-1}\}$.
In this section, we will determine the arrangement of the curves of
$\sigma(\X)$ in $\Sigma$ while identifying the special curves and
the normal curves. In other words, we will be able to describe the
graphs $\Gamma(\Sigma\,;\,\sigma_s(\X))$ (cf. Definition
\ref{defi:SubASigma}) and $\Gamma(\Sigma\,;\,\sigma(\X))$. Our main
tool will be the action of the subgroup $\H$ of $\B_n$ on the subset
$\X$ of $\G_0$.
\bigskip


\noindent \point\mbox{} In Subsection \ref{par:H_colorations}, 
we state a result concerning the action of $\H$ on $\X$. We define
notably the \emph{$\H$-colorations}. They are $\H$-equivariant
functions which will help us to express the constraints coming from
the structure of $\B_n$ on the characteristic elements of the
mapping classes of $\G_0$. \smallskip

\noindent \point\mbox{} In Subsection \ref{par:embedding_sigmas(X)_in_Sigma} 
we determine the graph $\Gamma(\Sigma,\,\sigma_s(\X))$ (cf.  Proposition \ref{prop:description_sigmasX}).
\smallskip

\noindent \point\mbox{}  In Subsection
\ref{par:embedding_sigma(X)_in_Sigma} we describe of the set
$\sigma(\X)=\sigma_s(\X)\cup\sigma_n(\X)$ with Propositions
\ref{prop:g_minoré} and \ref{prop:graph_curves_normales_1}.
\smallskip


\subsection{Action of $\H$ on $\X$ and $\H$-colorations}
          \label{par:H_colorations}
\medskip

Let us recall the definitions of the group $\H$ and
of the subset $\X$ of $\G_0$ (cf. Subsection \ref{par:H_and_X}). For any nonzero integer
$m$, let $\Imp(m)$ be the set of the first odd integers smaller than
or equal to $m$. Let $\X$ be the subset $\{A_i,\ i\in\Imp(n)\}$ of
$\G_0$ consisting in $\frac{n}{2}$ elements which pairwise commute. We
denote by $\sigma(\X)$ the curve simplex
$\cup_{i\in\Imp(n)}\sigma(A_i)$. We set
$\sigma_s(\X)=\sigma(\X)\cap\sigma_s(\G_0)$. The group $\H$ is the
subgroup $\langle\,\gamma_i\,,\ i\in\Imp(n)\,\rangle$ of $\B_n$
where for all $i\in\Imp(n)$, the element $\gamma_i$ is the product
$\tau_i\tau_{i+1}\tau_i\tau_{i+2}\tau_{i+1}\tau_i$.

Let us also recall the main properties of $\H$ (see Proposition
\ref{prop:properties_of_H}):
\smallskip

\begin{itemize}
\item[\;\;i)] The action of $\H$ on $\PMod(\Sigma)$ via $\rho$ preserves
$\X$. Indeed, for all $i,j\in\Imp(n)$, we have:
\smallskip

\centrer{$\gamma_i.A_j=\rho(\gamma_i)A_j\rho(\gamma_i)^{-1}=\left\{\begin{array}{lll}
  A_{j} &\mbox{ if }& i\not\in\{j,\,j-2\}\;\\
  A_{j-2} &\mbox{ if }& i=j-2\;\\
  A_{j+2} &\mbox{ if }& i=j.
\end{array}\right.$}
\smallskip
\item[\;(ii)] The homomorphism $\H\to\Ss(\X)$ describing the action of $\H$ on $\X$, where
$\Ss(\X)$ is the symmetric group on the elements of $\X$, is
surjective. In particular, this action is $\frac{n}{2}$ times
transitive.
\item[(iii)] The action of $\H$ on $\Courb(\Sigma)$ preserves $\sigma_s(\X)$.
\end{itemize}
\bigskip

We have already given the definition of $\J$-coloration in
Subsection \ref{par:action_of_J}. We define what is a
\emph{$\H$-coloration} in the same way.
\smallskip

\TITRE{$\H$-colorations on $\X$.}\mbox{}\\
Let $\E$ be an $\H$-set (i.e. a set together with an action of $\H$)
and $\P(\E)$ the power set of $\E$. An \emph{$\H$-coloration} is a
function $\col_\X\;:\; \X\longto\P(\E)$ that is
\emph{$\H$-equivariant}, which means that for all $\xi\in\H$ and all
$A\in\X$, we have:
\smallskip

\centrer{$\xi.\col_\X(A)=\col_\X(\xi.A)$.}
\smallskip

\noindent The integers of $\Imp(n)$, in bijection with $\X$, are
called \emph{colors}. We will say that an element $e\in\E$ \emph{is
of color $i$} if $e\in\col_\X(A_i)$. An element $e\in\E$ can be of
several colors in the meantime or possibly of none color.
\smallskip

Conversely, starting from a $\H$-coloration $\col_\X\;:\;
\X\longto\P(\E)$, let us define the map called \emph{$\H$-spectrum}:
\smallskip

\centrer{$\DEF{\sp_\X}{\E}{\P(\X)}{e}{\{A\in\X\ |\
e\in\col_\X(A)\}}$.}

\smallskip

\noindent The map $\sp_\X$ is $\H$-equivariant: for all $\xi\in\H$
and all $e\in\E$, we have:
\smallskip

\centrer{$\xi.\sp_\X(e)=\sp_\X(\xi.e)$.}
\bigskip

\begin{prop} \label{prop:H-coloration}
The map $\sigma_s\;:\; \X\to\Courb(\Sigma)$ is an $\H$-coloration.
\end{prop}
\smallskip

\DEM According to Proposition \ref{prop:properties_sigma}.(i), we
have $\sigma(\gamma.A)=\gamma.\sigma(A)$ for any $\gamma\in\B_n$ and
any $A\in\PMod(\Sigma)$, so the function $\Sigma$ is an
$\H$-coloration on $\sigma(\X)$. Moreover, $\sigma_s(\X)$ is
$\H$-stable, according to Proposition \ref{prop:properties_of_H}).
Then for all $i,j\in\Imp(n)$, we have:
\smallskip

\centrer{$\begin{array}{lll}
\sigma_s(\gamma_i.A_j) &=& \sigma(\gamma_i.A_j)\,\cap\,\sigma_s(\X)\\
            &=& \gamma_i.\sigma(A_j)\,\cap\,\sigma_s(\X)\\
            &=& \gamma_i.\sigma(A_j)\,\cap\,\gamma_i.\sigma_s(\X)\\
            &=& \gamma_i.\big(\sigma(A_j)\,\cap\,\sigma_s(\X)\big)\\
            &=& \gamma_i.\sigma_s(A_j).
\end{array}$}

\noindent So the restriction of the function $\sigma$ on $\X$ is an
$\H$-coloration.\fin
\bigskip

\begin{lem}
        \label{lem:Ckg}
Let $\E$ be an $\H$-set. Let $\col_\X\;:\; \X\to\P(\E)$ be an
$\H$-coloration and $\sp_\X\;:\; \E\to\P(\X)$ the associated
$\H$-spectrum. Let $e$ be an element of $\E$ such that the
$\H$-spectrum of $e$ contains $k$ mapping classes $(0\leqslant
k\leqslant \frac{n}{2})$. Then there exists an integer
$\ell\geqslant 1$ such that $|\H.e|=\ell C_r^k$ where
$r=\frac{n}{2}$. \end{lem}
\bigskip

\DEM This is an easy application of general principles about group
actions. In a general way, if $\F$ is a finite set on which $\H$
acts, if we choose an element $f_0\in\F$ and if we denote by
$S=\Stab_\H(f_0)$ the subgroup of $\H$ that fixes $f_0$, and by $\Z$
a transversal of $\quot{\H}{S}$, then,
\begin{itemize}
\item we have the disjoint union $\displaystyle \H=\dessous{\bigsqcup}{\gamma\in\Z}\gamma. S$
\item if the action of $\H$ on $\F$ is transitive, we have
$|\Z|=\big|\quot{\H}{S}\big|=|\F|$. \end{itemize}
\smallskip

Given an element $e$ of $\E$ such that $\sp_\X(e)$ is a set of $k$
mapping classes of $\X$ with $0\leqslant k\leqslant r$ where
$r=\frac{n}{2}$, let us consider the set $\P_k(\X)$ of the subsets
of $\X$ containing $k$ elements. We replace $\F$ by $\P_k(\X)$ the
and $f_0$ by $\sp_\X(e)$. Notice that the action of $\H$ on
$\P_k(\X)$ is transitive for the action of $\H$ on $\X$ is $k$ times
transitive. We get then:\smallskip

\centrer[1]{$\displaystyle
\H=\dessous{\bigsqcup}{\gamma\in\Z}\gamma. \Stab_\H(\sp_\X(e))$,}

\centrer[2]{and
$|\Z|=\big|\quot{\H}{\Stab_\H(\sp_\X(e))}\big|=|\P_k(\X)|=\comb{r}{k}$.}

\noindent Now for all distinct $\gamma'$ and $\gamma''$ belonging to
$\Z$, the elements $\gamma'.\sp_\X(e)$ and $\gamma''.\sp_\X(e)$ are
different, so the elements of
$\gamma'.\big(\,\Stab_\H(\sp_\X(e)).e\,\big)$, that all have the
same spectrum, which is different from the spectrum of the elements
of $\gamma''.\big(\,\Stab_\H(\sp_\X(e)).e\,\big)$. Hence the two
sets $\gamma'.\big(\,\Stab_\H(\sp_\X(e)).e\,\big)$ and
$\gamma''.\big(\,\Stab_\H(\sp_\X(e)).e\,\big)$ are disjoint.
Therefore, the assertion (1) implies:
\smallskip

\centrer[3]{$\H.e=\dessous{\bigsqcup}{\gamma\in\Z}\gamma.\big(\,
\Stab_\H(\sp_\X(e)).e\,\big)$,}

\smallskip

\noindent We set $\ell=|\big(\,\Stab_\H(\sp_\X(e))\,\big).e|$. Then,
we deduce from (2) and (3) the following equality:
\smallskip

\centrer[4]{$|\H.e|=\ell\, \comb{r}{k}$.} \fin
\bigskip


\subsection{Description of the embedding of $\sigma_s(\X)$ in $\Sigma$}
        \label{par:embedding_sigmas(X)_in_Sigma}
\medskip

This subsection is devoted to the proof of the following proposition.
\bigskip

\begin{prop}[Arrangement of the curves of $\sigma_s(\X)$ in $\Sigma$]\mbox{}
        \label{prop:description_sigmasX}
\begin{itemize}
\item[\;(i)]
The set $\sigma_s(\G_0)$ contains $n$ curves (hence for all $A\in\G_0$,
we have $|\sigma_s(A)|=1$).
\smallskip

\item[(ii)]
The set $\sigma_s(\X)$ is non-separating in $\Sigma$, or it is
separating but for any $a\in\sigma_s(\X)$, the set of curves
$\sigma_s(\X)\smallsetminus\{a\}$ is non-separating. In other words,
the graph $\Gamma(\Sigma,\,\sigma_s(\X))$ is one of those depicted
in Figure \ref{fig:lesDeuxGraphesPossibles}.
\begin{figure}[!h]
 \Includegraphics{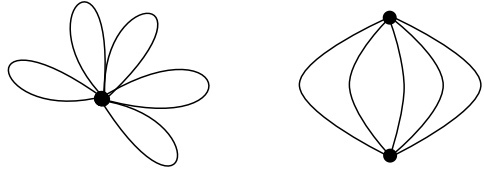}
 \caption{The two possible embeddings of $\sigma_s(\X)$ in $\Sigma$ (here $n=10$).}
 \label{fig:lesDeuxGraphesPossibles}
\end{figure}
\end{itemize}
\end{prop}
\smallskip

\DEM According to Proposition \ref{prop:cardinal_de_sigmas},
$\sigma_s(\G_0)$ contains $n$ or $2n$ curves, depending on whether
$\sigma_s(A_1)$ contains one or two curves. Since $\X$ contains only
the mapping classes with odd indices, $\sigma_s(\X)$ contains
$\frac{n}{2}$ or $n$ curves. While showing item (ii) we will show
that $\sigma_s(\X)$ contains $\frac{n}{2}$ curves, which will prove
item (i).
\smallskip

Let us show that if we prove this proposition in the case where
$\bord\Sigma$ is empty, then the proposition in the case where
$\bord\Sigma$ is not empty can be deduced easily. So, we assume that
$\bord\Sigma$ is not empty. Let us consider the graph
$\Gamma(\Sigma,\,\sigma_s(\X))$ and the subgroup $\H$ preserving the
simplex $\sigma_s(\X)$. Since the cardinality of each orbit of
curves of $\sigma_s(\X)$ under $\H$ is greater than or equal to
$\frac{n}{2}\geqslant 3$, we can apply Proposition
\ref{prop:as_if_surface_without_boundary}. In other words,
$\Gamma(\Sigma,\,\sigma_s(\X))$ is canonically isomorphic to
$\Gamma\big(\sq(\Sigma),\,\sq(\sigma_s(\X))\big)$ that is the graph
associated to the homomorphism $\sq^*\rond\rho\;:\;
\B_n\to\Mod\big(\sq(\Sigma)\big)$, in the same way as the graph
$\Gamma(\Sigma,\,\sigma_s(\X))$ is associated to the homomorphism
$\rho\;:\; \B_n\to\PMod(\Sigma)$. Thanks to this isomorphism, in
order to prove Proposition \ref{prop:description_sigmasX}, it is
enough to show the parts (i) and (ii) when $\bord\Sigma$ is empty.
\smallskip

From now on, we assume that $\bord\Sigma$ is empty. We are going to
use the action of $\H$ on the graph $\Gamma(\Sigma,\,\sigma_s(\X))$.
In order to make the action of $\H$ on the graph
$\Gamma(\Sigma,\,\sigma_s(\X))$ more obvious, we ``color in'' the
different edges and the different vertices as follows: \bigskip

\TITRE{The graph $\Gamma=\Gamma(\Sigma,\,\sigma_s(\X))$.}
\smallskip

\noindent We denote by $\Gamma$ the graph
$\Gamma(\Sigma,\,\sigma_s(\X))$. Its vertices are in bijection with
$\Sub_{\sigma_s(\X)}(\Sigma)$ and its edges are in bijection with
$\sigma_s(\X)$. The action of $\H$ on $\Sub_{\sigma_s(\X)}(\Sigma)$
induces an action of $\H$ on the set of vertices of $\Gamma$ and the
action of $\H$ on $\sigma_s(\X)$ induces an action of $\H$ on the
set of edges of $\Gamma$. Moreover, these two actions are compatible
with the graph structure of $\Gamma$. In addition, we also have two
$\H$-colorations:
\begin{itemize}
\item[\point] the $\H$-coloration $\sigma_s\;:\;
\X\to\sigma_s(\X)$, that can be seen as an $\H$-coloration on the
edges of $\Gamma$. In this way, the $\H$-spectrum of each edge of
$\Gamma$ contains a unique color $i\in\Imp(n)$;

\item[\point]the $\H$-coloration
$\DEF{\omega_\X}{\X}{\vblanc\Sub_{\sigma_s(\X)}(\Sigma)}{A}{\{S\in
\Sub_{\sigma_s(\X)}(\Sigma)\ |\ \Bord(S)\cap
\sigma_s(A)\not=\vide\}}$, \vrblanc[0,7cm] that can be seen as an
$\H$-coloration on the vertices of $\Gamma$. Thus the $\H$-spectrum
of each vertex contains exactly the colors of the edges incident to
this vertex. \end{itemize}
\smallskip

\noindent Since $\sigma_s$ and $\omega_X$ are $\H$-colorations, the
action of $\H$ on the vertices and on the edges of $\Gamma$ is
compatible with the action of $\H$ on their colors. Finally, the
group $\H$ acts on the graph $\Gamma$ together with its colors.
\bigskip

\TITRE{The cycles and the degrees in $\Gamma$.}
\smallskip

\noindent Let us denote by $c$ the number of  independent cycles
which exist in $\Gamma$. Then, we have of course:
\smallskip

\TITRE{Fact 1:} \emph{$c\leqslant g$.}
\smallskip

The vertices of degree 1 correspond to the subsurfaces bounded by
some separating curves. However, according to Proposition
\ref{prop:curves_spéciales_separating}, there does not exist in
$\sigma_s(\X)$ any separating curve in $\Sigma$, so:
\smallskip

\TITRE{Fact 2:} \emph{The degree of each vertex is at least 2.}
\smallskip

Hence the graph $\Gamma$ contains some cycles:
\smallskip

\TITRE{Fact 3:} \emph{$c\geqslant 1$.}
\smallskip

Notice that the vertices of degree 2 correspond to the connected
components having only two boundary components. Remember that
$\bord\Sigma$ is empty, so such connected components must be of
nonzero genus. We deduce:
\smallskip

\TITRE{Fact 4}: \emph{The number of vertices of degree 2 is bounded
by $g-c$.} \smallskip

Moreover, if there exists a vertex $p$ of degree 2, such that $i$
and $j$ in $\Imp(n)$ are the colors of the two edges having an
extremity in $p$, the orbit of $p$ under the action of $\H$ is of
cardinality at least $\frac{n}{2}$ if $i=j$, and of cardinality at
least $\comb{\frac{n}{2}}{2}\geqslant \frac{n}{2}$ if $i\not=j$,
according to Lemma \ref{lem:Ckg}. But according to Facts 3 and 4,
this number should be bounded by $g-1$. However, $g-1<\frac{n}{2}$,
so finally:
\smallskip

\TITRE{Fact 5}: \emph{There does not exist any vertex of degree 2.}
\bigskip

We are going to treat separately the graphs containing some edges
whose both extremities are equal, from the graphs where both
extremities of each edges are distinct.
\smallskip

\TITRE{Petals.}
            \label{defi:petals}
Let us call \emph{petal} an edge whose both extremities are equal.
\smallskip

\TITRE{Graphs with petals.} If there exist some petal in $\Gamma$,
then there exist some in each color, so there exist at least
$\frac{n}{2}$ petals. But each petal produces an independent cycle
of $\Gamma$, and the number of independent cycles is bounded by $g$,
and so by $\frac{n}{2}$. So we have $c=g=\frac{n}{2}$. Therefore, if
there exist some petals in $\Gamma$, then there exist exactly
$\frac{n}{2}$ petals and $g=\frac{n}{2}$. Moreover, since the
maximal number of independent cycles is achieved just because of the
petals, it follows that if we remove these $\frac{n}{2}$ petals from
$\Gamma$, we get a tree. Since each edge is \emph{separating} in a
tree (we say that an edge is \emph{separating} in a connected graph
if removing this edge from the graph makes it disconnected), then
all the edges in $\Gamma$ that are not a petal are separating. But
according to Proposition \ref{prop:curves_spéciales_separating},
$\sigma_s(\X)$ does not contain any separating curves, so $\Gamma$
contains no other edge but the $\frac{n}{2}$ petals. Therefore, the
graph $\Gamma$ is a rose, that is, a graph with only one vertex (cf.
graph on the left-hand side on Figure \ref{fig:lesDeuxGraphesPossibles}).
\bigskip

\TITRE{Graphs without petals.} Let us assume now that $\Gamma$ is a
graph without petals.
\smallskip

According to Proposition \ref{prop:cardinal_de_sigmas}, for each
$A\in\G_0$, the cardinality $|\sigma_s(A)|$ equals 1 or 2. In other
words, for each $i\in\{0,\,1,\dots,\,n-1\}$, there are only one or
two curves in $\sigma_s(\G_0)$ that are in $\sigma_s(\A_i)$. In the
graph $\Gamma$, this implies that:
\medskip

\TITRE{Fact 6}: \emph{There exist at most two edges of a same
color.}
\bigskip

\TITRE{If there is only one edge per color:}
\smallskip

According to Facts 2 and 5, all the vertices are at least of degree
3. Since there is no petal, at least three distinct edges are
incident to each vertex.
\smallskip

\begin{itemize}
\item
When $n=6$, since there is only one edge per color, $\Gamma$
contains only three edges. Hence in the case $n=6$, there exists
exactly 2 vertices of degree 3 and the graph is drawn in Figure
\ref{fig:lesDeuxGraphesPossibles}, left-hand side.
\smallskip

\item
Suppose now that $n\geqslant 8$. We take $j=n-1$ and we set
$\H_{n-2}=\langle\gamma_i,\,i\in\Imp(n-2)\rangle$ so that $\H_{n-2}$
is a subgroup of $\H$ that fixes the color $n-1$ and acts
$\frac{n}{2}-1$ times transitively on the set of the $\frac{n}{2}-1$
other colors of $\Imp(n)$. Given a vertex $q$ of degree $v$ and an
edge of color $n-1$ incident to $q$, according to Lemma
\ref{lem:Ckg}, there should exist at least
$\comb{(\frac{n}{2}-1)}{v-1}$ distinct vertices in the orbit of $q$
under the action of $\H_{n-2}$. But there is only one edge of color
$n-1$ in $\Gamma$, so there should be at most two vertices in the
orbit of $q$ under $\H_{n-2}$, including $q$. So
$\comb{(\frac{n}{2}-1)}{v-1}$ must be equal to 1 or 2. Since
$\frac{n}{2}-1\geqslant 3$ and $v-1\geqslant 2$, $v$ must be equal
to $\frac{n}{2}$. So the graph $\Gamma$ is fully determined: there
are exactly two vertices and $\frac{n}{2}$ edges incident to these
two vertices: the graph is drawn in Figure
\ref{fig:lesDeuxGraphesPossibles}, right-hand side.
\end{itemize}
\bigskip


\TITRE{If there are exactly two edges per color:}
\smallskip

\centrer{\textbf{We suppose from now on that there exist two curves
in each color.}}
\smallskip

We want to show that this case is absurd.
\smallskip

\noindent Notice that we have the following: \bigskip

\TITRE{Fact 7}: \emph{Two edges of the same color cannot have the
same extremities.}
\smallskip

\begin{figure}[!h]
 \Includegraphics{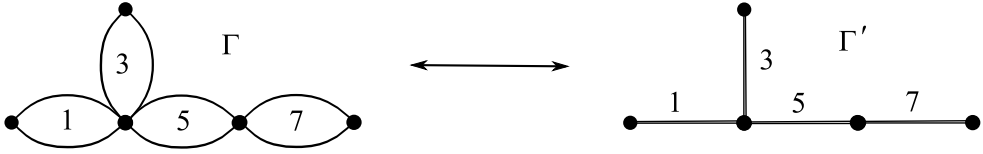}
 \caption{Case $g\geqslant4$, when two edges of a same color have the same extremities.}
 \label{fig:grapheSigmaGg4}
\end{figure}

\DEM Assume that two edges of a same color have the same
extremities. Then they would form a cycle, and we can associate
their color to this cycle. We would obtain $\frac{n}{2}$ independent
cycles of distinct colors, hence the equalities $c=g=\frac{n}{2}$
would hold and it would not exist other independent cycles in
$\Gamma$, according to Fact 1. It would not exist either other edges
than the $n$ edges constituting these $\frac{n}{2}$ cycles,
according to Fact 6. Let us identify in $\Gamma$ two edges if they
have the same color and let us call $\Gamma'$ this new graph (cf.
Figure \ref{fig:grapheSigmaGg4}). The graph $\Gamma'$ is a tree, so
it contains leaves (vertices of degree 1). These leaves correspond
in $\Gamma$ to vertices of degree 2. But this is forbidden,
according to Fact 5. Hence Fact 7 is shown.\fin
\bigskip

Let us go further:
\smallskip

\TITRE{Fact 8}: \emph{Two edges of the same color cannot share an
extremity in common.}
\smallskip

\DEM We argue again by contradiction. We suppose that two edges of
the same color share an extremity in common (if this is true for one
color, this must be true for each color). Let us recall that there
exist only two edges of the same color, according to Fact 6. Let us
also recall that two such edges cannot have the same extremities,
according to Fact 7. So for all $i\in\Imp(n)$, there exists only one
vertex that is a common extremity the two edges of color $i$. We
will call it $p_i$. By symmetry of the action of $\H$ on $\X$, the
$p_i$ are all distinct or all equal.
\smallskip

\begin{itemize}
\item[\textbf{a)}] \textbf{The vertices $p_i$, $i\in\Imp(n)$, coincide.}
We keep in mind Figure \ref{fig:graphePiConfondus}. Let us call $p$
the vertex $p_1=p_3=\dots=p_{n-1}$. The vertex $p$ is at least of
degree $n$. It cannot be of a greater degree, for there is no petal
in $\Gamma$, by assumption. Thus, each on the $n$ edges in $\Gamma$
has an extremity in $p$. The number of the other extremities is $n$,
and they are incident to vertices distinct from $p$. But since these
vertices are of degree greater than or equal to 3 according to Facts
2 and 5, their number is at most $\frac{n}{3}$. So, $\Gamma$
contains $n$ edges and at most $1+\frac{n}{3}$ vertices, so its
number of independent cycles is at least of
$1+n-(\frac{n}{3}+1)=\frac{2n}{3}$. But this is absurd for
$\frac{2n}{3}>\frac{n}{2}\geqslant g$.
\begin{figure}[!h]
 \Includegraphics{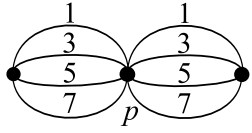}
 \caption{Example of graph where all the $p_i$ coincide.}
 \label{fig:graphePiConfondus}
\end{figure}

\item[\textbf{b)}] \textbf{The vertices $p_i$, $i\in\Imp(n)$, are pairwise distinct.}
We keep in mind Figure \ref{fig:graphePiDistincts}. As the $p_i$ are
not of degree 2, according to Fact 5, each $p_i$ is the extremity of
the two edges of color $i$ and of at least an edge of color
$j\not=i$. But the subgroup of $\H$ that fixes the color $i$ acts
transitively on the other colors, so $p_i$ is also the extremity of
an edge of color $k$ for all $k\in\Imp(n)\smallsetminus\{i\}$. Thus
$p_i$ is of degree at least $\frac{n}{2}+1$ (actually, this must be
an equality with our assumptions). As we know that there is at least
$\frac{n}{2}$ vertices in $\Gamma$ (think of the $p_i$,
$i\in\Imp(n)$), we deduce that the degrees of all the vertices in
$\Gamma$ sum to at least $\frac{n}{2}(\frac{n}{2}+1)$. Now, this sum
should be equal to two times the number of edges. Since there are
exactly $n$ edges in $\Gamma$, we have the equality:\smallskip

\centrer[$*$]{$\frac{n}{2}(\frac{n}{2}+1)\leqslant 2n$.}
\smallskip

\noindent We get then $n^2+2n\leqslant 8n$, hence $n\leqslant 6$. By
hypothesis, $n\geqslant 6$, so $n=6$ and $(*)$ is a equality and
becomes: $\frac{n}{2}(\frac{n}{2}+1)=12$. That means that the
vertices in $\Gamma$ is reduced to the set $\{p_1,\,p_3,\,p_5\}$.
The graph now is perfectly determined, and drawn in Figure
\ref{fig:graphePiDistincts}. But this graph contains $c=1+6-3=4$
independent cycles, which is absurd for $c\leqslant g\leqslant
\frac{n}{2}=3$.
\begin{figure}[!h]
 \Includegraphics{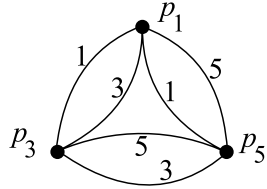}
 \caption{Example of a graph where the $p_i$ are distinct.}
 \label{fig:graphePiDistincts}
\end{figure}
\end{itemize}
\smallskip

\noindent Thus, Fact 8 is shown.\fin
\bigskip

Now, with these eight facts, we can terminate the proof of
Proposition \ref{prop:description_sigmasX}. Remember that we assume
that $|\sigma_s(\X)|=n$. We prove separately the cases $n=6$ and
$n\geqslant 8$. We start by the case $n\geqslant 8$ which is the
easiest one.
\smallskip

\TITRE{Let us show that the conditions $|\sigma_s(\X)|=n$ and
$n\geqslant 8$ lead to a contradiction.}

\smallskip

\noindent Let us denote by $S$ the number of vertices. Since
$\Gamma$ has $n$ edges, the number of independent cycles must
satisfy $c=1+n-S\leqslant\frac{n}{2}$. Hence:
\smallskip

\centrer{$S\geqslant \frac{n}{2}+1$.}

\noindent Let $v$ be the minimum degree among the vertices of
$\Gamma$. We have seen, according to Facts 2 and 5 that
$v\geqslant3$. But the sum of the degrees of all the vertices of
$\Gamma$ is at least equal to $vS$. On the other hand, it must be
equal to two times the number of vertices. Hence we have
$2n\geqslant vS$, so:\smallskip

\centrer{$S\leqslant\frac{2n}{v}$.}
\smallskip

\noindent It is then absurd that $v\geqslant4$ for we would have
$\frac{n}{2}+1\leqslant S\leqslant \frac{n}{2}$. Now, suppose that
$v=3$. Thanks to Fact 8, we know that all the edges incident to a
vertex are of pairwise distinct colors, so when $v=3$, to each
vertex corresponds a choice of three colors among $\frac{n}{2}$.
According to Lemma \ref{lem:Ckg}, we have then at least
$\comb{(\frac{n}{2})}{3}$ vertices in $\Gamma$ . When $n\geqslant
10$, we check that $\comb{(\frac{n}{2})}{3}\geqslant n$, so this
leads to the following contradiction:\smallskip

\centrer{$n\leqslant \comb{(\frac{n}{2})}{3}\leqslant S\leqslant
\frac{2n}{3}$.}
\smallskip

\begin{figure}[!h]
 \Includegraphics{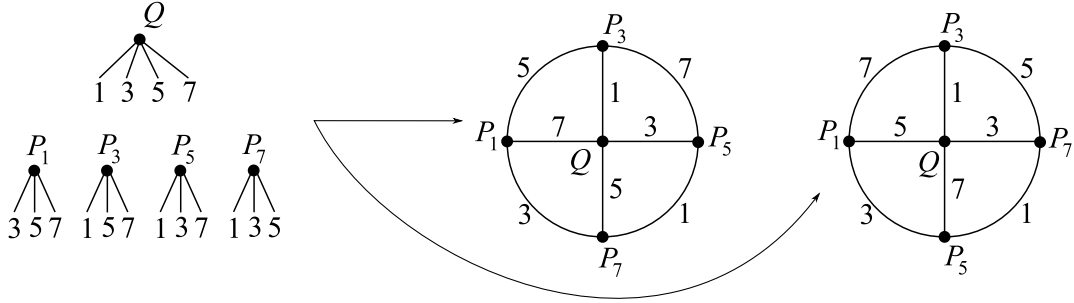}
 \caption{Two examples of \emph{a priori} possible graphs with $n=8$ and $|\sigma_s(\X)|=8$.
 \emph{(They
 are different, for $Q$ and $P_1$ are joined
 by an edge of color 7 in one case, and 5 in the other.
 But, none of them is $\H$-stable since the action of $\gamma_5$ swap them.)}}
 \label{fig:grapheCasN4A8}
\end{figure}

\noindent When $n=8$, since $v=3$, there exists an orbit of vertices
of degree 3. The cardinality of this orbit is a multiple of 4
vertices, according to Lemma \ref{lem:Ckg}. But as there exists 8
edges in $\Gamma$ and consequently 16 half-edges, there must exist
in $\Gamma$ exactly 4 vertices of degree 3 and one vertex of degree
4 (replace the vertex of degree 4 by two vertices of degree 2 is
forbidden by Fact 5) to satisfy the equality: $16 = 3\times
4+1\times 4$. We thus get a graph similar to those depicted in
Figure \ref{fig:grapheCasN4A8}. Let us call $P_1$, $P_3$, $P_5$ and
$P_7$ the four vertices of degree 3 and $Q$ the vertex of degree 4,
as in Figure \ref{fig:grapheCasN4A8}. The subgroup of $\H$ that
fixes the color 1 does not preserve such a graph, for it fixes the
vertices $P_1$ and $Q$ but modifies the color of the unique edge
that joint the vertices $P_1$ and $Q$: this edge can be of color 3,
5 or 7, and gives rise to different graphs (in Figure
\ref{fig:grapheCasN4A8}, we give two examples with 7 and 5). Hence
this case $n=8$ is absurd. Finally, the conditions
$|\sigma_s(\X)|=n$ and $n\geqslant 8$ lead only to contradictions.
\bigskip

\TITRE{Let us show that the conditions $|\sigma_s(\X)|=n$ and $n=6$
lead to a contradiction.}
\smallskip

\begin{figure}[!h]
 \Includegraphics{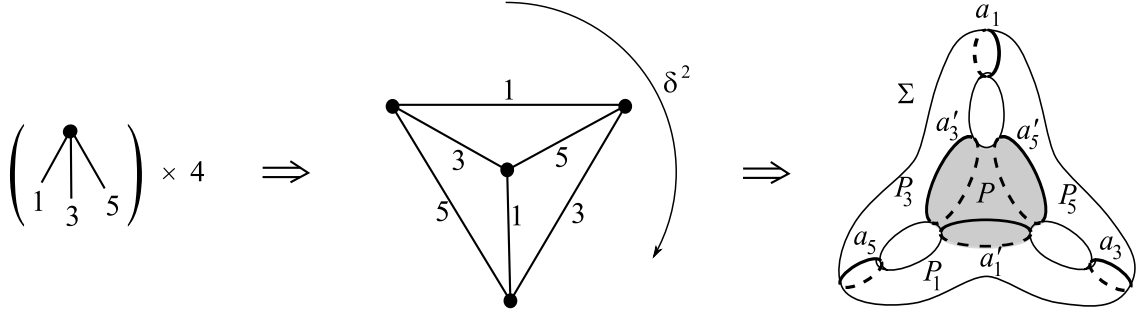}
 \caption{Case $n=6$: the surface $\Sigma_{3,\,0}$ together with the curves of $\sigma_s(\X)$ when $|\sigma_s(\X)|=6$.}
 \label{fig:Surface30Avec6Courbes}
\end{figure}

\noindent According to Facts 2, 5 and 8, all the vertices are of
degree 3, each being the extremity of three edges of colors 1, 3 and
5. We get the graph represented on Figure
\ref{fig:Surface30Avec6Courbes}. Let us recall that according to
Proposition \ref{prop:spéciales_normales} (5d) describing the
equivalent properties of the special curves, each orbit of edges of
$\sigma_s(\X)$ under the action of $\J=\langle\delta\rangle$ is of
cardinality $n=6$, and each orbit of edges of $\sigma_s(\X)$ under
the action of $\langle\delta^2\rangle$ is of cardinality
$\frac{n}{2}=3$. So we distinguish two orbits among the 6 curves of
$\sigma_s(\X)$, and the action of $\delta^2$ on the graph $\Gamma$
is periodic of order three. Then one of the vertices is preserved
whereas the three others are cyclicly permuted. Let us call $P$ the
preserved vertex. Its boundary components form an orbit under the
action of $\delta^2$, we will call them $a'_1$, $a'_3$, $a'_5$,
where the indices respect the color of each curve. Let us call the
other curves $a_1$, $a_3$, $a_5$, where the indices respect the
color of each curve. For all $i\in\{1,\,3,\,5\}$, let us denote by
$P_i$ the subsurface different from $P$ that contains the curve
$a'_i$ in its boundary, see Figure \ref{fig:Surface30Avec6Courbes}.
\smallskip

We are going to show, by means of the special curves of
$\Y=\G_0\smallsetminus\X=\{A_0,\,A_2,\,A_4\}$, that such a
configuration cannot happen. Notice that $\sigma_s(\G_0)$ is a set of
curves without triple intersection, so according to Proposition
\ref{prop:existence_uniqueness_representatives_curves}, there exists
a unique system of representatives of the curves of $\sigma_s(\G_0)$
in tight position, up to isotopy. Moreover, any representative of
the mapping class $\rho(\delta)$ sends such a system of
representatives on itself, up to isotopy. Moreover, we have,
$\sigma_s(\Y)=\delta.\sigma_s(\X)$, and we denote by $a_0$, $a_2$,
$a_4$, $a'_0$, $a'_2$, $a'_4$ the curves of $\sigma_s(\Y)$, so as to
be coherent with their colors, as we did with the curves of
$\sigma_s(\X)$. Let us consider the subsurface $P_1$. Then:
\smallskip

\begin{itemize}
  \item[\point]
$a'_2=\delta.a'_1$ and $a'_0=\delta^{-1}.a'_1$, so according to
Proposition \ref{prop:propiétés_normales_spéciales},
$I(a'_1,\,a'_2)=I(a'_1,\,a'_0)\not=0$, hence $a'_0\cap P_1$ and
$a'_2\cap P_1$ are not empty.
\smallskip

  \item[\point]
Similarly if we consider the curves $a_3$ and $a_5$, it follows that
$I(a_3,\,a_4)=I(a_3,\,a_2)\not=0$ and
$I(a_5,\,a_0)=I(a_5,\,a_4)\not=0$, so $a_0\cap P_1$, $a_2\cap P_1$
and $a_4\cap P_1$ are not empty.
\smallskip

  \item[\point]
Let $x$ be a curve of $\sigma_s(\Y)$ whose restriction to $P_1$ is a
path with extremities in a same boundary component, and $y$ another
curve of $\sigma_s(\Y)$ whose restriction to $P_1$ is a path with
extremities in the two other boundary components of $P_1$. Then $x$
and $y$ must intersect, as illustrated in Figure
\ref{fig:CourbeExtremitesMemeBord} (two cases are to be considered,
depending on whether the extremities of $y$ belong to the same
boundary of $P_1$ or not). But $\sigma_s(\Y)$ is a curve simplex, so
this situation cannot happen.
\begin{figure}[!h]
 \Includegraphics{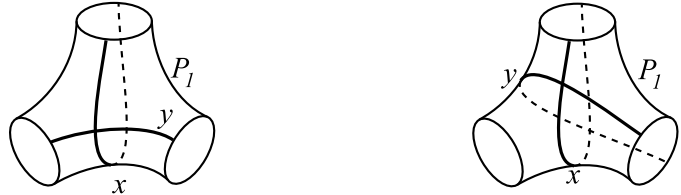}
 \caption{Case $n=6$: the curve $x$ and the curve $y$ have to intersect.}
 \label{fig:CourbeExtremitesMemeBord}
\end{figure}
\smallskip

  \item[\point]
Let us apply the last point to the curves $x=a'_0$ and $y=a_4$.
Suppose that the extremities of one of the connected components of
$a'_0\cap P_1$ lie in $a'_1$. Notice that $a_4$ cannot intersect
$a'_1$, so the extremities of any connected component of $a_4\cap
P_1$ lie in $a_3\cup a_5$. This leads to a contradiction as we just
have seen it above. Therefore $a'_0\cap P_1$ cannot contain any path
with extremities in $a'_1$. Hence $a'_0\cap P_1$ consists in paths
joining the boundary components $a'_1$ and $a_5$. Similarly,
$a'_2\cap P_1$ consists in paths joining the boundary components
$a'_1$ and $a_3$.
\smallskip

  \item[\point]
Let us now apply the last-but-one point to the curves $x=a_2$ and
$y=a_0$. We conclude just as before that $a_2\cap P_1$ cannot
contain any path with extremities in $a_3$. Similarly, if we take
the curves $x=a_4$ and $y=a_0$, we see that $a_4\cap P_1$ cannot
contain any path with extremities in $a_3$. The reader can check
that in fact, no curve of $\sigma_s(Y)$ can contain some path in
$P_1$ with extremities in a same boundary of $P_1$. The situation is
summed up in Figure \ref{fig:CourbeExtremitesBordDifferent}.
\begin{figure}[!h]
 \Includegraphics{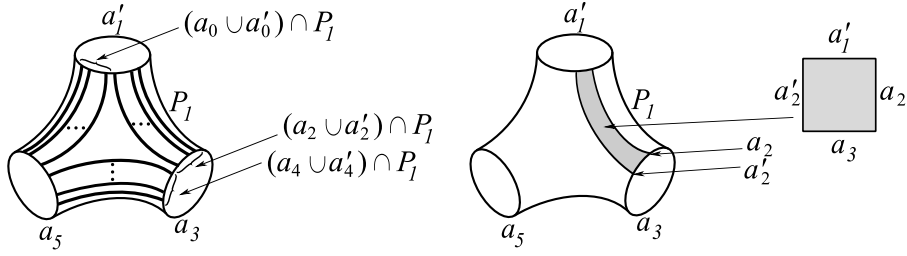}
 \caption{Case $n=6$: the intersection of the curves of $\sigma_s(\Y)$ and $P_1$.}
 \label{fig:CourbeExtremitesBordDifferent}
\end{figure}
\smallskip

  \item[\point]
Then, there exists in $\Sigma_{\sigma_s(\G_0)}$ a connected component
homeomorphic to a disk whose boundary consists in four arcs of
curves, each arc being included in  $a'_1$, $a_2$, $a'_2$ and $a_3$
respectively. We then deduce that there exists in
$\Sub_{\sigma_s(\Y)}(\Sigma)$ a subsurface bounded altogether by
$a_2$ and $a'_2$. Then, taking the image of this situation by
$\rho(\delta^{-1})$, we deduce that there should exist in
$\Sub_{\sigma_s(\X)}(\Sigma)$ a subsurface bounded altogether by
$a_1$ and $a'_1$. However, this is not the case. This is the
expected contradiction.
\end{itemize}
\smallskip

\noindent Finally, the conditions $|\sigma_s(\X)|=n$ and $n=6$ lead
to a contradiction. This concludes the proof of Proposition
\ref{prop:description_sigmasX}.\fin
\bigskip




\subsection{Description of $\sigma_n(\G_0)$ and of $\sigma(\X)$ in $\Sigma$}
        \label{par:embedding_sigma(X)_in_Sigma}
\medskip

According to Proposition \ref{prop:description_sigmasX}, for all
integers $i\in\{0,\dots,\,n-1\}$, the set of curves $\sigma_s(A_i)$
is reduced to a unique curve. We denote it by $a_i$. Since we know
that the set $\sigma_s(\X)$ is equal to
$\{a_1,\,a_3,\dots,\,a_{n-1}\}$  and since we know how these curves
are arranged in $\Sigma$ (cf. Proposition
\ref{prop:description_sigmasX}), we turn now to the simplex
$\sigma(\X)$. Let us recall that according to Proposition
\ref{prop:propiétés_normales_spéciales}, we have
$\sigma_n(A_1)=\sigma_n(\G_0)$, whence the equality
\smallskip

\centrer{$\sigma(\X)=\sigma_s(\X)\cup\sigma_n(\G_0)$.}
\smallskip

\noindent Recall that $\sigma(\X)$ is a simplex since $\X$ span an
abelian group. Consequently, $\sigma_s(\X)$ and $\sigma_n(\G_0)$ are
also simplexes. Moreover, $I(\sigma_s(\X),\,\sigma_n(\G_0))=0$.
Actually, we have the following, which is stronger, and which comes
from the properties of the normal curves (see Proposition
\ref{prop:propiétés_normales_spéciales}.(ii)):
\smallskip

\centrer{$I(\sigma_s(\G_0),\,\sigma_n(\G_0))=0$.}
\bigskip

\begin{prop}[Existence of the surface $\widehat\Sigma$]
        \label{prop:g_minoré}
\mbox{}
\begin{itemize}
\item[\;\;(i)] There exists a unique subsurface $\widehat\Sigma$ in
$\Sub_{\sigma_n(\G_0)}(\Sigma)$ that contains the curves of
$\sigma_s(\G_0)$.
\item[\;(ii)] The boundary of each subsurface belonging to $\Sub_{\sigma_s(\X)}(\widehat\Sigma)$
contains $\sigma_s(\X)$.
\item[(iii)] The surface $\widehat\Sigma$ is of genus $\widehat g\in\{\frac{n}{2}-1,\,\frac{n}{2}\}$.
\item[(iv)] The surface $\Sigma$ is of genus $g\in\{\frac{n}{2}-1,\,\frac{n}{2}\}$.
\item[\;(v)] The set $\Sub_{\sigma_s(\X)}(\widehat\Sigma)$
is reduced to a connected component of genus zero, or to two
connected components, one of them is of genus zero and the other is
of genus zero or one.
\end{itemize}
\end{prop}
\medskip

\DEM
\medskip

\emph{Let us show item (i).}
\smallskip

Let us recall that the normal curves do not intersect any curve of
$\sigma(\G_0)$ and hence any special curve. For all
$i\in\{0,\,1,\dots,\,n-1\}$, let $\widehat\Sigma_i$ be the
subsurface of $\Sub_{\sigma_n(\G_0)}(\Sigma)$ that contains the curve
$a_i$. Let us recall that for all $i\in\{0,\,1,\dots,\,n-1\}$, the
curve $\delta.a_i$ is a special curve of $\sigma_s(A_{i+1})$, so
$\delta.a_i=a_{i+1}$. But $I(a_i,\,\delta.a_i)\not=0$ according to
Proposition \ref{prop:propiétés_normales_spéciales} for $a_i$ is
special, so $I(a_i,\,a_{i+1})\not=0$, hence
$\widehat\Sigma_i=\widehat\Sigma_{i+1}$. Thus all the special curves
are included in a unique subsurface $\widehat\Sigma$ of
$\Sub_{\sigma_n(\G_0)}(\Sigma)$.
\medskip

\emph{Let us show item (ii).}
\smallskip

Let $S$ be a subsurface belonging to
$\Sub_{\sigma_s(\X)}(\widehat\Sigma)$. If none of the boundary
components of $S$ is a normal curve, then $S$ belongs to
$\Sub_{\sigma_s(\X)}(\Sigma)$ and according to Proposition
\ref{prop:description_sigmasX}, the boundary of $S$ contains
$\sigma_s(\X)$ and item (ii) is proved in this case.

Let us then assume that $\sigma_n(\G_0)$ is nonempty and let us focus
on the surfaces of $\Sub_{\sigma_s(\X)}(\widehat\Sigma)$ such that
at least one of the boundary components is a normal curve. We define
a subgroup $\H^*$ of $\H$ ($\H$ was defined in Subsection \ref{par:H_and_X}) by:
\smallskip

\centrer{$\H^*:=\langle\,\gamma_i\gamma_j^{-1}\,,\
i,j\in\Imp(n)\,\rangle$}
\medskip

\noindent It acts via $\rho$ on $\Courb(\Sigma)$ and notably on
$\sigma(\X)$ so that: \smallskip

\begin{itemize}
\item \emph{Action of $\H^*$ via $\rho$ on $\sigma_s(\X)$}.
We have seen that the action of $\H$ via $\rho$ on $\Courb(\Sigma)$
preserves $\sigma_s(\X)$, hence so does the action of $\H^*$ on
$\Courb(\Sigma)$, since $\H^*$ is included in $\H$. Moreover, the
homomorphism $\H\to\Ss(\X)$ sends $\H^*$ on $\Aa(\X)$, the alternating
group on the elements of $\X$. Indeed, for all $i\in\Imp(n)$, this
homomorphism sends $\gamma_i\gamma_{i+2}^{-1}$ on the circular
permutation on the three elements $A_i$, $A_{i+2}$, $A_{i+4}$. But
$\Aa(\X)$ acts transitively on $\X$ (recall that $\X$ contains at
least three elements). Hence the action of $\H^*$ on $\PMod(\Sigma)$
preserves $\X$ and is transitive on $\X$. Therefore $\H^*$ acts
transitively on $\sigma_s(\X)$.
\smallskip

\item \emph{Action of $\H^*$ via $\rho$ on $\sigma_n(\G_0)$}.
Let us recall that according to Proposition
\ref{prop:stability_des_curves_normales}, the action of $\B_n$ on
$\sigma_n(\G_0)$ and on $\Bord(\Sigma_{\sigma_n(\G_0)})$ is cyclic. Then
the action of $\F_n=\Ker(\{\lambda\ :\ \B_n\to\ZZ,\ \tau_1\mapsto1\})$ on
$\sigma_n(\G_0)$ and on $\Bord(\Sigma_{\sigma_n(\G_0)})$ is trivial
according to Lemma \ref{lem:cyclic_and_F_n}.(iii). Since $\H^*$ is a subgroup
of $\F_n$, the action of $\H^*$ on $\Courb(\Sigma)$ fixes each curve
of $\sigma_n(\G_0)$ and each boundary component of
$\Sigma_{\sigma_n(\G_0)}$.
\smallskip

\item \emph{Action of $\H^*$ via $\rho$ on $\Sub_{\sigma_n(\G_0)}(\widehat\Sigma)$}.
According to the action of $\H^*$ on $\sigma_n(\G_0)$, the action of
$\H^*$ via $\rho$ on $\Sub(\Sigma)$ preserves $\widehat\Sigma$ and
preserves the set $\Sub_{\sigma_s(\X)}(\widehat\Sigma)$ of
subsurfaces. Let $S$ be a surface belonging to
$\Sub_{\sigma_s(\X)}(\widehat\Sigma)$ such that at least one of the
boundary components is a normal curve. According to the action of
$\H^*$ via $\rho$ on $\sigma_n(\G_0)$, the action of $\H^*$ via $\rho$
on $\Sub(\Sigma)$ preserves the surface $S$.
\end{itemize}
\medskip

\noindent Let us exploit this. For any subsurface $S$ belonging to
$\Sub_{\sigma_s(\X)}(\widehat\Sigma)$ such that $\Bord(S)$ contains
a normal curve, $S$ must be stable by $\H^*$. Hence $\Bord(S)$ is
$\H^*$-stable,  and so is $\Bord(S)\cap\sigma_s(\X)$. So $\Bord(S)$
must contain all $\sigma_s(\X)$ since $\H^*$ acts transitively on
$\sigma_s(\X)$.
\bigskip

\emph{Let us show items (iii) and (iv).}
\smallskip

According to item (ii), each subsurface belonging to
$\Sub_{\sigma_s(\X)}(\widehat\Sigma)$ contains at least
$\frac{n}{2}$ special curves in its boundary. Since $\sigma_s(\X)$
contains $\frac{n}{2}$ curves, $\Sub_{\sigma_s(\X)}(\widehat\Sigma)$
contains one only connected component having $n$ special boundary
components, or $\Sub_{\sigma_s(\X)}(\widehat\Sigma)$ contains two
connected components having each $\frac{n}{2}$ special boundary
components. In the first case $\widehat\Sigma$ contains a
non-separating set of $\frac{n}{2}$ curves, so $\widehat\Sigma$ is
of genus $\widehat g\geqslant\frac{n}{2}$. In the second case
$\widehat\Sigma$ contains a non-separating set of $\frac{n}{2}-1$
curves, so $\widehat\Sigma$ is of genus $\widehat
g\geqslant\frac{n}{2}-1$. Let us recall that by hypothesis,
$g\leqslant \frac{n}{2}$. Since $\Sigma$ contains $\widehat\Sigma$,
we have $\widehat g\leqslant g$. So finally, we have:
\smallskip

\centrer{$\frac{n}{2}-1\leqslant\widehat g\leqslant
g\leqslant\frac{n}{2}$}
\bigskip

\emph{Let us show item (v).}
\smallskip

We have just seen in items (iii) and (iv) that
$\Sub_{\sigma_s(\X)}(\widehat\Sigma)$ is a set of exactly one or two
subsurfaces.
\begin{itemize}
\item
If $\Sub_{\sigma_s(\X)}(\widehat\Sigma)$ contains only one
subsurface, let us denote it by $S$. The surface $\widehat\Sigma$ is
the gluing of $S$ on itself by making coincide the $n$ boundary
components between them. The difference between the genera of
$\widehat\Sigma$ and $S$ is then $\frac{n}{2}$, but $\widehat\Sigma$
is already of maximal genus $\frac{n}{2}$, so $S$ is a genus-0
surface.

\item If $\Sub_{\sigma_s(\X)}(\widehat\Sigma)$ contains two
subsurfaces, let us denote them by $S_1$ and $S_2$. The gluing of
$S_1$ on $S_2$ along the $n$ special boundary components brings a
contribution of $\frac{n}{2}-1$ to the genus of $\widehat\Sigma$.
Since the genus of $\widehat\Sigma$ is at most $\frac{n}{2}$, the
sum of the genera of $S_1$ and of $S_2$ is at most 1.
\end{itemize}


\fin
\bigskip

Proposition \ref{prop:g_minoré} allows us to set the following
definitions that will still be useful in Section \ref{sec:end_grande_proof}.
Actually, the full description of $\sigma_n(\G_0)$ is useless to prove Theorem \ref{thm:Theorem_principal}.
Focussing on the subset $\U$ (defined just below) of $\sigma_n(\G_0)$ will be enough.
\smallskip

\begin{defi}[$\widehat\Sigma$, $\check\Sigma$ and $\U$]
            \label{defi:U}
\mbox{} \begin{itemize}
  \item[\point]
Let $\widehat\Sigma$ be the subsurface of
$\Sub_{\sigma_n(\G_0)}(\Sigma)$ that contains the special curves.

  \item[\point]
Let us set $\U=\Bord(\widehat\Sigma)\cap\sigma_n(\G_0)$.

  \item[\point]
Let $\check\Sigma$ be the union of the subsurfaces of
$\Sub_{\U}(\Sigma)$ different from $\widehat\Sigma$. If
$\widehat\Sigma$ is the only subsurface of $\Sub_{\U}(\Sigma)$, we
will say that $\check\Sigma$ is empty. \end{itemize} \end{defi}
\medskip

\TITRE{Attention:} {\em Let us recall that we have supposed in this
section that $\sigma_p(\G_0)$ was empty, in other words, the
separating curves belonging to $\sigma_n(\G_0)$ separate $\Sigma$ in
two surfaces of nonzero genus. Without this hypothesis, Proposition
\ref{prop:graph_curves_normales_1} would be false.}
\bigskip

\REM If $\check\Sigma$ is nonempty, we have
$\Bord(\widehat\Sigma)\cap\Bord(\check\Sigma)=\U$. In the following
section, the only information about $\sigma_n(\G_0)$ that will help us
concerns $\U$. That is why we focus only on $\U$. According to this
remark, the next proposition deals only with $\U$ instead of
$\sigma_n(\G_0)$. All the same, we terminates this section by giving
all the possible graphs of $\Gamma(\Sigma,\,\sigma(\X))$.
\bigskip

\begin{prop}[Description of $\Sub_\U(\Sigma)$]
          \label{prop:graph_curves_normales_1}
\mbox{}\\ We have $|\U|\leqslant 2$, and $\Sub_\U(\Sigma)$ satisfies
the following properties:
\begin{itemize}
  \item
if $\U$ is empty, then $\Sigma=\widehat\Sigma$; if $\widehat\Sigma$
is of genus $\frac{n}{2}$, then $\U$ is empty,
  \item
if $\U$ is reduced to a non-separating curve $u$, then
$\check\Sigma$ is empty and $\{u\}$ is the mark of $\widehat\Sigma$,
  \item
if $\U$ is reduced to a separating curve $u$, then $\check\Sigma$ is
a connected subsurface of genus 1,
  \item
if $\U$ contains two curves, then they are non-separating and
$\check\Sigma$ is connected, of genus zero.
\end{itemize}

\noindent These assertions can be summed up as follows: the graph
$\Gamma(\Sigma,\,\U)$ is one of the four graphs depicted in Figure
\ref{fig:GraphesCourbesNormales1} where:
\begin{itemize}
\item the circled vertices represent the subsurfaces of nonzero genus of $\Sub_{\U}(\Sigma)$,
\item the integer placed beside the circled vertices indicates the genus.
\end{itemize}

\begin{figure}[!h]
 \Includegraphics{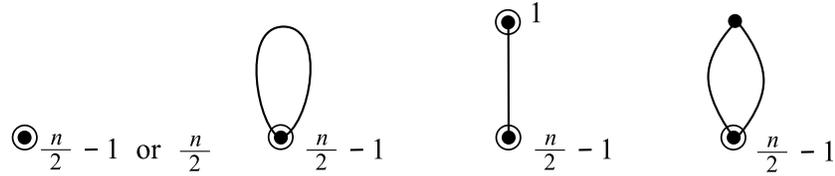}
 \caption{The four possible graphs for
 $\Gamma(\Sigma,\,\U)$ (cf. Proposition \ref{prop:graph_curves_normales_1}).}
 \label{fig:GraphesCourbesNormales1}
\end{figure}
\end{prop}
\smallskip

\DEM The proof comes from a genus computation. Let us recall that
$g$ is the genus of $\Sigma$. Let us denote by $\widehat g$ the
genus of $\widehat\Sigma$ and by $\check g$ the genus of
$\check\Sigma$ (which is zero by convention when $\check\Sigma$ is
empty). According to Proposition \ref{prop:g_minoré}, we have:
\smallskip

\centrer{$g\in\{\frac{n}{2}-1,\,\frac{n}{2}\}$, \ $\widehat
g\in\{\frac{n}{2}-1,\,\frac{n}{2}\}$, \ $\check g\in\{0,\,1\}$.}
\smallskip

\noindent Notice that each (nonempty) connected component of
$\check\Sigma$, that is separated from $\widehat\Sigma$ by a
separating curve of $\U$, is of nonzero genus, for we have assumed
that there exists no peripheral curve in $\sigma(\G_0)$. Hence it can
exist only one nonempty connected component of $\check\Sigma$
separated from $\widehat\Sigma$ by a separating curve. Besides, if a
(nonempty) connected component of $\check\Sigma$ is of genus zero,
it is separated from $\widehat\Sigma$ by at least two curves of $\U$
(again because the curves of $\U$ are not peripheral curves). So the
gluing of $\widehat\Sigma$ and of a genus-0 connected component of
$\check\Sigma$ brings also a nonzero contribution to the genus of
$\Sigma$. Hence $\check\Sigma$ contains only one connected
component. In other words, $\check\Sigma$ is connected. Hence the
following formula holds as soon as $\check\Sigma$ is nonempty:
\smallskip

\centrer{$g-\widehat g-\check g=|\U|-1$.}
\smallskip

\noindent This number must be equal to zero or one. Hence
$|\U|\in\{1,\,2\}$ (still under the hypothesis: $\check\Sigma$ is
nonempty). If $|\U|=2$, then $\check g=0$. We have seen that
conversely, if $\check g=0$ whereas $\check\Sigma$ is nonempty, we
have $|\U|\geqslant 2$, so $|\U|=2$. Similarly, if $\check\Sigma$ is
nonempty, $|\U|=1$ if and only if $\check g=1$. These two last cases
correspond to the two graphs drawn in Figure
\ref{fig:GraphesCourbesNormales1}, right-hand side. When
$\check\Sigma$ is empty, we get obviously the two graphs in Figure
\ref{fig:GraphesCourbesNormales1}, left-hand side. \fin
\bigskip

\section{Expression of the mapping classes of $\G_0$}
         \label{sec:end_grande_proof}
         \label{sec:théorèmes}
\bigskip

\TITRE{Hypotheses.}\\
Let $n\geqslant 6$ an even number, let $\Sigma=\Sigma_{g,\,b}$ with $g\leqslant\frac{n}{2}$, and let
$\rho\ :\ \B_n\to\PMod(\Sigma)$ such that:

\begin{tabular}{lll}
\point & $\rho$ is non-cyclic      & by assumption,\\
\point & $\sigma_p(\G_0)=\vide$      & by assumption, inspired by Proposition \ref{prop:pas_of_curve_peripheral},\\
\point & $\sigma_s(A_i)=\{a_i\}$ & according to \ref{prop:description_sigmasX}\\
\point & $a_i$ is not separating & according to Proposition \ref{prop:curves_spéciales_separating}\\
\end{tabular}
\bigskip

We want to show that $\rho$ is a transvection of monodromy homomorphism.
To do so, we first show in Subsection \ref{par:description_complète_sigmaG} that $\rho$
induces on $\widehat\Sigma$ (cf. Definition \ref{defi:U})
a transvection of monodromy homomorphism. Then we show
in Subsection \ref{par:end_of_demo} that $\rho$ is itself a transvection of monodromy homomorphism.
Finally, in Subsections \ref{par:théorème_cas_pair} and \ref{par:théorèmes},
we show Theorem \ref{thm:Theorem_principal} when $n$ is even and when $n$ is odd, respectively.
\bigskip


\subsection{The homomorphism $\widehat \rho$ induced by $\rho$ on $\Mod(\widehat\Sigma)$}
    \label{par:description_complète_sigmaG}
\medskip

We then define the homomorphism
\smallskip

\centrer{$\DEF{\wedge}{\PMod_{\sigma_n(\G_0)}(\Sigma)}
  {\Mod(\widehat\Sigma)}{A}{\widehat A}$.}
\smallskip

We denote by $\widehat\rho$ the composition
$\wedge\rond\rho\;:\; \B_n\to\Mod(\widehat\Sigma)$. This homomorphism is
well-defined for $\rho(\B_n)$ preserves $\sigma_n(\G_0)$ according to
Proposition \ref{prop:stability_des_curves_normales}. We will prove in this subsection
that $\widehat\rho$ is a transvection of monodromy homomorphism.
\smallskip

We denote by $\widehat \G_0$ and $\widehat\X$ the images of $\G_0$ and $\X$ by
$\wedge$, where
$\X=\{A_i,\,i\in\Imp(n)\}=\{A_1,\,A_3,\dots,\,A_{n-1}\}$. In order
to study the homomorphism $\widehat\rho$, we focus on the mapping
classes induced by those belonging to $\widehat\X$, in
$\Mod(\tildeS)$, where $\tildeS$ is the surface
$\widehat\Sigma_{\sigma_s(\widehat\X)}$. Notice that, according to
Lemma \ref{lem:syst_red_included}, we have:
\smallskip

\centrer{$\left\{\begin{array}{l}
\sigma_s(\widehat\X)=\sigma_s(\X)\cap\Courb(\widehat\Sigma)=\{a_i,\,i\in\Imp(n)\},\\
\sigma_n(\widehat\X)=\sigma_n(\X)\cap\Courb(\widehat\Sigma)=\vide.
\end{array}\right.$}
\smallskip

\medskip

\TITRE{Cutting $\widehat\Sigma$ along the curves of
$\sigma_s(\widehat\X)$}
\smallskip

\begin{itemize}

\item[\point]
  Let $\widetilde\Sigma$ be the surface $\widehat\Sigma_{\sigma_s(\widehat\X)}$.
  According to Proposition \ref{prop:g_minoré}.(v), $\widetilde\Sigma$
  is a connected genus-0 surface, or contains two connected components
  such that one is of genus zero and the other is of genus at most one.
\smallskip

\item[\point]
  Let $\sim$ be the canonical homomorphism $\DEF{\sim}{\Mod_{\sigma_s(\widehat\X)}(\widehat\Sigma)}
  {\Mod(\tildeS)}{\widehat A}{\tildeA}$ , where
  $\Mod_{\sigma_s(\widehat\X)}(\widehat\Sigma)$ is the group of the mapping
  classes of $\Mod(\widehat\Sigma)$ preserving
  $\sigma_s(\widehat\X)$. We will denote by $\tildeX$ the image of
  $\widehat\X$.
\smallskip

\item[\point]
  For all $i\in\Imp(n)$, let us denote by $a_i^+$ and $a_i^-$ the two boundary components
  coming from cutting $\widehat\Sigma$ along
  $a_i$. We set $\Bord^+(\tildeS)=\{a_i^+,\,i\in\Imp(n)\}$ and
  $\Bord^-(\tildeS)=\{a_i^-,\,i\in\Imp(n)\}$. When $\tildeS$ is not
  connected, the boundary components $a_i^+$ and
  $a_i^-$ are such that $\Bord^+(\tildeS)$ is included in the boundary of one connected components of
  $\tildeS$ and $\Bord^-(\tildeS)$ is included in the boundary of the other connected components of
  $\tildeS$. For these definitions, see Figure
  \ref{fig:decoupageUltime} that represents the case where
  $\tildeS$ is connected.
\begin{figure}[!h]
 \Includegraphics{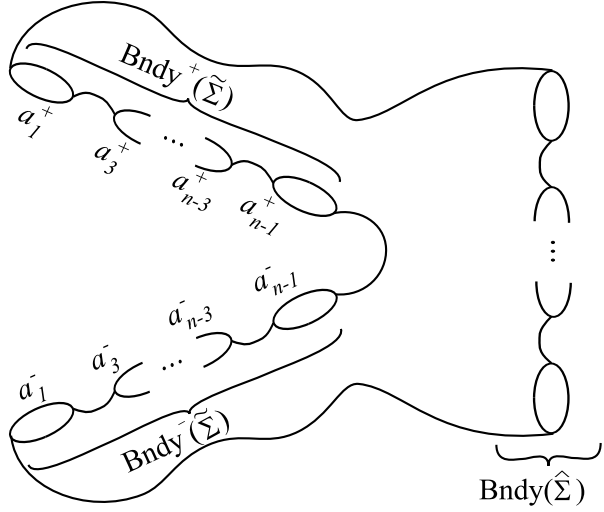}
 \caption{The surface $\widetilde \Sigma$ (case where $\widetilde\Sigma$ is connected).}
 \label{fig:decoupageUltime}
\end{figure}
\smallskip


\item[\point]
Recall that the subgroup $\H$ of $\B_n$ (see Subsection
\ref{par:H_and_X}) is the subgroup
$\langle\,\gamma_i\,,\ i\in\Imp(n)\,\rangle$ of $\B_n$ where for all
$i\in\Imp(n)$, the element $\gamma_i$ is the product
$\tau_i\tau_{i+1}\tau_i\tau_{i+2}\tau_{i+1}\tau_i$.

\item[\point]
Let us also recall the main properties of $\H$ (see Proposition
\ref{prop:properties_of_H}):
\smallskip

\begin{itemize}
\item[\;\;(i)] The action of $\H$ on $\PMod(\Sigma)$ via $\rho$ preserves
$\X$. Indeed, for all $i,j\in\Imp(n)$, we have:
\smallskip

\centrer{$\gamma_i.A_j=\rho(\gamma_i)A_j\rho(\gamma_i)^{-1}=\left\{\begin{array}{lll}
  A_{j} &\mbox{ if }& i\not\in\{j,\,j-2\}\;\\
  A_{j-2} &\mbox{ if }& i=j-2\;\\
  A_{j+2} &\mbox{ if }& i=j.
\end{array}\right.$}
\smallskip
\item[\;(ii)] The homomorphism $\H\to\Ss(\X)$ describing the action of $\H$ on $\X$, where
$\Ss(\X)$ is the symmetric group on the elements of $\X$, is
surjective. In particular, this action is $\frac{n}{2}$ times
transitive.
\item[(iii)] The action of $\H$ on $\Courb(\Sigma)$ preserves $\sigma_s(\X)$.
\end{itemize}
\smallskip

\item[\point]
  For all $\xi\in\H$, the mapping class $\rho(\xi)$ preserves
  $\sigma_s(\X)$, so the element $\sim\rond\wedge(\rho(\xi))$ is well-defined.
  Thus, we can define an action of $\H$ on $\tildeX$ as
  follows. For all $\xi\in\H$ and all $A\in\X$, we set:
  \smallskip

  \centrer{$\xi.\tildeA=\sim\rond\wedge\big(\rho(\xi)\,A\,\rho(\xi)^{-1}\big)$.}
  \smallskip

  \noindent Notice that the action of $\H$ on $\tildeX$ can be deduced from the action of
  $\H$ on $\X$. Then, for all $i,j\in\Imp(n)$, we have:
  \smallskip

\centrer{$\gamma_i.\tildeA_j=\left\{\begin{array}{lll}
  \tildeA_{j} &\mbox{ if }& i\not\in\{j,\,j-2\}\;\\
  \tildeA_{j-2} &\mbox{ if }& i=j-2\;\\
  \tildeA_{j+2} &\mbox{ if }& i=j.
\end{array}\right.$} \end{itemize}

\bigskip


\begin{lem}[The mapping classes $\tildeA_i$, $i\in\Imp(n)$]
              \label{lem:widehat_Ai}
\mbox{}\\All the mapping classes $\tildeA_i$, $i\in\Imp(n)$,
coincide. Either they are trivial, or there are periodic of order 2
and swap $a_j^+$ and $a_j^-$ for all $j\in\Imp(n)$.
\end{lem}
\bigskip

\DEM
\smallskip

\TITRE{Action of $\H$ on $\Bord^+(\tildeS)\sqcup\Bord^-(\tildeS)$:}
\smallskip

\noindent Notice that each mapping class of $\widehat \X$ preserves
$\sigma_s(\widehat A_i)$ for all $i\in\Imp(n)$, so each mapping
class of $\tildeX$ preserves $\{a_i^+,\,a_i^-\}$ for all
$i\in\Imp(n)$. If $\tildeA_1$ fixes both boundary components $a_3^+$
and $a_3^-$, then $\tildeA_1$ fixes both boundary components $a_i^+$
and $a_i^-$ for all $i\in\Imp(n)\smallsetminus\{1\}$, as we are
going to show it right now. Let us recall that according to
Proposition \ref{prop:properties_of_H}, $\H$ acts $\frac{n}{2}$
times transitively on $\X$, hence on $\tildeX$. So, for all
$i\in\Imp(n)\smallsetminus\{1\}$, there exists $\xi\in\H$ such that
$\xi.\tildeA_1=\tildeA_1$ and $\xi.\tildeA_3=\tildeA_i$. Then if
$\tildeA_1$ fixes the boundary components $a_3^+$ and $a_3^-$,
$\xi.\tildeA_1$ fixes the boundary components $\xi.a_3^+$ and
$\xi.a_3^-$ ($\H$ acts canonically via $\rho$ on
$\Bord(\widetilde\Sigma)$), so $\tildeA_1$ fixes both boundary
components $a_i^+$ and $a_i^-$. So there exist \emph{a priori} four
possible actions of $\tildeX$ on
$\Bord^+(\tildeS)\sqcup\Bord^-(\tildeS)$:\smallskip

\begin{itemize}
      \item[a)] the mapping classes of $\tildeX$ fixes the curves of
  $\Bord^+(\tildeS)\sqcup\Bord^-(\tildeS)$;
      \item[b)] the mapping classes of $\tildeX$ swap $a_i^+$ and $a_i^-$
  for all $i\in\Imp(n)$;
      \item[c)] for all $i\in\Imp(n)$, $\tildeA_i$ swaps $a_i^+$ and
  $a_i^-$, and fixes $a_j^+$ and $a_j^-$ for all
  $j\in\Imp(n)\smallsetminus\{i\}$;
      \item[d)] for all $i\in\Imp(n)$, $\tildeA_i$ fixes $a_i^+$ and
  $a_i^-$, and swaps $a_j^+$ and $a_j^-$ for all
  $j\in\Imp(n)\smallsetminus\{i\}$.
\end{itemize}
\medskip

\noindent Both cases a) and b) correspond to the situations
described in the statement of Lemma \ref{lem:widehat_Ai}. In the
remainder of this proof, we assume to be in the case c) or d) and we
expect to find some contradiction. Let us show first the below
assertion (1), stating that $\tildeS$ must be connected. Notice that
in the cases c) and d), $\tildeA_1$ sends some boundary components
of $\Bord^+(\tildeS)$ in $\Bord^+(\tildeS)$ and some other boundary
components of $\Bord^+(\tildeS)$ in $\Bord^-(\tildeS)$. Yet, on one
hand the boundary components of $\Bord^+(\tildeS)$ belong to a same
connected component of $\tildeS$ and so do the boundary components
of $\Bord^-(\tildeS)$. On the other hand, $\tildeA_1$ sends the
connected components of $\tildeS$ on the connected components of
$\tildeS$. Hence there exists a connected component in $\tildeS$
that contains some boundary components of $\Bord^+(\tildeS)$ and
some other boundary components of $\Bord^-(\tildeS)$. By a
transitivity argument, the boundary components of $\Bord^+(\tildeS)$
and of $\Bord^-(\tildeS)$ belong to a same connected component,
hence:
\smallskip

\centrer[1]{the surface $\tildeS$ is connected.}
\medskip

We pursue the proof of Lemma \ref{lem:widehat_Ai} by studying the
nature of the mapping classes of $\tildeX$, (reducible, periodic or
pseudo-Anosov). We will show that they are periodic and we will
determine their order.
\medskip

\TITRE{The nature of the mapping classes of $\tildeX$:}
\smallskip

\noindent The mapping class $\tildeA_1$ is pseudo-Anosov or
periodic, for $\sigma(\tildeA_1)=\vide$. Remember that $\tildeS$ is
the cut of $\widehat\Sigma$ along the curves $\sigma_s(\widehat X)$.
So, if $\tildeA_1$ was pseudo-Anosov in $\Mod(\tildeS)$, the mapping
class $\widehat A_1$ in $\Mod(\widehat\Sigma)$ would satisfy
$\sigma(\widehat A_1)=\sigma_s(\widehat X)$, according to Lemma
\ref{lem:curves_red_ess}. However $\sigma_s(\widehat A_1)=\{a_1\}$,
whence a contradiction. So $\tildeA_1$ is periodic. Now, notice that
$\tildeA_1^{\,2}$ fixes each curve of
$\Bord^+(\tildeS)\sqcup\Bord^-(\tildeS)$, so $\tildeA_1^{\,2}$ is a
periodic  mapping class that fixes more than three boundary
components in a genus-0 surface. Hence, according to Corollary
\ref{cor:periodic_we_a_sphere}, $\tildeA_1^{\,2}$ is the identity.
Hence $\tildeA_1$ is the identity or is periodic of order two. The
same argument can be held for $\tildeA_i$ for all $i\in\Imp(n)$.
Hence:
\smallskip

\centrer[2]{
\parbox{9cm}{either the mapping classes of $\tildeX$ all are the
identity,\\
or they all are periodic of order two.}}
\medskip

We are now ready to focus on the case c) and then on the case d), in
order to find some contradiction.
\medskip

\TITRE{Refutation of the case c):}
\smallskip

\noindent In the case c), $\tildeA_1$ fixes at least $n-2$ boundary
components of $\tildeS$. But as we saw it in (1), $\tildeS$ is a
connected genus-0 surface. Furthermore, $\tildeA_1$ is periodic
according to (2), so we can apply Corollary
\ref{cor:periodic_we_a_sphere} and conclude that $\tildeA_1$ is the
identity. But then $\tildeA_1$ must fix $a_1^+$ and $a_1^-$, which
contradicts the hypotheses of the case c).
\medskip

\TITRE{Study and refutation of the case d):}
\smallskip

\noindent A following simple argument allows us to reject the case
d) when $n\geqslant 8$. Let us consider the mapping class
$\tildeA_1\tildeA_3$. It is periodic of order 2, for $\tildeA_1$ and
$\tildeA_3$ commute. According to Corollary
\ref{cor:periodic_we_a_sphere}, a periodic mapping class that fixes
three or more boundary components of a genus-0 surface is the
identity. Here, $\tildeA_1\tildeA_3$ fixes all the special boundary
components except the boundary components $a_1^+$, $a_1^-$, $a_3^+$
and $a_3^-$. So, when $n\geqslant 8$, the mapping class
$\tildeA_1\tildeA_3$ fixes at least four boundary components. Then
$\tildeA_1\tildeA_3$ must be the identity. But this is absurd, for
$\tildeA_1\tildeA_3$ does not fix the boundary components $a_1^+$,
$a_1^-$, $a_3^+$ and $a_3^-$. Hence the case d) can \emph{a priori}
happen only if $n=6$.
\smallskip

Let us consider the case $n=6$ and let us describe the situation.
Recall that according to (1), $\tildeS$ is connected. Let
$\Bord^0(\tildeS)$ be the set of boundary components of $\tildeS$
that do not belong to $\Bord^+(\tildeS)\sqcup\Bord^-(\tildeS)$.
Again, the mapping class $\tildeA_1\tildeA_3$ is periodic and fixes
the boundary components $a_5^+$, $a_5^-$, whereas it permutes the
boundary components $a_1^+$, $a_1^-$, $a_3^+$ and $a_3^-$
non-trivially. According to Corollary
\ref{cor:periodic_we_a_sphere}, such a mapping class fixes at most
two boundary components. Hence $\tildeA_1\tildeA_3$ fixes no
boundary component of $\Bord^0(\tildeS)$. Yet, on one hand,
$\Bord^0(\tildeS)$ is included in $\Bord(\Sigma)\cup\U$, and on the
other hand the mapping classes of $\tildeX$ fixes the curves of
$\Bord(\Sigma)\cap\Bord^0(\tildeS)$ and have all the same action on
the curves of $\U\cap\Bord^0(\tildeS)$. So
$\Bord(\Sigma)\cap\Bord^0(\tildeS)$ must be empty. Concerning $\U$,
notice that in each of the three cases $|\U|\in\{0,\,1,\,2\}$ which
are authorized by Proposition \ref{prop:graph_curves_normales_1},
the mapping class $\tildeA_1\tildeA_3$ fixes the curves of $\U$.
Thus $\U\cap\Bord^0(\tildeS)$ and
$\Bord(\Sigma)\cap\Bord^0(\tildeS)$ are finally empty sets. Hence
$\U$, $\Bord(\Sigma)$ and $\Bord^0(\tildeS)$ are all empty sets.
Hence:
\smallskip

\centrer{$\widehat\Sigma=\Sigma=\Sigma_{3,\,0}$, and
$\tildeS=\Sigma_{0,\,6}$.}
\smallskip

\begin{figure}[!h]
 \Includegraphics{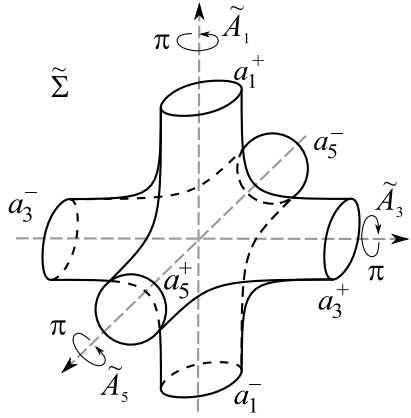}
 \caption{The case $n=6$ where $\widetilde A_1$ fixes $a_1^+$ and $a_1^-$, and permute $a_3^+$ with $a_3^-$, and $a_5^+$ with $a_5^-$.
 We see $\widetilde\Sigma$ in $\RR^3$ and the mapping classes $\widetilde A_1$, $\widetilde A_3$ and $\widetilde A_5$ are isotopy classes of rotations.}
 \label{fig:casN6B0}
\end{figure}

\noindent We have represented the surface $\tildeS$ in Figure
\ref{fig:casN6B0} and we have described in it an example of an
\emph{a priori} possible set of the three periodic mapping classes
$\widetilde A_1$, $\widetilde A_3$, and $\widetilde A_5$, such that
their action on $\Bord(\tildeS)$ is coherent with the case d).

\medskip

Now, since $\U$ is empty (and hence since $\sigma_n(\G_0)$ is empty),
the homomorphism $\wedge$ is trivial so we can forget it. The
contradiction we aim will comes from a study on the mapping class
$Z=(A_3A_4A_5)^2$: we will see that it induces on a subsurface of
$\Sigma$ a periodic mapping class of order 4 and that the square of
$Z$ coincides with $A_1$ modulo a power of Dehn twist along the
curve $a_1$. But as we will see it, such a mapping class has no
square root. This is absurd since $Z$ is a square itself.
\medskip

\noindent 1. \emph{First attempt to describe $\sigma(Z)$.}
\smallskip

Let us recall that $\sigma(Z)=\sigma(Z^2)$. Since $A_1$ commutes
with $A_1$, $A_3$, $A_4$ and $A_5$, the mapping classes $Z$ and
$A_1$ commute. Hence $Z$ fixes the curve $a_1$. Moreover,
$ZA_3Z^{-1}=A_5$ and $ZA_5Z^{-1}=A_3$, so $Z$ swaps the curves $a_3$
and $a_5$. Hence $Z$ induces a mapping class $\tildeZ$ in
$\Mod(\tildeS)$. Moreover, let us justify that $Z^2$ commutes with
$A_3$, $A_4$ and $A_5$. In $\B_4$, the element
$\delta=\tau_1\tau_2\tau_3$ acts cyclically on the set
$\{\tau_1,\,\tau_2,\,\tau_3,\,\tau_0\}$ by conjugation, so
$(\tau_1\tau_2\tau_3)^4$ acts trivially on this set by conjugation.
Now, the fact that $Z^2$ commutes with $A_3$, $A_4$ and $A_5$ comes
from an obvious homomorphism from $\B_4$ to $\langle
A_3,\,A_4,\,A_5\rangle$ that sends $\tau_1$, $\tau_2$ and $\tau_3$
respectively on $A_3$, $A_4$ and $A_5$, and that sends
$(\tau_1\tau_2\tau_3)^4$ on $Z^2$. Hence:
\smallskip

\centrer[3]{$Z^2$ commutes with $A_1$, $A_3$, $A_4$, $A_5$.}
\smallskip

\noindent So,
\smallskip

\centrer[4]{$Z^2$ fixes the curves $a_1$, $a_3$, $a_4$ and $a_5$.}
\smallskip

\noindent Hence the curves $a_1$, $a_3$, $a_4$, $a_5$ are some
reduction curves of $Z^2$, so we have \linebreak
$I\big(\,{\sigma(Z^2)}\,,\,\{a_1,\,a_3,\,a_4,\,a_5\}\,\big)=0$.
Hence, if we see $\tildeS$ as a subsurface of $\Sigma$, we have:
\smallskip

\centrer[5]{$\sigma(Z^2)\subset\Courb(\tildeS)\cup\{a_1,\,a_3,\,a_5\}$.}
\smallskip

\noindent But $a_4$, which is a reduction curve of $Z^2$, intersects
$a_3$ and $a_5$, so neither $a_3$ nor $a_5$ can belong to
$\sigma(Z^2)$. So (5) can be replaced by (6):
\smallskip

\centrer[6]{$\sigma(Z^2)\subset\Courb(\tildeS)\cup\{a_1\}$.}
\smallskip

\noindent We are interested in $\sigma(Z)$ which is equal to
$\sigma(Z^2)$. According to (6), we should investigate the set
$\sigma(Z)\cap\Courb(\tildeS)$. To do so, we focus on $\tildeZ$,
which we define as being $\sim(Z)$, the induced mapping class  by
$Z$ in $\Mod(\tildeS)$. Notice that $\sim(Z)$ is well-defined, since
$Z$ preserves the set $\{a_1,\,a_3,\,a_5\}$. According to Lemma
\ref{lem:syst_red_included}, we have:
\smallskip

\centrer[7]{$\sigma(\tildeZ)=\sigma(Z)\cap\Courb(\tildeS)$.}
\smallskip

\begin{figure}[!h]
 \Includegraphics{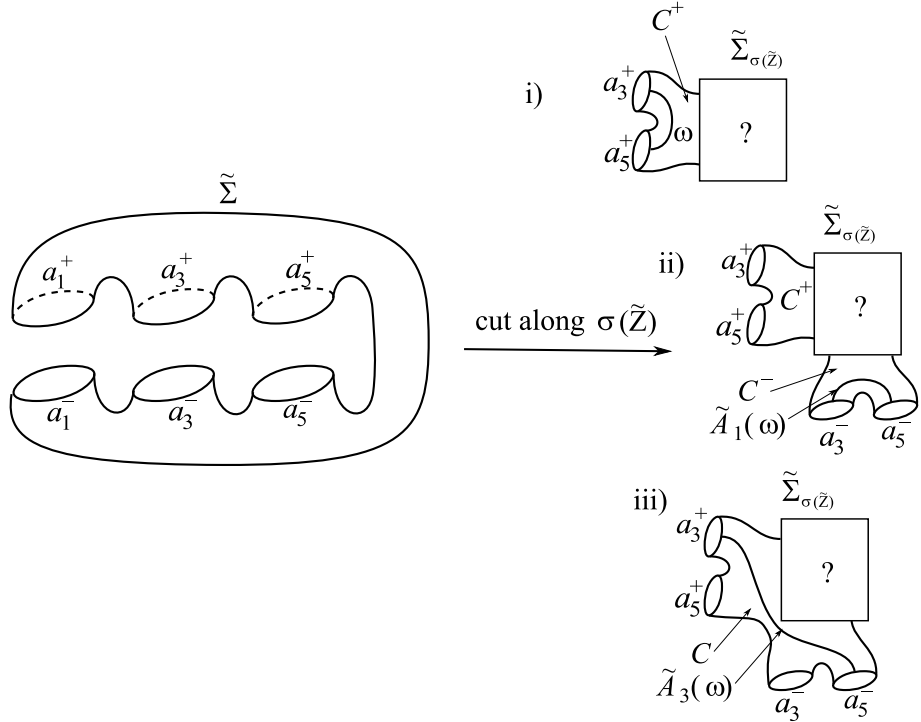}
 \caption{Following the argument of step 1. of the proof of Lemma \ref{lem:widehat_Ai}.}
 \label{fig:sigmaTildeZ}
\end{figure}

\noindent Let us then describe the surface
$\tildeS_{\sigma(\tildez)}$. In the remainder of step 1., we refer
to Figure \ref{fig:sigmaTildeZ}.
\smallskip

\begin{itemize}
\item[\;\;(i)]
Since $a_4$ intersects $a_3$ and $a_5$ but does not intersect the
curves of $\sigma(Z)$, there exists in $\tildeS_{\sigma(\tildez)}$ a
path $\omega$ included in $a_4$ such that $\omega$ has one of its
extremities in $a_3^+\cup a_3^-$ and the other extremity in
$a_5^+\cup a_5^-$. Even if it means renaming the curves, we can
assume that this path is with extremities in $a_3^+$ and $a_5^+$. So
we can assume without loss of generality that $a_3^+$ and $a_5^+$
belong to a same connected component of $\tildeS_{\sigma(\tildez)}$.
Let us denote this connected component by $C^+$.
\smallskip

\item[\;(ii)]
Since $A_1$ commutes with $Z^2$ according to (3), it preserves the
curves of $\sigma(Z)$ and induces an action on the boundary
components of $\tildeS_{\sigma(\tildez)}$. According to the
hypotheses of the case d), $\tildeA_1$ sends $a_3^+$ and $a_5^+$
respectively on $a_3^-$ and $a_5^-$. So $\tildeA_1$ sends $\omega$
on a path joining $a_3^-$ and $a_5^-$. Therefore $a_3^-$ and $a_5^-$
belong to a same connected component of $\tildeS_{\sigma(\tildez)}$.
Let us denote this connected component by $C^-$.
\smallskip

\item[(iii)]
Since $A_3$ commutes with $Z^2$ according to (3), it preserves the
curves of $\sigma(Z)$ and induces an action on the boundary
components of $\tildeS_{\sigma(\tildez)}$. According to the
hypotheses of the case d), $\tildeA_3$ sends $a_3^+$ and $a_5^+$
respectively on $a_3^+$ and $a_5^-$. So $\tildeA_3$ sends $\omega$
on a path joining $a_3^+$ and $a_5^-$. Therefore $a_3^+$ and $a_5^-$
belong to a same connected component of $\tildeS_{\sigma(\tildez)}$.
So $C^+=C^-$. Let us denote by $C$ this connected component.
\end{itemize}
\bigskip

\noindent 2. \emph{We show that the set $\sigma(Z)$ is included in
$\{a_1\}$.}
\smallskip

We argue by contradiction. Let us assume that there exists a curve
$x$ in $\sigma(Z)$ different from $a_1$. According to (6), the curve
$x$ lies in $\sigma(\tildeZ)$. It is a separating curve of $\tildeS$
for $\tildeS$ is a genus-0 surface. Let $x^+$ and $x^-$ be the two
boundary components of $\tildeS_{\sigma(\tildez)}$ coming from the
cut along $x$. We can assume that $x^-$ is a boundary component of
$C$ whereas $x^+$ belongs to another connected component of
$\tildeS_{\sigma(\tildez)}$ that we call $P$. Notice that
$\chi(\tildeS_{\sigma(\tildez)})=\chi(\tildeS)=-4$, whereas
$\chi(C)\leqslant -3$ since it has at least five boundary
components: $a_3^+$, $a_3^-$, $a_5^+$, $a_5^-$ and $x^-$, and
$\chi(P)\leqslant -1$. Since $-4=-3-1$, the connected component $C$
is a sphere with five holes, $P$ is a sphere with three holes, and
$C$ and $P$ are the only connected components of
$\tildeS_{\sigma(\tildez)}$. So
\smallskip

\centrer[8]{$x$ is the only curve in $\sigma(\tildeZ)$.}
\smallskip

\begin{figure}[!h]
 \Includegraphics{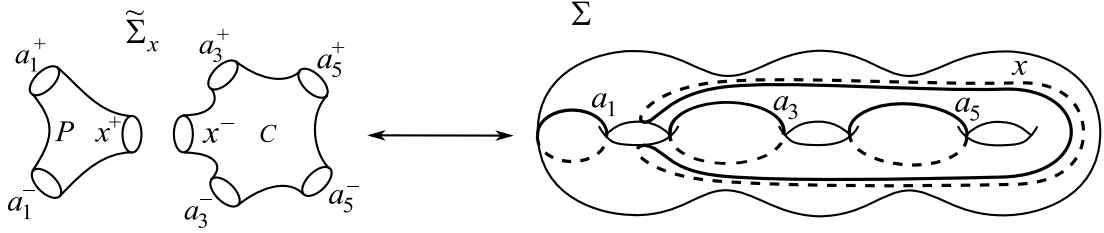}
 \caption{The curve $x$ such that $\sigma(\widetilde Z)=\{x\}$.}
 \label{fig:laCourbex}
\end{figure}

\noindent Notice that $\Bord(P)$ is equal to $\{x^+,
a_1^+,\,a_1^-\}$. This situation is summed up in Figure
\ref{fig:laCourbex}.
\smallskip

Now, $\tildeA_3$ commutes with $\tildeZ^2$, so $\tildeA_3$ fixes the
curve $x$. Then $\tildeA_3$ induces a mapping class in
$\Mod(\tildeS_x)$. But $\tildeS_x$ contains two non-homeomorphic
connected components , so $\tildeA_3$ induces a mapping class on
each connected component. Let $F$ be the so induced mapping class by
$\tildeA_3$ on $C$. It is periodic of order two on the surface $C$
(for $\tildeA_3$ was periodic of order two, according to (2)).
However, $F$ fixes $x^-$, $a_3^+$ and $a_3^-$. This is in
contradiction with Corollary \ref{cor:periodic_we_a_sphere}, since
$C$ is a genus-0 surface, and any periodic mapping class on a
genus-0 surface that fixes three boundary components is the
identity. This is the expected contradiction. We have then shown
that $\sigma(Z)\subset\{a_1\}$.
\bigskip

\noindent 3. \emph{Let us show that $Z^2$ and $A_1$ coincide, up to
a power of a Dehn twist along the curve $a_1$}.
\smallskip

Let us recall that, according to (2), $\tildeA_1$ is periodic of
order two. The mapping class $\tildeZ$ satisfies
$\sigma(\tildeZ)=\vide$, so it is either pseudo-Anosov, or periodic.
If it was pseudo-Anosov, the curves $a_1$, $a_3$ and $a_5$ would be
essential reduction curves of $Z$, but we saw in step 1. that $a_3$
and $a_5$ were not. Hence $\tildeZ$ is periodic. Then the mapping
classes $\tildeZ$ and $\tildeA_1$ are periodic and commute, so they
span a finite group. Notice that according to the hypotheses of the
case d), $\tildeA_1$ fixes the curves $a_1^+$ and $a_1^-$.
Concerning $\tildeZ$, remember that $A_3$, $A_4$ and $A_5$ fix the
curve $a_1$, hence so does the mapping class $A_3A_4A_5$. So the
mapping class $Z=(A_3A_4A_5)^2$ fixes the curve $a_1$ and does not
permute the two connected components of $\V\smallsetminus a_1$ where
$\V$ is a tubular neighbourhood of $a_1$, so $\tildeZ$ fixes the
curves $a_1^+$ and $a_1^-$.  Hence $\tildeZ$ and $\tildeA_1$ are
periodic mapping classes and both fix the curves $a_1^+$ and
$a_1^-$, so according to Lemma \ref{lem:group_cyclic}, $\tildeZ$ and
$\tildeA_1$ span a cyclic group. Since the order of $\tildeZ$ is
even (for $\tildeZ$ swaps $\{a_3^+,\,a_3^-\}$ with
$\{a_5^+,\,a_5^-\}$) and since the order of $\tildeA_1$ is two, then
there exists an integer $k\geqslant 1$ such that
$\tildeA_1=\tildeZ^k$. But $\tildeZ^4$ has a trivial action on the
boundary components of $\tildeS$, so it is the identity, hence
$k\in\{1,\,2,\,3\}$. But $\tildeZ$ and $\tildeZ^{3}$ swap
$\{a_3^+,\,a_3^-\}$ and $\{a_5^+,\,a_5^-\}$ whereas $\tildeA_1$
preserves each of these sets, so $k\not\in\{1,\,3\}$. Finally,
\smallskip

\centrer{$\tildeA_1=\tildeZ^2$.}
\smallskip

\noindent Let us consider the product $Z^2A_1^{-1}$. We just have
seen that $\sim\!(Z^2A_1^{-1})$ is trivial in $\Mod(\tildeS)$, so
according to the following exact sequence:
\smallskip

\centrer{$1\to\langle
T_{a_1},\,T_{a_3},\,T_{a_5}\rangle\to\Mod_{\sigma_s(\X)}(\Sigma)\xrightarrow{\sim}\Mod(\tildeS)\to1$,}
\smallskip

\noindent the mapping class $Z^2A_1^{-1}$ is a multitwist along
the curves $a_1$, $a_3$ and $a_5$. But $\sigma(A_1)$ is equal
to $\{a_1\}$ and according to step 2., $\sigma(Z^2)$ is
included in $\{a_1\}$, so according to Proposition
\ref{prop:properties_sigma}.(v), $\sigma(Z^2A_1^{-1})$ is
included in $\{a_1\}$.
\medskip

3. We get a contradiction when we examine on the mapping class
$Y=A_3A_4A_5$.
\smallskip

\begin{figure}[!h]
 \Includegraphics{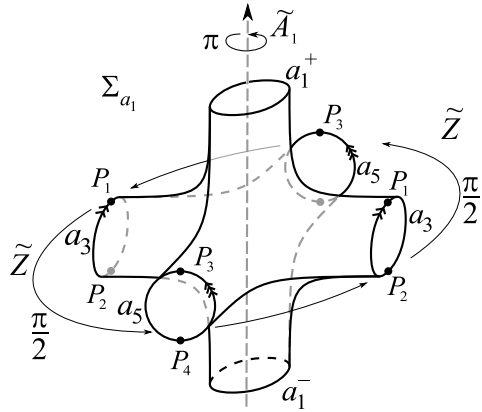}
 \caption{The periodic mapping class $\widetilde A_1$ induced by $A_1$ in $\Sigma_{a_1}$
 (the arrows along the curves $a_3$ and $a_5$ stand for a gluing along these curves, so that $\Sigma_{a_1}$
 is indeed a genus-2 surface with two boundary components: $a_1^+$ and $a_1^-$). }
 \label{fig:quatrePointsFixes}
\end{figure}

We can see $\tildeS$ as a punctured sphere where each boundary
component has been replaced by a puncture. The mapping class
$\tildeZ$ of $\Mod(\tildeS)$ is periodic so according to Kerckhoff's
Theorem, $\tildeZ$ can be represented by a periodic diffeomorphism
of $\tildeS$. According to Kerékj\`art\`o's Theorem, this periodic
diffeomorphism is conjugate to a rotation on the sphere $\tildeS$ by
a diffeomorphism isotopic to the identity. Since $\tildeZ$ is of
order 4, $\tildeS$ is the isotopy class of an angle $\pm\pi/2$
rotation over an axis containing two punctures that correspond to
the boundary components $a_1^+$ and $a_1^-$. The square of this
rotation is in the isotopy class of $\tildeA_1$. This justifies
Figure \ref{fig:quatrePointsFixes} in which we have represented the
genus-2 surface $\Sigma_{a_1}$ with two boundary components $a_1^+$
and $a_1^-$, and the periodic mapping classes $A'_1$ and $Z'$
induced by $A_1$ and $Z$ in $\PMod(\Sigma_{a_1})$. We can see that
$A'_1$ has exactly four fixed points, namely $P_1$, $P_2$, $P_3$ and
$P_4$, in Figure \ref{fig:quatrePointsFixes}. Therefore:
\smallskip

\begin{itemize}
\item[(1)] $Z'$ has no fixed points,
\item[(2)] $Z'^2=A'_1$ has four fixed points.
\end{itemize}

\smallskip

\noindent Let us set $Y:=A_3A_4A_5$. According to Proposition
\ref{prop:stability_des_curves_spéciales}, $A_3$, $A_4$ and $A_5$
fix $\sigma_s(A_1)$, hence fix $a_1$. So $Y$ induces a mapping class
in $\Mod(\Sigma_{a_1})$, which we will call $Y'$. Since $Z'$ is
periodic of order 4, the mapping class $Y'$ is periodic of order 8.
Assertions (1) and (2) imply assertions (3)-(5) below. Only
assertion (5) needs to be proved.
\smallskip

\begin{itemize}
  \item[(3)] $Y'$ is periodic of order 8 without fixed points,
  \item[(4)] $Y'^2=Z'$ is periodic of order 4 without fixed points,
  \item[(5)] $Y'^4$ has four fixed points belonging to a same orbit under $Y'$.
\end{itemize}
\smallskip

\noindent Let us justify assertion (5). If $Y'$ preserved
$\{P_1,\,P_3\}$, then $Y'^2$ would fix $P_1$ and $P_3$. But it does
not, so the cardinality of the orbit of $P_1$ under $Y'$ is at least
4. But $Y'^2$ preserves $\{P_1,\,P_3\}$, so the cardinality of the
orbit of $P_1$ under $Y'$ is  at most 4. Hence the orbit of $P_1$
and $P_3$ and the orbit of $P_2$ and of $P_4$ contain exactly 4
points. Since $Y'^4$ contains 4 fixed points instead of 8, the orbit
of $P_1$ and the orbit of $P_2$ coincide. Thus assertion (5) is
proved. We are then ready to apply Lemma \ref{lem:periodic_Euler},
linking the Euler characteristics of $\Sigma_{a_1}$ and of the
quotient surface $\Sigma_{a_1}/\langle Y'\rangle$. To do this, let
us compute $\chi(\,\Sigma_{a_1}/\langle Y'\rangle\,)$. The mapping
classes $A'_3$, $A'_4$ and $A'_5$ induced by $A_3$, $A_4$ and $A_5$
in $\Mod(\Sigma_{a_1})$ swap $a_1^+$ and $a_1^-$ (we knew it already
for $A_3$ and $A_5$ by hypothesis, we deduce it easily for $A_4$).
So $Y'$ swaps $a_1^+$ and $a_1^-$, so the surface
$\Sigma_{a_1}/\langle Y'\rangle$ must have a unique boundary
component. Since $Y'$ preserves the orientation,
$\Sigma_{a_1}/\langle Y'\rangle$ is a disk, a torus with one hole,
or a genus-2 surface with one hole. So $\chi(\,\Sigma_{a_1}/\langle
Y'\rangle\,)\in\{1,\,-1,\,-3\}$. The ramification points  of the
covering $\Sigma_{a_1}\to\Sigma_{a_1}/\langle Y'\rangle$ come also
into account. According to (5), there is only one ramification point
$Q$ that has 4 preimages $P_1$, $P_2$, $P_3$ and $P_4$ in
$\Sigma_{a_1}$. So, with the notation of Lemma
\ref{lem:periodic_Euler}, we have $o(Q)=4$. Let us recall that in
Lemma \ref{lem:periodic_Euler}, the cardinality of $\langle
Y'\rangle$ is denoted by $m$, then here, $m=8$. We get:
\smallskip

\centrer{$
\dessous{\underbrace{\chi(\Sigma_{a_1})}}{-4}+\dessous{\underbrace{
(m-o(Q))}}{8-4}=\dessous{\underbrace{m}}{8}.\dessous{\underbrace{\chi(\,\Sigma_{a_1}/\langle
Y'\rangle\,)}}{1,\,-1\mbox{\small \ or } -3}\ .$}
\smallskip

\noindent This equality cannot be satisfied, since the left-hand
side is zero whereas the right-hand side is nonzero. This is the
expected contradiction and terminates the proof of Lemma
\ref{lem:widehat_Ai}.\fin
\bigskip


We now can prove the following proposition:
\medskip

\begin{prop}[Description of $\widehat\rho$]
          \label{prop:description_widehat_rho}
\mbox{}\\There exist an integer $\varepsilon\in\{\pm1\}$ and a
mapping class $\alpha$ of $\Mod(\widehat\Sigma)$ that is either the
identity or the periodic mapping class of order two that fixes the
curves $a_i$, $i\in\{1,\,2,\dots,\,n-1\}$ such that, for all
$i\in\{0,\,1,\dots,\,n-1\}$, the homomorphism $\widehat\rho$ satisfies:
\smallskip

\centrer{$\widehat\rho(\tau_i)=T_{a_i}^{\,\varepsilon}\,\alpha$.}
\end{prop}
\bigskip

\DEM
\smallskip

1. Let us begin to deal with the case where $\tildeA_1$ is trivial.
In this case, $\widehat A_1$ is a multitwist along the curves of
$\sigma(\widehat \X)$. But $\sigma(\widehat
A_1)=\sigma_s(A_1)\cap\Courb(\widehat\Sigma)=\{a_1\}$, so there
exists a nonzero integer $\varepsilon$ such that $\widehat
A_1=T_{a_1}^\varepsilon$. Then, by applying the conjugation by the
elements of $\widehat\rho(\J)$, for all $i\in\{0,\,1,\dots,\,n-1\}$,
we get:
\smallskip

\centrer[1]{$\widehat \rho(\tau_i)=\widehat
A_i=T_{a_i}^{\,\varepsilon}$.}
\smallskip

\noindent Now, for all $i\in\{0,\,1,\dots,\,n-1\}$, $\widehat A_i$
and $\widehat A_{i+1}$ satisfy a braid relation, so according to
Proposition \ref{prop:properties_twists_Dehn}, $\varepsilon$ belongs
to $\{\pm1\}$. Thus, in the case where $\tildeA_1$ is trivial, the
proof is over.
\bigskip

2. When $\tildeA_1$ is not the identity, we have
$\tildeA_1=\tildeA_3=\dots=\tildeA_{n-1}$ according to Lemma
\ref{lem:widehat_Ai}. According to Corollary
\ref{cor:reconstruction_periodic}, these periodic mapping
classes of order 2 all induce a unique mapping class $\alpha$
of order 2 on $\Mod(\widehat\Sigma)$. According to Proposition
\ref{prop:properties_sigma}.(v), for all $i\in\Imp(n)$, we have
the inclusions $\sigma(\widehat
A_i\,\alpha)\subset\sigma(\widehat A_i)\cup\sigma(\alpha)$ and
$\sigma(\widehat A_i)\subset\sigma(\widehat
A_i\,\alpha)\cup\sigma(\alpha^{-1})$, but
$\sigma(\alpha)=\sigma(\alpha^{-1})=\vide$, so we get the
equality $\sigma(\widehat A_i\,\alpha)=\sigma(\widehat
A_i)=\{a_i\}$. By definition of $\alpha$, the mapping class
$\widehat A_i\,\alpha$ induces a trivial mapping class in
$\Mod(\tildeS)$. Hence, for all $i\in\Imp(n)$, there exists an
integer $k_i$ such that:
\smallskip

\centrer[2]{$\widehat A_i\,\alpha=T_{a_i}^{\,k_i}$.}
\smallskip

\noindent Notice that for all $\xi\in\H$, the mapping class
$\widehat\rho(\xi)\,\alpha\,\widehat\rho(\xi)^{-1}$ is periodic of
order 2 and induces in $\Mod(\tildeS)$ a mapping class that
coincides with $\xi.\tildeA_1$. But $\xi.\tildeA_1=\tildeA_1$
according to Lemma \ref{lem:widehat_Ai}, so by uniqueness of the
construction of $\alpha$ (cf. Corollary
\ref{cor:reconstruction_periodic}), we have the equality
$\widehat\rho(\xi)\,\alpha\,\widehat\rho(\xi)^{-1}=\alpha$. Then, by
conjugation by the elements of $\widehat\rho(\H)$, the $k_i$ are all
equal to an integer which we denote by $\varepsilon$. Hence for all
$i\in\Imp(n)$, we have:
\smallskip

\centrer[3]{$\widehat A_i\,\alpha=T_{a_i}^{\,\varepsilon}$.}
\smallskip

\noindent Since $\widehat A_3\widehat A_4\widehat
A_3(\sigma_s(\widehat A_3))=\sigma_s(\widehat A_4)$, then $\widehat
A_3\widehat A_4\widehat A_3(a_3)=a_4$. Hence the product $\widehat
A_3\widehat A_4\widehat A_3$ sends by conjugation $T_{a_3}$ on
$T_{a_4}$. Moreover, $\widehat A_3\widehat A_4\widehat A_3$ sends
also by conjugation $\widehat A_3$ on $\widehat A_4$. Hence, by
conjugating the equality $\widehat A_1T_{a_1}^{-\varepsilon
}=\widehat A_3T_{a_3}^{-\varepsilon }$  by $\widehat A_3\widehat
A_4\widehat A_3$, we get $\widehat A_1T_{a_1}^{-\varepsilon
}=\widehat A_4T_{a_4}^{-\varepsilon }$. Hence $\widehat
A_4T_{a_4}^{-\varepsilon }=\widehat A_3T_{a_3}^{-\varepsilon }$.
Now, let us make $\delta$ act on this last equality. We get (4):
\smallskip

\centrer[4]{$\widehat A_1T_{a_1}^{-\varepsilon }=\widehat
A_2T_{a_2}^{-\varepsilon }=\dots=\widehat
A_{n-1}T_{a_{n-1}}^{-\varepsilon }=\widehat A_0T_{a_0}^{-\varepsilon
}$.}
\smallskip

\noindent Therefore $\alpha$ (equal to $\widehat
A_1T_{a_1}^{-\varepsilon }$) is stable by the action of $\delta$. In
other words, $\alpha$ commutes with $\widehat\rho(\delta)$. Then,
since $\alpha(a_1)=a_1$, it follows that for all
$i\in\{0,\,1,\dots,\,n-1\}$, we have $\alpha(a_i)=a_i$. Hence, for
all $i\in\{0,\,1,\dots,\,n-1\}$, the mapping class $\alpha$ commutes
with $T_{a_i}$, hence with $T_{a_i}^{\varepsilon}\,\alpha$. That is,
$\alpha$ commutes with $\widehat A_i$. Therefore, the transvection
of $\widehat \rho$ with direction $\alpha$ is well-defined. Let us
denote it by $L_\alpha(\widehat \rho)$. Thus, for all
$i\in\{0,\,1,\dots,\,n-1\}$, we have:
\smallskip

\centrer[5]{$L_\alpha(\widehat \rho)(\tau_i)=\widehat
A_i\,\alpha=T_{a_i}^{\,\varepsilon }$.}
\smallskip

\noindent This kind of equality (5) is very similar to (1). We can
then prove, exactly as in the case where $\tildeA_1$ is trivial,
that $\varepsilon \in\{\pm1\}$.\fin
\bigskip
\bigskip


\subsection{The homomorphism $\rho$ is a transvection of monodromy homomorphism}
        \label{par:end_of_demo}
\medskip

At last, we prove in this subsection Proposition
\ref{prop:phase_finale}: we know that the homomorphism $\widehat\rho$
induced by $\rho$ is a transvection of monodromy homomorphism. We have
now to check that the homomorphism $\rho$ itself is a transvection of
monodromy homomorphism.
\bigskip

Let us gather the main information on $\widehat\Sigma$,
$\check\Sigma$ and $\U$ contained in Definition \ref{defi:U} and
Proposition \ref{prop:graph_curves_normales_1}.
\medskip

\TITRE{Recalls} (Definition \ref{defi:U} and Propositions
\ref{prop:g_minoré} and
\ref{prop:graph_curves_normales_1})\textbf{.}
\begin{itemize}
  \item[\point]
Let $\widehat\Sigma$ be the subsurface of
$\Sub_{\sigma_n(\G_0)}(\Sigma)$ containing the special curves. The
surface $\widehat\Sigma$ is of genus $\frac{n}{2}-1$ or
$\frac{n}{2}$.

  \item[\point]
We set $\U:=\Bord(\widehat\Sigma)\cap\sigma_n(\G_0)$. The set of
curves $\U$ can be empty and contains at most two curves.

  \item[\point]
Let $\check\Sigma$ be the subsurface of $\Sub_{\U}(\Sigma)$
different from $\widehat\Sigma$ (well-defined according to
Proposition \ref{prop:graph_curves_normales_1}). If $\widehat\Sigma$
is the only subsurface of $\Sub_{\U}(\Sigma)$, we will say that
$\check\Sigma$ is empty. \end{itemize}
\smallskip

\noindent The links between $\widehat\Sigma$, $\check\Sigma$ and
$\U$ are the following:
\begin{itemize}
  \item
if $\U$ is empty, $\Sigma=\widehat\Sigma$; if $\widehat\Sigma$ is of
genus $\frac{n}{2}$, then $\U$ is empty,
  \item
if $\U$ is reduced to a non-separating curve $u$, then
$\check\Sigma$ is empty and $\{u\}$ is the mark of $\widehat\Sigma$,
  \item
if $\U$ is reduced to a separating curve $u$, then $\check\Sigma$ is
of genus 1,
  \item
if $\U$ contains two curves, then they are non-separating and
$\check\Sigma$ is a nonempty genus-0 surface. \end{itemize}
\bigskip

\begin{prop}
    \label{prop:direction_of_rho}
There exists a mapping class $W\in\PMod(\Sigma)$ such that for all
\linebreak
$i\in\{1,\,2,\dots,\,n-1\}$, the following holds:
\smallskip

\centrer{
 \begin{tabular}{ll}
   \point & in $\Mod(\widehat\Sigma)$, $\wedge(W)$ commutes with $\wedge(T_{a_i})$ and $\wedge(A_i)$,\\
   \point & in $\Mod(\widehat\Sigma)$, $\wedge(A_i\,W^-1)=\wedge(T_{a_i}^{\,\varepsilon})$,\\
   \point & in $\Mod(\widehat\Sigma)$, $\vee(W)$ commutes with $\vee(A_i)$,\\
   \point & in $\Mod(\check\Sigma)$,  $\vee(A_i\,W^-1)=\vee(T_{a_i}^{\,\varepsilon})$.
 \end{tabular}
}
\end{prop}

\DEM Let us distinguish the case according to $\U$.
\smallskip

a) If $\U$ is empty, then $\widehat\Sigma=\Sigma$, so
$A_1T_{a_1}^{-1}$ coincides with $\alpha$, the mapping class defined
in Proposition \ref{prop:description_widehat_rho}. According to this
last proposition, $\alpha$ satisfies the four assertions that $W$
must satisfy. Then in this case, Proposition
\ref{prop:direction_of_rho} is proved.
\smallskip

b) If $\U$ is reduced to a non-separating curve $u$, then, if
$\alpha$ is the identity, we set $W=\Id$ in the group
$\PMod(\Sigma)$. And if $\alpha$ is not the identity, according to
Corollary \ref{cor:reconstruction_periodic}, there exists a unique
mapping class $W$ of $\PMod(\Sigma)$, periodic of order two, fixing
the curve $u$ and such that $\cut_u(W)=\alpha$. Then, again
according to Proposition \ref{prop:description_widehat_rho}, this
definition of $W$ suits.
\smallskip

c) If $\U$ is separating in $\Sigma$, which gathers all the cases
not treated by a) and b) above, then we are going to show that the
mapping class $W=A_1T_{a_1}^{-1}$ belonging to $\PMod(\Sigma)$
suits. Let us start by showing that the homomorphism $\vee\rond\rho$ is
cyclic.
\smallskip

The set of curves $\U$ is stable by the mapping classes of $\G_0$, for
the curves of $\U$ are topologically different from the other curves
of $\sigma_n(\G_0)$: they are the only ones that separate
$\widehat\Sigma$ from $\check\Sigma$. Let us distinguish two cases,
depending on whether $\check\Sigma$ is of genus 1, or of genus 0.

\begin{itemize}
  \item[\point]
   If $\check\Sigma$ is of genus 1, according to Proposition \ref{prop:graph_curves_normales_1},
   $\U$ is reduced to one only curve, that is hence stable by the action of
   $\rho(\B_n)$. Hence the image of the homomorphism $\vee\rond\rho$
   is included in $\PMod(\check\Sigma)$ (the boundary components are not permuted). Then we can apply
   Proposition \ref{prop:g_minoré}: since $\check\Sigma$ is
   of genus $1$, hence smaller than $\frac{n}{2}-1$, then the homomorphism $\vee\rond\rho$ is
   cyclic.
   \smallskip

 \item[\point]
   If $\check\Sigma$ is of genus zero, then $\U$ can contain two
   curves. They can be swapped by the elements of
   $\G_0$, but according to the proposition
   \ref{prop:stability_des_curves_normales}, the action of $\B_n$ on $\U$ is cyclic.
   Remember that we set in Subsection \ref{par:action_we_simplex_of_curves}:
   \smallskip

 \centrer{$\begin{array}{ccccc}
 \F_n
      & =& \big\langle\!\langle\,\tau_3\tau_{1}^{-1}\,\rangle\!\big\rangle_{\B_n}
      & =& \Ker(\lambda),\\
 \F_n^* &=& \big\langle\,\tau_i\tau_{1}^{-1}\,,\ 3\leqslant i\leqslant n-1\,\big\rangle_{\B_n}
      & \subsetnot& \Ker(\lambda).\vlblanc
 \end{array}$}
\smallskip

   \noindent where $\big\langle\!\langle\ \rangle\!\big\rangle_{\B_n}$ is the normal closure in $\B_n$
   and $\lambda\ :\ \B_n\to\ZZ$ is the unique homomorphism such that $\lambda(\tau_1)=1$.
   Since the action of $\B_n$ on $\U$ is cyclic, the induced action of $\F_n$ on $\U$ is trivial.
   Hence $\rho(\F_n)$ is included in $\PMod(\check\Sigma)$. Now, remember that
   $\F_n^*$ is isomorphic to $\B_{n-2}$, that $\check\Sigma$ is a genus-0 surface,
   so, according to Theorem \ref{thm:BndansMCG0b}, the homomorphism
   $\vee\rond\rho\ :\ \F_n^*\to \PMod(\check\Sigma)$ is cyclic.
   Hence again, the homomorphism $\vee\rond\rho\ :\ \B_n\to \PMod(\check\Sigma)$
   is cyclic after Lemma \ref{lem:cyclic_and_F_n}.(ii).
\end{itemize}
\bigskip

Then, according to Proposition \ref{prop:description_widehat_rho},
the mapping class $W=A_1T_{a_1}^{-\varepsilon}$ belonging to
$\PMod(\Sigma)$ satisfies for all $i\in\{1,\,2,\dots,\,n-1\}$ the
following facts (where $\alpha$ is the mapping class of
$\Mod(\widehat\Sigma)$ introduced in Proposition
\ref{prop:description_widehat_rho}):
\smallskip

\centrer{
 \begin{tabular}{ll}
   \point & in $\Mod(\widehat\Sigma)$, $\wedge(W)$ is equal to $\alpha$ and then commutes with $T_{a_i}$,\\
   \point & in $\Mod(\widehat\Sigma)$, $\wedge(T_{a_i}^{\,\varepsilon}\,W)=\widehat A_i$,\\
   \point & in $\Mod(\check\Sigma)$, $\vee(A_iT_{a_i}^{\,-\varepsilon})=\vee(A_i)=\vee(W)$.
 \end{tabular}
} \smallskip

This terminates the proof of Proposition
\ref{prop:direction_of_rho}.\fin
\bigskip

\begin{prop}
        \label{prop:phase_finale}
Under the hypotheses of this section, the homomorphism
$\rho\;:\; \B_n\to\PMod(\Sigma)$ is a transvection of monodromy homomorphism.
\end{prop}
\medskip

\DEM Let $W$ be the mapping class of $\PMod(\Sigma)$ defined in
Proposition \ref{prop:direction_of_rho}. Notice that according to
this same proposition, $\wedge(A_i\,W^{-1})$ belongs to
$\PMod(\widehat\Sigma)$ for all $i\in\{1,\,\dots,\,n-1\}$, since
$\wedge(A_i\,W^{-1})=\wedge(T_{a_i}^{\varepsilon})$. Similarly, when
$\check\Sigma$ is nonempty, $\vee(A_i\,W^{-1})$ belongs to
$\PMod(\check\Sigma)$ since
$\vee(A_i\,W^{-1})=\vee(T_{a_i}^{\varepsilon})$. Hence $A_iW^{-1}$
and $T_{a_i}^{\varepsilon}$ belong to $\P_\U\Mod(\Sigma)$. Let us
consider then the following central exact sequence:
\smallskip

    \centrer[1]{$1\to\langle
    T_u,\,u\in\U\rangle\to\P_\U\Mod(\Sigma)\xrightarrow{\cut_u}\PMod(\Sigma_\U)\to1$}
    \smallskip

\noindent where $\cut_\U$ is the canonical homomorphism. For any
$i\in\{1,\,\dots,\,n-1\}$, both mapping classes $A_iW^{-1}$ and
$T_{a_i}^{\varepsilon}$ induce the same mapping classes in
$\PMod(\Sigma_\U)$. So $A_iW^{-1}$ and $T_{a_i}^{\varepsilon}$
differ by a central element. Consequently, $W$ commutes with $A_i$.
We can then define a homomorphism $\rho'\;:\; \B_n\to\PMod(\Sigma)$ by
setting:
\smallskip

\centrer{$\rho'(\tau_i)=A_i\,W^{-1}=\rho(\tau_i)\,W^{-1}$}
\smallskip

\noindent for all integers $i\in\{1,\dots,\,n-1\}$. Let $\rho''\;:\;
\B_n\to\PMod(\Sigma)$ be the monodromy homomorphism defined by
\smallskip

\centrer{$\rho'(\tau_i)=T_{a_i}^{\,\varepsilon}$}
\smallskip

\noindent for all $i\in\{1,\,\dots,\,n-1\}$. The homomorphisms $\rho'$
and $\rho''$ satisfy $\cut_\U(\rho')=\cut_\U(\rho'')$ according to
Proposition \ref{prop:direction_of_rho}. Let us apply Lemma
\ref{lem:transvections} to the central exact sequence (1). It
follows that $\rho'$ and $\rho''$ are of the same nature, hence
$\rho'$ is a transvection of monodromy homomorphism. So by construction,
$\rho$ is also a transvection of monodromy homomorphism.
\bigskip

\subsection{Proof of Theorem \ref{thm:Theorem_principal} when $n$ is even}
            \label{par:théorème_cas_pair}
\medskip

Let $n$ be an even integer greater than or equal to 6, and
$\Sigma$ a surface $\Sigma_{g,\,b}$ where $g\leqslant
\frac{n}{2}$. We already know (cf. Lemma
\ref{lem:existence_des_homomorphisms_of_monodromie}) that monodromy
homomorphisms exist if and only if $g\geqslant \frac{n}{2}-1$. Let
$\rho$ be a noncyclic homomorphism from $\B_n$ to $\PMod(\Sigma)$.
Let $\Sigma'$ be the connected component of
$\Sigma_{\sigma_p(\G_0)}$ of genus $g$ (recall that
$\sigma_p(\G_0)$ is the set of peripheral curves which where
introduced in Section \ref{sec:curves_peripheral}). According
to Proposition \ref{prop:pas_of_curve_peripheral}.(i) and (ii),
the homomorphism $\rho$ induces a noncyclic homomorphism $\rho'\;:\;
\B_n\to\PMod(\Sigma')$. Then according to Proposition
\ref{prop:phase_finale}, $\zeta'$ is a transvection of
monodromy homomorphism. Then, according to Proposition
\ref{prop:pas_of_curve_peripheral}.(ii), the homomorphism $\rho$
was a transvection of monodromy homomorphism. Hence Theorem
\ref{thm:Theorem_principal} when $n$ is even is proved.
\bigskip


\subsection{Proof of Theorem \ref{thm:Theorem_principal} when $n$ is odd}
\label{par:théorèmes}
\medskip

Let $n$ be an odd integer greater than or equal to 7.
Notice that in this case, the
condition $g\leqslant \frac{n}{2}$ is equivalent to $g\leqslant
\frac{n-1}{2}$. Let $\Sigma$
be a surface $\Sigma_{g,\,b}$ where $g\leqslant \frac{n-1}{2}$ and
$b\geqslant 0$. We already
know (cf. Lemma \ref{lem:existence_des_homomorphisms_of_monodromie}) that monodromy homomorphisms
exist if and only if $g\geqslant \frac{n}{2}-1$.
Let $\rho$ be a noncyclic homomorphism from $\B_n$ to $\PMod(\Sigma)$. We are going to show that $\rho$ is a transvection of monodromy homomorphism.
\medskip

We adopt the following notation:
\smallskip

\centrer{$\B_{n-1}^{(1)}=\big\langle\,\tau_1,\,\tau_2,\dots,\,\tau_{n-2}\big\rangle_{\B_n}$,}

\centrer{$\B_{n-1}^{(2)}=\big\langle\,\tau_2,\,\tau_3,\dots,\,\tau_{n-1}\big\rangle_{\B_n}$.}
\smallskip

The homomorphism $\rho$ from $\B_n$ to $\PMod(\Sigma)$ induces by
restriction to $\B_{n-1}^{(1)}$ and $\B_{n-1}^{(2)}$ the homomorphisms
$\rho^{(1)}\,:\,\B_{n-1}^{(1)}\to\PMod(\Sigma)$ and
$\rho^{(2)}\,:\,\B_{n-1}^{(2)}\to\PMod(\Sigma)$. The homomorphism $\rho$
is not cyclic, so the mapping classes $\rho(\tau_2)$ and
$\rho(\tau_3)$ are distinct according to Lemma
\ref{lem:homomorphism_cyclic}. So the homomorphisms $\rho^{(1)}$ and
$\rho^{(2)}$ are not cyclic either. Then, according to Theorem
\ref{thm:Theorem_principal} when $n$ is even, $\rho^{(1)}$ and $\rho^{(2)}$ are
transvections of monodromy homomorphisms. So there exist two $(n-2)$-chains:
$(a_i,\ 1\leqslant i\leqslant n-2)$ and $(c_i,\ 2\leqslant
i\leqslant n-1)$; two mapping classes: $V$ belonging to the
centralizer of $\big\langle T_{a_i},\,1\leqslant i\leqslant
n-2\big\rangle$ in $\PMod(\Sigma)$ and $W$ belonging to the
centralizer of $\big\langle T_{c_i},\,2\leqslant i\leqslant
n-1\big\rangle$ in $\PMod(\Sigma)$; and two integers $\varepsilon$
and $\eta$ belonging to $\{\pm1\}$, such that for all
$i\in\{1\,\dots,\,n-2\}$:
\smallskip

\noindent \centrer{
    $\left\{\begin{array}{l}
      \rho^{(1)}(\tau_i)=T_{a_i}^\varepsilon\, V,\vrblanc[1em]\\
      \rho^{(2)}(\tau_{i+1})=T_{c_{i+1}}^\eta\, W.
    \end{array}\right.$
} \smallskip

\noindent The homomorphisms $\rho^{(1)}$ and $\rho^{(2)}$ coincide on
the standard generators $\tau_i$ with $2\leqslant i\leqslant n-2$,
so they coincide on at least four consecutive standard generators
(four when $n=7$). Then according to Lemma
\ref{lem:uniqueness_triple}, we have $V=W$, $\varepsilon=\eta$, and
$a_i=c_i$ for all $i\in\{2,\,3,\dots,\,n-2\}$. Let us denote by
$a_{n-1}$ the curve $c_{n-1}$. Then, we have for all
$i\in\{1,\,\dots,\,n-1\}$:

\noindent \centrer{$\rho(\tau_i)=T_{a_i}^\varepsilon\, V.$}
\smallskip

\noindent We just have to check that this is a transvection of
monodromy homomorphism.
\begin{itemize}
\item  On one hand, the mapping class $V$ is in the centralizer of
$\big\langle T_{a_i},\,1\leqslant i\leqslant n-1\big\rangle$, since
the equality $V=W$ implies that $V$ commutes with $T_{a_{n-1}}$ as
well;

\item On the other hand, the ordered list of curves $(a_i,\ 1\leqslant i\leqslant n-1)$,
is an $(n-1)$-chain, since the curve $a_{n-1}$ intersects $a_{n-2}$
in one point, does not intersect the curves $a_i$ for
$i\in\{2,\dots,\,n-3\}$, and does not intersect either $a_1$. Let us
justify this last point: $\tau_1$ and $\tau_{n-1}$ commute, so
$\rho(\tau_1)$ and $\rho(\tau_{n-1})$ commute, so $T_{a_1}V$ and
$T_{a_{n-1}}V$ commute. But $V$ commutes with $T_{a_1}$ and
$T_{a_{n-1}}$, so finally, $T_{a_1}$ and $T_{a_{n-1}}$ commute, so
we have $I(a_1,\,a_{n-1})=0$.
\end{itemize}
Thus, $\rho$ is a transvection of monodromy homomorphism. Hence Theorem \ref{thm:Theorem_principal}
is proved.

\begin{appendix}

\pageimpaire\addcontentsline{toc}{part}{\protect\textsc{\textbf{Appendix}}}

\part*{Appendix}


\bigskip

\section{Miscellaneous on the mapping class group}
\bigskip

\subsection{Parabolic subgroups of the mapping class group}
\medskip

The concept of \emph{parabolic subgroup of the  mapping class
group} has been introduced by L. Paris and D. Rolfsen in
[PaRo]: given a surface $\Sigma$, the \emph{parabolic subgroups
of the mapping class group of $\Sigma$} are the subgroups
induced by the inclusion of subsurfaces in $\Sigma$. The next
Theorem, due to Paris and Rolfsen, deals with some kernels
associated with parabolic subgroups of the mapping class group.
\bigskip

\begin{thm}[Paris and Rolfsen, {[PaRo]}]
            \label{thm:inclusionMCG}
\mbox{}\\ Let $\Sigma$ be a surface and $\Sigma'$ a subsurface
in $\Sigma$ such that $\bord\Sigma'$ and $\bord\Sigma$ are
disjoint. We denote by $a_1,\,a_2,\dots,\,a_r$ the boundary
components of $\Sigma'$ that bound a disk in $\Sigma$; we
denote by $b_j,\,b'_j$ for $1\leqslant j\leqslant s$ the pairs
of boundary components of $\Sigma'$ that cobound an annulus in
$\Sigma$. Then the inclusion $\iota\;:\; \Sigma'\to\Sigma$
induces a homomorphism $\iota_*\;:\;
\Mod(\Sigma',\,\bord\Sigma')\to\Mod(\Sigma,\,\bord\Sigma)$
whose kernel is the abelian group $\T$ of rank $r+s$ spanned by
$T_{a_i}$, $1\leqslant i\leqslant r$ and by
$T_{b_j}^{-1}T_{b'_j}$, $1\leqslant j\leqslant s$. \end{thm}
\bigskip


\subsection{Interactions between Dehn twists and the braid group}
\label{par:twists_of_Dehn}
\medskip

The properties of the Dehn twists are well-known. We simply
mention them.
\medskip

\begin{prop}[Dehn twists' properties, N.V. Ivanov, J.D. McCarthy \mbox{[Mc1]}]
   \label{prop:properties_twists_Dehn}
\mbox{}
\begin{itemize}
 \item For all $F$ in $\Mod(\Sigma)$, we have $F\,T_a =
     T_a\,F$ if and only if $F(a)=a$.
 \item Let $T_a$ and $T_b$ be two Dehn twists, and $i$ and
     $j$ two nonzero integers. The relation
     ${T_a}^i\,{T_b}^j = {T_b}^j\,{T_a}^i$ holds if and
     only if $I(a,\,b)=0$.
 \item Let $T_a$ and $T_b$ be two Dehn twists, and $i$ and
     $j$ two nonzero integers. The braid relation
     ${T_a}^i\,{T_b}^j\,{T_a}^i =
     {T_b}^j\,{T_a}^i\,{T_b}^j$ holds if and only if
     $I(a,\,b)=1$ and $i=j=\pm1$. \fin
 \end{itemize}
 \end{prop}
\bigskip


\begin{prop}[C. Labruère, L. Paris, cf. \mbox{[LaPa]}]
          \label{prop:relations_chaine_of_curves}
\mbox{}\\Let $(c_1,\,c_2,\dots,\,c_{k})$ be a $k$-chain of
curves where $k$ is an integer greater than or equal to 2.
Then,
\begin{itemize}
\item if $k$ is even, the tubular neighbourhood of $c_1\cup
    c_2\cup\dots\cup c_{k}$ is a surface $S$ (in grey on Figure \ref{fig:alphaOrdre2casPair}) of genus
    $\frac{k}{2}$ with one boundary component which we call
    $d$, and the product
    $\big(T_{c_1}(T_{c_2}T_{c_1})\dots(T_{c_{k}}\dots
    T_{c_2}T_{c_1})\big)^2$ is the mapping class $\alpha$
    that preserves each curve $c_i$, $1\leqslant i
    \leqslant k$, whose restriction outside of $S$
    coincides with the identity, and such that
    $\alpha^2=T_d$ (cf. Figure
    \ref{fig:alphaOrdre2casPair}). Notice that after having
    given orientations to the curves $c_i$, $1\leqslant
    i\leqslant k$, the mapping class $\alpha$ inverse them.
  \begin{figure}[!h]
   \Includegraphics{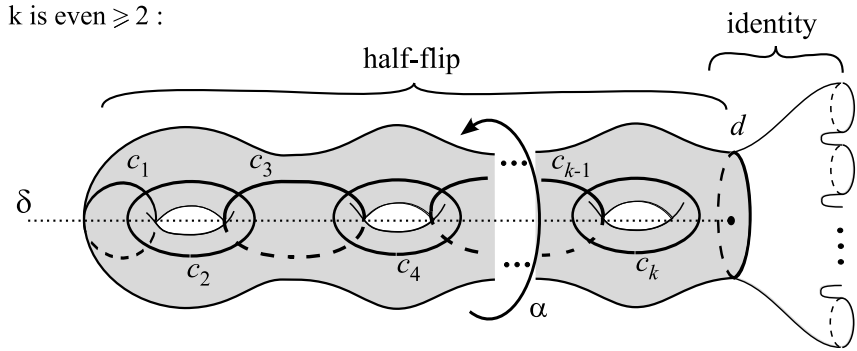}
   \caption{The product $\big(T_{c_1}(T_{c_2}T_{c_1})\dots(T_{c_{k}}\dots T_{c_2}T_{c_1})\big)^2$ when $k$ is even.}
   \label{fig:alphaOrdre2casPair}
  \end{figure}
\item if $k$ is odd, the tubular neighbourhood of $c_1\cup
    c_2\cup\dots\cup c_{k}$ is a surface $S$ (in grey on Figure \ref{fig:alphaOrdre2casImpair}) of genus
    $\frac{k-1}{2}$ with two boundary components $d_1$ and
    $d_2$, and the product \linebreak
    $\big(T_{c_1}(T_{c_2}T_{c_1})\dots(T_{c_{k}}\dots
    T_{c_2}T_{c_1})\big)^2$ equals the product $T_{d_1}\,
    T_{d_2}$ (cf. Figure
    \ref{fig:alphaOrdre2casImpair}).\fin
  \begin{figure}[!h]
   \Includegraphics{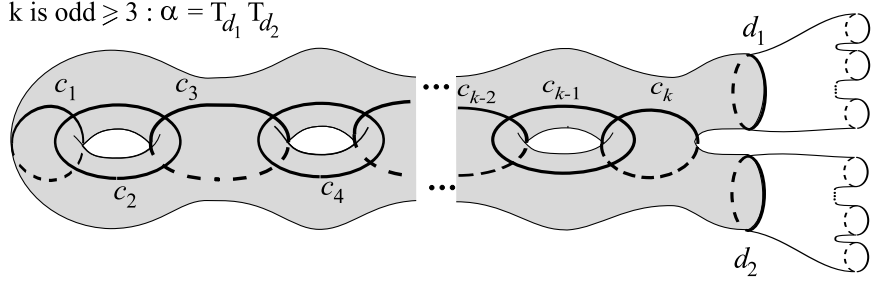}
   \caption{The product $T_{c_1}(T_{c_2}T_{c_1})\dots(T_{c_{k}}\dots T_{c_2}T_{c_1})$ when $k$ is odd.}
   \label{fig:alphaOrdre2casImpair}
  \end{figure}
\end{itemize}
\end{prop}
\bigskip

Let $k$ be an integer and let $\Sigma$ be the surface in grey
in Figure \ref{fig:alphaOrdre2casPair} if $k$ in odd or
in Figure \ref{fig:alphaOrdre2casImpair} if $k$ in even.
Using Proposition \ref{prop:properties_twists_Dehn},
there is an obvious homomorphism from $\B_{k+1}$ to $\Mod(\Sigma,\,\bord\Sigma)$
sending $\tau_i$ on $T_{c_i}$ for all $i\leqslant k+1$.
Let us denote it by $\rhor$. Then:
\medskip

\begin{thm}[Geometric embedding the braid group, Birman \& Hilden, \mbox{[BiHi]}]
            \label{thm:BiHi}
            \label{thm:Perron_Vannier_1}
\mbox{}\\
The homomorphism $\rhor\ :\ \B_{k+1}\to\Mod(\Sigma,\,\bord\Sigma)$ is injective.
\end{thm}
\bigskip

The following Theorem is due independently to M. Dehn [D] and
W.B.R. Lickorish [Lk]. For the proof of Part 1, the reader is
referred to [Bi], Theorem 4.1. Part 2 is a consequence of Part 1,
using \emph{lantern relations} (cf. [FaMa] section 5.1) in
order to obtain the Dehn twists along the boundary curves.
\smallskip

\begin{figure}[!h]
 \Includegraphics{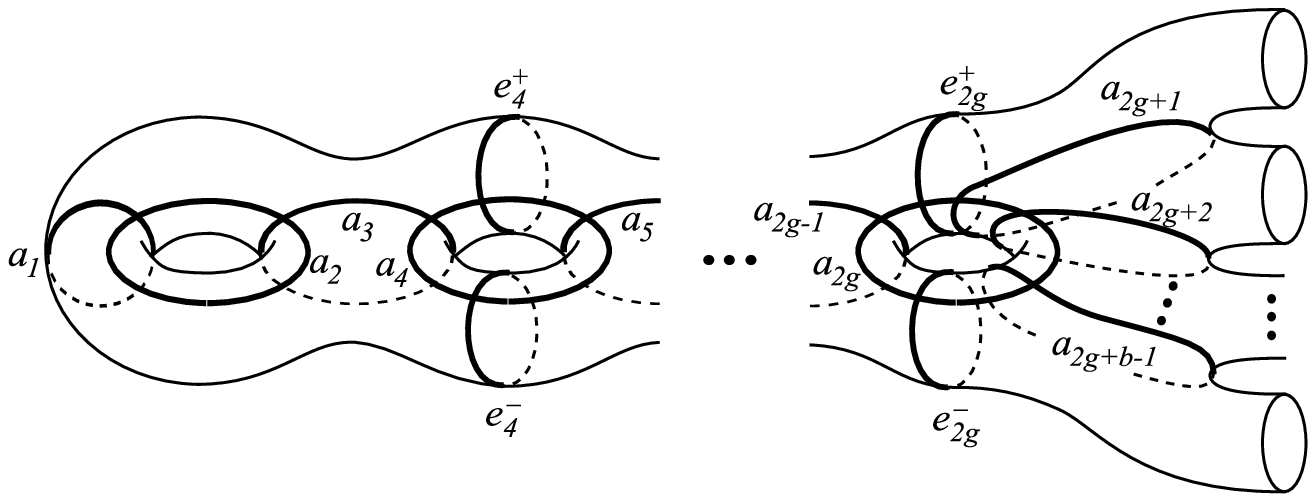}
 \caption{Curves in $\Sigma_{g,b}$ span $\PMod(\Sigma_{g,b})$.}
 \label{fig:genererMod}
\end{figure}

\begin{thm}[The mapping class group is spanned by Dehn twists]
   \label{thm:Mod_spanned_by_the_twists_of_Dehn}
\mbox{}

1. Let $\Sigma_{g,b}$ be a surface whose genus satisfies
$g\geqslant 1$. The group $\PMod(\Sigma_{g,b})$ is spanned by
the Dehn twists along the non-separating curves of
$\Sigma_{g,b}$. Actually, one can for instance consider only
the curves drawn in Figure \ref{fig:genererMod}.

2. The same Dehn twists can be seen as lying in
$\Mod(\Sigma_{g,b},\,\bord\Sigma_{g,b})$; in this case they
span $\Mod(\Sigma_{g,b},\,\bord\Sigma_{g,b})$. \fin \end{thm}

\bigskip


\subsection{Periodic Mapping classes}
\medskip

\subsubsection{Kerckhoff's Theorem, Riemann-Hurwitz formula, and related results}
\smallskip

According to the Nielsen realization theorem,
every periodic mapping class on a closed surface
is the isotopy class of an isometry of the surface with a
specific metric depending on the mapping class (and of course,
the isometry has the same order as the mapping class). Given a
connected surface $\Sigma$ together with a hyperbolic metric
$g$, we denote by $\Isom^+(\Sigma,\,g)$ the group of
positive isometries of $(\Sigma,\,g)$, which is a subgroup of
$\Diff^+(\Sigma)$.
This Theorem has been
generalized by Kerckhoff to finite subgroups of $\Mod(\Sigma)$
where $\bord\Sigma$ is possibly nonempty (cf.[Ke1] and [Ke2]):
\bigskip

\begin{thm}[``Nielsen Realization Problem'', Kerckhoff, cf. \mbox{[Ke2]}]
             \label{thm:Kerckhoff}
\mbox{}\\Let $\Sigma$ be a surface with a possibly nonempty
boundary. Let $\Gamma$ be a finite subgroup of $\Mod(\Sigma)$.
Then, there exists a finite group $\bar\Gamma$ of
$\Diff^+(\Sigma)$ such that the natural homomorphism
$\Diff^+(\Sigma)\to\Mod(\Sigma)$ sends isomorphically
$\bar\Gamma$ on $\Gamma$. Moreover, we can choose $\bar\Gamma$
as a subgroup of the isometry group of $\Sigma$ equipped with a
metric of constant curvature, where the boundary components are
geodesics. \end{thm}
\bigskip

\begin{lem}[Riemann-Hurwitz Formula]
        \label{lem:periodic_Euler}
\mbox{}\\ Let $\Sigma$ be a surface with a possibly trivial
boundary and let $\Gamma$ be a finite subgroup of order $m$ of
$\Diff^+(\Sigma)$, such that $m\geqslant 1$. Then
$\Sigma/\Gamma$ is a surface and the quotient map $\pi\;:\;
\Sigma\to\Sigma/\Gamma$ is a ramified covering. Let $Q_1,\dots,
Q_\ell$ be the ramification points of $\pi$, and, for $1\leqslant
i \leqslant \ell$, let $o(Q_i)$ be the number of preimages of
$Q_i$ by $\pi$. Then the Euler characteristics of the surfaces $\Sigma$ and
$\Sigma/\Gamma$ are linked by the formula:
$$
\chi(\Sigma)+\sum_{i=1}^{\ell} (m-o(Q_i))=m.\chi(\Sigma/\Gamma).
$$
 \fin
\end{lem}
\smallskip

\TITRE{Ramification points and singular points.}
\mbox{}\\ The points of $\Sigma/\Gamma$ whose number of
preimages by $\pi$ is smaller than $|\Gamma|$ are called
\emph{ramification points}. Their preimages in $\Sigma$ are
called \emph{singular points}. Let us make some remarks:
\smallskip

\Point For all $i$, the group $\Gamma$ acts transitively on
    $\pi^{-1}(\{Q_i\})$. For this reason, the cardinality
    of $\pi^{-1}(\{Q_i\})$, denoted by $o(Q_i)$, divides
    $m$.
\smallskip

\Point When $\Gamma$ is spanned by a unique element $f$, for
    all $i\leqslant \ell$, the action of $f$ on
    $\pi^{-1}(\{Q_i\})$ is cyclic, so $\pi^{-1}(\{Q_i\})$
    belongs to $\Fix(f^{o(Q_i)})$, the set of fixed points
    of $f^{o(Q_i)}$. However, according to Lefschetz
    Theorem, the number of fixed points of $f$ and of its
    powers depends only on the isotopy class of $f$ (let us
    make clear that in Lefschetz Theorem, the number of
    fixed points takes into account the multiplicity of
    each fixed point, but in the case of nontrivial
    isometries, this integer always equals 1). Hence the
    data of $\ell$ and $\{o(Q_i)\ ,\ 1\leqslant
    i\leqslant\ell\}$ is an invariant of the isotopy class
    of $f$. Since $\chi(\Sigma/\langle f\rangle)$ and the
    number of boundary components of $\Sigma/\langle
    f\rangle$ are also invariant by isotopy on $f$, the
    surface $\Sigma/\langle f\rangle$ itself is an
    invariant of the isotopy class of $f$.
\smallskip

\Point This lemma together with Kerckhoff's Theorem have a
    lot of corollaries and here are some of them that we
    will use in this paper.
\bigskip

\begin{cor}[Fixed points and preserved boundary components preserved]
        \label{cor:points_fixes}
\mbox{}\\
Let $\Sigma$ be a surface $\Sigma_{g,\,b}$ and let $F$ be a
periodic mapping class of $\Mod(\Sigma)$ of order $m$. Then,
the sum of the number of boundary components preserved by $F$
and the number of fixed points of $F$ is bounded by
$2+\frac{2g}{m-1}$. \end{cor}
\bigskip

\DEM Let $\bar F$ be a diffeomorphism of order $m$ representing
$F$. Let $b'$ be the number of boundary components preserved by
$\bar F$ and let $\ell'$ be the number of fixed points of $\bar
F$. Let $\D$ be the set of boundary components of $\Sigma$ that
are not fixed by $\bar F$. Let $\Sigma'$ be the surface
obtained from $\Sigma$ where each boundary components of $\D$
have been filled. Let $\bar F'$ the diffeomorphism induced by
$\bar F$ on $\Sigma'$. Let us apply Lemma
\ref{lem:periodic_Euler} to $\bar F'$, and let us adopt the
notation of the statement. The order $m$ of $\bar F'$ (equal to
the order of $\bar F$) coincides with the order of the group
$\Gamma=\langle \bar F'\rangle$. In addition, notice that:
\begin{itemize}
\item[\point] $\chi(\Sigma')=2-2g-b'$;
\item[\point] for all fixed point $Q$, we have $o(Q)=1$;
\item[\point] $\sum_{i=1}^{\ell} (m-o(Q_i))\geqslant
    \ell'(m-1)$;
\item[\point] $\chi(\Sigma'/\Gamma)\leqslant 2-b'$ .
\end{itemize} Finally, we get:
\smallskip

\centrer{$(2-2g-b')+\ell'(m-1)\leqslant m.(2-b')$,}

\noindent whence:

\centrer{$(m-1)b'+\ell'(m-1)\leqslant 2(m-1)+2g$,}
\smallskip

\noindent whence the formula: $b'+\ell'\leqslant
2+\frac{2g}{m-1}$\ . \fin
\bigskip

The following corollary is the special case $g=0$.
\medskip

\begin{cor}[Periodic mapping classes on a sphere]
   \label{cor:periodic_we_a_sphere}
\mbox{}\\Let $S$ be a holed sphere and $F$ a periodic mapping
class of $\Mod(S)$. If there exist at least three boundary
components in $S$ preserved by $F$, then $F$ is the identity of
$\Mod(S)$.\fin
\end{cor}
\bigskip

Finally the following result is a technical corollary of Lemma
\ref{lem:periodic_Euler} which gives an upper bound to the cardinality
of a finite subgroup of the mapping class group of a surface
$\Sigma$ of genus $g$. It is due to Hurwitz. For a proof, see
[FaMa], section 6.2.
\smallskip

\begin{cor}[``84($g-1$)'' Theorem]
      \label{cor:borne_cardinality_Mod}
\mbox{}\\ Let $\Sigma$ be a connected surface of genus at least
2 without boundary . The order of a finite subgroup of the
mapping class group is bounded by $42|\chi(\Sigma)|$ (which is
equal to $84(g-1)$).\fin
\end{cor}
\medskip

\noindent When the finite subgroup of the mapping class group is cyclic,
we have much more precise upper bound:
\medskip

\begin{prop}[``$4g+2$'' Theorem]
 \mbox{}\\The order of a periodic mapping class in
 $\Mod(\Sigma_{g,\,0})$ is bounded by $4g+2$.
 \end{prop}
 \medskip

We need a more general statement when $g=1$.
\medskip

\begin{lem}[Order of periodic mapping classes in $\PMod(\Sigma_{1,\,b})$, $b\geqslant 0$]
  \label{cor:order_periodic_à_boundary_genre_1}
\mbox{}\\ Any mapping class of $\PMod(\Sigma_{1,\,b})$ with $b\geqslant 0$ is of order smaller than or equal to 6.
\end{lem}
\smallskip

\DEM When $b=0$, this is the ``$4g+2$'' theorem.
When $b=1$, it comes from the fact that $\Mod(\Sigma_{1,\,0})=\Mod(\Sigma_{1,\,1})$.
For the case $b=2$, we use the fact that $\PMod(\Sigma_{1,\,2})$ is isomorphic to the quotient of $\B_4$
by its center, and it is well-known that the order of the periodic elements of $\B_n$ modulo its center is
smaller than or equal to $n$. Finally, when $b\geqslant 3$,
according to Corollary \ref{cor:points_fixes}, for any periodic element of order $m$
we have the inequality $b\leqslant 2+\frac{2g}{m-1}$ where $g=1$,
hence $m\leqslant 1+\frac{2}{b-2}$. So $m\leqslant 3$.
In all cases, the order of a periodic mapping class is less than or equal to 6.
\fin\bigskip

\subsubsection{When periodic mapping classes preserve a boundary component}
\smallskip

In this subsubsection, we consider a surface $\Sigma$ with a
nonempty boundary, we choose a boundary component $d$, and look
only at the mapping class group of $\Sigma$ that preserve globally
$d$. All the following results are based on a classic
result coming from the Riemannian manifolds theory:
\medskip

\begin{lem} \label{lem:rotation}
Let $M$ a Riemannian manifold and $f$ an isometry on $M$. If
$f$ fixes a point $x$ and if its differential in $x$ is
the identity, then $f$ is the identity. We have the same
conclusion when $f$ fixes a boundary curve of $M$
pointwise.\fin \end{lem}
\medskip

\TITRE{Corollaries of Lemma \ref{lem:rotation}.} Lemma
\ref{lem:rotation} is fundamental to the comprehension of the
periodic mapping classes and induces many essential corollaries
in this paper:

\begin{itemize}
\item Lemma \ref{lem:without_torsion}: ``\emph{If
    $\bord\Sigma$ is nonempty, then
    $\Mod(\Sigma,\,\bord\Sigma)$ is torsion-free.}''

\item Corollary \ref{cor:reconstruction_periodic} of Lemma
    \ref{lem:without_torsion}, which allows us to deduce
    the existence of some periodic mapping classes from the
    existence of some other periodic mapping classes lying
    in ``\emph{smaller}'' subsurface (i.e. of smaller genus
    or of greater Euler characteristic).

\item Lemma \ref{lem:group_cyclic}: ``\emph{If
    $\bord\Sigma$ is nonempty and if $d$ belongs to
    $\Bord(\Sigma)$, then any finite subgroup of
    $\Mod_d(\Sigma)$ is cyclic.}''
\end{itemize}
\bigskip

\begin{lem}[Behaviour of periodic mapping classes in the neighbourhood of $\bord\Sigma$]
      \label{lem:without_torsion}
\mbox{}\\Let $\Sigma$ be a connected surface such that
$\bord\Sigma\not=\vide$. Let $d$ be a boundary curve of
$\Sigma$ and let $\Mod(\Sigma,\,d)$ be the group of the mapping
classes of $\Sigma$ that fix $d$ pointwise. Then the group
$\Mod(\Sigma,\,d)$ is torsion-free. Moreover, let $F$ be a
periodic mapping class of period $m\geqslant 2$ and belonging
to $\Mod(\Sigma,\,d)$. Then there exists an integer $\ell$
coprime to $m$ such that $F^m=T_{d}^{\,\ell}$ (in particular,
$\ell$ is nonzero). \end{lem}

\DEM Remember that a mapping class in $\Mod(\Sigma,\,d)$ is
said to be periodic if one of its nonzero power induces the
trivial mapping class in $\Mod(\Sigma)$. So any mapping class
that would belong to the torsion of $\Mod(\Sigma,\,d)$ is by
definition a periodic mapping class. However if we show that
any nonzero power of any periodic mapping class is nontrivial
in $\Mod(\Sigma,\,d)$, we will have shown that
$\Mod(\Sigma,\,d)$ is torsion-free. So the second part of Lemma
\ref{lem:without_torsion} implies the first part. Hence, by
showing the second part, we will be done.

\begin{figure}[!h]
 \Includegraphics{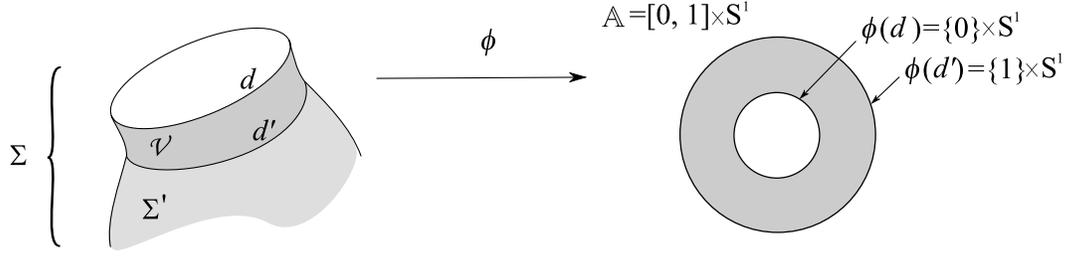}
 \caption{The situation described in proof of Lemma \ref{lem:without_torsion}.}
 \label{fig:carte2}
\end{figure}

Let $\bar F$ be a diffeomorphism in $\Diff^+(\Sigma,\,d)$ which
represents $F$ and fixes $d$ pointwise. Let $d'$ be a curve
isotopic to $d$ in $\Sigma$ that lies outside of $\bord\Sigma$.
Let us call $\V$ the compact cylinder included in $\Sigma$
whose both boundary components are $d$ and $d'$. Let $\Sigma'$
be the closure of the complement of $\V$ in $\Sigma$ (cf.
Figure \ref{fig:carte2}).

Let us apply the Nielsen - Kerckhoff realization theorem: there
exist a hyperbolic metric $g$ on $\Sigma'$ and an isometry
$\bar F_1$ of $(\Sigma',\,g)$ representing the restriction of
$F$ to $\Sigma'$. Let us denote by $f$ the restriction of $\bar
F_1$ to $d'$. Since $\bar F_1$ is a periodic isometry  of order
$m$, there exists an integer $k$ such that $f$ is a rotation of
angle $\frac{2k\pi}{m}$ on $d'$ equipped with the induced
metric. But for all $p\in\{1,\,2,\dots,\,m-1\}$, the mapping
class $(\for_d(F))^p$ (where $\for_d(F)$ is the mapping class
induced by $F$ in $\Mod(\Sigma)$) is different from the
identity, so according to Lemma \ref{lem:rotation}, $f^p$ is
different from the identity. Hence $k$ is coprime to $m$.

\begin{figure}[!h]
 \Includegraphics{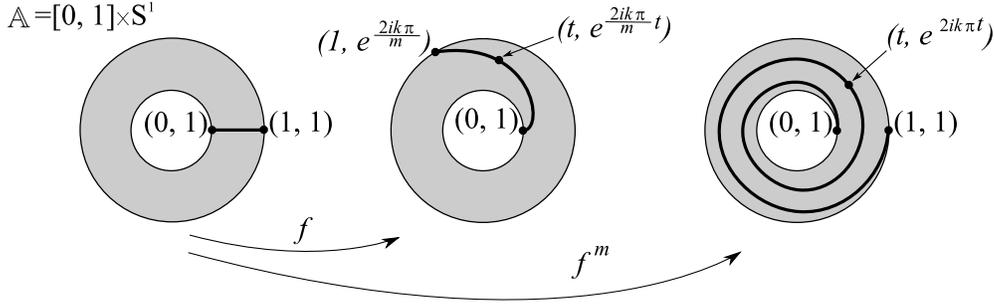}
 \caption{Image of the segment $[0,\,1]\times \{1\}$ by $f$ and $f^{\,m}$.}
 \label{fig:periodiqueAuxBords2}
\end{figure}

Let $\AA$ be the annulus $[0,\,1]\times S^1$ and $\phi$ a
positive diffeomorphism of $\V$ in $\AA$ such that
$\phi(d)=\{0\}\times S^1$ and $\phi(d')=\{1\}\times S^1$ (cf.
Figure \ref{fig:carte2}). Moreover, we can construct $\phi$ so
that the map $\phi\rond f\rond\phi^{-1}$ from $\{1\}\times S^1$
to $\{1\}\times S^1$ coincides with the function
$(1,e^{i\theta})\mapsto(1,e^{i(\theta+\frac{2k\pi}{m})})$, for
$f$ is an angle $\frac{2k\pi}{m}$ rotation on $d'$. We extend
$f$ on $\V$ (cf. Figure \ref{fig:periodiqueAuxBords2}) by
setting for all $t\in[0,\,1]$ and $\theta\in[0,\,2\pi[$:
\smallskip

\centrer{
$f\big(\,\phi^{-1}\big((t,e^{i\theta})\big)\,\big)=\phi^{-1}\big((t,e^{i(\theta+\frac{2k\pi}{m}t)})\big)$.
}
\smallskip

\noindent In this way, $f$ fixes $d$ pointwise and coincides
with $\bar F_1$ on $d'$. Let $\bar F_2$ be the homeomorphism of
$\Sigma$ fixing $d$ pointwise and coinciding with $f$ on $\V$
and with $\bar F_1$ on $\Sigma'$. By definition of a Dehn
twist, the isotopy class of the homeomorphism $\bar F_2^{\,m}$
is equal to $T_d^{\,k}\in\Mod(\Sigma,\,d)$. Let $G$ be the
mapping class of $\Mod(\Sigma,\,d)$ containing the
homeomorphism $\bar F\bar F_2^{-1}$. By construction, we have
$\for_d(\bar F)=\for_d(\bar F_2)$, so, according to the
following central exact sequence:
\smallskip

\centrer{$1\to\langle
T_d\rangle\to\Mod(\Sigma,\,d)\xrightarrow{\for_d}\Mod_d(\Sigma)\to
1$,}
\smallskip

\noindent the isotopy classes of $\bar F$ and $\bar F_2^{-1}$
in $\Mod(\Sigma,\,d)$ commute, so the homeomorphisms $(\bar
F\bar F_2^{-1})^m$ and $\bar F^m\,\bar F_2^{-m}$ are isotopic,
hence $G^m=F^mT_d^{-k}$. But $\for_d(\bar F)=\for_d(\bar F_2)$,
so, according to the same exact sequence, there exists an
integer $j$ such that $G=T_d^{\,j}$. These two last equalities
imply that $F^m=G^mT_d^{\,k}=T_d^{\,jm+k}$, and $jm+k$ is
coprime to $m$. \fin
\bigskip

Let us give a corollary of Lemma \ref{lem:without_torsion} that
will allow us to show the existence of some periodic mapping
classes from periodic mapping classes defined on
``\emph{smaller}'' surfaces (i.e. of smaller genus or of
greater Euler characteristic).
\smallskip

\begin{cor}[Rebuilding a periodic mapping class]
              \label{cor:reconstruction_periodic}
\mbox{}\\
Let $\Sigma$ be a connected surface. Let $I$ be a finite set
and let ${\A_I}=\{a_i,\,i\in I\}$ be a curve simplex in
$\Sigma$. For all $i\in{I}$, we denote by $a_i^+$ and $a_i^-$
the boundary components of $\Sigma_{\A_I}$ coming from the cut
of ${\A_I}$ along the curve $a_i$. Let $\tildeF$ be a periodic
mapping class of $\Mod(\Sigma_{\A_I})$ of order two that
preserves $\{a_i^+,\,a_i^-\}$ for all $i\in I$. Then there
exists a unique periodic mapping class $F\in\Mod(\Sigma)$ of
order two such that $F$ induces $\tildeF$ in
$\Mod(\Sigma_{\A_I})$. \end{cor}
\smallskip

\DEM
\smallskip

1. \emph{Notation}
\smallskip

\noindent Let $J$ and $K$ be two subsets of $I$ such that
$J\sqcup K$ forms a partition of $I$ and such that for all
$j\in J$, the boundary components $a_j^+$ and $a_j^-$ of
$\Sigma_{\A_I}$ are swapped by $\tildeF$, whereas for all $k\in
K$, the boundary components $a_k^+$ and $a_k^-$ of
$\Sigma_{\A_I}$ are not. Let $\A_J=\{a_j,\,j\in J\}$ and
$\A_K=\{a_k,\,k\in K\}$.
\medskip

2. \emph{Existence of $F$}
\smallskip

\noindent Let $\bar F_1$ be a diffeomorphism of
$\Diff^+(\Sigma_{\A_I})$ of order 2 representing $\tildeF$. We
construct a diffeomorphism $\bar F_2$ of
$\Diff^+(\Sigma_{\A_K})$ by identifying in $\Sigma_{\A_I}$
$a_j^+$ and $a_j^-$ thanks to the following relation in
$\Sigma_{\A_I}$: for all pairs of points
$(x,\,y)\in\Sigma_{\A_I}$,
\smallskip

\centrer{$x\sim y$ if and only if
  $\left\{\begin{array}{l}
    \blanc\! x=y,\\
    \mbox{or}\\
    \blanc\! x,\,y\in a_j^+\cup a_j^- \mbox{ and } y=\bar F_1(x).
  \end{array}\right.$}
\smallskip

\noindent The differential structure of $\Sigma_{\A_K}$ induced
by the one of $\Sigma_{\A_I}$ is well-defined up to conjugation
by a diffeomorphism (cf. [Hi] page 184: \emph{Gluing Manifolds
Together}); and by construction, $\bar F^2$ is a diffeomorphism
of $\Sigma_{\A_K}$. Let $F_2$ be the mapping class of
$\Mod(\Sigma_{\A_K})$ containing $\bar F_2$. It is clear that
$F_2$ is periodic of order two. Notice that $F_2$ preserves
each boundary $a_k^\varepsilon$, $k\in K$,
$\varepsilon\in\{+,\,-\}$. One can then define $F_3$ as being a
representative of $F_2$ in
$\Mod(\Sigma_{\A_K},\{a_k^+,\,a_k^-,\,k\in K\})$. According to
Lemma \ref{lem:without_torsion}, for all $k\in{K}$, there exist
two odd integers $m_k^+$ and $m_k^-$ such that:
\smallskip

\centrer{$F_3^{\,2}=\prod_{k\in{K}}
T_{a_k^+}^{\,m_k^+}\,T_{a_k^-}^{\,m_k^-}$.}
\smallskip

\noindent Let $F_4$ be the mapping class of $\Mod(\Sigma)$
obtained by identifying the boundary components $a_k^+$ and
$a_k^-$ for all $k\in K$. Thus $F_4$ satisfies the equality:
\smallskip

\centrer{$F_4^{\,2}=\prod_{k\in K} T_{a_k}^{\,(m_k^++m_k^-)}$.}
\smallskip

\noindent Let us set $F_5=F_4 \big(\prod_{k\in K}
T_{a_k}^{\,\ell_k}\big)^{-1}$, where, for all $k\in K$, the
rational $\ell_k$ defined by $\ell_k=\frac{m_k^++m_k^-}{2}$ is
an integer since $m_k^+$ and $m_k^-$ are both odd. Since this
product is commutative, we have $F_5^{\,2}=\Id$. By
construction, the mapping class $F_5$ fixes the curves of
$\A_I$ and induces in $\Mod(\Sigma_{\A_I})$ the mapping class
$\tildeF$. So the mapping class $F_5$ plays the role of $F$ in
the statement of Corollary \ref{cor:reconstruction_periodic}.
\bigskip

3. \emph{Uniqueness of $F$}
\smallskip

\noindent Let us assume that there exist two mapping classes
$F$ and $F'$ in $\Mod_{\A_I}(\Sigma)$ such that $F^2=F'^2=\Id$
and such that $F$ and $F'$ induce $\tildeF$ in
$\Mod(\Sigma_{\A_I})$. According to the following exact
sequence:
\smallskip

\centrer{$1\to\langle
T_{a_i},\,i\in{I}\rangle\to\Mod_{\A_I}(\Sigma)\to\Mod(\Sigma_{\A_I})\to
1$,}
\smallskip

\noindent there exist some integers $p_i$, $i\in I$ and a
mapping class $W=\prod_{i\in I} T_{a_i}^{\,p_i}$ such that
$F'=FW$. But $F$ and $F'$ fix the curves of $\A_I$, so these
two mapping classes commute with $W$ and we have
$(FF'^{-1})^2=W^2$. However, $FF'^{-1}$ is periodic of order
two as a commutative  product of two periodic mapping classes
of order two, so the mapping class $W^2=\prod_{i\in I}
T_{a_i}^{\,2p_i}$ is trivial. Then for all $i\in I$, we have
$2p_i=0$, so $W$ is trivial. Therefore $F=F'$. \fin
\bigskip

\begin{lem}
        \label{lem:group_cyclic}
Let $\Sigma$ be a connected surface with a nonempty boundary
and let $d$ be a boundary component. Any finite subgroup of
$\Mod_d(\Sigma)$ is cyclic. \end{lem}
\smallskip

\DEM Let $\Gamma$ be a finite group included in
$\Mod_d(\Sigma)$. According to Kerckhoff's Theorem (cf. Theorem
\ref{thm:Kerckhoff}) there exist a hyperbolic metric $g$ on
$\Sigma$ and a finite group $\bar\Gamma$ of
$\Isom^+(\Sigma,\,g)$ that is sent isomorphically on $\Gamma$
by the natural homomorphism $\Diff^+(\Sigma)\to\Mod(\Sigma)$. Let
us give an orientation to $d$ and define the map
\smallskip

\centrer{$\theta\;:\; \bar\Gamma\to\RR/2\pi\ZZ$}
\smallskip

\noindent that sends an isometry $\bar K$ to the angle of the
induced rotation by $\bar K$ on $d$. It is clear that $\theta$
is a homomorphism. According to Lemma \ref{lem:rotation}, any
isometry that fixes $d$ pointwise is the identity, so $\theta$
is an injective homomorphism. Thus $\Gamma$ is isomorphic to a
finite subgroup of $\RR/2\pi\ZZ$, hence $\Gamma$ is cyclic.\fin
\bigskip


\subsection{The theory of Pseudo-Anosov diffeomorphisms.}
\label{par:théorie_pseudo-Anosov}
\medskip

We will now recall the fundamental results of the theory of the
\emph{pseudo-Anosov diffeomorphisms} on surfaces without
boundary, in order to prove Proposition
\ref{prop:structure_centralizer}. There exists a very rich
literature concerning this theory. In this subsection, we will
lean mainly on the article of C. Bonatti and L. Paris [BoPa],
section 2. The classic references are [FLP], [Th], [BlCa],
[Iv2]. Let us begin by a remark on the terminology concerning
the ``pseudo-Anosov diffeomorphisms''. In this subsection,
$\Sigma$ is a surface without boundary.
\medskip

\TITRE{Pseudo-Anosov diffeomorphisms and homeomorphisms.} Let
us recall that the ``pseudo-Anosov diffeomorphisms'' on a
surface $\Sigma$ are actually both homeomorphisms on $\Sigma$,
and diffeomorphisms on the surface $\Sigma$ minus a finite
number of points. In the neighbourhood of these points, called
\emph{singular points}, the homeomorphisms are not
differentiable although they are perfectly known, cf. [BoPa].
For this reason C. Bonatti and L. Paris [BoPa] speak about
``pseudo-Anosov homeomorphisms'' when we speak about
``pseudo-Anosov diffeomorphisms'', but we deal with the same
objets. Let us recall that there exist true diffeomorphisms in
the isotopy class of a ``pseudo-Anosov diffeomorphism''
(isotopic means here linked by a path of homeomorphisms) and a
``pseudo-Anosov diffeomorphism'' defines always a unique
isotopy class of true diffeomorphisms, so it defines a unique
mapping class.
\medskip

\TITRE{Properties of pseudo-Anosov diffeomorphisms.} The
fundamental theorem of the pseudo-Anosov theory, due to
Thurston (cf. [BlCa], [Th]), is the following:
\medskip

\centrer{\parbox{14cm}{\emph{A mapping class $F\in\M(\Sigma)$
is pseudo-Anosov in the sense of Definition
\ref{défi:pseudo-Anosov}, if and only if there exists a
``\emph{pseudo-Anosov diffeomorphism}'' (cf. below) $\bar
F\in\Diff(\Sigma)$ representing $F$.}}}
\medskip

\noindent Instead of a definition of ``\emph{pseudo-Anosov
diffeomorphism}'', we describe some of their properties. This
will be enough for our purpose. Let $\bar F$ be any
``\emph{pseudo-Anosov diffeomorphism}'' in $\Diff(\Sigma)$. By
definition, $\bar F$ satisfies the following properties:
\bigskip

\TITRE{The foliations $\F^s$, $\F^u$, Singular points,
separatrices and indexes.} There exist a finite set $\S$ of
points in $\Sigma$ preserved by $\bar F$ and a pair of
transverse measured regular foliations on
$\Sigma\smallsetminus\S$. We say we have a unique pair of
transverse measured singular foliations $\F^s$ and $\F^u$
invariant by $\bar F$ on $\Sigma$ that are preserved by
$\Sigma$. These foliations are called the \emph{stable and
unstable foliations} (cf. [BoPa]). An integer $k\geqslant 3$ is
associated to each singular point $P\in\S$ and corresponds to
the number of leaves of $\F^s$ that end in $P$ (taking $\F^u$
instead of $\F^s$ would have lead to the same number). These
leaves ending in a singular point will be called
\emph{separatrices} (cf. [BoPa]). The integer $k$ associated to
a singular point $P$ will be denoted by
$\Ind(\F^s,\,\F^u\,:\,P)$ and called \emph{the index of
$(\F^s,\,\F^u)$ at the singular point $P$}. Let us put the
emphasis on the fact that for any $P\in\S$, we have
$\Ind(\F^s,\,\F^u\,:\,P)\geqslant 3$. Furthermore, each
separatrix contains a unique singular point (cf. [BoPa]
proposition 2.1, assertion (3)). The indices and the Euler
characteristic of $\Sigma$ satisfy the following  (see [FLP],
Exposé 5):
\smallskip

\begin{prop}
        \label{prop:indices_and_chi}
\mbox{}\smallskip

\centrer{$\displaystyle
\chi(\Sigma)=\sum_{P\in\S}\bigg(1-\frac{\Ind(\F^s,\,\F^u\,:\,P)}{2}\bigg)$.}
\smallskip

\end{prop}
\bigskip

\TITRE{The group $\Norm(\F^s,\,\F^u)$ and the homomorphism $L$.}
Let us denote by $\Norm(\F^s,\,\F^u)$ the set of
diffeomorphisms of $\Diff(\Sigma)$ that preserve $\F^s$ and
$\F^u$. The elements of $\Norm(\F^s,\,\F^u)$ send the singular
points on singular points of same index. They are either
pseudo-Anosov or periodic. Proposition
\ref{prop:properties_of_L} gives the main properties of the
\emph{dilatation coefficient} $\lambda^u(.)$ (cf. [BoPa]) of
the diffeomorphisms of $\Norm(\F^s,\,\F^u)$, using the homomorphism
\smallskip

\centrer{$\DEF{L}{\Norm(\F^s,\,\F^u)}{\mathbb{R}}{\bar
G}{\log\big(\, \lambda^u(\bar G)\,\big)}$.}
\smallskip

\noindent We won't use the definition of the dilatation
coefficient, but only the existence of the homomorphism $L$ and
Proposition \ref{prop:properties_of_L} below. \bigskip

\begin{prop}[Properties of $L$]
      \label{prop:properties_of_L}
\mbox{}\\
The homomorphism $L\;:\; \Norm(\F^s,\,\F^u)\to \mathbb{R}$
(depending on the pseudo-Anosov diffeomorphism $\bar F$)
satisfies the five following properties:

\begin{itemize}
\item[\;\;(i)] the real number $L(\bar F)$ satisfies
    $L(\bar F)>0$,
\item[\;(ii)] the image of $L$ is isomorphic to $\ZZ$,
\item[(iii)]  the kernel $\Ker(L)$ acts freely on the set
    of separatrices of $\F^u$,
\item[(iv)]  the kernel $\Ker(L)$ is a finite group.
\end{itemize} \end{prop}
\medskip

\DEM
\begin{itemize}
 \item[\;\;(i)] By definition of the dilatation
     coefficient, $L(F)=\log(\lambda^u(F))$ satisfies
     $L(F)>0$ for $\lambda^u(F)>1$ (cf. [FLP], Exposé 9).
\smallskip

 \item[\;(ii)] We refer to [ArYo] for this result.
\smallskip

 \item[(iii)] Since $\Ker(L)$ is included in
     $\Norm(\F^s,\,\F^u)$, the elements of $\Ker(L)$
     preserve $\F^u$ and permute the separatrices.
     According to Lemma 2.11 of [BoPa], the induced action
     by $\Ker(L)$ on the set of separatrices of $\F^u$ is
     free (the group $\Ker(L)$ is denoted by
     $\Sym(\F^s,\,\F^u)$ in [BoPa]).

  \item[(iv)] Since $\S$ is finite and since only a finite
      number of separatrices end in each singular point,
      the set of separatrices is finite. But $\Ker(L)$ acts
      freely on it, hence $\Ker(L)$ is a finite group.\fin
\end{itemize}
\bigskip

Let us recall that for all pseudo-Anosov mapping classes
$F\in\Mod(\Sigma)$, we have defined $\Centr(F)$ as being the
set of mapping classes $G\in\Mod(\Sigma)$ such that for some
nonzero integer $m$, the mapping classes $G$ and $F^m$ commute
(cf. Definition \ref{defi:centr}).
\smallskip

Let us make clear that the group $\Centr(F)$ contains the
centralizer of $F$, but is not equal to it in general. The aim
of this definition is to prepare the following proposition
which establishes an isomorphism between $\Centr(F)$ and
$\Norm(\F^s,\,\F^u)$.
\bigskip

\begin{prop}[Realization of $\Centr(F)$ in
$\Norm(\F^s,\,\F^u)$]
                          \label{prop:réalisation_pseudo-Anosov}
\mbox{}\\Let $F$ be a pseudo-Anosov mapping class of
$\Mod(\Sigma)$ and let $\F^s$ and $\F^u$ be the  stable and
unstable foliations of a pseudo-Anosov diffeomorphism $\bar F$
representing $F$. Then there exists a unique
isomorphism\smallskip

$\centrer{\DEF{\psi}{\Centr(F)}{\Norm(\F^s,\,\F^u)}{G}{\bar
G}}$
\smallskip

\noindent where for all $G\in\Centr(F)$, $\bar G$ is a
representative of $G$ in $\Norm(\F^s,\,\F^u)$. \end{prop}
\bigskip

\DEM We show actually the existence of the inverse isomorphism
of $\psi$ that we denote by $\phi$ in this proof:
\smallskip

$\centrer{\DEF{\phi}{\Norm(\F^s,\,\F^u)}{\Centr(F)}{\bar
G}{G},}$
\smallskip

\noindent where $G$ is the isotopy class of $\bar G$. There
exists a nonzero integer $k$ such that $\bar F^k$ fixes all the
separatrices.
\smallskip

Let us first show that $\phi$ is well-defined. For any
diffeomorphism $\bar G$ belonging to $\Norm(\F^s,\,\F^u)$,
$\bar G\bar F^k\bar G^{-1}\bar F^{-k}$ belongs to
$\Norm(\F^s,\,\F^u)$ and fixes all the separatrices. But since
$L$ is a homomorphism in $\RR$, $L(\bar G\bar F^k\bar G^{-1}\bar
F^{-k})=0$, so $\bar G\bar F^k\bar G^{-1}\bar F^{-k}=\Id$ Hence
$\bar G$ commutes with $\bar F^k$, hence $G$ belongs to
$\Centr(F)$.
\smallskip

Let us show that the homomorphism $\phi$ is injective. For any
diffeomorphism $\bar A\in\Ker(\phi)$, $\bar A$ is isotopic to
the identity. Then $\bar A$ preserves the isotopy class of
curves, hence $\bar A$ cannot be pseudo-Anosov, hence $\bar A$
is a diffeomorphism of finite order, isotopic to the identity.
According to [FLP], Exposé 12, this implies that $\bar A$ is
the identity.

At last, let us show that the homomorphism $\phi$ is surjective.
Let $G$ be an element of $\Centr(F)$: there exists a nonzero
integer $\ell$ such that $G$ and $F^\ell$ commute. Let
$\bar{G'}$ be a representative of $G$ in $\Diff(\Sigma',\,\P)$;
then $\bar{G'}\bar F^\ell\bar{G'}^{-1}$ is isotopic to $\bar
F^\ell$, so according to Theorem III, Exposé 12 in [FLP] (see
also Theorem 2.14 in [BoPa]), there exists a diffeomorphism
$\bar H$ isotopic to the identity such that $\bar H\bar{G'}\bar
F^\ell(\bar H\bar{G'})^{-1}=\bar F^\ell$. Let us set $\bar
G=\bar H\bar {G'}$. We get then $\bar{G}\bar
F^\ell\bar{G}^{-1}=\bar F^\ell$. Hence $\bar G$ preserves the
stable and unstable foliations of $\bar F^\ell$, hence $\bar G$
belongs to $\Norm(\F^s,\,\F^u)$. Thus any $G$ of $\Centr(F)$
has an preimage in $\Norm(\F^s,\,\F^u)$.\fin
\bigskip

We are now ready to prove Proposition
\ref{prop:structure_centralizer}.
\medskip

\noindent \textbf{Proposition \ref{prop:structure_centralizer}}
(Structure of $\Centr(F)$).\\
 {\em
Let $\Sigma$ be a connected surface without boundary and let
$F$ be a pseudo-Anosov mapping class in $\Mod(\Sigma)$. There
exists a surjective homomorphism $\ell_F\,:\ \Centr(F)\to\ZZ$
satisfying the following properties:
\begin{itemize}
\item[\;\;(i)] The kernel $\Ker(\ell_F)$ coincides with the
    set of all the finite order mapping classes of
    $\Centr(F)$.
\item[\;(ii)] The kernel $\Ker(\ell_F)$ is a group of order
    smaller than or equal to $6|\chi(\Sigma)|$.
\item[(iii)] The homomorphism $\ell_F$ does not depend on $F$,
    but only on $\Centr(F)$ up to the sign. If for two
    pseudo-Anosov $F$ and $F'$ there exist two positive
    integers $p$ and $q$ such that $F^p={F'}^q$, then
    $\Centr(F)=\Centr(F')$ and $\ell_F=\ell_{F'}$.
\item[(iv)] The group $\Centr(F)$ is a semi-direct product
    $\Per\rtimes\ZZ$ where $\Per$ is a finite group
    isomorphic to $\Ker(\ell_F)$. In particular,
    $\Centr(F)$ is virtually infinitely cyclic.
\end{itemize}}
\bigskip

\DEM Let $\Sigma$ be a surface without boundary. According to
Proposition \ref{prop:réalisation_pseudo-Anosov}, there exists
an isomorphism $\psi\;:\; \Centr(F)\to\Norm(\F^s,\,\F^u)$,
$G\mapsto \bar G$, such that $\bar G$ is a representative of
$G$. Let $\ell_F$ be the homomorphism $L\rond\psi$, that we
normalize, so that $\Im(\ell_F)=\ZZ$. Let us prove items (i),
(ii) and (iii).
\medskip

(i) It is clear that for all $G\in\Centr(F)$, if
$\ell_F(G)\not=0$, then for all integers $m$ different from 1,
$\ell_F(G^m)\not=\ell_F(G)$, so $G$ is not of finite order. In
the contrary, $\Ker(\ell_F)$ is equal to $\psi^{-1}(\Ker(L))$
which is finite, so the elements of $\Ker(\ell_F)$ are of
finite order.
\medskip

(ii) According to Proposition
\ref{prop:réalisation_pseudo-Anosov}, the homomorphism $\psi$ is an
isomorphism, so $\Ker(\ell)$, which is equal to
$\Ker(L\rond\psi)$, is isomorphic to $\Ker(L)$. Let us compute
the cardinality of $\Ker(L)$. According to Proposition
\ref{prop:properties_of_L}.(iv), we know it is finite.

Let us recall that all the indices of the singular points are
greater than or equal to 3. We lean on Proposition
\ref{prop:indices_and_chi}:
\smallskip

\centrer{$\displaystyle
\chi(\Sigma)=\sum_{P\in\S}\bigg(1-\frac{\Ind(\F^s,\,\F^u\,:\,P)}{2}\bigg)$.}
\smallskip

\noindent The surface $\Sigma$ is of negative Euler
characteristic, so the set of singular points is nonempty. Let
$X$ be the set of all the separatrices of $\F^u$. Let $P$ be a
singular point and let $k$ be the index
$\Ind(\F^s,\,\F^u\,:\,P)$. The $k$ separatrices ending at $P$
bring together a $\left(1-\frac{k}{2}\right)$-contribution to
the Euler characteristic of $\Sigma$, so each of them bring a
contribution of $\frac{1}{k}-\frac{1}{2}$. Since $k\geqslant
3$, the contribution per separatrix ending at $P$ to the Euler
characteristic of $\Sigma$ is smaller than or equal to
$-\frac{1}{6}$. This is also true for all the separatrices of
$\F^u$, so the cardinality of $X$, the set of separatrices of
$\F^u$, must equal at most $6|\chi(\Sigma)|$. Now the action of
$\Ker(L)$ on $X$ is free according to Proposition
\ref{prop:properties_of_L}.(iii), so the cardinality of
$\Ker(L)$ is smaller than or equal to the cardinality of $X$,
whence:
\smallskip

\centrer{$|\Ker(L)|\leqslant 6|\chi(\Sigma)|$.}

\smallskip

\noindent But $\Ker(L)$ is isomorphic to $\Ker(\ell_F)$ via
$\psi$, so:
\smallskip

\centrer{$|\Ker(\ell_F)|\leqslant6|\chi(\Sigma)|$.}
\bigskip

(iii) Let $G$ be a pseudo-Anosov mapping class such that
$\Centr(G)=\Centr(F)$ and let $\ell_G$ and $\ell_F$ be the two
homomorphisms associated to $G$ and $F$. The periodic elements of
$\Centr(G)$ and $\Centr(F)$ are the same, so
$\Ker(\ell_G)=\Ker(\ell_F)$. Hence we have the following
commutative diagram where the lines are exact and the vertical
full arrows are equalities. It follows that the dotted arrow is
an isomorphism.
\smallskip

\centrer{$\begin{diagram}
  \node{1}\arrow{e}\node{\Ker(\ell_G)}\arrow{e}\arrow{s,r}{=}
  \node{\Centr(G)}\arrow{e,t}{\ell_G}\arrow{s,r}{=}\node{\ZZ}\arrow{e}\arrow{s,..}\node{1}\\
  \node{1}\arrow{e}\node{\Ker(\ell_F)}\arrow{e}\node{\Centr(F)}\arrow{e,t}{\ell_F}\node{\ZZ}\arrow{e}\node{1}
\end{diagram}$}\smallskip

\noindent But $\ell_F$ (respectively $\ell_G$) is the quotient
map $\Centr(G)\to\Centr(G)/\Ker(\ell_G)\cong \ZZ$ (resp.
$\Centr(F)\to\Centr(F)/\Ker(\ell_F)\cong \ZZ$). Hence $\ell_G$
and $\ell_F$ coincide up to an isomorphism of $\ZZ$, so
$\ell_F$ and $\ell_G$ are equal up to multiplication by $-1$.

Let $F$ and $F'$ be two pseudo-Anosov mapping classes such that
there exist two positive integers $p$ and $q$ such that
$F^p={F'}^q$, let $\ell_F$ and $\ell_{F'}$ be the two homomorphisms
associated to $F$ and $F'$. Since $F^p={F'}^q$, we have
$\Centr(F)=\Centr(F')$, so according to what precedes,
$\ell=\ell'$ or $\ell=-\ell'$. Let $u$ and $v$ be two positive
integers such that $\ell_F(F)=u$ and $\ell_{F'}(F')=v$. Then
$\ell_F(F^p)=pu$ and $\ell_{F'}(F^p)=\ell_{F'}({F'}^q)=qv$, so
$\ell_F(F^p)$ and $\ell_{F'}(F^p)$ have the same sign, hence
$\ell_F=\ell_{F'}$.
\bigskip

(iv) Let $M$ be a mapping class of $\Centr(F)$ such that
$\ell(M)=1$. Then the quotient set $\Centr(F)/\langle M\rangle$
is in bijection with $\Ker(\ell)$, so we have
$[\Centr(F):\langle M\rangle]=|\Ker(\ell)|$. Thus $\langle
M\rangle$ is an infinite cyclic subgroup of finite index of
$\Centr(F)$. Moreover, the exact sequence
\smallskip

\centrer{$1\to\Per\xrightarrow{\ incl.\
}\Centr(F)\xrightarrow{\ \ell_F\ }\ZZ\to1$}
\smallskip

\noindent is split since the last but one term of the sequence
is $\ZZ$. So $\Centr(F)$ is a semi-direct product. This is a
direct product if and only if there exists a pseudo-Anosov
mapping class belonging to $\ell_F^{-1}(\{1\})$, whose action
on the separatrices is trivial.\fin
\bigskip


\subsection{Reducible mapping classes, Thurston's theory and the canonical reduction system}
\medskip

\begin{lem}[Canonical reduction system on subsurfaces]
      \label{lem:syst_red_included}
\mbox{}\\ Let $\Sigma$ be a surface and $\Sigma'$ be a
subsurface of $\Sigma$ non-homeomorphic to a pair of pants. Let
$F$ be a mapping class of $\Mod(\Sigma)$ that preserves
$\Sigma'$ and let $F'$ in $\Mod(\Sigma')$ be the restriction of
$F$ to $\Sigma'$.

\begin{itemize}
\item[\;\;(i)] If there exists a reduction curve of $F$ in
    $\Courb(\Sigma)$ that is not included in
    $\Sigma\smallsetminus\Sigma'$, then there exists a
    reduction curve of a nonzero power of $F'$ in
    $\Courb(\Sigma')$.

\item[\;(ii)] Let $x$ be a curve belonging to
    $\Courb(\Sigma')$. If there exists a reduction curve of
    $F$ in $\Courb(\Sigma)$ that intersects $x$, then there
    exists a reduction curve of a nonzero power of $F'$ in
    $\Courb(\Sigma')$ that intersects $x$.

\item[(iii)] Any curve of $\sigma(F)$ non-isotopic to a
    boundary component of $\Sigma'$ is included either in
    $\Sigma'$ or in $\Sigma\smallsetminus\Sigma'$.

\item[(iv)] Moreover,
    $\sigma(F')=\sigma(F)\cap\Courb(\Sigma')$.
\end{itemize} \end{lem}
\medskip

\DEM Item (ii) implies item (i): it is indeed enough to choose
a curve $x$ in $\Sigma'$ that intersects a reduction curve of
$F$. Then, according to item (ii), there exists a reduction
curve of $F'$ included in $\Sigma'$.

Let us show item (ii). Let $x$ be a reduction curve of
$\Courb(\Sigma')$, let $c$ be a reduction curve of $F$ that
intersects $x$. We are going to show that there exists a
reduction curve $c'$ of $F'$ that intersects $x$.
\smallskip

If $c$ is included in $\Sigma'$, there is nothing to be shown.
Let us now assume that $c$ is not included $\Sigma'$. As the
curve $c$ intersects the curve $x$, which is included in
$\Sigma'$, then $c$ must intersect $\bord\Sigma'$.
Consequently, the intersection $c\cap\Sigma'$ consists in a
finite nonzero number of paths with extremities in
$\bord\Sigma'$. Then there exists a nonzero integer $m$ such
that $F^m$ preserves each boundary of $\Sigma'$ and preserves
each connected component of $c\cap\Sigma'$. Since
$I(c,\,x)\not=0$, at least one of these paths cuts the curve
$x$. Let us choose one such path which we call $d$. The
extremities of $d$ lie in one or two boundary components of
$\Sigma'$. In both cases, let $P$ be the pair of pants
corresponding to the intersection of $\Sigma'$ with the tubular
neighbourhood of the union of $d$ and the boundary components
of $\Sigma'$ containing the extremities of $d$ (darkened in
Figures \ref{fig:courbeRedEss1Bord} and
\ref{fig:courbeRedEss2Bords}). Then $F^{2m}$ preserves each
boundary component of $P$. Notice that no boundary component of
$P$ can bound any disk:
\begin{itemize}
  \item if $P$ has only one boundary component $a$ in
      $\bord\Sigma'$, then both of the other boundary
      components are isotopic to the union of a path
      included in $a$ and a path included in $d$, but $a$
      and $d$ do not cobound any bigon, so both boundary
      components of $P$ different from $a$ cannot bound any
      disk;
\begin{figure}[!h]
 \Includegraphics{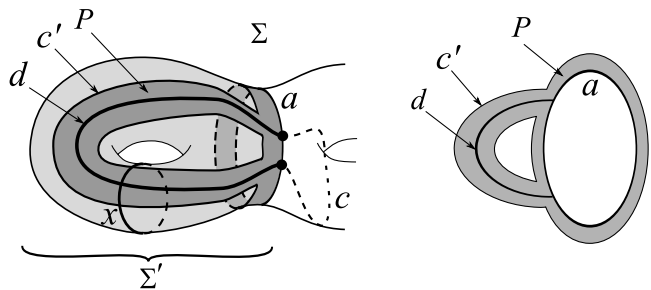}
 \caption{The pair of pants $P$ in $\Sigma'$, when both extremities of $d$ belong to the same boundary $a$ of $\Sigma'$.}
 \label{fig:courbeRedEss1Bord}
\end{figure}
  \item and if $P$ has two boundary components in
      $\bord\Sigma'$, both of these boundary components do
      not bound any disk. Moreover, as they are not
      isotopic in $\Sigma$, they do not cobound any
      cylinder. So the third boundary component of $P$
      cannot bound any disk.
\begin{figure}[!h]
 \Includegraphics{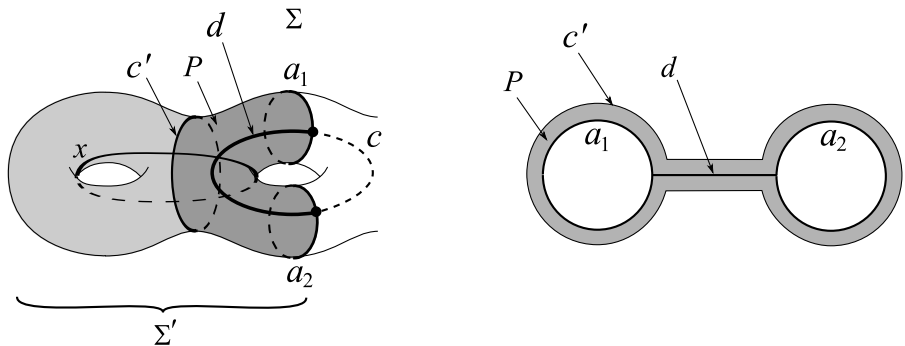}
 \caption{The pair of pants $P$ in $\Sigma'$ when both extremities of $d$ belong to two distinct boundary components $a_1$ and $a_2$ of $\Sigma'$.}
 \label{fig:courbeRedEss2Bords}
\end{figure} \end{itemize}
\smallskip

\noindent Notice also that the curve $x$ cannot be one of the
boundary components of $P$, for $x$ and $d$ intersect. The
curve $x$ cannot be included in $P$ for any curve that is
included in a pair of pants is isotopic to one of its boundary
components. Consequently, there exists a boundary of $P$, say
$c'$, that intersects $x$. This boundary $c'$ cannot be a
boundary component of $\Sigma'$ because the boundary components
of $\Sigma'$ do not intersect $x$. Hence $c'$ belongs to
$\Courb(\Sigma')$, is stable by $F^{2m}$, so is stable by
$(F')^{2m}$, and intersects the curve $x$. We have proved item
(ii).
\medskip

Let us show item (iii). The mapping class $F$ preserves
$\Sigma'$, so $F(\Bord(\Sigma'))=\Bord(\Sigma')$, so the curves
of $\Bord(\Sigma')$ are reduction curves of $F$. But no curve
of $\sigma(F)$ intersects any reduction curve of $F$, hence any
curve of $\sigma(F)$ that is not isotopic to a boundary
component of $\Sigma'$ is included either in $\Sigma'$ or in
$\Sigma\smallsetminus\Sigma'$.
\medskip

Let us show item (iv). Let $x$ be a reduction curve of $F'$. If
$x$ is not an essential reduction curve of $F$, then there
exists a reduction curve $c$ of $F$ in $\Courb(\Sigma)$ that
intersects $x$. Then according to item (ii), there exists a
nonzero integer $m$ and a reduction curve $c'$ of $(F')^m$ in
$\Courb(\Sigma')$ that intersects $x$. Hence $x$ is not an
essential reduction curve of $(F')^m$. But
$\sigma(F')=\sigma((F')^m)$, so $x$ is not an essential
reduction curve of $F'$. In other words, we have
$\sigma(F')\subset \sigma(F)\cap\Courb(\Sigma')$. The converse
inclusion is obvious: if a curve belongs to $\sigma(F)$, it
belongs \emph{a fortiori} to $\sigma(F')$. \fin
\bigskip

\begin{lem}[Characterization of the essential reduction curves]
        \label{lem:curves_red_ess}
\mbox{}\\Let $\Sigma$ be a surface and let $F$ be a mapping
class in $\Mod(\Sigma)$. Let $a$ be an oriented curve such that
$F$ preserves $a$ and its orientation. Then, there exist two
connected subsurfaces $S_1$ and $S_2$ (they may be equal) in
$\Sigma$ such that:
\begin{itemize}
  \item the curve $a$ bounds $S_1$ (respectively $S_2$) on
      the left (rep. on the right),
  \item both surfaces $S_1$ and $S_2$ are stable by $F$,
  \item for all $i\in\{1,\,2\}$, the mapping class induced
      by $F$ in $\Mod(S_i)$, denoted by $\restr{F}{S_i}$,
      is either periodic or pseudo-Anosov.
\end{itemize}
\smallskip

\noindent Let us denote by $S_{12}$ the union of $S_1$ and
$S_2$ along  $a$. Let us denote by $\restr{F}{S_{12}}$ the
mapping class induced by $F$ in $\Mod(S_{12})$. Then $a$
belongs to $\sigma(F)$ if and only if we are in one of the
three following cases:
\begin{itemize}
  \item[a)] $\restr{F}{S_1}$ or $\restr{F}{S_2}$ is
      pseudo-Anosov;
  \item[b)] $\restr{F}{S_1}$ and $\restr{F}{S_2}$ are both
      periodic of the same order $m\geqslant 1$, and
      $\restr{F^m}{S_{12}}$ is a nontrivial power of a
      Dehn twist along the curve $a$;
  \item[c)] $\restr{F}{S_1}$ and $\restr{F}{S_2}$ are
      periodic with orders $m_1$ and $m_2$ respectively
      such that $m_1\not=m_2$.
\end{itemize} \end{lem}
\bigskip

\DEM Let $\Gamma$ be the set of curves $\sigma(F)\cup\{a\}$.
Let us give to $a$ an orientation and let us denote by $S_1$
(respectively $S_2$) the connected component of $\Sigma_\Gamma$
bounding $a$ on the left (resp. on the right). For all
$i\in\{1,\,2\}$, notice that $\sigma(\restr{F}{S_i})=\vide$, so
$\restr{F}{S_i}$ is either pseudo-Anosov, or periodic.
Moreover, the surface $S_{12}$ cannot be a pair of pants, so
there exist some curves in $\Courb(S_{12})$ that intersect $a$.
\smallskip

Let us show first that if $a$ belongs to $\sigma(F)$, then one
of the cases a), b) or c) is satisfied. The cases a), b), c)
describe all the possible cases except the one where
$\restr{F}{S_1}$ and $\restr{F}{S_2}$ are both periodic of
order $m\geqslant 1$ and where $\restr{F^m}{S_{12}}$ coincides
with the identity. Let $c$ be a curve of $\Courb(S_{12})$ which
intersects $a$. But this curve $c$ is preserved by $F^m$, so
$c$ is a reduction curve of $F^m$. Therefore $a$ cannot belong
to $\sigma(F^m)$. So $a$ does not belong to $\sigma(F)$.
\smallskip

Conversely, let us show now that each case a), b), or c)
implies that $a\in\sigma(F)$ (or that there exists a nonzero
integer $p$ such that $a\in\sigma(F^p)$, which is equivalent).
\smallskip

\begin{figure}[!h]
 \Includegraphics{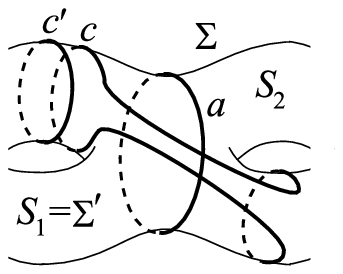}
 \caption{Description of the case a) in the proof of Lemma \ref{lem:curves_red_ess}.}
 \label{fig:caracterisationCourbeEssentielle}
\end{figure}

In the case a), let us denote by $\Sigma'$ the subsurface $S_1$
or $S_2$ on which $F$ induces a pseudo-Anosov mapping class.
Then $\Sigma'$ is not a pair of pants. If there existed a
reduction curve $c$ of $F$ that intersected $a$, it would not
be included in $\Sigma\smallsetminus\Sigma'$, so we could apply
item (i) of Lemma \ref{lem:syst_red_included}: there would
exist a reduction curve $c'$ of a nonzero power of $F$ that
would be included in $\Sigma'$ (cf. Figure
\ref{fig:caracterisationCourbeEssentielle}). This is absurd for
the restriction of $F$ (and of its nonzero powers) to $\Sigma'$
is pseudo-Anosov. Hence any curve $c$ of $\Courb(\Sigma)$ such
that $I(a,\,c)\not=0$ is not a reduction curve of $F$, so the
reduction curve $a$ of $F$ is
 an essential reduction curve of $F$.

\smallskip

In the case b), there exists a nonzero integer $\ell$ such that
$(\restr{F}{S_{12}})^m=T_a^{\,\ell}$. But then for any integer
$p$ and any curve $c$ in $\Courb(S_{12})$ intersecting $a$,
applying the famous intersection formula between a curve and
its image by the power of a Dehn twist (cf. \mbox{[FLP]}], we
have:
\smallskip

\centrer{$I((\restr{F}{S_{12}})^{mp}(c),\,c)=I\big(\,T_a^{\,\ell
p}(c),\ c\,\big)=|\ell p|\ I(a,\ c)^2\not=0$.}
\smallskip

\noindent Hence no curve $c$ that intersect $a$ can be
preserved by some power of $(\restr{F}{S_{12}})^m$, hence $a$
belongs to $\sigma((\restr{F}{S_{12}})^m)$. Hence $a$ belongs
to $\sigma(F^m)$ according to item (iv) of Lemma
\ref{lem:syst_red_included}.
\smallskip

In the case c), $S_1$ and $S_2$ cannot be equal. For all
$i\in\{1,\,2\}$, let $\Mod(S_i,\,a)$ be the set of isotopy
classes of diffeomorphisms of $S_i$ fixing (a representative
of) the curve $a$ pointwise, and let $F_i$ be a mapping class
in $\Mod(S_i,\,a)$ such that $F_i$ induces the mapping class
$\restr{F}{S_i}$ in $\Mod(S_i)$. Let us introduce some
notation:
\begin{itemize}
  \item For two integers $k$ and $\ell$, let us denote
      their greatest common divisor by $k\wedge\ell$ and
      their least common multiple by $k\vee\ell$.
  \item For all $i\in\{1,\,2\}$, let $d_i$ be an integer
      such that $F_i^{m_i}=T_a^{\,d_i}$ with $d_i\wedge
      m_i=1$ according to Lemma \ref{lem:without_torsion}.
  \item Let us set $m'_1=\frac{m_1}{m_1\wedge
      m_2}=\frac{m_1\vee m_2}{m_2}$.
  \item Let us set $m'_2=\frac{m_2}{m_1\wedge
      m_2}=\frac{m_1\vee m_2}{m_1}$; thus, $m'_1$ and
      $m'_2$ are coprime.
\end{itemize} Then for all $i\in\{1,\,2\}$, we have
\smallskip

\centrer{$F_i^{(m_1\vee m_2)}=T_a^{\,(\frac{m_1\vee
m_2}{m_i}d_i)}=T_a^{\,m'_{(3-i)}d_i}.$}
\smallskip

\noindent Hence there exists an integer $q$ such that:
\smallskip

\centrer{$(\restr{F}{S_{12}})^{(m_1\vee m_2)}=T_a^{(m'_2d_1 +
m'_1d_2)}\,T_a^{\,q(m_1\vee m_2)}. $}
\smallskip

\noindent Since $m_1$ and $m_2$ are different, one of them, by
instance $m_2$, is not equal to $m_1\vee m_2$, so $m'_1$ is
different from 1. But $d_1$ and $m'_2$ are both coprime with
$m'_1$, so $m'_2d_1$ is coprime with $m'_1$. Now $m'_1$ divides
$m'_1d_2+q(m_1\vee m_2)$, so $m'_1$ is coprime with
$m'_2d_1+m'_1d_2+q(m_1\vee m_2)$, so $m'_2d_1+m'_1d_2+q(m_1\vee
m_2)$ is nonzero. Hence $(\restr{F}{S_{12}})^{(m_1\vee m_2)}$
is a nonzero power of $T_a$. Thus, if we set $p=(m_1\vee m_2)$
and if we consider $F^p$ instead of $F$, then we are back to
the case b) with $m=1$. \fin
\bigskip

As a corollary, we establish a link between the canonical reduction systems of $F$, $G$ and $FG$,
assuming that $F$ and $G$ are two commuting mapping classes.
\bigskip

\begin{prop}
          \label{prop:sigma(FG)}
Let $F$ and $G$ be two commuting mapping classes. Then,
\smallskip

\centrer{$\sigma(FG)\subset\sigma(F)\cup\sigma(G)$.}
\end{prop}
\medskip

\DEM Let $a$ be a curve belonging to $\sigma(FG)$. Since $FG$
commutes with $F$ and with $G$, we have
$I(a,\,\sigma(F)\cup\sigma(G))=0$. Let us set
$\A=\sigma(F)\cup\sigma(G)$. Let us assume that $a$ does not
belong to $\A$ and let us show that this is absurd.

Let $m$ be a nonzero integer such that $F^m $ and $G^m$ preserve
each subsurface in $\Sub_\A(\Sigma)$. According to Lemma
\ref{lem:syst_red_included}.(iv), the restrictions of the
mapping classes induced by $F$ and $G$ on each subsurface in
$\Sub_\A(\Sigma)$ have empty canonical reduction systems, so
they are either pseudo-Anosov or periodic. We can even assume
that they are either pseudo-Anosov or the identity mapping
class, provided that the integer $m$ is large enough.

Since $a$ does not belong to $\A$ although $I(a,\,\A)=0$, there
exists a subsurface $S$ belonging to $\Sub_\A(\Sigma)$ such
that $a$ is included in $S$. Since $F$ and $G$ commute with
$FG$, they preserve $\sigma(FG)$, so we can assume that $F^m$
and $G^m$ preserve $a$, even if it means multiplying $m$ by
some positive integer. Therefore $a$ is a reduction curve of
$F^m$ and $G^m$. So the restrictions of $F^m$ and $G^m$ to $S$
cannot be pseudo-Anosov, so they are the identity mapping
class. Hence the restriction of $(FG)^m$ to $S$ is equal to the
identity mapping class, hence according to Lemma
\ref{lem:curves_red_ess}, the reduction curve $a$ of  $FG$
cannot be essential. This is the expected contradiction. \fin

\end{appendix}

\addcontentsline{toc}{part}{\protect\textsc{\textbf{References}}}

\section*{References}
\bigskip

\begin{supertabular}{lp{14cm}l}

{[ArYo]} &
 P. Arnoux, J.-C. Yoccoz. \emph{Construction de difféomorphismes
 pseudo-Anosov} C. R. Acad. Sci. Paris Sér. I Math. (1) \textbf{292} (1981), 75-78.
 \\

{[At1]} &
 E. Artin, \emph{Theorie der Zöpfe}, Abhandlungen Hamburg \textbf{4} (1925),
 47-72,  (JMF 51.0450.01).
 \\

{[At2]} &
 E. Artin. \emph{Theory of braids}, Annals of Math. (2) \textbf{48} (1947), 101-126.
 \\

{[At3]} &
 E. Artin, \emph{Braids and permutations}, Annals of Math. (2)
 \textbf{48} (1947),
 643-649.
 \\

{[BeMa]} &
 R. Bell, D. Margalit, \emph{Braid groups and the
 co-Hopfian property}, Journal of Algebra \textbf{303} (2006), 275-294.
 \\

{[BkMt]} &
 J. Berrick, M. Matthey \emph{Stable classical groups and strongly torsion generated groups},
 Comment. Math. Helv. \textbf{84} (2009), 909-928.
 \\

{[Bi]} &
 J. Birman, \emph{Braids, Links and Mapping Class Groups},
 Annals of Mathematics Studies, Princeton University Press
 \textbf{82} (1975), Princeton.
 \\

{[BiBr]} & J. Birman, T. Brendel, \emph{Braids: A survey}, in
Handbook of Knot Theory (Ed. W. Menasko and T. Thistlethwaite)
 (2005), Elsevier.
 \\

{[BiHi]} &
 J. Birman, H. Hilden, \emph{On isotopies of homeomorphisms of Riemann surfaces},
 Annals of Math. (2)  \textbf{97} (1973), 424-439.
 \\

{[BiLuMc]} &
 J. Birman, A. Lubotzky, J. McCarthy, \emph{Abelian
 and Solvable Subgroups of the Mapping Class Group}, Duke
 Mathematical Journal (4) \textbf{50} (1983),  1007-1020.
 \\

{[Bk]} &
 N. Bourbaki, \emph{Groupes et alg\`ebres de Lie}, chapters 4, 5 et 6, Hermann, Paris, 1968.
 \\

{[BlCa]} &
 S.A. Bleiler, A.J. Casson, \emph{Automorphisms of surfaces after
 Nielsen and Thurston}, London Mathematical Society Students Texts
 \textbf{9} (1988), Cambridge University Press, Cambridge.
 \\

{[BoPa]} &
 C. Bonatti, L. Paris, \emph{Roots in the mapping class groups},
 Proc. Lond. Math. Soc. (Third Series) (2) \textbf{98} (2009), 471-503.
 \\

{[Bri]} &
 E. Brieskorn, \emph{Sur les groupes de tresses (d'apr\`es V.I. Arnold)},
 Séminaire Bourbaki, 24\`eme
ann´ee (1971/72), Exp. \textbf{401},
 Lecture Notes Math. \textbf{317},
 Springer, Berlin Heidelberg New York, (1973)  1-24.
 \\

{[BriSa]} &
 E. Brieskorn, K. Saito, \emph{Artin Gruppen und Coxeter-Gruppen},
 Invent. Math. \textbf{17} (1972), 245-271.
 \\

{[Bro]} &
 L.E.J. Brouwer, \emph{Über die periodischen Transformationen
 der Kugel}, Math. Ann. \textbf{80} (1919), 39-41.
 \\


{[Ch1]} &
 R. Charney, \emph{Artin groups of finite type are biautomatic}, Math. Ann. \textbf{292} (1992), 671-683.
 \\

{[Ch2]} &
 R. Charney, \emph{Geodesic automation and growth functions for Artin groups of finite type},
 Math. Ann. \textbf{301} (1995), 307-324.
 \\

{[ChCr]} &
 R. Charney, J. Crisp,
 \emph{Automorphism groups of some affine and finite type Artin groups},
 Mathematical Research Letters \textbf{12} (2005), 321-333.
 \\

{[CoPa]} &
 A.M. Cohen, L. Paris, \emph{On a theorem of Artin},
 J. Group Theory \textbf{6} (2003), 421-441.
 \\

{[CrPa1]} &
 J. Crisp, L. Paris,
 \emph{The solution to a conjecture of Tits on subgroups
 generated by the squares of the generators of an Artin group},
 Inventiones mathematicae \textbf{145} (2001), 19-36.
 \\

{[CrPa2]} &
 J. Crisp, L. Paris,
 \emph{Artin groups of type $B$ and $D$},
 Advances in Geometry \textbf{5} (2005), 607-636.
 \\

{[Ct]} &
 F. Castel, \emph{Geometric representation of the Braid groups},
 translation of the author's PhD. thesis in English,
 available on Archiv.\\

{[Cx]} &
 H.S.M. Coxeter. \emph{The complete enumeration of finite groups of the form
 ${R_i}^2=(R_i R_j )^{k_{i j}} = 1$},
 J. London Math. Soc. \textbf{10} (1935), 21-25.
 \\

{[De1]} &
 M. Dehn, \emph{On curve  systems on two-sided surfaces, with
 application to the mapping problem}, Lecture (supplemented) to the math.
 colloquium,  Breslau 11-2-1922, translated in English in [De3], 234-252.
 \\

{[De2]} &
 M. Dehn, \emph{Die Gruppe der Abbildungsklassen}, Acta Math. \textbf{69}, 135-206;
 translated in English in [De3], 256-362.
 \\

{[De3]} &
 M. Dehn, \emph{Papers in Group Theory and Topology}, Springer-Verlag, New York (1987),
 translated from the German and introduced by John Stillwell.
 \\

{[Dl]} &
  P. Deligne, \emph{Les immeubles des groupes de tresses g´en´eralis´es},
  Invent. Math. \textbf{17} (1972), 273-302.
 \\

{[DyGr]} &
 J.L.Dyer, E.K. Grossman. \emph{The automorphism groups of the braid
 groups}, Amer. J. Math. (6) \textbf{103} (1981), 1151-1169.
 \\

{[Ei]} &
 S. Eilenberg, \emph{Sur les transformations périodiques de la
 surface de sphère}, Fund. Math. \textbf{22} (1934), 28-41.
 \\

{[Ep]} &
 D.B.A. Epstein, \emph{Curves on $2$-manifolds and
 isotopies}, Acta Math. \textbf{115} (1966), 83-107.
 \\

{[Fa]} &
 B. Farb, \emph{Problems on Mapping Class Groups and Related Topics},
 Proceedings of Symposia in Pure Mathematics \textbf{74} (2006), American Mathematical
 Society, AMS Bookstore.
 \\

{[FaMa]} &
 B. Farb, D. Margalit \emph{A Primer on Mapping Class Group},
 version 3.1 (June 2009) is available on the web page:
 http://www.math.uchicago.edu/~margalit/mcg/mcgv31.pdf
 \\

{[FdNe]} &
 E. Fadell, L. Neuwirth, \emph{Configuration spaces},
 Math. Scand. \textbf{10} (1962), 111-118.
 \\

{[FLP]} &
 A. Fathi, F. Laudenbach, V. Poénaru, \emph{Travaux de
 Thurston sur les surfaces}, Séminaire Orsay, Astérisque \textbf{66} et \textbf{67} (1979),
 Société Mathématique de France, Paris.
 \\

{[Hi]} & Morris W. Hirsch, \emph{Differential topology},
 Graduate texts in Mathematics \textbf{33} (1976), Springer-Verlag, New York Inc.
\\

{[Hu]} &
 J.E. Humphreys, \emph{Reflection groups and Coxeter groups},
 Cambridge Studies in Advanced  Mathematics \textbf{29} (1990),
 Cambridge University Press, Cambridge.
 \\

{[Hv1]} &
 W.J. Harvey, \emph{Geometric structure of surface mapping class groups},
 Homological group theory (Proc. Durham, 1977),
 London Math.Soc. Lecture Note Ser. \textbf{36} (1979), 255-269,
 Cambridge Univ. Press, Cambridge New York.\\

{[Hv2]} &
 W.J. Harvey, \emph{Boundary structure of the modular group}, in
 \emph{Riemann surfaces and related topics: Proceedings of the 1978 Stony Brook
 Conference}, edited by I. Kra and B. Maskit, Annals of Math.
 Studies \textbf{97} (1981), Princeton University Press, 245-251.
 \\

{[HvKo]} &
 W.J. Harvey, M. Korkmaz \emph{Homomorphisms from mapping class groups},
 Bull. London Math. Soc. \textbf{37} (2005), 275-284.
 \\

{[HtTh]} &
 A. Hatcher, W. Thurston, \emph{A presentation for the mapping class group of a closed orientable surface},
  Topology (3) \textbf{19} (1980), 221-237.
 \\

{[Iv1]} &
 N.V. Ivanov, \emph{Automorphisms  of  Teichmu¨ller  modular  groups},  Lecture
Notes in Math., \textbf{1346} (1988), Springer-Verlag, Berlin and New York, 199–270.
 \\

{[Iv2]} &
 N.V. Ivanov, \emph{Subgroups of Teichmüller Modular
 Groups}, translated from the Russian by E.J.F. Primrose,
 Translations of Mathematical Monographs
 \textbf{115} (1992), American Mathematical Society, Providence, RI.
 \\

{[Iv3]} &
 N.V. Ivanov, \emph{Automorphisms of complexes of curves and of Teichmüller Spaces},
 Math. Res. Notices \textbf{14} (1997), 651-666; also in:
 \emph{Progress in knot theory and related topics}, Travaux en
 cours \textbf{56} (1997), Hermann, Paris,  113-120.
 \\

{[IvMc]} &
 N.V. Ivanov, J. McCarthy, \emph{On injective
 homomorphims between Teichmüller Modular Groups I}, Inventiones
 mathematicae \textbf{135} (1999), 425-486, Springer-Verlag.
 \\

{[J]} &
 D. Johnson, \emph{The structure of Torelli group I: A finite set of generators for $\I$},
 Annals of Mathematics \textbf{118} (1983), 423-442.
 \\

{[KaTu]} &
 C. Kassel, V. Turaev, \emph{Braid groups}, Graduate Texts in Mathematics, \textbf{247} (2008),
 Springer, New York.\\

{[Ke1]} &
 S.P. Kerckhoff, \emph{The Nielsen realisation problem}, Bull. AMS (3)
 \textbf{2} (1980),  452-454.
 \\

{[Ke2]} &
 S.P. Kerckhoff, \emph{The Nielsen realisation problem}, Ann. of math. (2),
 \textbf{117} (1983), 235-265.
 \\

{[Kj]} &
 B. de Kerékj\`art\`o. \emph{Über die periodischen
Tranformationen der Kreisscheibe und der Kugelfläche}, Math. Annalen
\textbf{80} (1919), 3-7.
 \\

{[Ko1]} &
 M. Korkmaz, \emph{Automorphisms of complexes of curves on punctured spheres
 and on punctured tori}, Topology and its Applications (2) \textbf{95} (1999),
 85-111.
 \\

{[Ko2]} &
 M. Korkmaz, \emph{Low-dimensional homology groups of mapping class groups: a survey},
 Turkish Journal of Mathematics (1) \textbf{26},
 (2002).
 \\

{[LaPa]} &
 C. Labruère, L. Paris, \emph{Presentations for the punctured mapping class groups
 in terms of Artin groups},
 Algebr. Geom. Topol. \textbf{1} (2001), 73-114.
 \\

{[Le]} &
 H. Van der Leck, \emph{The Homotopy Type of Complex Hyperplane
 Complements}, Ph.D. Thesis, Utrecht (1994).
 \\

{[Lk]} &
 W.B.R. Lickorish, \emph{A finite set of generators for the
 homeotopy group of a 2-manifold}, Proc. Cambridge Philos. Soc.,
 \textbf{60} (1964), 769-778.
 \\

{[Ln1]} &
 V. Lin, \emph{Artin braids and the groups and spaces connected with them},
 Plenum Publishing Corporation (1982), 736-788.
 \\

{[Ln2]} &
 V. Lin, \emph{Braids and Permutations},
 arXiv:math.GR/0404528 v1, 29 Apr 2004.
 \\

{[Lu]} &
 F. Luo, \emph{Automorphisms of complexes of curves},
 Topology (2) \textbf{39} (2000), 283-298.
 \\

{[Mc1]} &
 J.D. McCarthy, \emph{Automorphisms of surface mapping class
 groups. A recent theorem of N. Ivanov}, Invent. Math. \textbf{84} (1986),  49-71.
 \\

{[Mc2]} &
 J.D. McCarthy, \emph{Normalizers and Centralizers of
 pseudo-Anosov mapping classes}, Preprint, 1982
 (This constitutes a part of the author's Ph. D thesis in Columbia University).
 \\

{[McPp]} &
 J.D. McCarthy, A. Papadopoulos \emph{Automorphisms of the complex of domains},
 preprint, Max-Plank Institut, Bonn (2006).
 \\

{[Ni1]} &
 J. Nielsen, \emph{Untersuchungen zur Theorie der geschlossenen zweiseitigen Flächen I},
  Acta Math. \textbf{50} (1927), 189-358.
 \\

{[Ni2]} &
 J. Nielsen, \emph{Untersuchungen zur Theorie der geschlossenen zweiseitigen Flächen II},
 Acta Math. \textbf{53} (1929), 1-76.
 \\

{[Ni3]} &
 J. Nielsen, \emph{Untersuchungen zur Theorie der geschlossenen zweiseitigen Flächen III},
 Acta Math. V. \textbf{58} (1932), 87-167.
 \\

{[Ni4]} &
 J. Nielsen, \emph{Abbildungsklassen endlicher Ordnung},
 {[J]} Acta Math. \textbf{73} (1942), 23-115, {[ISSN 0001-5962]}.
 \\

{[Pa1]} &
 L. Paris, \emph{Parabolic subgroups of Artin groups}, J. Algebra \textbf{196} (1997), 369-399.
 \\

{[Pa2]} &
 L. Paris, \emph{From braid groups to mapping class groups},
 Proceedings of Symposia in Pure Mathematics \textbf{74} (2006), Amer. Math. Soc., Providence, RI, 355-371
 \\

{[PaRo]} &
 L. Paris, D. Rolfsen,
 \emph{Geometric subgroups of mapping class groups},
 Journal für die reine und angewandte Mathematik \textbf{521} (2000), 47-83,
 Walter de Gruyter, Berlin-New York.
 \\

{[PeV]} &
 B. Perron, J.P. Vannier,
 \emph{Groupe de monodromie géométrique des singularités simples},
 Mathematische Annalen \textbf{306} (1996), 231-245,
 Springer-Verlag.
 \\

{[Ro]} &
 D. Rolfsen, \emph{Knots and links}, Publish or Perish,
 Berkeley 1976.
 \\

{[Se]} &
 V. Sergiescu, \emph{Graphes planaires et présentations des groupes de tresses},
 Mathematische Zeitschrift, vol. \textbf{214} (1993), 477-490,
 Springer-Verlag.
 \\

{[Th]} &
 W.P. Thurston, \emph{On the geometry and dynamics of surfaces.},
 Bull. Amer. Math. Soc. (N.S.) (2) \textbf{19} (1988), 417-431.
 \\

{[Ti]} &
 J. Tits, \emph{Normalisateurs de tores. I. Groupes de Coxeter
 étendus}, J. Algebra vol \textbf{4} (1966), 96-116.
 \\

{[VL]} &
 H. Van der Lek, \emph{The homotopy type of complex hyperplane complements}, Ph. D. Thesis,
 Nijmegen, 1983.
 \\

{[W]} &
 B. Wajnryb, \emph{Artin groups and geometric monodromy}, Invent. Math. \textbf{138} (1999), 563-571.
\\

\end{supertabular}

\end{document}